\documentclass[review]{elsarticle} 
\usepackage{lineno,hyperref}
\usepackage{graphicx} 
\usepackage[toc,page]{appendix} 
\usepackage{color} 
\usepackage{caption}
\usepackage{amsfonts}
\usepackage{array}
\usepackage{float}

 \modulolinenumbers[5] 
\newcommand{\fr}{\displaystyle \frac} 
\newcommand{\del}{\partial} 
\newcommand{\be}{\begin{equation}} 
\newcommand{\ee}{\end{equation}} 
\newcommand{\bea}{\begin{eqnarray}} 
\newcommand{\eea}{\end{eqnarray}} 
\newcommand{\ba}{\begin{array}} 
\newcommand{\ea}{\end{array}}  
\newcommand{\sgn}{\mathop{\mathrm{sgn}}}

\usepackage{tikz}
\tikzset{
main node/.style={inner sep=0,outer sep=0},
label node/.style={inner sep=0,outer ysep=.2em,outer xsep=.4em,font=\scriptsize,overlay},
strike out/.style={shorten <=-.2em,shorten >=-.5em,overlay}
}

\newcommand{\bcancelto}[3][]{\tikz[baseline=(N.base)]{
  \node[main node](N){$#2$};
  \node[label node,#1, anchor=north west] at (N.south east){$#3$};
  \draw[strike out,-latex,#1]  (N.north west) -- (N.south east);
}}

\begin{document}  

%
%
%
%
%

\title{A Flexible Velocity Boltzmann Scheme for Convection-Diffusion Equations} 

\author[Aero-IISc,UNISA]{S.V.Raghurama Rao} \fnref{SVR} \ead[S.V. Raghurama Rao]{raghu@iisc.ac.in}
\author[K.S.Shrinath,HAL]{K.S.Shrinath}  \ead[K.S.Shrinath]{shrinath.k.s@gmail.com} 
\author[RBI]{Ankit Ruhi} \ead[Ankit Ruhi] {ankruh@iisc.ac.in} 
\author[UNISA]{Veeredhi Vasudeva Rao}  \ead[Veeredhi Vasudeva Rao]{vasudvr@unisa.ac.za} 
 \address[Aero-IISc]{Department of Aerospace Engineering, Indian Institute of Science, Bangalore, India.}
 \address[UNISA]{Department of Mechanical and Industrial Engineering, University of South Africa (UNISA), Johannesburg}
\address[K.S.Shrinath]{Research Scholar, Department of Aerospace Engineering, Indian Institute of Science, Bangalore, India} 
\address[HAL]{Hindustan Aeronautics Limited, Bangalore, India} 
\address[RBI]{Consultant Applied Mathematician, Reserve Bank of India, India} 
\fntext[SVR]{Visiting faculty at UNISA; present and permanent affiliation: Indian Institute of Science, Bangalore} 

\begin{abstract} 
A framework of finite-velocity model based Boltzmann equation has been developed for convection-diffusion equations.  These velocities are kept flexible and adjusted to control numerical diffusion. A flux difference splitting based kinetic scheme is then introduced for solving a wide variety of nonlinear convection-diffusion equations numerically. Based on this framework, a generalized kinetic Lax-Wendroff scheme is also derived, recovering the classical Lax-Wendroff method as one of the choices. Further, a total variation diminishing version of this kinetic flux difference splitting scheme is presented, combining it with the kinetic Lax-Wendroff scheme using a limiter function. The numerical scheme has been extensively tested and the results for benchmark test cases, for 1D and 2D nonlinear convection and convection-diffusion equations, are presented.  
\end{abstract} 
 
\maketitle 

\section{Introduction}  
     Conservation equations of fluid and gas dynamics are of importance in several branches of engineering, applied physics and applied mathematics.  Typically physical principles representing conservation of mass, momentum, energy and species lead to such conservation equations.  The well-known Navier-stokes equations, representing conservation principles, model the nonlinear convection and diffusion processes occuring in fluid flows.  The essential features of the Navier-Stokes equations are easily represented by a scalar model, the convection-diffusion equation, with a nonlinear convection term.  In this paper, a new finite velocity Boltzmann scheme is presented for simulating nonlinear convection-diffusion equations.  The kinetic schemes are some of the alternative approaches to traditional macroscopic description based algorithms for modeling compressible fluid flows, pursued from the 1970s.  Some of the kinetic schemes developed earlier can be found in the following papers \cite{Sanders_Prendergast,Pullin,Rietz,Deshpande,Mandal_Deshpande,Kaniel,Perthame,SVRRao_Deshpande,Prendergast_KunXu}.  This research is further continued by pursuing {flexible } velocity based Boltzmann schemes as in \cite{Driollet_Natalini,SVRRao_Balakrishna,SVRRao_SubbaRao,Arun_Lukacova_Prasad_SVRRao,SVRRao_Venkat} which are also related to the relaxation schemes \cite{Jin_Xin}. In this paper, this research is advanced further by making the velocities flexible to control numerical diffusion.  A variant of the current work, which has been designed specifically for simulating compressible flows, together with an alternative idea to fix entropy violation, has been presented in \cite{SVRR_VVR_SKS}.  In this paper, the focus is on developing kinetic schemes for nonlinear convection equations and convection-diffusion equations. The basic foundation of the number of finite velocities is also different in this work, with the focus being on the minimal number of velocities required to obtain an efficient scheme.  The strategy used to obtain second order accuracy is also different. Further, a specific strategy is also introduced to deal with the source terms that arise naturally in some convection systems.


   
In the current work, first a flexible velocity Boltzmann equation (FVBE) is introduced, with two velocities in 1-D and four velocities in 2-D.  This model resembles the discrete velocity Boltzmann system in structure, though the velocities are flexible and not constants.  A flux difference splitting numerical scheme is then developed based on this framework.     The flexible velocities at the level of designing the  numerical algorithm are fixed based on sub-characteristic condition, derived from the Chapman-Enskog expansion (named {\em Kinetic Flux Difference Splitting} KFDS) or by using the Rankine-Hugoniot(RH) jump conditions (named KFDS+).  The use of RH condition to fix the numerical diffusion results in a low diffusive algorithm capable of exact capturing of steady discontinuities.  To obtain higher order accuracy, a second order Kinetic Lax-Wendroff type scheme (KLW) has been developed based on the FVBE.  Further,  a total variation diminishing model of KFDS scheme (TVD-KFDS) has been constructed using limiters.  This results in deploying optimal numerical diffusion based on Rankine-Hugoniot condition in resolving the discontinuities and sub-characteristic condition based numerical diffusion in smooth regions.  Further, the formulation is extended to systems of hyperbolic PDEs representing the shallow water flows, together with a special way of treating the source terms.  Several typical bench-mark problems are solved using the new algorithm to demonstrate its capability in resolving the propagating and diffusing nonlinear waves in numerical simulations.  

 The new scheme is first derived for pure convection equations in the following sections. 
 
\section{Kinetic Schemes for Hyperbolic Conservation Equations} 
The kinetic or Boltzmann schemes are based on the fact that conservation equations of gas dynamics can be obtained as moments of the classical Boltzmann equation.  This strategy can be extended the to scalar conservation laws, with appropriate definition of the equilibrium distribution functions. For example, 1-D Burgers equation, given by 
\be \label{1D_Burgers}
\fr{\del u}{\del t} + \fr{\del g(u)}{\del x} = 0, \ \textrm{where} \ g(u) = \fr{1}{2} u^{2}   
\ee 
can be written in the moment form 
\be \label{BE_as_moments}
\left\langle \left( \underbrace{\fr{\del f}{\del t} +{v} \cdot \fr{\del f}{\del {x}} = J(f,f)}_{\rm Boltzmann \ equation} \right) \right\rangle 
\ee 
where the moments are defined by 
\be 
\left \langle  f \right \rangle = \int_{-\infty}^{\infty} f \ dv 
\ee 
Here, $f$ is the molecular velocity distribution function, $v$ is the molecular velocity and $J(f,f)$ is the collision term. A simpler model for the collision term, $J(f,f)$, in the Boltzmann equation is given by the popular B-G-K model \cite{Bhatnagar_Gross_Krook}, which reduces the otherwise integro-differential equation to a PDE with a simple 
(algebraic) source term.     
\be 
J(f,f) = - \fr{1}{\epsilon} \left[f - f^{eq} \right] 
\ee 
Here, $\epsilon$ is the relaxation time and $f^{eq}$ is the distribution representing local thermodynamic equilibrium, given by a simpe Gaussian as 
\be 
f^{eq} = u \fr{1}{\sqrt{2 \pi}} e^{- \frac{1}{2} \left( v - \tilde{u} \right)^{2} } 
\ee  
Utilizing the B-G-K model and operator splitting, the solution of the Boltzmann equation can be split into two steps as 
\be   
\textrm{\bf Convection Step:} \ \fr{\del f}{\del t} + \fr{\del {h}}{\del {x}} = 0 
\ee  
\be   
\textrm{\bf Collision Step:} \ \fr{d f}{d t} = - \fr{1}{\epsilon} \left[f - f^{eq} \right]  
\ee 
where $ {h} =  {v}  f$ is the flux, with the Boltzmann equation being written in conservation form.   
Choosing an {\em instantaneous relaxation to equilibrium} ($\epsilon=0$), the collision step becomes a simple relaxation step as $f = f^{eq}$.  Thus, a conservation equation can be written in an intriguing form as 
\be 
\left \langle \left( \underbrace{\fr{\del f}{\del t} +  \fr{\del \vec{h}}{\del \vec{x}}=0}_{\rm convection}, \ \underbrace{f = f^{eq}}_{\rm relaxation} \right) \right \rangle
\ee 
The essential advantage of using this representation for macroscopic Burgers equation is the linearity of the convection terms in the Boltzmann equation.  This makes introducing upwinding easier in a kinetic or Boltzmann scheme, easily leading to an upwind scheme for a conservation equation.  While this strategy is explained for Burgers equation here, it is possible to extend the strategy to any conservation equation or a system of conservation laws.

\subsection{Strategy for flexible velocity kinetic schemes} 
Here, we utilize the above formulation and introduce new versions of the equilibrium based on flexible velocities. These discrete distributions are then utilized to introduce novel kinetic schemes for solving various conservation equations.  Compared to traditional kinetic schemes, the expressions here are much simpler, as integrations are replaced by summations and Maxwellians are replaced by simpler algebraic expressions.  It is worth noting here that the numerical schemes developed in this work begin with the conservation form of the Boltzmann equation.  Further, the finite volume framework and the flux splitting form of numerical diffusion ensure the preservation of the conservation form.  These frameworks ensure conservation even if we make the velocities functions of conserved variables. Further, the essential numerical diffusion is fixed based on Chapman Engskog type expansion (KFDS) or {\em by enforcing Rankine-Hugoniot (R-H) conditions} (KFDS+). The former approach (KFDS) ensures boundedness to the maximum eigenvalue and ensures conservation, while the latter (KFDS+) utilises R-H conditions which represent the quintessential conservation of fluxes across discontinuities.  In the next section, the new flexible velocity Boltzmann equation is introduced.

\section{Flexible  Velocity Boltzmann Scheme for 1-D convection equation}  
\subsection{1-D convection equation in moment form}   
     Consider the 1-D nonlinear convection equation as 
\be \label{1D_scalar_PDE} 
\fr{\del u}{\del t} + \fr{\del g(u)}{\del x} = 0 
\ee 
where $u$ is the dependent variable and $g(u)$ its nonlinear flux.  If we consider the inviscid Burgers equation, $g(u) = \fr{1}{2}u^{2}$.  The framework which will be constructed in the following sections will be valid for any non-linear flux $g(u)$.  We construct here a kinetic framework to derive this nonlinear hyperbolic convection equation as a moment of a Boltzmann equation given by 
\be \label{1D_BE_for_scalar}
\fr{\del f}{\del t} +  \fr{\del (vf)}{\del x} = - \fr{1}{\epsilon} \left[ f - f^{eq} \right] 
\ee
Here, $f$ is the molecular velocity distribution function, $v$ is the molecular velocity, $\epsilon$ is a small  relaxation time and $f^{eq}$ is the equilibrium distribution function, which is typically the Maxwellian.  The right hand side is the popular BGK relaxation model \cite{Bhatnagar_Gross_Krook} in which the distribution function relaxes to its equilibrium within a small relaxation time, $\epsilon$.   
The moments to obtain the macroscopic variables are defined by 
\be \label{1D_moments_scalar} 
u = \int_{-\infty}^{\infty} f dv \ \textrm{and} \ g(u) = \int_{-\infty}^{\infty} v f d v 
\ee
so that the taking a moment of the Boltzmann equation (\ref{1D_BE_for_scalar}) as 
\be 
\int_{-\infty}^{\infty} \left[ \fr{\del f}{\del t} +  \fr{\del (vf)}{\del x} = 
- \fr{1}{\epsilon} \left\{ f - f^{eq} \right\}  \right] dv  
\ee  
yields, using the moment relations (\ref{1D_moments_scalar}), the original scalar convection equation (\ref{1D_scalar_PDE}).  The equilibrium distribution function, $f^{eq}$ is typically the Maxwellian in the traditional kinetic theory of gases.  To construct $f^{eq}$, let us start with the Gaussian with unit variance as 
\be 
f^{eq} = u \fr{1}{\sqrt{2 \pi}} e^{- \frac{1}{2} \left( v - \tilde{u} \right)^{2} } 
\ee
Its moment is given by 
$$ \ba{rcl}  
\displaystyle \int_{-\infty}^{\infty} f^{eq} dv & = & \displaystyle \int_{-\infty}^{\infty} 
u \fr{1}{\sqrt{2 \pi}} e^{- \frac{1}{2} \left( v - \tilde{u} \right)^{2} } d v \\ 
 & = & u \fr{1}{\sqrt{2 \pi}} \int_{-\infty}^{\infty} e^{- \frac{1}{2} \left( v - \tilde{u} \right)^{2} } d v 
\ea $$ 
Let us define a new variable $w$ as 
\be 
w = \fr{1}{\sqrt{2}} \left( v - \tilde{u} \right)  
\ee
The limits change as follows. 
\bea 
\ba{c} 
\textrm{As} \ v\rightarrow -\infty, \ w \rightarrow -\infty \\ 
\textrm{As} \ v\rightarrow \infty, \ w \rightarrow \infty \\ 
\textrm{As} \ v\rightarrow 0, \ w \rightarrow - \fr{1}{\sqrt{2}} \tilde{u} \\  
\ea 
\eea 
Note also that $dv = dw \sqrt{2}$.  Therefore, we obtain 
$$  
\displaystyle \int_{-\infty}^{\infty} f^{eq} dv =  u \fr{1}{\sqrt{2 \pi}} \int_{-\infty}^{\infty} e^{-w^{2}} dw \sqrt{2} = u \fr{1}{\sqrt{\pi}} \int_{-\infty}^{\infty} e^{-w^{2}} dw = u \fr{1}{\sqrt{\pi}} \sqrt{\pi} = u $$ 
thus recovering the first of the moment conditions (\ref{1D_moments_scalar}). Enforcing the second of the moment conditions (\ref{1D_moments_scalar}) leads to 
$$ \ba{rcl} 
g(u) & = & \displaystyle \int_{-\infty}^{\infty} v f^{eq} dv \\ 
 & = & \displaystyle \int_{-\infty}^{\infty} v u \fr{1}{\sqrt{2 \pi}} 
e^{- \frac{1}{2} \left( v - \tilde{u} \right)^{2} } d v \\ 
& = & u \fr{1}{\sqrt{2 \pi}} \int_{-\infty}^{\infty} \left( w \sqrt{2} + \tilde{u} \right) e^{-w^{2}} dw \sqrt{2} \\
& = & u \fr{1}{\sqrt{\pi}} \left[ 0 + \tilde{u} \sqrt{\pi} \right] \\ 
& = & u \tilde{u} 
\ea $$ 
Therefore, we obtain 
\be 
\tilde{u} = \fr{1}{u} g(u) 
\ee 
For the inviscid Burgers equation, we get 
\be
\tilde{u} = \fr{1}{u} g(u) = \fr{1}{u} \fr{1}{2} u^{2} = \fr{1}{2} u 
\ee
\subsection{Dirac delta based equilibrium distribution and flexible velocities} 
     Let us replace the equilibrium distribution by a combination of Dirac delta functions as 
\be
f^{eq} = f^{eq}_{+} \delta (v - \lambda_{+}) + f^{eq}_{-} \delta (v - \lambda_{-})  
\ee
Let us further assume, for simplicity, that the flexible velocities, $\lambda_{+}$ and $\lambda_{-}$ 
are given by 
\be \label{simple_DVs_1D_scalar}
\lambda_{+} = \lambda \ \textrm{and} \  \lambda_{-} = - \lambda 
\ee
Thus, we have three unknowns, namely, $f^{eq}_{+}$, $f^{eq}_{-}$ and $\lambda$ to be fixed for the equilibrium distribution.  We therefore first use the two moment relations 
\be 
u = \int_{-\infty}^{\infty} f^{eq} dv \ \textrm{and} \ g(u) = \int_{-\infty}^{\infty} vf^{eq} dv 
\ee
together with either the Chapman-Enskog type expansion or the Rankine-Hugoniot jump condition (or both, with each being used in a different part of the domain) as the third relation for fixing the unknowns. The R-H condition is given by   
\be \label{RHC_1D_scalar}
\Delta g(u) = s \Delta u 
\ee
Here $s$ is the shock speed and we utilize this condition in constructing the cell-interface flux in a finite volume method, with the hope of designing a scheme which can capture discontinuities like shock waves accurately.  

    Applying the moment conditions (\ref{1D_moments_scalar}), we obtain the following.  
$$ \ba{rcl}  
u & = & \displaystyle \int_{-\infty}^{\infty} f^{eq} dv \\ 
 & = & \displaystyle \int_{-\infty}^{\infty} \left\{ f^{eq}_{+} \delta (v - \lambda_{+}) 
+ f^{eq}_{-} \delta (v - \lambda_{-}) \right\} dv \\ 
 & = & f^{eq}_{+} (1) + f^{eq}_{-} (1)   
\ea $$ 
or 
\be 
f^{eq}_{+} + f^{eq}_{-} = u 
\ee 
$$ \ba{rcl} 
g(u) & = & \displaystyle \int_{-\infty}^{\infty} v f^{eq} dv \\ [3mm]
 & = & \displaystyle \int_{-\infty}^{\infty} v 
\left\{ f^{eq}_{+} \delta (v - \lambda_{+}) + f^{eq}_{-} \delta (v - \lambda_{-}) \right\} dv \\ [2mm]
 & = &  f^{eq}_{+} \displaystyle \int_{-\infty}^{\infty} \phi(v) \delta (v - \lambda_{+}) dv + 
        f^{eq}_{-} \displaystyle \int_{-\infty}^{\infty} \phi(v) \delta (v - \lambda_{-}) dv  
				\ \textrm{where} \ \phi(v) = v \\ 
 & = & f^{eq}_{+} \lambda_{+} + f^{eq}_{-} \lambda_{-} 
\ea $$ 
or 
\be 
f^{eq}_{+} \lambda_{+} + f^{eq}_{-} \lambda_{-} = g(u) 
\ee 
Solving the above two equations, we obtain 
\be 
f^{eq}_{+} = \fr{- \lambda_{-} u + g(u) }{\lambda_{+} - \lambda_{-} } \ \textrm{and} \ 
f^{eq}_{-} = \fr{ \lambda_{+} u - g(u) }{\lambda_{+} - \lambda_{-} }  
\ee
With our simplification of flexible velocities (\ref{simple_DVs_1D_scalar}), these components of the equilibrium distribution become 
\be 
f^{eq}_{+} = \fr{1}{2} u + \fr{1}{2 \lambda} g(u) \ \textrm{and} \ f^{eq}_{-} = \fr{1}{2} u - \fr{1}{2 \lambda} g(u) 
\ee 
with $\lambda$ to be fixed while applying the finite volume method, using Chapman-Enskog expansion or the Rankine-Hugoniot condition 
(\ref{RHC_1D_scalar}).  

   We can also rewrite the expressions derived above as 
\be  
u = \sum_{i=1}^{2} f^{eq}_{i} \ \textrm{and} \ g(u) = \sum_{i=1}^{2} \lambda_{i} f^{eq}_{i}
\ee
where 
\be \label{swap_variable}
f^{eq}_{1} = f^{eq}_{+}, \ \lambda_{1} = \lambda_{+}, \ f^{eq}_{2} = f^{eq}_{-} \ \textrm{and} \ 
\lambda_{2} = \lambda_{-}, 
\ee 
We now define 
\be \label{1D_scalar_moments_in_summation_form}
u = \sum_{i=1}^{2} f_{i} = \sum_{i=1}^{2} f^{eq}_{i} \ \textrm{and} \ g(u) = \sum_{i=1}^{2} \lambda_{i} f^{eq}_{i} 
\ee 
so that the first moment is true for any distribution function.  Comparing these moments with the original moments 
(\ref{1D_moments_scalar}) clearly shows that the continuous molecular velocity $v$ is replaced by two flexible velocities ($\pm \lambda$) and the molecular velocity distribution function $f$ is replaced by two components: $f_{1}$ and $f_{2}$.  
We thus obtain the framework of {\em flexible velocity Boltzmann equation} as 
\be 
\fr{\del f_{i}}{\del t} 
+  \fr{\del (\lambda_{i}f_{i})}{\del x} = - \fr{1}{\epsilon} 
\left[ f_{i} - f^{eq}_{i} \right], \ i=1,2 
\ee 
with the moment relations (\ref{1D_scalar_moments_in_summation_form}) to be used for recovering the nonlinear convection equation 
\be 
\fr{\del u}{\del t} + \fr{\del g(u)}{\del x} = 0 
\ee 
The essential advantage of this formulation is the linearity of the convection terms in the flexible velocity Boltzmann equation, which makes the job of developing a numerical method for the nonlinear convection equation simpler.   
 
\subsection{Flexible Velocity Boltzmann Equation and Moments} 

The {\em Flexible Velocity Boltzmann Equation} (FVBE) derived before can be written in matrix form as    
\be \label{DVBE_1D_scalar} 
\fr{\del \bf f}{\del t} +  \fr{\del \bf h}{\del x} = - \fr{1}{\epsilon} \left[ \bf f - \bf f^{eq} \right] 
\ee
where 
\be \label{disc_defns}
\bf f = \left[ \ba{c} f_{+} \\ f_{-} \ea \right], \  \
\bf h = \Lambda\bf f ,\
\Lambda = \left[ \ba{cc} \lambda_{+} & 0 \\ 0 & \lambda_{-} \ea \right] \ \textrm{and} \ 
\bf f^{eq} = \left[ \ba{c} f^{eq}_{+} \\ f^{eq}_{-} \ea \right] 
= \left[ \ba{c} \fr{1}{2} u + \fr{1}{2 \lambda} g(u) \\[4mm]  
\fr{1}{2} u - \fr{1}{2 \lambda} g(u) \ea \right]  
\ee
To obtain the moments (which are now based on summations instead of integrals), we introduce  
\be  \label{P_moment_scalar}
{\bf P} = [ 1 \ 1] 
\ee
Therefore, the moments (\ref{1D_scalar_moments_in_summation_form}) are now defined by 
\be 
u = {\bf P} {\bf f^{eq}} \ \textrm{and} \ g(u) = {\bf P} \Lambda {\bf f^{eq}} 
\ee
Thus, the moment of the FVBE can be written as 
\be \label{1D_scalar_eqn_as_DVBE_moment} 
{\bf P} \left[ \fr{\del \bf f}{\del t} +  \fr{\del \Lambda\bf f}{\del x} 
= - \fr{1}{\epsilon} \left[ \bf f - \bf f^{eq} \right] 
 \right] 
\ee
which is equivalent, in the limit $\epsilon \rightarrow 0$, to the original conservation law 
\be 
\fr{\del u}{\del t} + \fr{\del (g(u))}{\del x} = 0 
\ee
as 
\be 
{\bf P} {\bf f} = {\bf P} {\bf f^{eq}} = u 
\ee 
A further simplification can be introduced by using operator splitting, by which the solution of the flexible velocity Boltzmann equation can be split into two steps: a {\em convection step} and a {\em relaxation step}.  Thus 
$$ \fr{\del \bf f}{\del t} +  \fr{\del \left(\Lambda \bf f\right)}{\del x} 
= - \fr{1}{\epsilon} \left[ \bf f - \bf f^{eq} \right] $$ 
can be written as 
$$ \fr{\del \bf f}{\del t} = -  \fr{\del \left(\Lambda \bf f\right)}{\del x} 
- \fr{1}{\epsilon} \left[ \bf f - \bf f^{eq} \right]  $$ 
or 
$$ \fr{\del \bf f}{\del t} = O_{1}({\bf f}) + O_{2}({\bf f}) $$ 
where $O_{1}({\bf f})$ and $O_{2}({\bf f})$ are mathematical operators of ${\bf f}$.   
The above equation can be split, using operator splitting, as 
$$ \fr{\del \bf f}{\del t} = O_{1}({\bf f}) \ \textrm{and} \ 
\fr{\del \bf f}{\del t} = O_{2}({\bf f}) $$ 
or 
$$ \fr{\del \bf f}{\del t} +\fr{\del \left(\Lambda \bf f\right)}{\del x} = 0 \ \textrm{and} \ 
\fr{d \bf f}{d t} = - \fr{1}{\epsilon} \left[ \bf f - \bf f^{eq} \right] $$ 
leading to a {\em convection step} and a {\em relaxation step}.  The relaxation step can be further simplified into an {\em instantaneous relaxation step} by assuming that $\epsilon = 0$. This leads to instantaneous relaxation of the distribution function to its equilibrium, ${\bf f}$ = ${\bf f^{eq}}$.  Thus, the solution of the flexible velocity Boltzmann equation can be written as 
\be 
\fr{\del {\bf f}}{\del t} +  \fr{\del \left(\Lambda \bf f\right)}{\del x} = 0, \ 
{\bf f} = {\bf f^{eq}} 
\ee 
Now taking moments, we obtain 
\be 
{\bf P} \left[ 
\fr{\del \bf f}{\del t} +  \fr{\del \left(\Lambda \bf f\right)}{\del x} = 0, \ 
{\bf f} = {\bf f^{eq}} 
\right] 
\ee 
which is the moment form of the original nonlinear convection equation given by 
\be 
\fr{\del u}{\del t} + \fr{\del g(u)}{\del x} = 0 
\ee 
Note that our strategy of deriving a flexible velocity Boltzmann equation based on approximating the Maxwellian with a combination of Dirac delta functions has led to recovering the diagonal relaxation model of Aregba-Driollet \& Natalini \cite{Driollet_Natalini} in 1-D, which is also equivalent to the relaxation system of Jin \& Xin \cite{Jin_Xin}.  However, our discretization strategy leads to a totally different numerical scheme, which is presented in the following.   
\subsection{Kinetic Flux Splitting Based on Flexible Velocities} 

Now consider the flexible velocity Boltzmann formulation given by (\ref{1D_scalar_eqn_as_DVBE_moment}) in a finite volume framework. 
\begin{figure}[!h]  
\begin{center} 
\includegraphics[height=3.2cm]{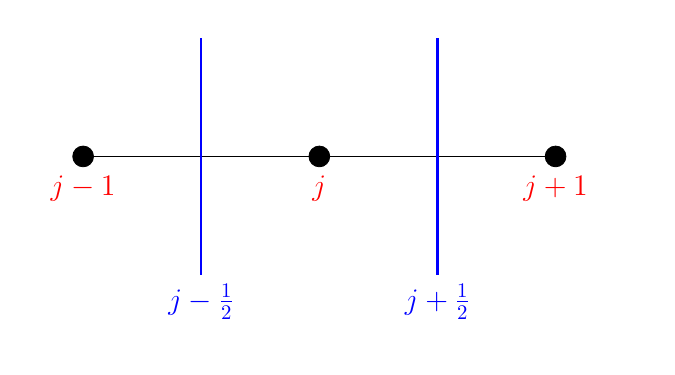} 
\caption{3-point stencil for a FVM} 
\label{3_point_stencil}
\end{center} 
\end{figure} 
The update formula can be written as 
\be \label{K_update_FVM} 
 {\bf f}_{j}^{n+1} = {\bf f}_{j}^{n} - \fr{\Delta t}{\Delta x} \left[ {\bf h}_{j+\frac{1}{2}}^{n} - {\bf h}_{j-\frac{1}{2}}^{n}  \right] 
\ee
where $ {\bf h} = {\bf \Lambda} {\bf f}$.  
To introduce upwinding in two velocity model based FVBE, let us split the flexible velocity matrix ${\bf \Lambda}$ in to two parts, separating the positive and negative velocities.  
\bea \label{DV_splitting_1D_2VM} 
{ \Lambda} = \left[ \ba{cc} \lambda^{+} & 0 \\ 0 & \lambda^{-} \ea \right]
= \left[ \ba{cc} \lambda^{+} & 0 \\ 0 & 0 \ea \right] + \left[ \ba{cc} 0 & 0 \\ 0 & \lambda^{-} \ea \right]
= { \Lambda}^{+} + { \Lambda}^{-}
\eea
It is possible to write
\be \label{mod_lambda} 
 {|\Lambda|}={\Lambda}^{+} - { \Lambda}^{-}
\ee
Thus, the upwind fluxes  applied on a three point stencil can be written as 
\bea 
 {\bf h}_{j+\frac{1}{2}}^{n}= [{ \Lambda}^{+} {\bf f}_{eq}]_{j} 
+  [ {\Lambda}^{-} {\bf f}_{eq}]_{j+1}\\
{\bf h}_{j-\frac{1}{2}}^{n}= [ { \Lambda}^{+} {\bf f}_{eq}]_{j-1} 
+  [ { \Lambda}^{-} {\bf f}_{eq}]_{j}
\eea
Using (\ref {DV_splitting_1D_2VM}) and (\ref{mod_lambda}) we can write 
 
\bea 
 {\bf h}_{j+\frac{1}{2}}^{n}= \underbrace{\fr{1}{2}[ {\bf h}_{j+1}^{n}+ {\bf h}_{j}^{n}]}_{\rm average \ flux} 
- \underbrace{\fr{1}{2}|\Lambda|[{\bf f}^{eq}_{j+1}-[{\bf f}^{eq}_{j}]}_{\rm diffusive \ flux} \\ 
{\bf h}_{j-\frac{1}{2}}^{n}= \fr{1}{2}[ {\bf h}_{j}^{n}+ {\bf h}_{j-1}^{n}] 
- \fr{1}{2}|\Lambda|[{\bf f}^{eq}_{j}-[{\bf f}^{eq}_{j-1}]  
\eea 

Clearly, $\Lambda$ represents the coefficients of numerical diffusion.  Choosing the values of $\lambda$ for numerical considerations is an efficient strategy to control the numerical diffusion without losing conservation, as the above expressions for cell-interface fluxes in the flux difference splitting form enforce conservation in the basic finite volume framework.  
The above expressions can further be rewritten in flux difference splitting form as 

\bea \label{2vel_KineticUpdate} 
 {\bf h}_{j+\frac{1}{2}}^{n}= \fr{1}{2}[ {\bf h}_{j+1}^{n}+ {\bf h}_{j}^{n}] 
- \fr{1}{2}[\Delta {\bf h}^{+}_{j+\frac{1}{2}}-\Delta {\bf h}^{-}_{j+\frac{1}{2}}]  \\ 
 {\bf h}_{j-\frac{1}{2}}^{n}= \fr{1}{2}[ {\bf h}_{j}^{n}+ {\bf h}_{j-1}^{n}] 
- \fr{1}{2}[\Delta {\bf h}^{+}_{j-\frac{1}{2}}-\Delta {\bf h}^{-}_{j-\frac{1}{2}}] 
\eea
where $ \Delta {\bf h}_{j\pm\frac{1}{2}}^{\pm} = [{\bf \Lambda}^{\pm}\Delta {\bf f}_{eq}]_{j\pm\frac{1}{2}}$.  In the next sub-section, the coefficient of numerical diffusion, which corresponds to $|\bf \Lambda|$, is chosen as a function of both the left and right states.  Thus, the flux difference splitting is the appropriate choice rather than flux vector splitting.  

To recover the macroscopic update formula, let is take moments by multiplying with $P$ as given in (\ref{P_moment_scalar}).  Therefore the macroscopic update formula for the Kinetic Flux Difference Splitting (KFDS) scheme thus developed using 2-velocity model based FVBE can be written as 
\be 
 u_{j}^{n+1} = u_{j}^{n} - \fr{\Delta t}{\Delta x}\left[ g(u)_{j+\frac{1}{2}}^{n} - g(u)_{j-\frac{1}{2}}^{n}  \right] 
\ee
where the interface fluxes are given by
\bea 
 g(u)_{j+\frac{1}{2}}^{n}= \fr{1}{2}[  g(u)_{j+1}^{n}+ g(u)_{j}^{n}]- \fr{1}{2}[\Delta g(u)^{+,n}_{j+\frac{1}{2}}-\Delta g(u)^{-,n}_{j+\frac{1}{2}}]\\ 
g(u)_{j-\frac{1}{2}}^{n}= \fr{1}{2}[  g(u)_{j}^{n}+  g(u)_{j-1}^{n}]- \fr{1}{2}[\Delta  g(u)^{+,n}_{j-\frac{1}{2}}-\Delta  g(u)^{-,n}_{j-\frac{1}{2}}] 
\eea
\bea \label{2Vel_KFDS} 
\Delta  g(u)_{j+\frac{1}{2}}^{\pm} = \fr{1}{2}[   g(u)_{j+1}-  g(u)_{j}]\pm \fr{1}{2}|\lambda|[u_{j+1}-u_{j}] \\ 
\Delta  g(u))_{j-\frac{1}{2}}^{\pm} = \fr{1}{2}[   g(u)_{j}-  g(u)_{j-1}]\pm \fr{1}{2}|\lambda|[u_{j}-u_{j-1}] 
\eea  

Note that all the variables of the flexible velocity Boltzmann equation have disappeared in the final scheme, in the spirit of the kinetic or Boltzmann schemes.  All the variables in the final scheme are purely macroscopic variables, as $\lambda$, fixing of which is explained in the next sections based on numerical considerations, also takes the values of macroscopic variables.  The kinetic framework is used only the conceptual level, exploiting the linearity of the convection terms in the flexible velocity Boltzmann equation, leading to an upwind scheme for the nonlinear convection equation.  
\subsection{Fixing \texorpdfstring{$\lambda$}{Lambda}} 
     Two different ways of fixing $\lambda$ are proposed here.  The first way is based on Chapman-Enskog type asymptotic expansion of the FVBE, which leads to a diffusive scheme.  The second way is based on utilizing the Rankine-Hugoniot jump condition, which leads to an accurate shock capturing scheme.  By substituting the expressions for the split fluxes into the finite volume update formula, one can easily see that $\lambda$ also represents the coefficient of numerical diffusion and therefore the different choices for $\lambda$ can be justified to modify the numerical diffusion in the scheme.       
\subsubsection{Chapman-Enskog type expansion of the flexible velocity Boltzmann equation} 
One simple way of fixing $\lambda$ is to use the Chapman-Enskog type expansion of the flexible velocity Boltzmann equation.  This expansion, which is a special form of asymptotic expansion for the flexible velocity Boltzmann equation, is given in the Appendix C.  The Chapman-Enskog type expansion shows that the flexible velocity Boltzmann equation leads to the original convection equation augmented by a diffusion term, given by 
\be
\fr{\del u}{\del t} + \fr{\del g(u)}{\del x} = \epsilon \fr{\del}{\del x} \left[ 
\fr{\del u}{\del x} \left\{ \lambda^{2} - \left( a\left(u\right) \right)^{2} \right\}
\right] + \mathcal{O}\left( \epsilon^{2} \right)
\ee
where $a(u)$ is the wave speed for the original convection equation, given by $a(u)=\fr{\del g(u)}{\del u}$.  The right hand side in the above equation contains a second derivative of the conserved variable $u$ and is multiplied by a small parameter ($\epsilon \rightarrow 0$).  Thus, the moment of the flexible velocity Boltzmann equation represents a vanishing diffusion approximation to the original convection equation.  For the approximation to be stable, the diffusion term must be non-negative.  Therefore, we obtain 
\be
\lambda^{2} - \left( a\left(u\right)\right)^{2} \ge 0 
\ee 
which can approximated to 
\be 
\lambda = \max \left[ a(u)_{j} \right], \ \forall j  
\ee 
where $j$ is the index representing the grid points.  The scheme derived using this approach is named as {\em KFDS} method, as it leads to a kinetic flux difference splitting scheme, based on flexible velocities.     

\subsubsection{Fixing \texorpdfstring{$\lambda$}{Lambda} based on R-H (jump) condition} \label{RHC}
\begin{figure}[!h]  
\begin{center} 
\includegraphics[height=3.2cm]{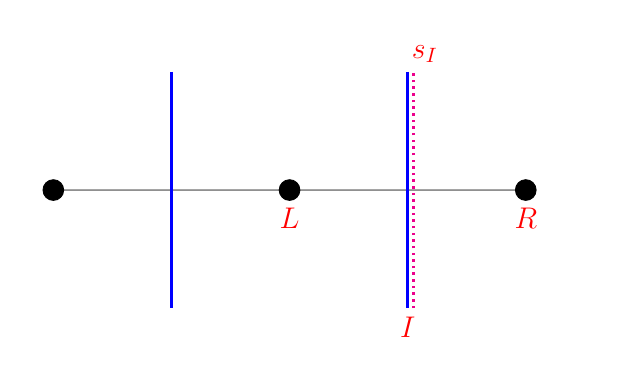} 
\includegraphics[height=3.2cm]{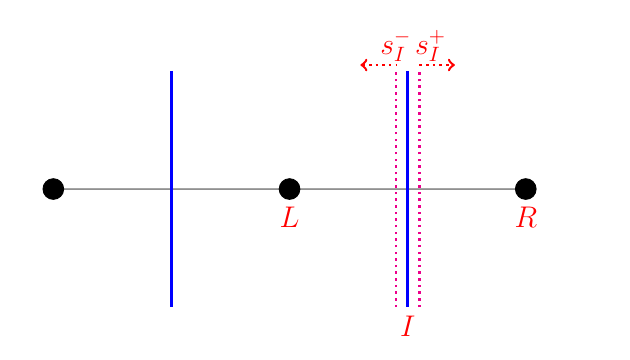} \\ 
\caption{Shock at a cell-interface in FVM} 
\label{shock_at_interface} 
\end{center} 
\end{figure}

Let us fix $\lambda$ based on enforcing Rankine-Hugoniot jump condition (R-H condition) at the cell-interface.  Let us assume that a shock discontinuity is present at the cell-interface $I$, moving with shock speed $s$, as shown in the figure (\ref{shock_at_interface}).  

The Rankine-Hugoniot jump condition connecting the left state $L$ and the right state $R$ is given by 
\be 
\Delta g(u) = s_{I} \Delta u 
\ee
where 
\be 
\Delta u = u_{R} - u_{L} \ \textrm{and} \ 
\Delta g(u)  = g_{R}(u)  - g_{L}(u)  
\ee 
Splitting the shock speed into a positive part (moving to the right) and a negative part (moving to the left), we can write 
(see figure(\ref{shock_at_interface}))  
\be 
s_{I} = s^{+}_{I} + s^{-}_{I} = \fr{s_{I} + |s_{I}|}{2} + \fr{s_{I} - |s_{I}|}{2}  
\ee 
Using the above shock speed splitting, we now split the R-H condition into two parts, one to the right of the cell-interface and the other to the left of the cell-interface (see figure \ref{3_point_stencil} ).  Obviously, out of the three variables ($u$, $g$, $s$) in the R-H equation, as we are modifying $u$ and $g$, the third part, $\Delta u$ needs to be kept unmodified.  Therefore, we obtain 
\be 
g_{R}(u) - g_{I}(u) = s^{+}_{I} \Delta u \ \textrm{and} \ 
g_{I}(u) - g_{L}(u) = s^{-}_{I} \Delta u 
\ee
Subtracting one from the other, we obtain 
$$ g_{R}(u) - g_{I}(u) - g_{I}(u) + g_{L}(u) = \left( s^{+}_{I} - s^{-}_{I}\right) \Delta u $$ 
or 
$$ g_{I}(u) = \fr{g_{L}(u) + g_{R}(u)}{2} 
- \fr{1}{2} \left( \fr{s_{I} + |s_{I}|}{2} - \fr{s_{I} - |s_{I}|}{2} \right) \Delta u $$ 
or 
\be \label{Interface_flux_RHC_1D_convection} 
g_{I}(u) = \fr{1}{2} \left[ g_{L}(u) + g_{R}(u) \right] - \fr{1}{2} |s_{I}| \Delta u 
\ee
This is the cell-interface flux derived from the R-H condition.  The cell-interface flux derived for the flexible velocity based upwind kinetic scheme is 
$$ \ba{rcl} g(u)_{I} 
& = & g^{+}(u)_{L} + g^{-}(u)_{R} \\ [3mm]
& = & \fr{1}{2} \left[  g(u)_{L} +  g(u)_{R} \right] - \fr{1}{2} \lambda \left[ u_{R} - u_{L} \right]  
\ea $$   
Therefore 
\be \label{Interface_flux_1D_KFVSDV} 
g(u)_{I} = \fr{1}{2} \left[  g(u)_{L} +  g(u)_{R} \right] - \fr{1}{2} \lambda \Delta u 
\ee
Comparing the cell-interface flux derived from the flexible velocity based upwind kinetic scheme 
(\ref{Interface_flux_1D_KFVSDV}) with that derived based on the R-H condition 
(\ref{Interface_flux_RHC_1D_convection}), we obtain 
\be 
\lambda = |s_{I}| 
\ee
For the nonlinear convection equation, as $g(u) = \fr{1}{2} u^{2}$, we obtain  
$$ s = \fr{\Delta g(u)}{\Delta u} = \fr{g(u)_{R} - g(u)_{L}}{u_{R} - u_{L}} 
= \fr{\fr{1}{2}u_{R}^{2} - \fr{1}{2} u_{L}^{2}}{u_{R} - u_{L}} 
= \fr{1}{2} \left( u_{L} + u_{R} \right) $$ 
Therefore 
\be 
\lambda = \left|\fr{1}{2} \left(u_{L} + u_{R} \right)\right|
\ee
The method derived thus is named as {\em KFDS+} scheme, which is an accurate shock capturing variation of {\em KFDS} scheme.  

\subsection {Stability Analysis}  
Consider a 1D Burgers equation system as given in [\ref{1D_scalar_PDE}] whose equivalent flexible velocity kinetic system of equations is given in [\ref{DVBE_1D_scalar}].  The finite volume version of the KFDS scheme for the above system takes the form
\be \label{start_eqn}
{\bf f}^{n+1}_{j} = {\bf f}^{n}_{j} - \fr{\Delta t}{\Delta x} \left[ {\bf h}_{j+\frac{1}{2}} - {\bf h}_{j-\frac{1}{2}} \right] ^{n}
\ee
The interface fluxes $\bf {{h}_{j+\frac{1}{2}}^{n}}$ and $\bf {{h}_{j-\frac{1}{2}}^{n}}$ are defined by 
\bea
{\bf h}_{j+\frac{1}{2}}^{n} = \left [ \bf {\Lambda}^{+} {f}^{eq} \right]^{n}_{j} + \left [ \bf {\Lambda}^{-} {f}^{eq } \right]^{n}_{j+1} \\ [2mm]
{\bf h}_{j-\frac{1}{2}}^{n} = \left [ \bf {\Lambda}^{+} {f}^{eq} \right]^{n}_{j-1} + \left [ \bf {\Lambda}^{-} {f}^{eq } \right]^{n}_{j}
\eea
Expanding the distribution by substituting [\ref{disc_defns}] and rewriting, we get 
\bea
\bf {{h}_{j+\frac{1}{2}}^{n}} = {\left[ \ba{c} {\lambda}_{+}f^{eq}_{+} \\ 0 \ea \right]}_{j}^{n}  \ + \  {\left[ \ba{c} 0 \\ {\lambda}_{-}f^{eq}_{-} \ea \right]}_{j+1}^{n}\\[4mm]
\bf {{h}_{j-\frac{1}{2}}^{n}} = {\left[ \ba{c} {\lambda}_{+}f^{eq}_{+} \\ 0 \ea \right]}_{j-1}^{n} + {\left[ \ba{c} 0 \\ {\lambda}_{-}f^{eq}_{-} \ea \right]}_{j}^{n}.
\eea 
Substituting this in the equation [\ref{start_eqn}], we get
\be
{\left[ \ba{c} f_{+} \\[3mm] f_{-} \ea \right]}^{n+1}_{j} =  {\left[ \ba{c} f_{+} \\[3mm] f_{-} \ea \right]}^{n}_{j}-\fr{\Delta t}{\Delta x}\left [ \ba{c} {\lambda}_{+}{f^{eq}_{+}}_{j} - {\lambda}_{+}{f^{eq}_{+}}_{j-1}  \\[4mm] {\lambda}_{-}{f^{eq}_{-}}_{j+1} - {\lambda}_{-}{f^{eq}_{-}}_{j} \ea  \right]
\ee
Substituting for ${\lambda}_{+}$ and ${\lambda}_{-}$, we obtain  
\be
{\left[ \ba{c} f_{+} \\[3mm] f_{-} \ea \right]}^{n+1}_{j} =  {\left[ \ba{c} f_{+} \\[3mm] f_{-} \ea \right]}^{n}_{j}-\fr{\Delta t}{\Delta x}\left [ \ba{c} |{\lambda}|{f^{eq}_{+}}_{j} - |{\lambda}|{f^{eq}_{+}}_{j-1}  \\[4mm] -|{\lambda}|{f^{eq}_{-}}_{j+1} + |{\lambda}|{f^{eq}_{-}}_{j} \ea  \right]
\ee 
Using expressions in [\ref{swap_variable}], the above equation can be written in simpler form as
\be \label{2DVBE_Simple}
{f}_{\alpha , j}^{n+1} = {f}_{\alpha , j}^{n} - \fr{|{\lambda}| \Delta t}{\Delta x} \left[ {f}^{eq}_{\alpha , j} - {f}^{eq}_{\alpha , (j+\sgn(\frac{3}{2}-\alpha))} \right]^{n}, \ \alpha = 1,2  
\ee 
with the order of $\alpha$ being changed.  Here the standard signum function is defined by 
\bea 
\sgn (\omega) = \left\{ \ba{c} - 1 \ \textrm{for} \ \omega < 1 \\ 0 \ \textrm{for} \ \omega = 0 \\ 1 \ \textrm{for} \ \omega > 1 \ea \right. 
\eea 
In order to evaluate the stability characteristics of the above system, we perform a von Neumann stability analysis. We introduce the following Fourier transforms. 
\bea
  {\mathcal{F}} ({f}_{j}^{n}) =  {\hat f}^{n}_{j}{e}^{i k j \Delta x} \\
  {\mathcal{F}} ({{f}^{eq , n}_{j}}) = {{\hat f}_{j}^{eq , n}} {e}^{i k j \Delta x}
\eea
Also from Chapman-Engkog expansion we have
\be
f = {f}^{eq} + {\epsilon {f}^{\prime}} + \cdots
\ee
It is interesting to note in the above equation that 
\be
 { {f}^{\prime}} \geq 0, \ \ \epsilon \geq 0, \ \  \cdots.
\ee
which implies  $ f \geq {f}^{eq} $ or $ \fr{{f}^{eq}}{f} \leq 1 $. 
Applying Fourier transform to [\ref{2DVBE_Simple}] and substituting the above relations we get
\be
{\hat f}_{\alpha,j}^{n+1}{e}^{i k j \Delta x} = {\hat f}_{\alpha,j}^{n}{e}^{i k j \Delta x} -\fr{|{\lambda}| \Delta t}{\Delta x} \left[ {\hat f}_{\alpha, j}^{eq , n}-{\hat f}_{\alpha}^{eq , n}{e}^{ik\left( \sgn \left ( \frac{3}{2}- \alpha \right ) \right) \Delta x}  \right]{e}^{i k j \Delta x} .
\ee
Rearranging the terms and substituting $C =\fr{|{\lambda}| \Delta t}{\Delta x}$, we obtain 
\be
\fr{{\hat f}_{\alpha,j}^{n+1}}{{\hat f}_{\alpha,j}^{n}} = 1 -  C \fr{{\hat f}_{\alpha, j}^{eq , n}}{{\hat f}_{\alpha,j}^{n}} \left[ 1 - {e}^{ik\left( \sgn \left ( \frac{3}{2}- \alpha \right ) \right) \Delta x}  \right]
\ee  
The exponential term can be further expanded using Euler's formula ${e}^{i \theta} = \cos \theta + i \sin \theta$. we get
\be
\fr{{\hat f}_{\alpha,j}^{n+1}}{{\hat f}_{\alpha,j}^{n}} = 1 -  C \fr{{\hat f}_{\alpha, j}^{eq , n}}{{\hat f}_{\alpha,j}^{n}}\left[ 1 - \left( \cos\left(k \Delta x\right) + \sgn \left ( \frac{3}{2}- \alpha \right ) i  \sin\left(k \Delta x\right) \right) \right]
\ee
Rearranging,  we obtain  
\be
\fr{{\hat f}_{\alpha,j}^{n+1}}{{\hat f}_{\alpha,j}^{n}} = 1 +  C \fr{{\hat f}_{\alpha, j}^{eq , n}}{{\hat f}_{\alpha,j}^{n}}\left[  \left( \cos\left(k \Delta x\right) + \sgn \left ( \frac{3}{2}- \alpha \right ) i  \sin\left(k \Delta x\right) \right) -1 \right].
\ee
or
\be
\fr{{\hat f}_{\alpha,j}^{n+1}}{{\hat f}_{\alpha,j}^{n}} \leq 1 +  C \left[  \left( \cos\left(k \Delta x\right) + \sgn \left ( \frac{3}{2}- \alpha \right ) i  \sin\left(k \Delta x\right) \right) -1 \right].
\ee 
Here, we evaluate the RHS part by multiplying it with its complex conjugate ${RHS}^{*}$.  
\bea
RHS = \left[ \left( 1 + C\cos \left( k \Delta x \right) - C\right) + i \sgn \left ( \frac{3}{2}- \alpha \right )C \sin \left( k \Delta x \right) \right] \\
{RHS}^{*} = \left[ \left( 1 + C\cos \left( k \Delta x \right) - C\right) + i \sgn \left ( \frac{3}{2}- \alpha \right )C \sin \left( k \Delta x \right) \right]
\eea
\be
{|RHS|}^{2} = RHS. {RHS}^{*}
\ee
\be
 {|RHS|}^{2}=  1 + 2C\left( \cos \left( k \Delta x \right) -1\right) + {C}^{2}\left( {\cos}^{2} \left( k \Delta x \right) +1 - 2\cos \left( k \Delta x \right) \right) + {C}^{2} {\sin}^{2} \left( k \Delta x \right)
\ee
which simplifies to  
\be
 {|RHS|}^{2}= 1 -2C\left( 1-C\right) \left(1-\cos \left( k \Delta x \right)  \right)
\ee
For the numerical scheme to be stable, we need to find the conditions at which $ {|RHS|}^{2} \leq1$. Therefore the above expression can be written as
\be
2C\left( 1-C\right) \left(1-\cos \left( k \Delta x \right)  \right) \leq 0.
\ee  
We know that $\left[1-\cos \left( k \Delta x \right)  \right] \geq 0 \ \forall \  k,\Delta x \ \in  \mathbb{R}$. Therefore for the scheme to be stable, we get 
$C\left( C-1\right) \leq 0$, which upon solving yields $ 0\leq \ C \ \leq 1$. Thus the condition for the scheme to be stable is 
\be
0 \leq \fr{|\lambda| \Delta t}{\Delta x} \leq 1
\ee


\section{Second order accuracy} \label{KLW_Start}
    In this section, the Lax-Wendroff method is first introduced in the flexible velocity Boltzmann framework. Following that, a high resolution TVD scheme in which KFDS method (correspondingly also KFDS+ method) switches over to the {\em Kinetic Lax-Wendroff Method} in the smooth regions, with the help of a limiter function, is introduced.  
\subsection{Kinetic Lax-Wendroff (KLW) Method}  
     While the second order accuracy can be obtained by several possibilities, including upwinding on extended stencils, the Lax-Wendroff approach is attractive because of its operation in a compact stencil of three points.  Here, the Lax-Wendroff method in the framework of flexible velocity Boltzmann equation is introduced and then moments are taken, a process which yields a {\em generalized Lax-Wendroff method}.    

     The flexible velocity Boltzmann framework we start with, to obtain the convection equation, is given by  
\be 
{\bf P} \left[ \fr{\del {\bf f}}{\del t} +  \fr{\del \left(\Lambda \bf f\right)}{\del x} = 0, 
{\bf f} = {\bf f}^{eq} \right] \Longrightarrow \fr{\del u}{\del t} + \fr{\del g(u)}{\del x} = 0  
\ee
Using Taylor series expansion, we can write 
$$ {\bf f}(x, t+\Delta t) = {\bf f}(x,t) 
+ \fr{\Delta t}{1!} \fr{\del {\bf f}(x,t)}{\del t} 
+ \fr{\Delta t^{2}}{2!} \fr{\del^{2} {\bf f}(x,t)}{\del t^{2}} 
+ \mathcal{O}(\Delta t^{3})  $$ 
Using the Boltzmann equation in the convection step, we can write 
$$ \fr{\del {\bf f}}{\del t} = - \fr{\del \left(\Lambda \bf f\right)}{\del x} $$ 
and further 
$$ \fr{\del^{2} {\bf f}}{\del t^{2}} 
= \Lambda^{2} \fr{\del^{2} {\bf f}}{\del x^{2}} $$  
Thus, we get 
$$ {\bf f}^{n+1}_{j} = {\bf f}^{n}_{j} 
- \Delta t \Lambda \fr{\del {\bf f}^{n}}{\del x}|_{j} 
+ \fr{\Delta t^{2}}{2} \Lambda^{2} \fr{\del {\bf f}^{n}}{\del x^{2}}|_{j} 
+ \mathcal{O}(\Delta t^{3})$$  
Using central differences, we obtain 
$$ {\bf f}^{n+1}_{j} = {\bf f}^{n}_{j} 
- \Delta t \Lambda \left[\fr{{\bf f}^{n}_{j+1} - {\bf f}^{n}_{j-1}}{2 \Delta x}\right] 
+ \fr{\Delta t^{2}}{2} \Lambda^{2} \left[ 
\fr{{\bf f}^{n}_{j+1} - 2 {\bf f}^{n}_{j} + {\bf f}^{n}_{j-1} }{\Delta x^{2}} 
\right] + \mathcal{O}(\Delta t^{3}, \Delta x^{2} \Delta t) $$ 
or 
$$ {\bf f}^{n+1}_{j} = {\bf f}^{n}_{j} 
- \fr{\Delta t}{2 \Delta x} \Lambda 
\left[ {\bf f}^{n}_{j+1} + {\bf f}^{n}_{j} - {\bf f}^{n}_{j} - {\bf f}^{n}_{j-1} \right]  + \fr{\Delta t^{2}}{2 \Delta x^{2}} \Lambda^{2} \left[ 
{\bf f}^{n}_{j+1} - 2 {\bf f}^{n}_{j} + {\bf f}^{n}_{j-1} \right] 
+ \mathcal{O}(\Delta t^{3}, \Delta x^{3}) $$   
assuming $\Delta t$ and $\Delta x$ to be of the same order.  With a little algebraic manipulation, we obtain 
$$ {\bf f}^{n+1}_{j} = {\bf f}^{n}_{j} 
- \fr{\Delta t}{\Delta x} 
\left[ {\bf h}^{n}_{j+\frac{1}{2}} - {\bf h}^{n}_{j-\frac{1}{2}}\right] + 
\mathcal{O}(\Delta t^{2}, \Delta x^{2}) $$ 
where 
$$ {\bf h}^{n}_{j+\frac{1}{2}} 
= \fr{1}{2} \left[ \Lambda {\bf f}_{j} + \Lambda {\bf f}_{j+1} \right] 
- \fr{1}{2} \Lambda^{2} \fr{\Delta t}{\Delta x} 
\left[ {\bf f}_{j+1} - {\bf f}_{j} \right]  $$  
and 
$$ {\bf h}^{n}_{j-\frac{1}{2}} 
= \fr{1}{2} \left[ \Lambda {\bf f}_{j-1} + \Lambda {\bf f}_{j} \right] 
- \fr{1}{2} \Lambda^{2} \fr{\Delta t}{\Delta x} 
\left[ {\bf f}_{j} - {\bf f}_{j-1} \right]  $$ 
If we now take moments, after using the instantaneous relaxation step as 
${\bf f}^{n}_{j} = {\bf f}^{eq,n}_{j}$  we obtain the finite volume update formula as 
\be 
{\bf P} \left[ {\bf f}^{n+1}_{j} = {\bf f}^{n}_{j} 
- \fr{\Delta t}{\Delta x} 
\left[ {\bf h}^{n}_{j+\frac{1}{2}} - {\bf h}^{n}_{j-\frac{1}{2}} \right] \right] 
\ee 
or 
\be \label{KLW_update_formula}
u^{n+1}_{j} = u^{n}_{j} 
- \fr{\Delta t}{\Delta x} 
\left[ g(u)^{n}_{j+\frac{1}{2}} - g(u)^{n}_{j-\frac{1}{2}} \right] 
\ee
where 
\be \label{KLW_interface_flux}
g(u)_{j+\frac{1}{2}} = \fr{1}{2} \left[ g(u)_{j} + g(u)_{j+1} \right] 
- \fr{1}{2} \lambda^{2} \fr{\Delta t}{\Delta x} \left[u_{j+1} - u_{j} \right] 
\ee 
Here we have used the expressions 
\be 
{\bf P} {\bf f} = u, \ {\bf P} \Lambda {\bf f^{eq}} = g(u) \ \textrm{and} \ 
{\bf P} \Lambda^{2} {\bf f^{eq}} = \lambda^{2} u 
\ee 
as 
$$ \ba{rcl} {\bf P} \Lambda^{2} {\bf f}^{eq} & = & 
\left[ 1 \ 1 \right] \left[ 
\ba{cc} \left( \lambda \right)^{2} & 0 \\ 0 & \left( - \lambda \right)^{2} \ea \right] \left[ \ba{c} \fr{1}{2} u + \fr{1}{2 \lambda} g(u) \\[5mm] 
\fr{1}{2} u - \fr{1}{2 \lambda} g(u) \ea \right] \\[5mm]  
& = &  \lambda^{2} u 
\ea $$ 
It is worth noting that if we take $\lambda = a(u)_{j+\frac{1}{2}}$, where 
$a(u)_{j+\frac{1}{2}}$ is the wave speed of the original convection equation at the cell interface, we obtain the classical Lax-Wendroff method.  Since other choices are also possible for $\lambda$, the scheme we obtained, which is the {\em Kinetic Lax-Wendroff (KLW) method} (\ref{KLW_update_formula},\ref{KLW_interface_flux}), can be thought of a {\em generalized Lax-Wendroff method}.  As a first attempt, we use the Chapman-Enskog type expansion to fix $\lambda$ in the KLW method.  
\subsection{High resolution TVD KFDS Method}  
    Since KLW method, being second order accurate, produces oscillations near large gradients like shock waves, a Total Variation Diminishing (TVD) Kinetic Lax-Wendroff method, as a high resolution KFDS method, is introduced in this section to obtain a non-oscillatory high resolution flexible velocity Boltzmann scheme.  The KLW method is written as DVKS method augmented by additional anti-diffusive terms and these anti-diffusive terms are switched off near large gradients using a limiter function, to obtain a non-oscillatory scheme.   

    Consider the finite volume update formula for the flexible velocity Boltzmann equation, given by 
\be \label{FVM_update_formulat_for_DVBE}
{\bf f}^{n+1}_{j} = {\bf f}^{n}_{j} - \fr{\Delta t}{\Delta x} \left[ 
{\bf h}^{n}_{j+\frac{1}{2}} - {\bf h}^{n}_{j+\frac{1}{2}} \right] 
\ee
with the cell-interface flux for deriving KFDS method being given by 
\be 
{\bf h}^{KFDS}_{j+\frac{1}{2}} 
= \fr{1}{2} \left[ \Lambda {\bf f}^{eq}_{j} + \Lambda {\bf f}^{eq}_{j+1} \right] 
- \fr{1}{2} |\Lambda| \left[ {\bf f}^{eq}_{j+1} - {\bf f}^{eq}_{j} \right] 
\ee 
Similarly, for deriving the KLW method, cell-inteface flux can be written as 
\be 
{\bf h}^{KLW}_{j+\frac{1}{2}} 
= \fr{1}{2} \left[ \Lambda {\bf f}^{eq}_{j} + \Lambda {\bf f}^{eq}_{j+1} \right] 
 - \fr{1}{2} \Lambda^{2} \fr{\Delta t}{\Delta x} 
\left[ {\bf f}^{eq}_{j+1} - {\bf f}^{eq}_{j} \right] 
\ee
For obtaining a TVD flux, we write the cell-interface flux as 
\be 
{\bf h}^{TVD-KFDS}_{j+\frac{1}{2}} = {\bf h}^{KFDS}_{j+\frac{1}{2}} 
+ \phi \left[ {\bf h}^{KLW}_{j+\frac{1}{2}} - {\bf h}^{KFDS}_{j+\frac{1}{2}} \right] 
\ee 
where $\phi = \phi(r)$ is a limiter function, to be chosen appropriately.  Introducing 
\be 
\Lambda^{\pm} = \fr{\Lambda \pm |\Lambda|}{2} 
\ee
we obtain 
\be 
\Lambda^{+} = \left[ \ba{cc} \lambda & 0 \\ 0 & 0 \ea \right], \  
\Lambda^{-} = \left[ \ba{cc} 0 & 0 \\ 0 & - \lambda \ea \right] \ \textrm{and} \ 
|\Lambda| = \left[ \ba{cc} \lambda & 0 \\ 0 & \lambda \ea \right] 
= \Lambda^{+} - \Lambda^{-} 
\ee    
We now write 
\bea \ba{rcl} 
{\bf h}^{TVD-KFDS}_{j+\frac{1}{2}} = {\bf h}^{KFDS}_{j+\frac{1}{2}}  - \phi({\bf r}^{+}_{j}) \left[ \fr{1}{2} \left( \fr{1}{2} \Lambda^{2} \fr{\Delta t}{\Delta x} - \Lambda^{+} \right) \left( {\bf f}^{eq}_{j+1} - {\bf f}^{eq}_{j} \right)
\right] \\[5mm]  - \phi({\bf r}^{-}_{j+1}) \left[ \fr{1}{2} \left( \fr{1}{2} \Lambda^{2} \fr{\Delta t}{\Delta x} + \Lambda^{-} \right) 
\left( {\bf f}^{eq}_{j+1} - {\bf f}^{eq}_{j} \right)
\right] \ea
\eea 
where 
\be 
{\bf r}^{+}_{j} = \fr{ {\bf f}^{eq}_{j} - {\bf f}^{eq}_{j-1} } 
{{\bf f}^{eq}_{j+1} - {\bf f}^{eq}_{j}} \ \textrm{and} \ 
{\bf r}^{-}_{j+1} = \fr{ {\bf f}^{eq}_{j+2} - {\bf f}^{eq}_{j+1} } 
{{\bf f}^{eq}_{j+1} - {\bf f}^{eq}_{j}} 
\ee 
following \cite{Thomas_2}.  
We choose the limiter function as the {\em minmod} function and define 
\be 
\phi ({\bf r}) = minmod (1,{\bf r}) 
\ee 
where 
\be 
minmod(a,b) = \left\{ \ba{l} 
a \ \textrm{if} \ |a| < |b| \ \& \ ab > 0 \\ 
b \ \textrm{if} \ |b| < |a| \ \& \ ab > 0 \\ 
0 \ \textrm{if} \ ab \le 0 
\ea \right. 
\ee
Therefore 
$$ \ba{rcl} 
{\bf h}^{TVD-KFDS}_{j+\frac{1}{2}} 
& = &{\bf h}^{KFDS}_{j+\frac{1}{2}} \\[5mm]  
& & - minmod \left(1, \fr{ {\bf f}^{eq}_{j} - {\bf f}^{eq}_{j-1} } 
{{\bf f}^{eq}_{j+1} - {\bf f}^{eq}_{j}} \right) \left[ 
\fr{1}{2} \left( \fr{1}{2} \Lambda^{2} \fr{\Delta t}{\Delta x} - \Lambda^{+} \right) 
\left( {\bf f}^{eq}_{j+1} - {\bf f}^{eq}_{j} \right)
\right] \\[5mm] 
& & - minmod \left(1, \fr{ {\bf f}^{eq}_{j+2} - {\bf f}^{eq}_{j+1} } 
{{\bf f}^{eq}_{j+1} - {\bf f}^{eq}_{j}}  \right) \left[ 
\fr{1}{2} \left( \fr{1}{2} \Lambda^{2} \fr{\Delta t}{\Delta x} + \Lambda^{-} \right) 
\left( {\bf f}^{eq}_{j+1} - {\bf f}^{eq}_{j} \right)
\right] \ea $$ 
Using the relation 
\be 
c \ minmod (a,b) = minmod (ca, cb) 
\ee
we obtain 
\bea \ba{rcl} 
{\bf h}^{TVD-KFDS}_{j+\frac{1}{2}} 
& = & {\bf h}^{KFDS}_{j+\frac{1}{2}}  \\[5mm]  
& & - \fr{1}{2} minmod \left[  \left( \fr{1}{2} \Lambda^{2} \fr{\Delta t}{\Delta x} 
- \Lambda^{+} \right) 
\left( {\bf f}^{eq}_{j+1} - {\bf f}^{eq}_{j} \right), \ \left( \fr{1}{2} \Lambda^{2} \fr{\Delta t}{\Delta x} - \Lambda^{+} \right) 
\left( {\bf f}^{eq}_{j} - {\bf f}^{eq}_{j-1} \right)
\right]  \\[5mm] 
& & - \fr{1}{2} minmod \left[ \left( \fr{1}{2} \Lambda^{2} \fr{\Delta t}{\Delta x} 
+ \Lambda^{-} \right)  
\left( {\bf f}^{eq}_{j+1} - {\bf f}^{eq}_{j} \right), \  \left( \fr{1}{2} \Lambda^{2} \fr{\Delta t}{\Delta x} + \Lambda^{-} \right)  
\left( {\bf f}^{eq}_{j+2} - {\bf f}^{eq}_{j+1} \right)
\right] \ea 
\eea  
Now taking moment of the above cell-interface flux and using the relations 
\be 
{\bf P} \Lambda {\bf f} = g(u), \ {\bf P} \Lambda^{2} {\bf f} = \lambda^{2} u, \ 
{\bf P} \Lambda^{\pm} {\bf f} = \fr{1}{2} g(u) \pm \fr{1}{2} \lambda u  
\ee  
we obtain 
\bea \ba{ll} 
{g(u)_{j+\frac{1}{2}}}^{TVD-KFDS} & = {g(u)_{j+\frac{1}{2}}}^{KFDS} \\[2mm]  
& - \fr{1}{2} minmod \left[ \ba{c}  
\left\{ \fr{1}{2} \fr{\Delta t}{\Delta x} \lambda^{2} \left( u_{j+1} - u_{j} \right) 
- \fr{1}{2} \lambda \left( u_{j+1} - u_{j} \right) 
- \fr{1}{2} \left( g(u)_{j+1} - g(u)_{j} \right) 
\right\}, \\[2mm]  
\left\{ \fr{1}{2} \fr{\Delta t}{\Delta x} \lambda^{2} \left( u_{j} - u_{j-1} \right) 
- \fr{1}{2} \lambda \left( u_{j} - u_{j-1} \right) 
- \fr{1}{2} \left( g(u)_{j} - g(u)_{j-1} \right) 
\right\} \ea 
\right] \\[5mm] 
& - \fr{1}{2} minmod \left[ \ba{c}  
\left\{ \fr{1}{2} \fr{\Delta t}{\Delta x} \lambda^{2} \left( u_{j+1} - u_{j} \right) 
- \fr{1}{2} \lambda \left( u_{j+1} - u_{j} \right) 
+ \fr{1}{2} \left( g(u)_{j+1} - g(u)_{j} \right) 
\right\}, \\[2mm]  
\left\{ \fr{1}{2} \fr{\Delta t}{\Delta x} \lambda^{2} \left( u_{j+2} - u_{j+1} \right) 
- \fr{1}{2} \lambda \left( u_{j+2} - u_{j+1} \right) 
+ \fr{1}{2} \left( g(u)_{j+2} - g(u)_{j+1} \right) 
\right\} \ea 
\right] \ea   
\eea 
  
In a similar way, starting with 
\bea \ba{rcl} 
{\bf h}^{TVD-KFDS}_{j-\frac{1}{2}} 
& = & {\bf h}^{KFDS}_{j-\frac{1}{2}}  - \phi({\bf r}^{+}_{j-1}) \left[ \fr{1}{2} \left( \fr{1}{2} \Lambda^{2} \fr{\Delta t}{\Delta x} - \Lambda^{+} \right) \left( {\bf f}^{eq}_{j} - {\bf f}^{eq}_{j-1} \right) \right] \\[5mm] 
& & - \phi({\bf r}^{-}_{j}) \left[ 
\fr{1}{2} \left( \fr{1}{2} \Lambda^{2} \fr{\Delta t}{\Delta x} + \Lambda^{-} \right) 
\left( {\bf f}^{eq}_{j} - {\bf f}^{eq}_{j-1} \right)
\right] \ea
\eea 
where 
\be 
{\bf r}^{-}_{j} = \fr{ {\bf f}^{eq}_{j+1} - {\bf f}^{eq}_{j} } 
{{\bf f}^{eq}_{j} - {\bf f}^{eq}_{j-1}} \ \textrm{and} \ 
{\bf r}^{+}_{j-1} = \fr{ {\bf f}^{eq}_{j-1} - {\bf f}^{eq}_{j-2} } 
{{\bf f}^{eq}_{j} - {\bf f}^{eq}_{j-1}} 
\ee
we obtain 
\bea \ba{ll} 
{g(u)_{j-\frac{1}{2}}}^{TVD-KFDS} & = {g(u)_{j-\frac{1}{2}}}^{KFDS} \\[2mm]  
& - \fr{1}{2} minmod \left[ \ba{c}  
\left\{ \fr{1}{2} \fr{\Delta t}{\Delta x} \lambda^{2} \left( u_{j} - u_{j-1} \right) 
- \fr{1}{2} \lambda \left( u_{j} - u_{j-1} \right) 
- \fr{1}{2} \left( g(u)_{j} - g(u)_{j-1} \right) 
\right\}, \\[2mm]  
\left\{ \fr{1}{2} \fr{\Delta t}{\Delta x} \lambda^{2} \left( u_{j-1} - u_{j-2} \right) 
- \fr{1}{2} \lambda \left( u_{j-1} - u_{j-2} \right) 
- \fr{1}{2} \left( g(u)_{j-1} - g(u)_{j-2} \right) 
\right\} \ea 
\right] \\[5mm] 
& - \fr{1}{2} minmod \left[ \ba{c}  
\left\{ \fr{1}{2} \fr{\Delta t}{\Delta x} \lambda^{2} \left( u_{j} - u_{j-1} \right) 
- \fr{1}{2} \lambda \left( u_{j} - u_{j-1} \right) 
+ \fr{1}{2} \left( g(u)_{j} - g(u)_{j-1} \right) 
\right\}, \\[2mm]  
\left\{ \fr{1}{2} \fr{\Delta t}{\Delta x} \lambda^{2} \left( u_{j+1} - u_{j} \right) 
- \fr{1}{2} \lambda \left( u_{j+1} - u_{j} \right) 
+ \fr{1}{2} \left( g(u)_{j+1} - g(u)_{j} \right) 
\right\} \ea 
\right] \ea   
\eea 
These expressions for $g(u)_{j+\frac{1}{2}}$ and $g(u)_{j-\frac{1}{2}}$ are used in the update formula given by 
\be 
u^{n+1}_{j} = u^{n}_{j} - \fr{\Delta t}{\Delta x} \left[ 
g(u)_{j+\frac{1}{2}} - g(u)_{j-\frac{1}{2}} \right] 
\ee
which is the moment of the update formula (\ref{FVM_update_formulat_for_DVBE}) for the flexible velocity Boltzmann equation. An advantage of the above formulation is that the contributions from the second order terms (Lax-Wendroff part) and the first order terms 
(KFDS part) can be seen separately in the expressions for the cell-interface fluxes and thus different values of $\lambda$ can be chosen in different regions, leading to a scheme in which numerical diffusion is chosen smartly, appropriate to smooth and discontinuous regions. 

\section{KFDS Method for 1-D Convection-Diffusion Equation} 
The flexible velocity Boltzmann equation (FVBE), with instantaneous relaxation to Chapman-Enskog distribution function (obtained after operator splitting), is given by
\bea \label{DVBE} 
\fr{\del {\bf f}}{\del t} + \Lambda \fr{\del {\bf f}}{\del x} 
= 0, \ {\bf f} ={\bf f}^{CE}
\eea 
where the Chapman-Enskog distribution is given by (see Appendix C) 
\be \label{fCE}
{\bf f}^{CE} = {\bf f}^{eq} - {\bf f}^{eq}_{v} \ \textrm{with} \ 
{\bf f}^{eq}_{v} = \fr{1}{2} \nu \fr{\del u}{\del x} 
\left[ \ba{c} \fr{1}{\lambda} \\[4mm] - \fr{1}{\lambda} \ea \right] 
\ee
which upon taking moments yield the nonlinear convection-diffusion equation, which is the viscous Burgers equation (see Appendix C).  
\be
\frac { \partial u } { \partial t } + \frac { \partial g(u) } { \partial x } 
= \frac { \partial  g_{v}(u) } { \partial x }
\ee 

\subsection {1-D KFDS method for convection-diffusion equation}
Now consider the flexible Boltzmann formulation of the convection-diffusion equation given by
\be
{\bf P} \left[\fr{\del {\bf f}}{\del t} + \Lambda \fr{\del {\bf f}}{\del x} 
= 0, \ {\bf f} ={\bf f}^{CE}\right]
\ee
where
\be 
{\bf f}^{CE} = {\bf f}^{eq} - {\bf f}_{v}^{eq}
\ee 
\be 
{\bf f}^{eq} = \left[ \ba{c} \fr{1}{2} u + \fr{1}{2\lambda} g(u) \\[4mm] 
\fr{1}{2} u - \fr{1}{2\lambda} g(u) \ea \right] \ \textrm{and} \ 
{\bf f}_{v}^{eq} = \fr{1}{2} \nu \fr{\del u}{\del x} 
\left[ \ba{c} \fr{1}{\lambda} \\[4mm] - \fr{1}{\lambda} \ea \right] 
\ee 
Applying finite volume method for the above flexible velocity Boltzmann equation (with instantaneous relaxation to discrete Chapman-Enskog distribution) based on upwinding of the flexible velocities and then taking moments, we obtain 
\be 
\fr{\left({\bf P f}\right)^{n+1}_{j} - \left({\bf P f}\right)^{n}_{j}} {\Delta t} 
+ \fr{1}{\Delta x} \left[ 
\left( {\bf P h}^{n}_{j+\frac{1}{2}} \right) - \left( {\bf P h}^{n}_{j-\frac{1}{2}} \right) \right] = 0, {\bf f}^{n}_{j} = {\bf f^{CE}}_{j}^{n} 
\ee
where 
\be 
{\bf h} = \Lambda {\bf f}  
\ee
and the interface fluxes based on upwinding are given by 
\be 
{\bf h}_{j+\frac{1}{2}} = {\bf h}^{+}_{j} + {\bf h}^{-}_{j+1} \ \textrm{and} \ 
{\bf h}_{j-\frac{1}{2}} = {\bf h}^{+}_{j-1} + {\bf h}^{-}_{j}
\ee 
We therefore obtain, for the convection-diffusion equation at the macroscopic level, the finite volume update formula as 
\be 
\fr{u^{n+1}_{j} - u^{n}_{j}} {\Delta t} 
+ \fr{1}{\Delta x} \left[g^{n}_{T,j+\frac{1}{2}} - g^{n}_{T,j-\frac{1}{2}} \right] = 0  
\ee 
where $g_{T}$ is the total flux (convective flux plus the diffusive flux).  
The cell-interface fluxes are defined by 
\be 
g_{T,j+\frac{1}{2}} = g^{+}_{T,j} + g^{-}_{T,j+1} \ \textrm{and} \ 
g_{T,j-\frac{1}{2}} = g^{+}_{T,j-1} + g^{-}_{T,j}  
\ee
The split fluxes are given by (see Appendix C for details)  
\be 
g^{\pm}_{T} = {\bf P} \Lambda^{\pm} {\bf f^{CE}} = {\bf P} \Lambda^{\pm} {\bf f}^{eq} 
- {\bf P} \Lambda^{\pm} {\bf f}^{eq}_{v} = g^{\pm}(u) - \fr{1}{2} g_{v}(u) 
\ee

where 
\be 
g^{\pm} = \fr{1}{2} g(u) \pm \fr{1}{2} \lambda u 
\ee  
The final update formula for the finite volume method based flexible velocity kinetic flux splitting scheme for viscous Burgers equation, after some rearrangement, is given by 
\be 
u^{n+1}_{j} = u^{n}_{j} - \fr{\Delta t}{\Delta x} 
\left[ g^{n}(u)_{j+\frac{1}{2}} - g^{n}(u)_{j-\frac{1}{2}} \right]  + \fr{\Delta t}{\Delta x} 
\left[ g_{v}^{n}(u)_{j+\frac{1}{2}} - g_{v}^{n}(u)_{j-\frac{1}{2}} \right] 
\ee 
where the interface fluxes are given by
where the interface fluxes are given by
\bea 
 g(u)_{j+\frac{1}{2}}^{n}= \fr{1}{2}[  g(u)_{j+1}^{n}+ g(u)_{j}^{n}]- \fr{1}{2}[\Delta g(u)^{+,n}_{j+\frac{1}{2}}-\Delta g(u)^{-,n}_{j+\frac{1}{2}}]\\ 
g(u)_{j-\frac{1}{2}}^{n}= \fr{1}{2}[  g(u)_{j}^{n}+  g(u)_{j-1}^{n}]- \fr{1}{2}[\Delta  g(u)^{+,n}_{j-\frac{1}{2}}-\Delta  g(u)^{-,n}_{j-\frac{1}{2}}] 
\eea
\bea
\Delta  g(u)_{j+\frac{1}{2}}^{\pm} = \fr{1}{2}[   g(u)_{j+1}-  g(u)_{j}]\pm \fr{1}{2}|\lambda|[u_{j+1}-u_{j}] \\ 
\Delta  g(u))_{j-\frac{1}{2}}^{\pm} = \fr{1}{2}[   g(u)_{j}-  g(u)_{j-1}]\pm \fr{1}{2}|\lambda|[u_{j}-u_{j-1}] 
\eea  
and
\be 
g_{v}(u)_{j+\frac{1}{2}}^{n} = \frac{1}{2}\left(g_{v}(u)_{j} + g_{v}(u)_{j+1}\right)^{n} \ \& \ 
g_{v}(u)_{j-\frac{1}{2}}^{n} = \frac{1}{2}\left(g_{v}(u)_{j-1} + g_{v}(u)_{j}\right)^{n}
\ee
It is worth noting that in this strategy, a parabolic PDE (viscous Burgers equation) is treated completely based on discretrizing a hyperbolic PDE (the flexible velocity Boltzmann equation) and then taking moments.  However, due to the specific expressions for the moments (derived in Appendix C), the discretization of the viscous fluxes do not lead to any additional numerical diffusion and simple central discretization of them is obtained, with upwinding being retained only for the inviscid fluxes.  

\subsection {KLW with discrete Chapman-Enskog distribution}

   In this section we derive a second order Kinetic scheme for the viscous  Burgers equation. Here,  the Lax-Wendroff method is used in the framework of flexible velocity Boltzmann equation  and then moments are taken, to obtain a generalized Lax-Wendroff method for the viscous Burgers Equation.  The flexible velocity Boltzmann framework we start with, to obtain the convection equation, is given by  

\be 
{\bf P} \left[ \fr{\del {\bf f}}{\del t} + \Lambda \fr{\del {\bf f}}{\del x} 
= 0, \  
{\bf f} ={\bf f}^{CE} \right] \Longrightarrow \fr{\del u}{\del t} + \fr{\del g(u)}{\del x} = \frac { \partial g_{v}(u) } { \partial x }  
\ee
Using the similar derivation as for the KLW method, along with the understanding that ${\bf f}^{CE} = {\bf f}^{eq} - {\bf f}_{v}^{eq}$, with a little algebraic manipulation, we obtain 
$$ {\bf f}^{n+1}_{j} = {\bf f}^{n}_{j} 
- \fr{\Delta t}{\Delta x} 
\left[ {\bf h}^{n}_{j+\frac{1}{2}} - {\bf h}^{n}_{j-\frac{1}{2}}\right] + \fr{\Delta t}{\Delta x} 
\left[ {\bf h_{v}}^{n}_{j+\frac{1}{2}} - {\bf h_{v}}^{n}_{j-\frac{1}{2}}\right] + 
\mathcal{O}(\Delta t^{2}, \Delta x^{2}) $$ 
where 
$$ {\bf h}^{n}_{j+\frac{1}{2}} 
= \fr{1}{2} \left[ \Lambda {\bf f}^{eq}_{j} + \Lambda {\bf f}^{eq}_{j+1} \right] 
- \fr{1}{2} \Lambda^{2} \fr{\Delta t}{\Delta x} 
\left[ {\bf f}^{eq}_{j+1} - {\bf f}^{eq}_{j} \right]  $$  
and 
$$ {\bf h}^{n}_{j-\frac{1}{2}} 
= \fr{1}{2} \left[ \Lambda {\bf f}^{eq}_{j-1} + \Lambda {\bf f}^{eq}_{j} \right] 
- \fr{1}{2} \Lambda^{2} \fr{\Delta t}{\Delta x} 
\left[ {\bf f}^{eq}_{j} - {\bf f}^{eq}_{j-1} \right]  $$ 
and the viscous interface fluxes are given by
$$ {\bf h_{v}}^{n}_{j+\frac{1}{2}} 
= \fr{1}{2} \left[ \Lambda {\bf f_{v}}^{eq}_{j} + \Lambda {\bf f_{v}}^{eq}_{j+1} \right] 
- \fr{1}{2} \Lambda^{2} \fr{\Delta t}{\Delta x} 
\left[ {\bf f_{v}}^{eq}_{j+1} - {\bf f_{v}}^{eq}_{j} \right]  $$  
and 
$$ {\bf h_{v}}^{n}_{j-\frac{1}{2}} 
= \fr{1}{2} \left[ \Lambda {\bf f_{v}}^{eq}_{j-1} + \Lambda {\bf f_{v}}^{eq}_{j} \right] 
- \fr{1}{2} \Lambda^{2} \fr{\Delta t}{\Delta x} 
\left[ {\bf f_{v}}^{eq}_{j} - {\bf f_{v}}^{eq}_{j-1} \right]  $$
Upon taking moments we obtain the finite volume update formula as 
\be 
{\bf P} \left[ {\bf f}^{n+1}_{j} = {\bf f}^{n}_{j} 
- \fr{\Delta t}{\Delta x} 
\left[ {\bf h}^{n}_{j+\frac{1}{2}} - {\bf h}^{n}_{j-\frac{1}{2}} \right] + \fr{\Delta t}{\Delta x} 
\left[ {\bf h_{v}}^{n}_{j+\frac{1}{2}} - {\bf h_{v}}^{n}_{j-\frac{1}{2}}\right] \right]
\ee 
or 
\be \label{KLWVisc_update_formula}
u^{n+1}_{j} = u^{n}_{j} 
- \fr{\Delta t}{\Delta x} 
\left[ g(u)^{n}_{j+\frac{1}{2}} - g(u)^{n}_{j-\frac{1}{2}} \right] + \fr{\Delta t}{\Delta x} 
\left[ g_{v}(u)^{n}_{j+\frac{1}{2}} - g_{v}(u)^{n}_{j-\frac{1}{2}} \right] 
\ee
where 
\be \label{KLWVisc_interface_flux}
g(u)_{j+\frac{1}{2}} = \fr{1}{2} \left[ g(u)_{j} + g(u)_{j+1} \right] 
- \fr{1}{2} \lambda^{2} \fr{\Delta t}{\Delta x} \left[u_{j+1} - u_{j} \right] 
\ee 
\be
g_{v}(u)_{j+\frac{1}{2}} = \fr{1}{2} \left[ g_{v}(u)_{j} + g_{v}(u)_{j+1} \right] 
\ee
Here we have used the expressions 
\be 
{\bf P} {\bf f^{eq}} = u, \ {\bf P} \Lambda {\bf f^{eq}} = g(u) \  \ 
{\bf P} \Lambda^{2} {\bf f^{eq}} = \lambda^{2} u 
\ee 
\be 
{\bf P} {\bf f_{v}^{eq}} = 0, \ {\bf P} \Lambda {\bf f_{v}^{eq}} = g_{v}(u) \  \ 
{\bf P} \Lambda^{2} {\bf f_{v}^{eq}} = 0 
\ee 
It is evident in this case also that the viscous fluxes simply get central differenced, yielding the KLW method derived before for the inviscid fluxes.   
\subsection {TVD-KFDS method for 1-D convection-diffusion equation}
In this section we present a TVD version of the KFDS method for the viscous Burgers equation.  Consider the finite volume update formula for the flexible velocity Boltzmann equation, given by 
\be \label{FVM_update_formulat_for_DVBEvisc}
{\bf f}^{n+1}_{j} = {\bf f}^{n}_{j} - \fr{\Delta t}{\Delta x} \left[ 
{\bf h}^{n}_{j+\frac{1}{2}} - {\bf h}^{n}_{j+\frac{1}{2}} \right] + \fr{\Delta t}{\Delta x} \left[ 
{\bf h_{v}}^{n}_{j+\frac{1}{2}} - {\bf h_{v}}^{n}_{j+\frac{1}{2}} \right] 
\ee 
where 
\be 
{\bf h}^{n} = \Lambda {\bf f}^{eq,n} \ \textrm{and} \ {\bf h_{v}} = - \Lambda {\bf f_{v}}^{eq,n}    
\ee 
with the cell-interface flux for deriving KFDS method being given by 
\be 
{\bf h}^{KFDS}_{j+\frac{1}{2}} 
= \fr{1}{2} \left[ \Lambda {\bf f}^{eq}_{j} + \Lambda {\bf f}^{eq}_{j+1} \right] 
- \fr{1}{2} |\Lambda| \left[ {\bf f}^{eq}_{j+1} - {\bf f}^{eq}_{j} \right] 
\ee 
$${\bf h_{v}}^{KFDS}_{j+\frac{1}{2}} 
= \fr{1}{2} \left[ \Lambda {\bf f_{v}}^{eq}_{j} + \Lambda {\bf f_{v}}^{eq}_{j+1} \right]$$ 
Similarly, for deriving the KLW method, cell-inteface flux can be written as 
\be
{\bf h}^{KLW}_{j+\frac{1}{2}} 
= \fr{1}{2} \left[ \Lambda {\bf f}^{eq}_{j} + \Lambda {\bf f}^{eq}_{j+1} \right] - \fr{1}{2} \Lambda^{2} \fr{\Delta t}{\Delta x} 
\left[ {\bf f}^{eq}_{j+1} - {\bf f}^{eq}_{j} \right] 
\ee
$${\bf h_{v}}^{KLW}_{j+\frac{1}{2}} 
= \fr{1}{2} \left[ \Lambda {\bf f_{v}}^{eq}_{j} + \Lambda {\bf f_{v}}^{eq}_{j+1} \right]$$ 
It can be observed that the diffusion flux at the interface is simply the average diffusion flux at the respective interface for both KFDS and KLW. Therefore for obtaining a TVD flux it is sufficient to form a flux limited term for the convection terms. Thus the cell interface flux for the convection term will be identical to that for the inviscid Burgers equation. Hence we write the cell-interface flux as 
\be  \label{TVD_strat}
{\bf h}^{TVD-KFDS}_{j+\frac{1}{2}} = {\bf h}^{KFDS}_{j+\frac{1}{2}} 
+ \phi \left[ {\bf h}^{KLW}_{j+\frac{1}{2}} - {\bf h}^{KFDS}_{j+\frac{1}{2}} \right] 
\ee 
where $\phi = minmod (1,{\bf r})$ is a minmod limiter function. 
\be 
{\bf r}^{+}_{j} = \fr{ {\bf f}^{eq}_{j} - {\bf f}^{eq}_{j-1} } 
{{\bf f}^{eq}_{j+1} - {\bf f}^{eq}_{j}} \ \textrm{and} \ 
{\bf r}^{-}_{j+1} = \fr{ {\bf f}^{eq}_{j+2} - {\bf f}^{eq}_{j+1} } 
{{\bf f}^{eq}_{j+1} - {\bf f}^{eq}_{j}} 
\ee 
following the same steps as defined in the TVD-KFDS section we get the convective interface flux as  
\bea \ba{rcl} 
{\bf h}^{TVD-KFDS}_{j+\frac{1}{2}} 
& = & {\bf h}^{KFDS}_{j+\frac{1}{2}}  - \fr{1}{2} minmod \left[  \left( \fr{1}{2} \Lambda^{2} \fr{\Delta t}{\Delta x} 
- \Lambda^{+} \right) 
\left( {\bf f}^{eq}_{j+1} - {\bf f}^{eq}_{j} \right), \ \left( \fr{1}{2} \Lambda^{2} \fr{\Delta t}{\Delta x} - \Lambda^{+} \right) 
\left( {\bf f}^{eq}_{j} - {\bf f}^{eq}_{j-1} \right)
\right]  \\[5mm] 
& & - \fr{1}{2} minmod \left[ \left( \fr{1}{2} \Lambda^{2} \fr{\Delta t}{\Delta x} 
+ \Lambda^{-} \right)  
\left( {\bf f}^{eq}_{j+1} - {\bf f}^{eq}_{j} \right), \  \left( \fr{1}{2} \Lambda^{2} \fr{\Delta t}{\Delta x} + \Lambda^{-} \right)  
\left( {\bf f}^{eq}_{j+2} - {\bf f}^{eq}_{j+1} \right)
\right] \ea 
\eea  
Now taking moments, we obtain 
$$ \ba{rcl} 
{\bf P} {\bf h}^{TVD-KFDS}_{j+\frac{1}{2}} 
& = &  {\bf Ph}^{KFDS}_{j+\frac{1}{2}} - \fr{1}{2} minmod \left[ {\bf P}  \left( \fr{1}{2} \Lambda^{2} \fr{\Delta t}{\Delta x} 
- \Lambda^{+} \right) 
\left( {\bf f}^{eq}_{j+1} - {\bf f}^{eq}_{j} \right), \ {\bf P} \left( \fr{1}{2} \Lambda^{2} \fr{\Delta t}{\Delta x} - \Lambda^{+} \right) 
\left( {\bf f}^{eq}_{j} - {\bf f}^{eq}_{j-1} \right)
\right]  \\[5mm] 
& & - \fr{1}{2} minmod \left[ {\bf P} \left( \fr{1}{2} \Lambda^{2} \fr{\Delta t}{\Delta x} 
+ \Lambda^{-} \right)  
\left( {\bf f}^{eq}_{j+1} - {\bf f}^{eq}_{j} \right), \ {\bf P}  \left( \fr{1}{2}  \Lambda^{2} \fr{\Delta t}{\Delta x} + \Lambda^{-} \right)  
\left( {\bf f}^{eq}_{j+2} - {\bf f}^{eq}_{j+1} \right)
\right] \ea $$ 
Using the moments  
\be 
{\bf P} \Lambda {\bf f}^{eq} = g(u), \ {\bf P} \Lambda^{2} {\bf f}^{eq} = \lambda^{2} u, \ 
{\bf P} \Lambda^{\pm} {\bf f}^{eq} = \fr{1}{2} g(u) \pm \fr{1}{2} \lambda u  
\ee  
we obtain 
\bea \ba{ll} 
g(u)_{j+\frac{1}{2}}^{TVD-KFDS} & = g(u)_{j+\frac{1}{2}}^{KFDS} \\[2mm]   
& - \fr{1}{2} minmod \left[ \ba{c}  
\left\{ \fr{1}{2} \fr{\Delta t}{\Delta x} \lambda^{2} \left( u_{j+1} - u_{j} \right) 
- \fr{1}{2} \lambda \left( u_{j+1} - u_{j} \right) 
- \fr{1}{2} \left( g(u)_{j+1} - g(u)_{j} \right) 
\right\}, \\[2mm]  
\left\{ \fr{1}{2} \fr{\Delta t}{\Delta x} \lambda^{2} \left( u_{j} - u_{j-1} \right) 
- \fr{1}{2} \lambda \left( u_{j} - u_{j-1} \right) 
- \fr{1}{2} \left( g(u)_{j} - g(u)_{j-1} \right) 
\right\} \ea 
\right] \\[5mm] 
& - \fr{1}{2} minmod \left[ \ba{c}  
\left\{ \fr{1}{2} \fr{\Delta t}{\Delta x} \lambda^{2} \left( u_{j+1} - u_{j} \right) 
- \fr{1}{2} \lambda \left( u_{j+1} - u_{j} \right) 
+ \fr{1}{2} \left( g(u)_{j+1} - g(u)_{j} \right) 
\right\}, \\[2mm]  
\left\{ \fr{1}{2} \fr{\Delta t}{\Delta x} \lambda^{2} \left( u_{j+2} - u_{j+1} \right) 
- \fr{1}{2} \lambda \left( u_{j+2} - u_{j+1} \right) 
+ \fr{1}{2} \left( g(u)_{j+2} - g(u)_{j+1} \right) 
\right\} \ea 
\right] \ea   
\eea 
In a similar way, starting with 
\bea \ba{rcl} 
{\bf h}^{TVD-KFDS}_{j-\frac{1}{2}} 
& = & {\bf h}^{KFDS}_{j-\frac{1}{2}}  - \phi({\bf r}^{+}_{j-1}) \left[ 
\fr{1}{2} \left( \fr{1}{2} \Lambda^{2} \fr{\Delta t}{\Delta x} - \Lambda^{+} \right) 
\left( {\bf f}^{eq}_{j} - {\bf f}^{eq}_{j-1} \right)
\right] \\[5mm] 
& & - \phi({\bf r}^{-}_{j}) \left[ 
\fr{1}{2} \left( \fr{1}{2} \Lambda^{2} \fr{\Delta t}{\Delta x} + \Lambda^{-} \right) 
\left( {\bf f}^{eq}_{j} - {\bf f}^{eq}_{j-1} \right)
\right] \ea
\eea 
where 
\be 
{\bf r}^{-}_{j} = \fr{ {\bf f}^{eq}_{j+1} - {\bf f}^{eq}_{j} } 
{{\bf f}^{eq}_{j} - {\bf f}^{eq}_{j-1}} \ \textrm{and} \ 
{\bf r}^{+}_{j-1} = \fr{ {\bf f}^{eq}_{j-1} - {\bf f}^{eq}_{j-2} } 
{{\bf f}^{eq}_{j} - {\bf f}^{eq}_{j-1}} 
\ee
we obtain 
\bea \ba{ll} 
g(u)_{j-\frac{1}{2}}^{TVD-KFDS} & =g(u)_{j-\frac{1}{2}}^{KFDS}\\[2mm] 
& - \fr{1}{2} minmod \left[ \ba{c}  
\left\{ \fr{1}{2} \fr{\Delta t}{\Delta x} \lambda^{2} \left( u_{j} - u_{j-1} \right) 
- \fr{1}{2} \lambda \left( u_{j} - u_{j-1} \right) 
- \fr{1}{2} \left( g(u)_{j} - g(u)_{j-1} \right) 
\right\}, \\[2mm]  
\left\{ \fr{1}{2} \fr{\Delta t}{\Delta x} \lambda^{2} \left( u_{j-1} - u_{j-2} \right) 
- \fr{1}{2} \lambda \left( u_{j-1} - u_{j-2} \right) 
- \fr{1}{2} \left( g(u)_{j-1} - g(u)_{j-2} \right) 
\right\} \ea 
\right] \\[5mm] 
& - \fr{1}{2} minmod \left[ \ba{c}  
\left\{ \fr{1}{2} \fr{\Delta t}{\Delta x} \lambda^{2} \left( u_{j} - u_{j-1} \right) 
- \fr{1}{2} \lambda \left( u_{j} - u_{j-1} \right) 
+ \fr{1}{2} \left( g(u)_{j} - g(u)_{j-1} \right) 
\right\}, \\[2mm]  
\left\{ \fr{1}{2} \fr{\Delta t}{\Delta x} \lambda^{2} \left( u_{j+1} - u_{j} \right) 
- \fr{1}{2} \lambda \left( u_{j+1} - u_{j} \right) 
+ \fr{1}{2} \left( g(u)_{j+1} - g(u)_{j} \right) 
\right\} \ea 
\right] \ea   
\eea 
These expressions for $g(u)_{j+\frac{1}{2}}$ and $g(u)_{j-\frac{1}{2}}$ are used in the update formula given by 
\be 
u^{n+1}_{j} = u^{n}_{j} - \fr{\Delta t}{\Delta x} \left[ 
g(u)_{j+\frac{1}{2}} - g(u)_{j-\frac{1}{2}} \right] + \fr{\Delta t}{\Delta x} \left[ g_{v}(u)_{j+\frac{1}{2}} - g_{v}(u)_{j-\frac{1}{2}} \right] 
\ee
which is the moment of the update formula (\ref{FVM_update_formulat_for_DVBEvisc}) for the flexible velocity Boltzmann equation.

\section{KFDS for 1-D Shallow Water Equations}
In this section, we expand the scope of the KFDS scheme to a system of hyperbolic equations. The system of shallow water equations in one dimension with a source term is given by
\be \label{1D_SWE}
\fr{\partial U}{\partial t}+\fr{\partial G(U)}{\partial x}+ S(U) = 0 
\ee
where $U$ is the conserved variable vector, $G(U)$ the Flux vector and $S(U)$ is the source term vector. These are defined by
\be \label{1D_SWE_COMP}
U=\left[\begin{array}{c} h \\ h u \end{array}\right], \quad G(U)=\left[\begin{array}{c} h u \\ h u^{2}+\frac{1}{2} g h^{2} \end{array}\right] \quad and \quad S(U)=\left[\begin{array}{c} 0 \\ -g h b_{x} \end{array}\right] 
\ee 
where h is the depth of the water, $u$ the fluid velocity, $g$ the acceleration due to gravity, $b(x)$ is the geometric description of the water bed and $b_{x}$ is its derivative in $x$. 
\subsection{1-D Shallow Water Equations}
In this section we construct the 1D shallow water equations with source terms, first from the Boltzmann equation with continuous molecular velocity case, and thereby strategise the FVBE model for shallow water equations. The Boltzmann equation including the external forces acting on the fluid can be written as
\be \label{1D_BOLTZ}
\fr{\del  f}{\del t} +  \fr{\del \left(vf\right)}{\del x} +\mathbb{F}\fr{\del f}{\del v} = - \fr{1}{\epsilon} \left[  f - f^{eq} \right] 
\ee  
where $\mathbb{F} = \mathbb{F}(x)$ is the Force acting on the fluid. The equilibrium distribution function, $f^{eq}$ is defined as
 
\be 
f^{eq} = \phi \sqrt{\fr{\beta}{\pi}} e^{- \beta \left( v - u \right)^{2} } \ \ \  \textrm{where} \ \beta = \fr{1}{\phi} 
\ee 
which leads to  
\be 
f^{eq} = \sqrt{\fr{\phi}{\pi}} e^{- \fr{\left( v - u \right)^{2}}{\phi}}.  
\ee

To construct the 1D shallow water equations given in [\ref{1D_SWE}], we multiply [\ref{1D_BOLTZ}] by a vector $\psi$ as
\be  \label{1D_BOLTZ2}
\psi \left \{ \fr{\del f}{\del t} +\fr{\del \left(vf\right)}{\del x} +\mathbb{F}\fr{\del f}{\del v} = - \fr{1}{\epsilon} \left[  f -  f^{eq} \right] \right \} , \ \ \textrm{where} \ \psi = \left[ \begin{array}{c} 1 \\ v \end{array}\right ] .
\ee
The introduction of $\psi$ vector converts the scalar FVBE into the desired vector form and the inclusion of force component at kinetic level leads to the source terms at the macroscopic level. The macroscopic system can be recovered by following the  moment relations for each term in  [\ref{1D_BOLTZ2}], as shown in Annexure D. 

\be \label{SWE_CV}
 \int\limits_{- \infty}^{\infty}{\psi {f}^{eq}} dv = \int\limits_{- \infty}^{\infty}{\left[ \begin{array}{c} {f}^{eq} \\ {v} {f}^{eq} \end{array}\right ]} dv\ = \left[ \begin{array}{c} \phi \\ \phi u \end{array}\right ].
\ee

\be \label{SWE_CF}
 \int\limits_{- \infty}^{\infty}{\psi v {f}^{eq}} dv = \int\limits_{- \infty}^{\infty}{\left[ \begin{array}{c} {v}{f}^{eq} \\ {v}^{2} {f}^{eq} \end{array}\right ]} dv\ = \left[ \begin{array}{c} \phi u \\ \fr{{\phi}^{2}}{2} + {\phi}{ u }^{2} \end{array}\right ] .
\ee

\be \label{SWE_CFOR}
 \int\limits_{- \infty}^{\infty}{\psi \mathbb{F}\fr{\del {f}^{eq}}{\del v}} dv =  \int\limits_{- \infty}^{\infty}{\left[ \begin{array}{c} {f}^{eq}\mathbb{F}\fr{\del (1)}{\del v} \\  {f}^{eq}\mathbb{F}\fr{\del (v)}{\del v} \end{array}\right ]} dv\ = \left[ \begin{array}{c} 0\\ -\phi  \mathbb{F} \end{array}\right ] .
\ee

Here we assign, $\phi = h$ and $\mathbb{F}=g{b}_{x}$, where $g$ is the acceleration due to gravity, $h$ is the geopotential height and ${b}_{x}$ represents the rate of change of the bottom floor of the shallow water system. Therefore, the moments of the equations [\ref{SWE_CV}], [\ref{SWE_CF}] and [\ref{SWE_CFOR}] lead to  [\ref{1D_SWE_COMP}]. 
The two dimensional shallow water equations can be derived by following the procedures similar to that in one dimension. 

\subsection{Flexible velocity model for shallow water equations}

Introducing two flexible velocities $\lambda^{+}$ and $\lambda^{-}$ to replace the continuous variation of the molecular velocity $v$ and also two corresponding components of $\tilde f^{eq}_{+}$ and $ \tilde f^{eq}_{-}$, we write 
\be
\tilde{f}^{eq} = { \left\{\tilde f^{eq}_{+} \delta (v - \lambda^{+}) +  \tilde f^{eq}_{-} \delta (v - \lambda^{-}) \right\} } .
\ee 
To recover the  (vector of) conserved variables we multiply the scalar $\tilde f^{eq}$  with $\psi$. Thus, the conserved variable vector becomes 
\be
U_{i} ={ \left[ \int_{-\infty}^{\infty}{\psi}_{i} \tilde f^{eq}_{+} \delta (v - \lambda^{+})  dv 
+ \int_{-\infty}^{\infty} {\psi}_{i}\tilde f^{eq}_{+} \delta (v - \lambda^{-})  dv \right] } .
\ee 
We combine the $\psi \tilde f^{eq}$ to form $f^{eq}$. Therefore we get the conserved variable vector as 
\be
U_{i} ={ \left[ \int_{-\infty}^{\infty} f^{eq}_{+} \delta (v - \lambda^{+})  dv 
+ \int_{-\infty}^{\infty}  f^{eq}_{+} \delta (v - \lambda^{-})  dv \right] }_{i} .
\ee 
Let us further assume, for simplicity, that the flexible velocities, $\lambda^{+}$ and $\lambda^{-}$ for each $i$ are given by 
\be \label{simple_2DVs_1D}
\lambda^{+}_{i} = \lambda_{i} \ \textrm{ and } \  \lambda^{-}_{i} = - \lambda_{i} .   
\ee
Thus, we have three unknowns, namely, $f^{eq}_{+}$,  $f^{eq}_{-}$ and $\lambda$  for each $i$ to be fixed for the equilibrium distribution. 
$$ \ba{rcl}  
U_{i}  & = & \displaystyle \int_{-\infty}^{\infty} \left\{ f^{eq}_{+} \delta (v - \lambda^{+}) + f^{eq}_{-} \delta (v - \lambda^{-}) \right\} dv  \\ 
 & = & \left\{ f^{eq}_{+} (1) + f^{eq}_{-} (1) \right\}_{i} .  
\ea $$ 
Thus 
\be 
\left\{f^{eq}_{+} + f^{eq}_{-}\right\}_{i} = U_{i}. 
\ee 
$$ \ba{rcl} 
G(U)_{i} & = & \displaystyle \int_{-\infty}^{\infty} v {\psi}_{i} \tilde {f}^{eq} dv \\ [3mm]
 & = & \displaystyle \int_{-\infty}^{\infty} v 
\left\{ f^{eq}_{+} \delta (v - \lambda^{+}) + f^{eq}_{-} \delta (v - \lambda^{-}) \right\}_{i} dv \\ [2mm]
 & = & \left\{ f^{eq}_{+} \displaystyle \int_{-\infty}^{\infty} v \delta (v - \lambda^{+}) dv +  f^{eq}_{-} \displaystyle \int_{-\infty}^{\infty} v \delta (v - \lambda^{-}) dv  \right\}_{i}  \\
 & = & \left\{  f^{eq}_{+} \lambda^{+} +  f^{eq}_{-} \lambda^{-}  \right\}_{i}
\ea $$ 
Thus 
\be 
 \left\{ f^{eq}_{+} \lambda^{+} +   f^{eq}_{-} \lambda^{-} \right\}_{i}  = G(U)_{i}
\ee 
Solving the above two equations and simplifying  of flexible velocities for each $i$ using (\ref{simple_2DVs_1D}), we get
\be \label{2_Vel_equilibria}
{f^{eq}_{+}}_{i} = \fr{1}{2} U_{i} + \fr{1}{2 {\lambda}_{i}} G(U)_{i} \ \textrm{and} \ {f^{eq}_{-}}_{i} = \fr{1}{2} U_{i}- \fr{1}{2 {\lambda}_{i}} G(U)_{i}
\ee
To recover the source terms we substitute the flexible velocities into the source term as explained in Annexure D.
\bea
\int\limits_{- \infty}^{\infty} {\psi}_{i} \mathbb{F}\fr{\del f}{\del v}dv = - \int\limits_{- \infty}^{\infty}{\mathbb{F} \tilde{f}_{eq}\fr{\del{ \psi}_{i} }{\del v}}dv  \\
\psi = 1 :  \int\limits_{- \infty}^{\infty}{\mathbb{F} \tilde {f}_{eq}\fr{\del {\psi}_{i} }{\del v}}dv =  \int\limits_{- \infty}^{\infty}{\mathbb{F} \left\{ \tilde f^{eq}_{+} \delta (v - \lambda^{+}) +\tilde f^{eq}_{-} \delta (v - \lambda^{-}) \right\}\fr{\del (1) }{\del v}}dv = 0 \\
\psi = v :   \int\limits_{- \infty}^{\infty}{\mathbb{F}\tilde {f}_{eq}\fr{\del {\psi}_{i} }{\del v}}dv = \int\limits_{- \infty}^{\infty}{\mathbb{F} \left\{\tilde f^{eq}_{+} \delta (v - \lambda^{+}) +\tilde f^{eq}_{-} \delta (v - \lambda^{-}) \right\}\fr{\del (v) }{\del v}}dv = \mathbb{F} h
\eea

Thus the moment strategy recovers the source terms as defined in [\ref{1D_SWE_COMP}]. Therefore, the {\em Flexible Velocity Boltzmann Equation} (FVBE) based on the two flexible velocity model can be written as    
\be \label{DVBE_1D_2_Vel} 
\left\{\fr{\del \bf f}{\del t} + \fr{\del  \left(\bf {\Lambda}\bf f\right)}{\del x}  + (i-1){\bf F}  = - \fr{1}{\epsilon} \left[ \bf f - \bf f^{eq} \right] \right\}_{i}     \  \textrm{ with} \ \ i=1,2
\ee
where 
\be \label{1D_2Vel_Feq}
\bf f_{i} = \left[ \ba{c} f_{+} \\ f_{-} \ea \right]_{i}, \ 
\bf F_{i} =  \left[ \ba{c} \mathbb{F}h \\ 0 \ea \right]_{i}
\Lambda_{i} = \left[ \ba{cc} \lambda^{+} & 0 \\ 0 & \lambda^{-} \ea \right]_{i} \ \textrm{and} \ 
\bf {f^{eq}}_{i} = \left[ \ba{c} f^{eq}_{+} \\ f^{eq}_{-} \ea \right]_{i}
= \left[ \ba{c} \fr{1}{2} U + \fr{1}{2 \lambda} G(U) \\[4mm]  
\fr{1}{2} U - \fr{1}{2 \lambda} G(U) \ea \right]_{i}.
\ee

\subsection{KFDS, KLW and TVD Schemes for Shallow Water Equations}

We follow the same methodology used in the scalar framework to obtain the finite volume update formula for the shallow water system of equations. To handle the source terms in the momentum equation, we follow a well-balancing approach. 

\subsubsection{KFDS}
The discrete framework for 1D Shallow water equation with operator split source terms can be obtained as follows.   
Let us write a finite volume update formula  for the 1D system of shallow water equations with source terms as (\ref{DVBE_1D_2_Vel})     
\be 
{\bf P} {\left[ \fr{\del \bf f}{\del t} +  \fr{\del \Lambda \bf f}{\del x} + (i-1) {\bf F} = 0, \ {\bf f} = {\bf f^{eq}}   \right]}_{i}  \ \ i=1,2.
\ee 
The numerical scheme is evolved following the KFDS methodology. 

$$ {\left[ \fr{\left( {\bf P f} \right)^{n+1}_{j} - \left({\bf P f} \right)^{n}_{j} } {\Delta t} +  \fr{\left({\bf P} \Lambda {\bf f}\right)^{n}_{j+\frac{1}{2}} 
- \left({\bf P} \Lambda {\bf f}\right)^{n}_{j-\frac{1}{2}}} {\Delta x} -(i-1){\bf P F}_{j}  = 0, \ 
{\bf f}^{n}_{j} = {\bf f^{eq,}}^{n}_{j}  \right]}_{i}  \ \ i=1,2.$$ 
To treat the source terms, we introduce well-balancing by evaluating the source term using the average of variables at the interface and substituting for $\mathbb{F}$.  Thus we obtain 
$$ {\left[ \fr{\left( {\bf P f} \right)^{n+1}_{j} - \left({\bf P f} \right)^{n}_{j} } {\Delta t} +  \fr{\left({\bf P} \Lambda {\bf f}\right)^{n}_{j+\frac{1}{2}} 
- \left({\bf P} \Lambda {\bf f}\right)^{n}_{j-\frac{1}{2}}} {\Delta x} -(i-1)g \fr{\left ( b_{j+\frac{1}{2}} - b_{j-\frac{1}{2}}\right)}{\Delta x} \fr{\left (h ^{n}_{j+\frac{1}{2}} 
+ h ^{n}_{j-\frac{1}{2}} \right )}{2} = 0, \ 
{\bf f}^{n}_{j} = {\bf f^{eq,}}^{n}_{j}  \right]}_{i}  \ \ i=1,2.$$ 
The Final update formula thus obtained by following finite volume method based upwind flexible velocity Boltzmann scheme is given by 

\be 
\left[ {U}^{n+1}_{j} = {{U}}^{n}_{j} - \fr{\Delta t}{\Delta x} \left[ {{G(U)}}^{n}_{j+\frac{1}{2}} - {{G(U)}}^{n}_{j-\frac{1}{2}} \right] +  \fr{(i-1)}{2}\fr{\Delta t}{\Delta x}(g)({h}_{j+\frac{1}{2}} + h_{j-\frac{1}{2}})({b_{j+\frac{1}{2}} - b_{j-\frac{1}{2}} }) \right]_{i} 
\ee 
where the interface fluxes are given by
\bea 
\left \{ G(U)_{j+\frac{1}{2}}^{n}= \fr{1}{2}[  G(U)_{j+1}^{n}+ G(U)_{j}^{n}]- \fr{1}{2}[\Delta G(U)^{+,n}_{j+\frac{1}{2}}-\Delta G(U)^{-,n}_{j+\frac{1}{2}}]\right \}_{i}\\ 
 \left \{G(U)_{j-\frac{1}{2}}^{n}= \fr{1}{2}[  G(U)_{j}^{n}+ G(U)_{j-1}^{n}]- \fr{1}{2}[\Delta G(U)^{+,n}_{j-\frac{1}{2}}-\Delta G(U)^{-,n}_{j-\frac{1}{2}}] \right \}_{i}
\eea
\bea \label{2Vel_NKFDS} 
 \left \{\Delta G(U)_{j+\frac{1}{2}}^{\pm} = \fr{1}{2}[  G(U)_{j+1}- G(U)_{j}]\pm \fr{1}{2}|\lambda|[U_{j+1}-U_{j}] \right \}_{i}\\ 
 \left \{\Delta G(U)_{j-\frac{1}{2}}^{\pm} = \fr{1}{2}[  G(U)_{j}- G(U)_{j-1}]\pm \fr{1}{2}|\lambda|[U_{j}-U_{j-1}] \right \}_{i}
\eea 

\subsection{Fixing \texorpdfstring{$\lambda$}{Lambda}}

Consistent with the scalar conservation laws, two different ways of fixing $\lambda$ are proposed here. While the first way is based on Chapman-Enskog type asymptotic expansion of the FVBE,  the second way is based on utilizing the Rankine-Hugoniot jump condition.

\subsubsection{Chapman-Enskog type expansion of the flexible velocity Boltzmann equation} 
The methodology followed for Chapman-Enskog type expansion of the discrete velocity Boltzmann equation followed in Appendix C can be extended to the FVBE framework of shallow water equations. The resulting expression leads to the original convection equations (plus a source term), augmented by a diffusion term, given by 
\be
\left[ \fr{\del U}{\del t} + \fr{\del G(U)}{\del x} + S(U) = \epsilon \fr{\del}{\del x} \left[ 
\fr{\del u}{\del x} \left\{ \Lambda^{2} - \left( A\left(U\right) \right)^{2} \right\}
\right] + \mathcal{O}\left( \epsilon^{2} \right) \right]
\ee
where $A(U)$ is the flux Jacobian, whose eigenvalues form the wave speeds for the original system of convection equations, given by $A(U)=\fr{\del G(U)}{\del U}$.  For the approximation to be stable, the diffusion term must be non-negative.  Therefore, we obtain 
\be
\Lambda^{2} - \left( A\left(U\right)\right)^{2} \ge 0  \ where \ \Lambda = \left[ \ba{cc} \lambda_{1} & 0 \\ 0 & \lambda_{2} \ea \right] 
\ee 
which can approximated to 
\be 
\Lambda = \max \left[ A(U)_{j} \right], \ \forall j  
\ee 
which is the maximum eigenvalue of the system. For the shallow water system of equations, the eigenvalues are given by $(u-a)$ and $(u+a)$, where $a = \sqrt{gh}$ is the celerity. Therefore the $\lambda$ for the scheme is fixed as
\be
 \lambda_{i} = max(u_L-a_L , u_L+a_L , u_R-a_R, u_R+a_R) .
\ee

\subsection{Fixing \texorpdfstring{$\lambda$}{Lambda} based on R-H (jump) condition}
Here we fix $\lambda$ based on the methodology followed in section [\ref{RHC}]. For the shallow water system of equations we have 
\be
 A(U) = \left[ \ba{cc} 0 & 1 \\ {a}^{2}-{u}^{2} & 2u \ea \right] 
\ee
The eigenvalues of A(U) are $(u-a)$ and $(u+a)$. Therefore, the wave speed for the scheme is assigned  as 
\be
\Lambda = \left[ \ba{cc} |u_L - a_L| & 0 \\ 0 & |u_L+a_L| \ea \right] 
\ee
\subsubsection{KLW}
The discretisation of the  source term part is inherently second order accurate. Therefore it is sufficient to derive a KLW scheme for the convection part of the system of equations.  We obtain the finite volume update formula for the KLW scheme by following procedures as in section \ref{KLW_Start}.
\be 
{\bf P} \left[ {\bf f}^{n+1}_{j} = {\bf f}^{n}_{j} 
- \fr{\Delta t}{\Delta x} 
\left[ {\bf h}^{n}_{j+\frac{1}{2}} - {\bf h}^{n}_{j-\frac{1}{2}} \right] - \Delta t (i-1){\bf F}_{j}\right]_{i}
\ee 
or 
\be \label{KLW_SWE_update}
\left[ U^{n+1}_{j} = U^{n}_{j}  - \fr{\Delta t}{\Delta x} \left[ G(U)^{n}_{j+\frac{1}{2}} - G(U)^{n}_{j-\frac{1}{2}} \right] + \Delta t S(U) \right]_{i}
\ee
where 
\be \label{KLW_SWE_interface_fluxA}
\left[ G(U)_{j+\frac{1}{2}} = \fr{1}{2} \left[ G(U)_{j} + G(U)_{j+1} \right] 
- \fr{1}{2} \lambda^{2} \fr{\Delta t}{\Delta x} \left[U_{j+1} - U_{j} \right] \right]_{i} 
\ee 

\be \label{KLW_SWE_interface_fluxB}
\left[ G(U)_{j-\frac{1}{2}} = \fr{1}{2} \left[ G(U)_{j-1} + G(U)_{j} \right] 
- \fr{1}{2} \lambda^{2} \fr{\Delta t}{\Delta x} \left[U_{j} - U_{j-1} \right] \right]_{i} 
\ee 
Again we rewrite the source term using the well-balancing approach. The final form of [\ref{KLW_SWE_update}] will become 
\be \label{KLW_SWE_update2}
\left[ U^{n+1}_{j} = U^{n}_{j}  - \fr{\Delta t}{\Delta x} \left[ G(U)^{n}_{j+\frac{1}{2}} - G(U)^{n}_{j-\frac{1}{2}} \right] + \fr{(i-1)}{2}\fr{\Delta t}{\Delta x}(g)({h}_{j+\frac{1}{2}} + h_{j-\frac{1}{2}})({b_{j+\frac{1}{2}} - b_{j-\frac{1}{2}} }) \right]_{i}
\ee
where $G(U)_{j+\frac{1}{2}}$ \& $G(U)_{j-\frac{1}{2}}$ are given by [\ref{KLW_SWE_interface_fluxA}] \& [\ref{KLW_SWE_interface_fluxB}] respectively.

\subsubsection{TVD-KFDS}
The TVD scheme is derived by following the strategy introduced in  [\ref{TVD_strat}]. The finite volume scheme for the TVD approach thus derived using the approaches defined in earlier section can be written as
\be 
\left[ U^{n+1}_{j} = U^{n}_{j} - \fr{\Delta t}{\Delta x} \left[ G(U)_{j+\frac{1}{2}} - G(U)_{j-\frac{1}{2}} \right] +  \fr{(i-1)}{2}\fr{\Delta t}{\Delta x}(g)({h}_{j+\frac{1}{2}} + h_{j-\frac{1}{2}})({b_{j+\frac{1}{2}} - b_{j-\frac{1}{2}} })\right]_{i} 
\ee  
where the expressions for $G(U)_{j+\frac{1}{2}}$ and $G(U)_{j-\frac{1}{2}}$ are given in the TVD form duly regulated by the minmod limiter based on the flux gradient.

\bea \ba{ll} 
{G_{i}(U)}_{j+\frac{1}{2}}^{TVD-KFDS}  & =  {G_{i}(U)}_{j+\frac{1}{2}}^{KFDS} \\[2mm]  
& - \fr{1}{2} minmod \left[ \ba{c}  
\left\{ \fr{1}{2} \fr{\Delta t}{\Delta x} \lambda^{2} \left( U_{j+1} - U_{j} \right) 
- \fr{1}{2} \lambda \left( U_{j+1} - U_{j} \right) 
- \fr{1}{2} \left( G(U)_{j+1} - G(U)_{j} \right) 
\right\}, \\[2mm]  
\left\{ \fr{1}{2} \fr{\Delta t}{\Delta x} \lambda^{2} \left( U_{j} - U_{j-1} \right) 
- \fr{1}{2} \lambda \left( U_{j} - U_{j-1} \right) 
- \fr{1}{2} \left( G(U)_{j} - G(U)_{j-1} \right) 
\right\} \ea 
\right]_{i} \\[5mm] 
& - \fr{1}{2} minmod \left[ \ba{c}  
\left\{ \fr{1}{2} \fr{\Delta t}{\Delta x} \lambda^{2} \left( U_{j+1} - U_{j} \right) 
- \fr{1}{2} \lambda \left( U_{j+1} - U_{j} \right) 
+ \fr{1}{2} \left( G(U)_{j+1} - G(U)_{j} \right) 
\right\}, \\[2mm]  
\left\{ \fr{1}{2} \fr{\Delta t}{\Delta x} \lambda^{2} \left( U_{j+2} - U_{j+1} \right) 
- \fr{1}{2} \lambda \left( U_{j+2} - U_{j+1} \right) 
+ \fr{1}{2} \left( G(U)_{j+2} - G(U)_{j+1} \right) 
\right\} \ea 
\right]_{i} \ea   
\eea

\bea \ba{ll} 
{G_{i}(U)}_{j-\frac{1}{2}}^{TVD-KFDS} & = {G_{i}(U)}_{j-\frac{1}{2}}^{KFDS} \\[2mm]  
& - \fr{1}{2} minmod \left[ \ba{c}  
\left\{ \fr{1}{2} \fr{\Delta t}{\Delta x} \lambda^{2} \left( U_{j} - U_{j-1} \right) 
- \fr{1}{2} \lambda \left( U_{j} - U_{j-1} \right) 
- \fr{1}{2} \left( G(U)_{j} - G(U)_{j-1} \right) 
\right\}, \\[2mm]  
\left\{ \fr{1}{2} \fr{\Delta t}{\Delta x} \lambda^{2} \left( U_{j-1} - U_{j-2} \right) 
- \fr{1}{2} \lambda \left( U_{j-1} - U_{j-2} \right) 
- \fr{1}{2} \left( G(U)_{j-1} - G(U)_{j-2} \right) 
\right\} \ea 
\right]_{i} \\[5mm] 
& - \fr{1}{2} minmod \left[ \ba{c}  
\left\{ \fr{1}{2} \fr{\Delta t}{\Delta x} \lambda^{2} \left( U_{j} - U_{j-1} \right) 
- \fr{1}{2} \lambda \left( U_{j} - U_{j-1} \right) 
+ \fr{1}{2} \left( G(U)_{j} - G(U)_{j-1} \right) 
\right\}, \\[2mm]  
\left\{ \fr{1}{2} \fr{\Delta t}{\Delta x} \lambda^{2} \left( U_{j+1} - U_{j} \right) 
- \fr{1}{2} \lambda \left( U_{j+1} - U_{j} \right) 
+ \fr{1}{2} \left( G(U)_{j+1} - G(U)_{j} \right) 
\right\} \ea 
\right]_{i} \ea   
\eea 

\section{KFDS, KLW and TVD Schemes in Two Dimensions}
 The KFDS, KLW and TVD schemes discussed presented in the previous sections for one dimensional scalar convection-diffusion equation can be extended to multi-dimensions using the finite volume method easily.  Across a cell-interface in the finite volume method, the usual mid-point rule based quadrature leads to a locally one-dimensional approximation and hence the procedure presented in the previous sections simply gets extended.  In this section, we consider the two dimensional  Burgers equation given by 
\be
\frac { \partial u } { \partial t } 
+ \frac { \partial g_{1}(u) } { \partial x } 
+ \frac { \partial g_{2}(u) } { \partial y } 
= 0
\ee 
 using the standard finite volume method.  The terms of the inviscid and viscous fluxes in the above equation have similar definitions as in one dimension.  
A four-velocity based FVBE for the above system can be derived as 
\be  
\frac{\partial \mathbf{f}}{\partial t}+ \frac{\partial \mathbf{{h}_{1}}}{\partial x}+ \frac{\partial \mathbf{\bf {h}_{2}}}{\partial y}= - \fr{1}{\epsilon}\left[ \mathbf{f} -\mathbf{f}^{eq}\right] 
\ee
where $ \bf {h}_{1} = {\Lambda}_{1}f, \bf {h}_{2} = {\Lambda}_{2}f $ and the flexible velocity matrices are given by 
\be
{\Lambda}_{1}=\left[\begin{array}{cccc}{-\lambda} & {0} & {0} & {0}  \\ {0} & {\lambda} & {0} & {0}  \\ {0} & {0}  & {\lambda} & {0}  \\ {0} & {0} & {0} & {-\lambda} \end{array}\right]   \Lambda_{2}=\left[\begin{array}{cccc}{-\lambda} & {0} & {0} & {0}  \\ {0} & {-\lambda} & {0} & {0}  \\ {0} & {0}  & {\lambda} & {0}  \\ {0} & {0} & {0} & {\lambda} \end{array}\right]
\ee
The equilibria are defined by 
\be
{\bf f}^{eq}=\left[\begin{array}{l}{f^{eq}_{1}} \\ {f^{eq}_{2}} \\ {f^{eq}_{3}} \\ {f^{eq}_{4}} \end{array}\right]=\left[\begin{array}{c}{\frac{1}{4} U-\frac{1}{4 \lambda}G_{1}-\frac{1}{4 \lambda} G_{2}} \\ {\frac{1}{4}U+\frac{1}{4 \lambda} G_{1}-\frac{1}{4 \lambda} G_{2}}  \\ {\frac{1}{4} U+\frac{1}{4 \lambda} G_{1}+\frac{1}{4 \lambda} G_{2}} \\ {\frac{1}{4} U-\frac{1}{4 \lambda} G_{1}+\frac{1}{4} \lambda} G_{2} \end{array}\right]
\ee
Let us denote $\mathbf{h}_{x}=\mathbf{\Lambda}_{1} \mathbf{f}, \mathbf{h}_{y}=\mathbf{\Lambda}_{2} \mathbf{f}$ as the fluxes along the $x-$ and $y-$ directions.   
The cell-interface fluxes across a cell are constructed be normal to the interfaces and can be obtained as 
\bea
\mathbf{h}_{ L}=\left(\mathbf{h}_{x} \cos \theta+\mathbf{h}_{y} \sin \theta\right)_{L}\\
\mathbf{h}_{ R}=\left(\mathbf{h}_{x} \cos \theta+\mathbf{h}_{y} \sin \theta\right)_{R}
\eea
where the suffix L and R represent the left and the right states of the cell-interface.  
A finite volume update formula for the 2D Euler system can be derived as 
\be 
\mathbf{f}_{j,k}^{n+1}=\mathbf{f}_{j,k}^{n}-\frac{\Delta t}{A_{j,k}} \sum_{I_{c}=1}^{4} \mathbf{h}_{n I_{c}} \Delta s_{I_{c}}
\ee
where  $A_{j,k}$ is the area of the cell centered at $(j,k)$ and $\Delta s_{I_{c}}$ is the length of the cell interface $I_{c}$.    
Applying the {\em Kinetic Flux Difference Splitting} across a finite volume cell-interface, we get  
\be
\mathbf{h}_{n, I_{c}}=\frac{1}{2}\left[h_{ R}^{n}+h_{ L}^{n}\right]-\frac{1}{2}\left[\Delta h_{I,n}^{+}-\Delta h_{I,n}^{-}\right]
\ee 
where 
\be
\Delta h_{I,n}^{\pm} = \left(\Lambda_{1 } \cos \theta\right)^{\pm}\Delta \mathbf{f}_{I,n}^{e q}+\left(\Lambda_{2 } \sin \theta\right)^{\pm}\Delta \mathbf{f}_{In}^{e q}
\ee 
Upon taking moments of the above equations and extending the same techniques that we had applied in one dimensional system in determining the  diffusion coefficients, we recover the macroscopic update formula as 
\be 
\mathbf{U}_{j,k}^{ n+1}=\mathbf{U}_{j,k}^{ n}-\frac{\Delta t}{A_{j,k}} \sum_{I_{c}=1}^{4} \mathbf{G}_{n, I_{c}} \Delta s_{I_{c}}
\ee
where 
\bea
\mathbf{G}_{n, I_{c}} = \frac{1}{2}\left[G(U)_{n, R}^{n}+G(U)_{n, L}^{n}\right]-\frac{1}{2}\left[\Delta G(U)_{I_{c}}^{+}-\Delta G(U)_{I_{c}}^{-}\right]\\
\Delta G(U)_{I_{c}}^{\pm} = \left.\frac{1}{2}\left[G(U)_{n, R}^{n}-G(U)_{n, L}^{n}\right]\right] \pm \frac{1}{2}|\lambda|\left[U_{n, R}-U_{n, L}\right]
\eea

\section{Experimental Order of Convergence}
In this section, we present the results of the computational experiment performed to evaluate order of accuracy and the accumulation of error with increasing grid size. For the numerical experiment, we choose the 1D inviscid Burgers framework given by  $\fr{\del u}{\del t} + \fr{\del g(u)}{\del x} = 0, \ g(u) = \fr{1}{2} u^{2}$. The computational domain is defined by $x  \in [0,1]$ with the initial conditions $u(x,0) = sin(2\pi x) $. The test case is unsteady by nature wherein the initial velocity gradients evolve into steeper gradient, until $t = \fr{1}{2\pi}$ when the gradients evolve into a steady standing  shock at $x=0.5$. For conducting our numerical experiment we choose the final time as $t=\fr{0.4}{\pi}$ where the velocity gradients are strong, yet without any discontinuities.

The numerical simulations are carried out by changing the number of computational cells methodically as 10, 20, 40, $\cdots$, 5120 cells. The L1 and L2 errors represented by $\left\|\mathcal{E}_{K}\right\|_{L^{1}} and \left\|\mathcal{E}_{K}\right\|_{L^{2}}$  are calculated at $t^{n}=\fr{0.4}{\pi}$ seconds using equation [\ref{L1Form}] and [\ref{L2Form}] as given below.

\be
\label{L1Form}
\left\|\mathcal{E}_{K}\left(t^{n}\right)\right\|_{L^{1}}=\Delta x \sum_{j=1}^{K}\left|u_{j}^{numerical}-u_{j}^{analytical}\right| 
\ee
\be
\label{L2Form}
\left\|\mathcal{E}_{K}\left(t^{n}\right)\right\|_{L^{2}}=\sqrt{\Delta x \sum_{j=1}^{K}\left(u_{j}^{numerical}-u_{j}^{analytical}\right)^{2}} \\
\ee
where $K$ is the number of cells. 
The Experimental order of convergence is computed using [\ref{EOC}].  
\be
\label{EOC}
\mathrm{EOC}=\log _{2}\left(\frac{\left\|\mathcal{E}_{K / 2}\right\|}{\left\|\mathcal{E}_{K}\right\|}\right)
\ee

\begin{figure}[!h] 
\begin{center} 
\includegraphics[width=14cm,angle=0]{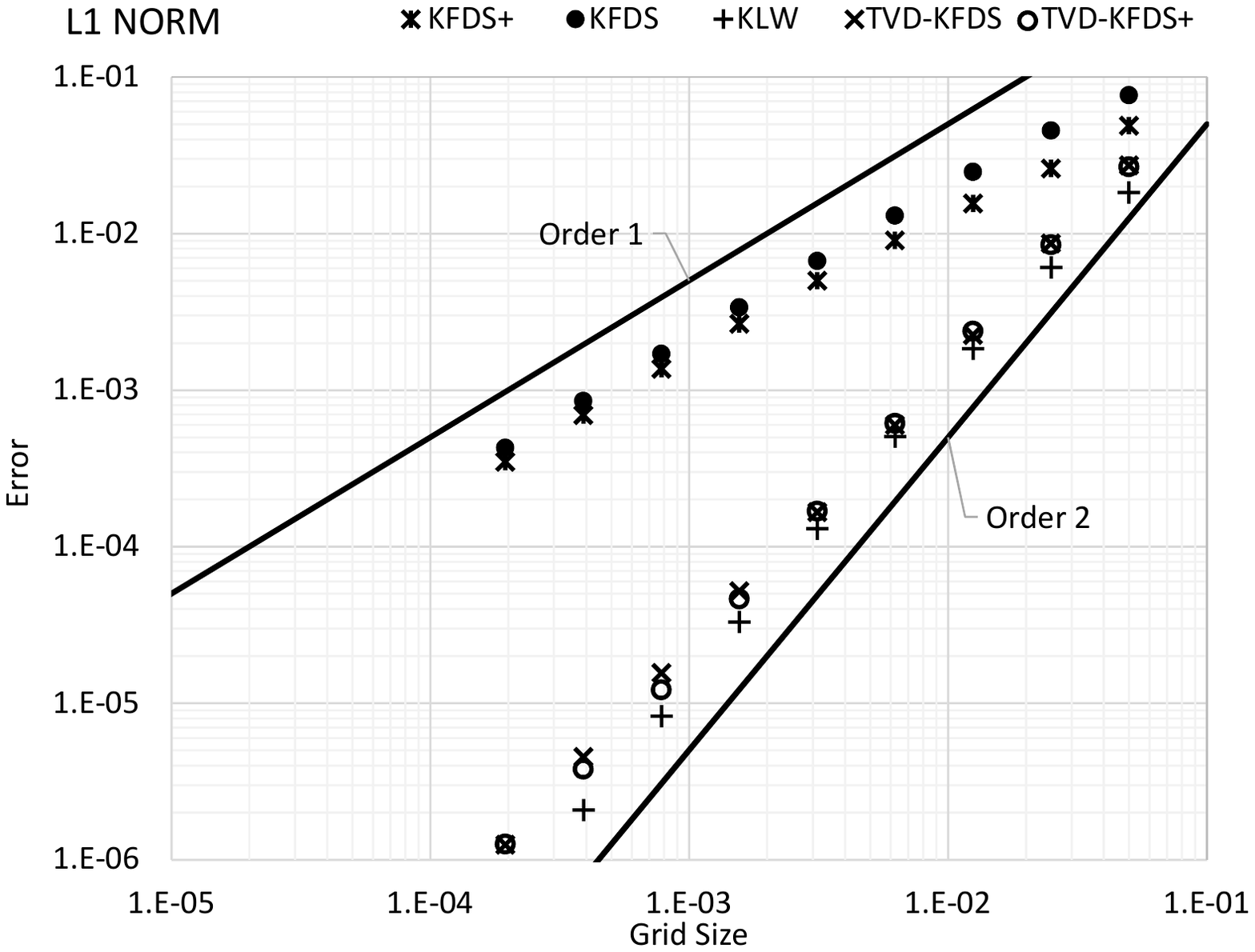}  
\caption{L1 Norm Errors for a) KFDS b) KFDS+ c) KLW d) TVD-KFDS \& e) TVD-KFDS+ schemes} 
\label{L1_Norm} 
\end{center} 
\end{figure}

\begin{figure}[!h] 
\begin{center} 
\includegraphics[width=14cm,angle=0]{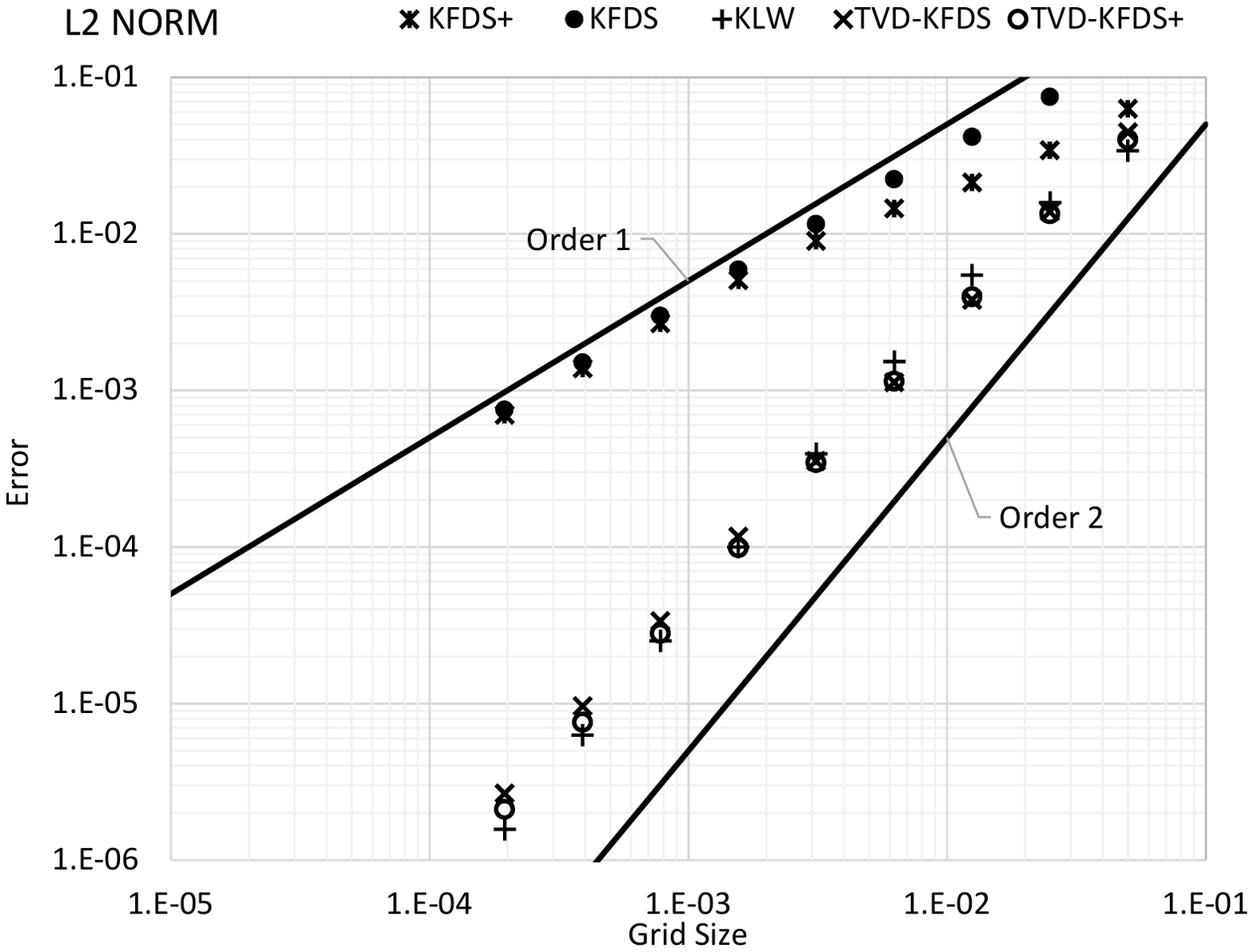}  
\caption{L2 Norm Errors for a) KFDS b) KFDS+ c) KLW d) TVD-KFDS \& e) TVD-KFDS+ schemes} 
\label{L2_Norm} 
\end{center} 
\end{figure}

 The L1 and L2 norms obtained for  the test case for each of the first order accurate and second order accurate versions of KFDS scheme and are tabulated in table 1 and 2 respectively. The log-log plots comparing the L1 and L2 norm errors with slope 1 and  slope 2 are shown in  figures [\ref{L1_Norm}] and [\ref{L2_Norm}] respectively. The KFDS schemes exhibit First order accuracy while the KLW scheme exhibits second order accuracy quite consistently. The TVD schemes which combine the positive aspects of KFDS \& KLW schemes exhibit mild loss in accuracy as compared to KLW and are comparable to second order accurate.

\begin{table} 
\scriptsize
\caption{EOC using L1  error norm for KFDS , KLW \& TVD-KFDS scheme with a smooth periodic test case.} 
\begin{tabular}{rrrrrrrrrrrrr}
\hline \\ GRID SIZE & \multicolumn{2}{c}{\bf KFDS } & \multicolumn{2}{c}{\bf KFDS $+$} & \multicolumn{2}{c}{ \bf KLW } & \multicolumn{2}{c}{\bf TVD-KFDS } & \multicolumn{2}{c}{\bf TVD-KFDS+ } \\
\multicolumn{1}{c}{$\Delta \mathbf{x}$} & ERROR & ORDER & ERROR & ORDER & ERROR & ORDER & ERROR & ORDER & ERROR & ORDER \\ \\
\hline \\ $0.05$ & $0.07672597$ & & $0.04876295$ & & $0.01832800$ & & $0.02733510$ & & $0.02678619$ & \\
$0.025$ & $0.04547810$ & $0.755$ & $0.02601692$ & $0.906$ & $0.00610302$ & $1.586$ & $0.00865670$ & $1.659$ & $0.00849075$ & $1.658$ \\
$0.0125$ & $0.02486742$ & $0.871$ & $0.01562377$ & $0.736$ & $0.00184254$ & $1.728$ & $0.00225355$ & $1.942$ & $0.00237682$ & $1.837$ \\
$0.00625$ & $0.01303956$ & $0.931$ & $0.00908054$ & $0.783$ & $0.00050626$ & $1.864$ & $0.00059727$ & $1.916$ & $0.00061311$ & $1.955$ \\
$0.003125$ & $0.00667554$ & $0.966$ & $0.00501462$ & $0.857$ & $0.00013075$ & $1.953$ & $0.00016656$ & $1.842$ & $0.00016857$ & $1.863$ \\
$0.0015625$ & $0.00338401$ & $0.980$ & $0.00265379$ & $0.918$ & $0.00003308$ & $1.983$ & $0.00005211$ & $1.677$ & $0.00004641$ & $1.861$ \\
$0.00078125$ & $0.00170426$ & $0.990$ & $0.00136685$ & $0.957$ & $0.00000830$ & $1.994$ & $0.00001568$ & $1.732$ & $0.00001221$ & $1.927$ \\
$0.000390625$ & $0.00085504$ & $0.995$ & $0.00069370$ & $0.978$ & $0.00000208$ & $1.999$ & $0.00000454$ & $1.790$ & $0.00000379$ & $1.687$ \\
$0.000195313$ & $0.00042833$ & $0.997$ & $0.00034948$ & $0.989$ & $0.00000052$ & $1.999$ & $0.00000125$ & $1.860$ & $0.00000125$ & $1.600$ \\ \hline \\
\end{tabular}
\caption{EOC using L2  error norm for KFDS , KLW \& TVD-KFDS scheme with a smooth periodic test case.} 
\begin{tabular}{rrrrrrrrrrrrr}
\hline \\ GRID SIZE & \multicolumn{2}{c}{\bf KFDS } & \multicolumn{2}{c}{\bf KFDS $+$} & \multicolumn{2}{c}{ \bf KLW } & \multicolumn{2}{c}{\bf TVD-KFDS } & \multicolumn{2}{c}{\bf TVD-KFDS+ } \\
\multicolumn{1}{c}{$\Delta \mathbf{x}$} & ERROR & ORDER & ERROR & ORDER & ERROR & ORDER & ERROR & ORDER & ERROR & ORDER \\ \\
\hline \\ $0.05$ & $0.12308970$ & & $0.06266486$ & & $0.03394767$ & & $0.04456875$ & & $0.03984737$ & \\
$0.025$ & $0.07466003$ & $0.721$ & $0.03410766$ & $0.878$ & $0.01582174$ & $1.101$ & $0.01379576$ & $1.692$ & $0.01341841$ & $1.570$ \\
$0.0125$ & $0.04161332$ & $0.843$ & $0.02127948$ & $0.681$ & $0.00543107$ & $1.543$ & $0.00378013$ & $1.868$ & $0.00394625$ & $1.766$ \\
$0.00625$ & $0.02222377$ & $0.905$ & $0.01448966$ & $0.554$ & $0.00152566$ & $1.832$ & $0.00112373$ & $1.750$ & $0.00114227$ & $1.789$ \\
$0.003125$ & $0.01151868$ & $0.948$ & $0.00897774$ & $0.691$ & $0.00039441$ & $1.952$ & $0.00035306$ & $1.670$ & $0.00034483$ & $1.728$ \\
$0.0015625$ & $0.00589079$ & $0.967$ & $0.00504854$ & $0.830$ & $0.00010002$ & $1.979$ & $0.00011610$ & $1.605$ & $0.00009851$ & $1.808$ \\
$0.00078125$ & $0.00298164$ & $0.982$ & $0.00267565$ & $0.916$ & $0.00002513$ & $1.993$ & $0.00003360$ & $1.789$ & $0.00002791$ & $1.819$ \\
$0.000390625$ & $0.00149956$ & $0.992$ & $0.00137548$ & $0.960$ & $0.00000629$ & $1.999$ & $0.00000961$ & $1.805$ & $0.00000756$ & $1.884$ \\
$0.000195313$ & $0.00075224$ & $0.995$ & $0.00069729$ & $0.980$ & $0.00000157$ & $1.999$ & $0.00000266$ & $1.854$ & $0.00000210$ & $1.849$ \\ \hline
\end{tabular}

\end{table}

\section{Results and Discussion} 
    First order accurate and second order accurate versions of both KFDS and KFDS+ schemes, together with KLW scheme, are tested on the following bench-mark test cases to test their capabilities in resolving the features like shock waves and expansion fans. Test case 1 is to chosen to test the scheme in the framework of a linear convection equation. Test case 2 is selected to study the behavior of the scheme in the framework of a linear convection-diffusion equation.  Test cases 3 to 5 have been been formulated for inviscid Burgers equation and test cases 6 to 8 have been formulated for viscous Burgers equation. Test case 9 explores the ability of the numerical schemes in  a system of equations such as shallow water equations with a source term. Test case 10 assesses the numerical scheme in two dimensions for linear convection equation whereas test cases 11 and 12 are chosen to check the numerical scheme in a 2D inviscid Burgers equation framework. Test cases 13 and 14 are chosen to evaluate the scheme for two dimensional viscous Burgers equation.  Test case 15 is chosen to test the numerical scheme in a two dimensional shallow water system.  These test cases will involve Riemann problems, which refer to the solution of nonlinear hyperbolic equations with initial discontinuities.  These solutions involve shock wave and expansion wave solutions, the analytical solutions of which, obtained using method of characteristics (see \cite{Leveque_1}), are given here.  For a non-linear convection equation as      
\be 
\fr{\del u}{\del t} + \fr{\del g(u)}{\del x} = 0, \ g(u) = \fr{1}{2} u^{2} 
\ee
with the initial condition (a discontinuity) given by  
\be 
u(x,t=0) = \left\{ \ba{c} u_{L} \ \textrm{for} \ x < 0 \\ u_{R} \ \textrm{for} \ x \ge 0 \ea \right. 
\ee 
the two possible solutions for this Riemann problem are as follows.  
\subsubsection*{For $u_{L} > u_{R}$: (shock wave)}      
\be 
u(x,t) = \left\{ \ba{c} u_{L} \ \textrm{for} \ x < s t \\ 
u_{R} \ \textrm{for} \ x \ge s t \ea \right. 
\ee
where 
\be 
s = \fr{u_{L} + u_{R}}{2}
\ee
is the shock speed. 
 
\subsubsection*{For $u_{L} < u_{R}$: (expansion wave)} 
\be 
u(x,t) = \left\{\ba{l} 
u_{L} \ \ \ \ \textrm{for} \ x < a(u_{L}) t \\[2mm]  
a^{-1}\left(\fr{x}{t}\right) \ \ \ \ \ \textrm{for} \ a(u_{L}) t \le x \le a(u_{R}) t \\[2mm]  
u_{R} \ \ \ \ \textrm{for} \ x > a(u_{R}) t 
\ea \right. 
\ee 
The results with the first order KFDS and KFDS+ schemes, second order KLW scheme and TVD high resolution KFDS and KFDS+ schemes are presented for each of the following test cases.  

\subsubsection*{Test case 1: Linear convection of an initial discontinuity}
Consider a one-dimensional linear convection equation as shown below.  
\be 
\fr{\del u}{\del t} + \fr{\del g(u)}{\del x} = 0, \ g(u) =  u
\ee
The initial conditions in the domain is given as 
\bea
u(x,0) = \left\{ \ba{l} 0 \ \textrm{for} \ x < -\fr{1}{3} \\ 
1 \ \textrm{for} \ -\fr{1}{3} \le x \le \fr{1}{3} \\    
-1 \ \textrm{for} \ x > \fr{1}{3} \ea \right. 
\eea
 The convection of the above mentioned discontinuities is to be obtained at t = 0.3 in a domain $[-1,1]$. The exact solution for the test case is given below.  
\bea
u(x,0) = \left\{ \ba{l} 0 \ \textrm{for} \ x < -\fr{1}{3} + t \\ 
1 \ \textrm{for} \ -\fr{1}{3} \le x \le \fr{1}{3} + t\\    
-1 \ \textrm{for} \ x > \fr{1}{3} + t \ea \right. 
\eea

\begin{figure}[!h] 
\begin{tabular}{ccc}
\includegraphics[height=3.5cm]{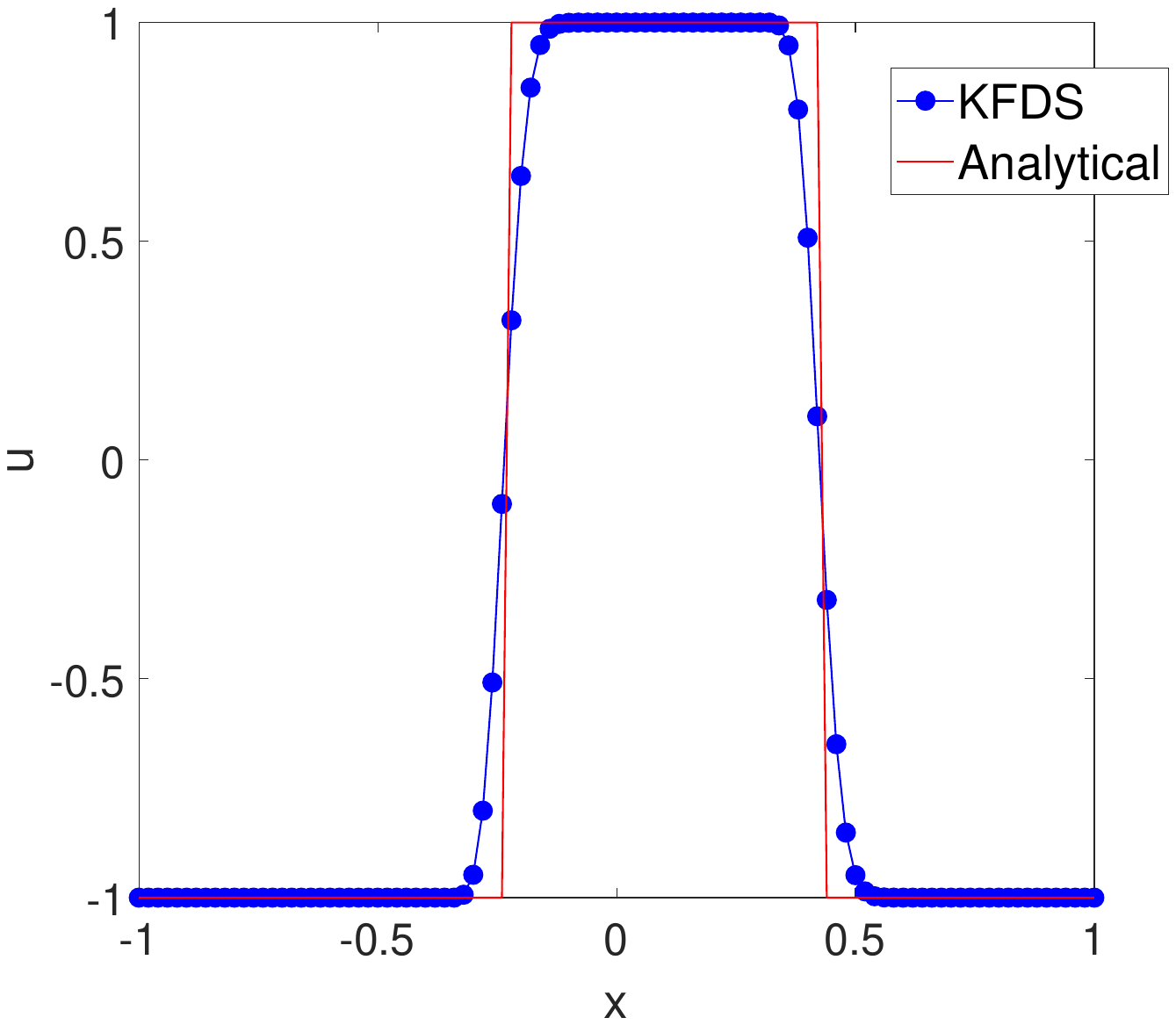} &
\includegraphics[height=3.5cm]{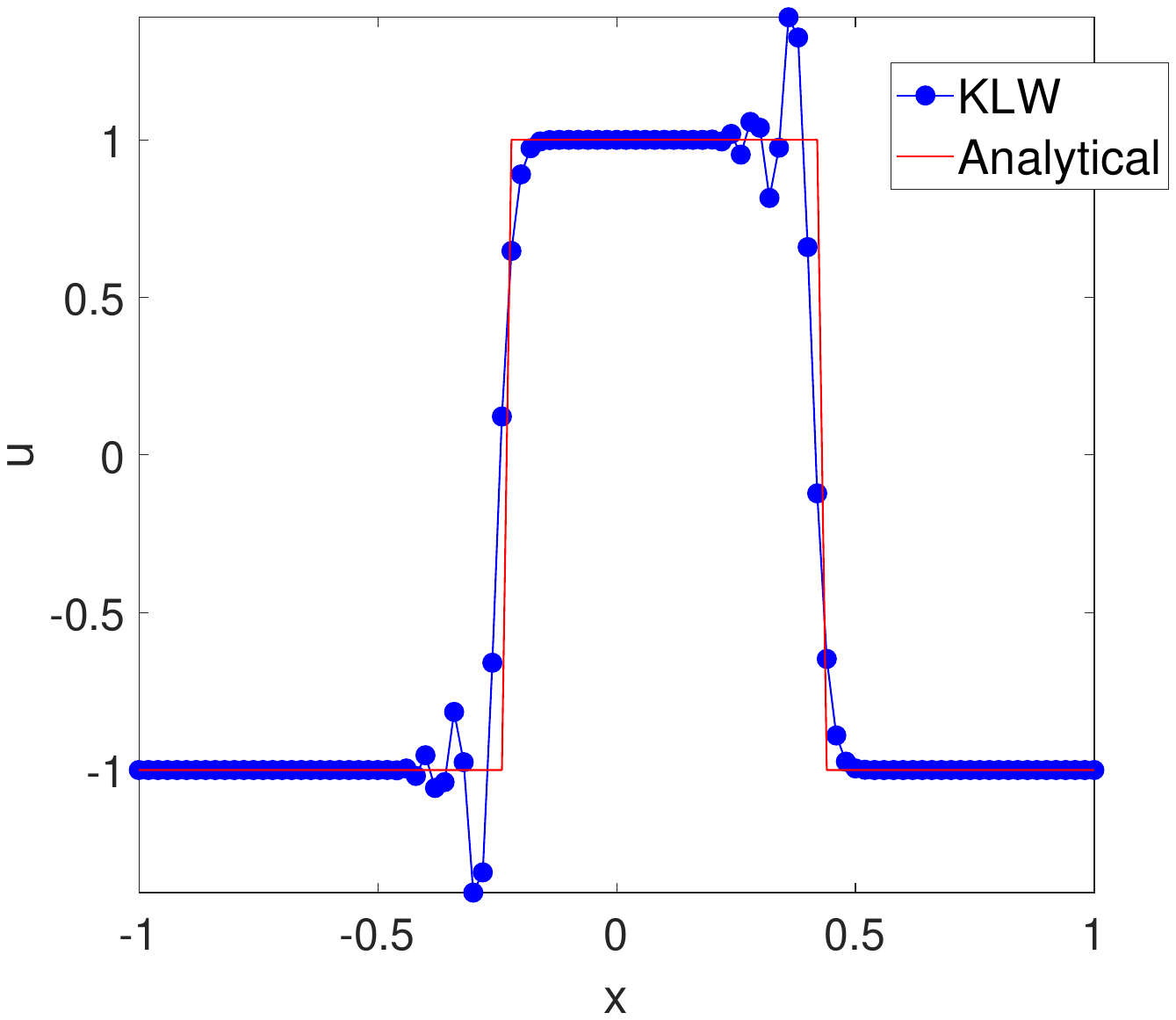} &
\includegraphics[height=3.5cm]{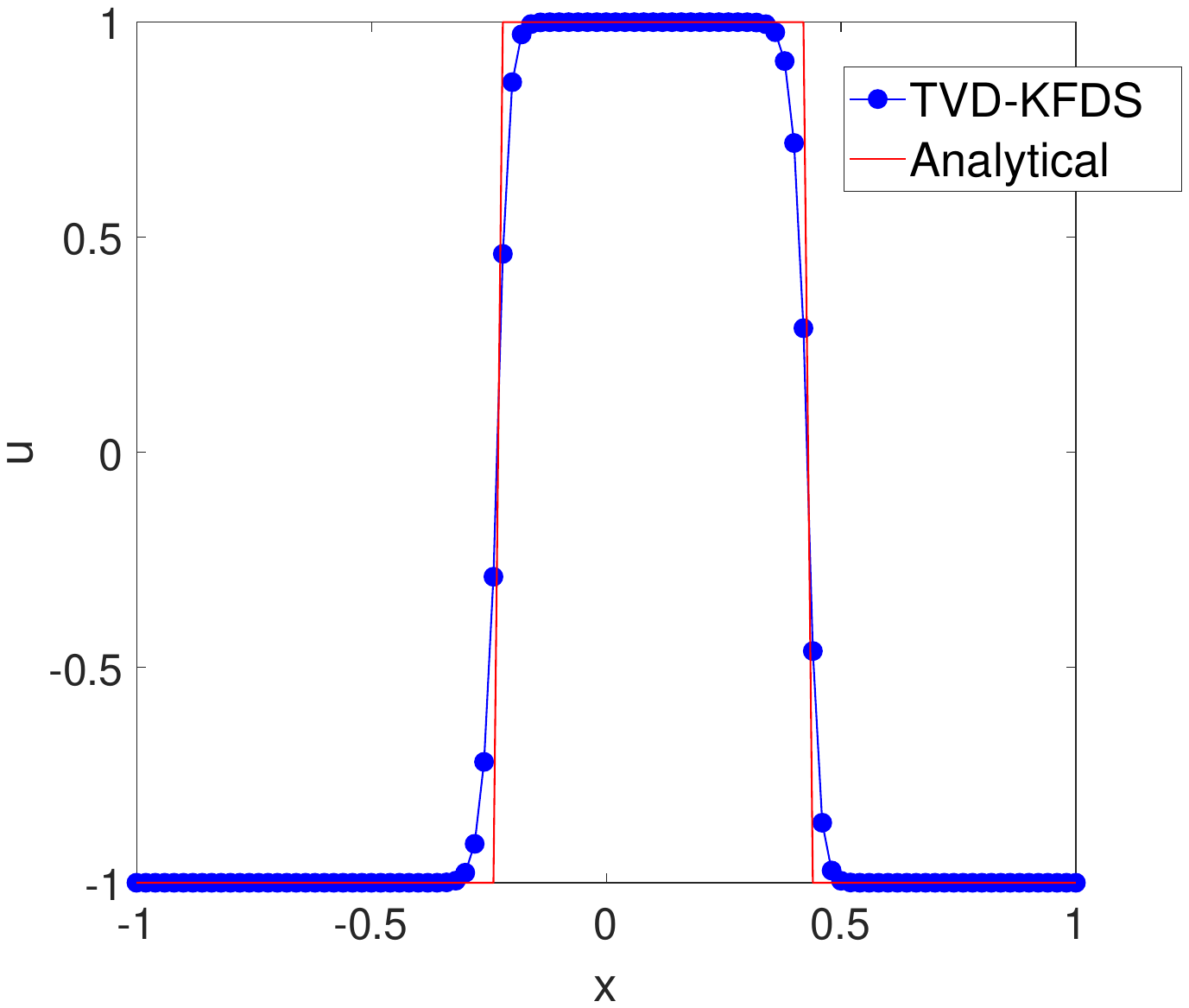} \\
 (a)First order KFDS & (b )Second order KLW & (c) Second order TVD-KFDS
\end{tabular}
\caption{Test Case 1 : KFDS, KLW \& TVD-KFDS schemes for linear convection}
\label{TC_1_1O_KFDS_LCE} 
\end{figure}

The results for this test case are given in figure (\ref{TC_1_1O_KFDS_LCE}).  Ideally, the linear convection equation convects the initial discontinuities in the system with a constant wave speed and the spread of the discontinuities in the domain remains undisturbed. It can be observed that the numerical schemes developed capture the position of the discontinuities with reasonable accuracy. However the schemes also introduce numerical diffusion by virtue of their  construction.  Both the First Order KFDS and TVD-KFDS schemes capture the discontinuities at the right positions without any undue oscillations.  However, the second order KLW scheme produces oscillations near both the discontinuities, as expected.     

\subsubsection*{Test case 2: Linear Convection-Diffusion Equation}    
Consider a one-dimensional linear convection equation with a viscous term as 
\be
\fr{\del u}{\del t} + \fr{\del u}{\del x} = Pe^{-1}\fr{\del^{2} u}{\del x^{2}}
\ee
where Pe is the Peclet Number. In this relation, the inverse of Peclet number resembles the viscosity coefficient and therefore demonstrates the effect of Peclet number on the diffusion. The solution is to be evolved to attain steady state condtions for Peclet Numbers 1 and 50. The domain for the test case is $[0,1]$ while  boundary conditions are $u(0) = 0$ and $u(1) =1$.  The exact solution for this test case is  
\be
{u}({x})=\frac{1-e^{xPe}}{1-e^{Pe}}
\ee  

\begin{figure}[!h] 
\begin{center}  
\begin{tabular}{ccc}
\includegraphics[height=3.5cm]{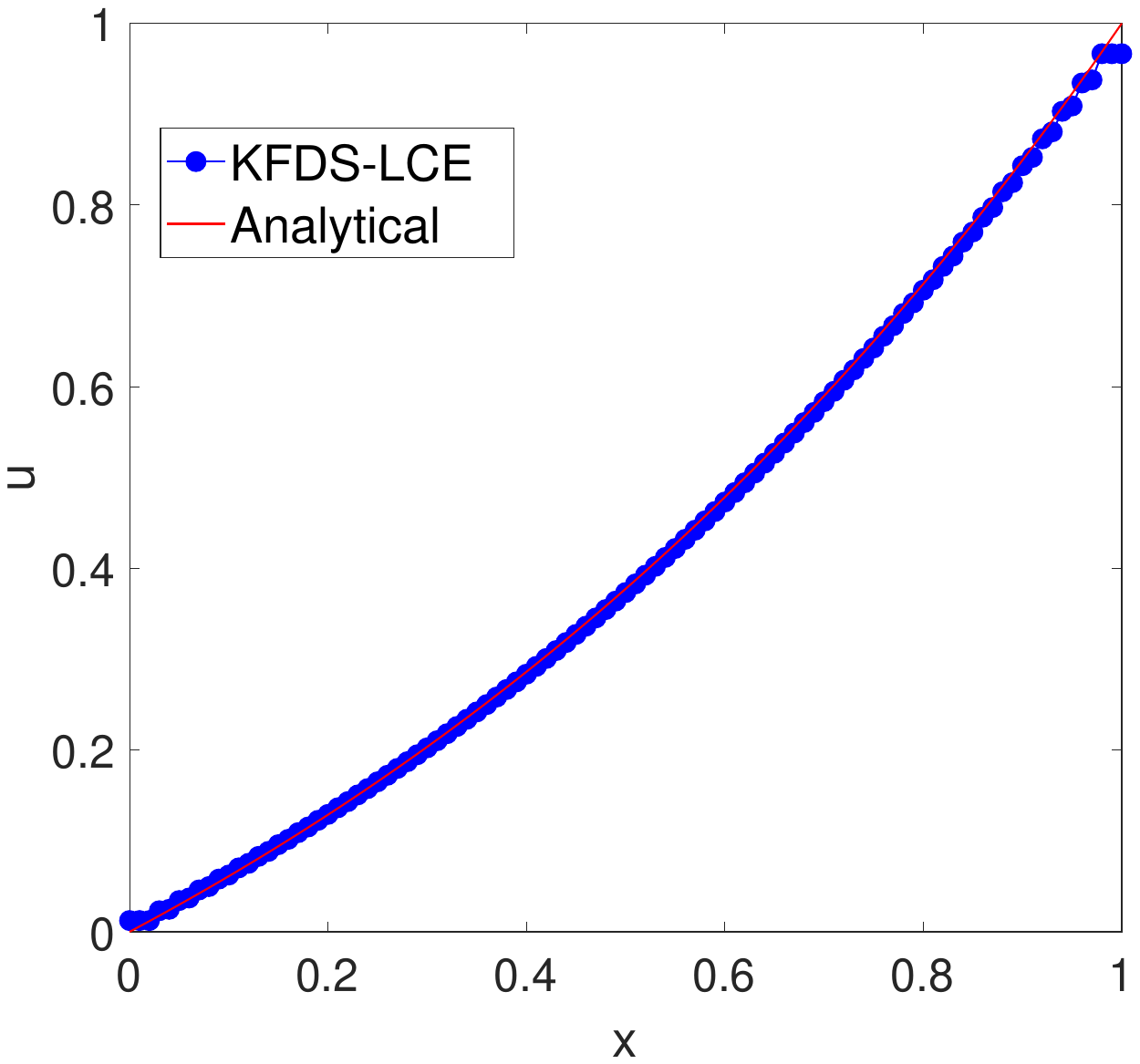} &
\includegraphics[height=3.5cm]{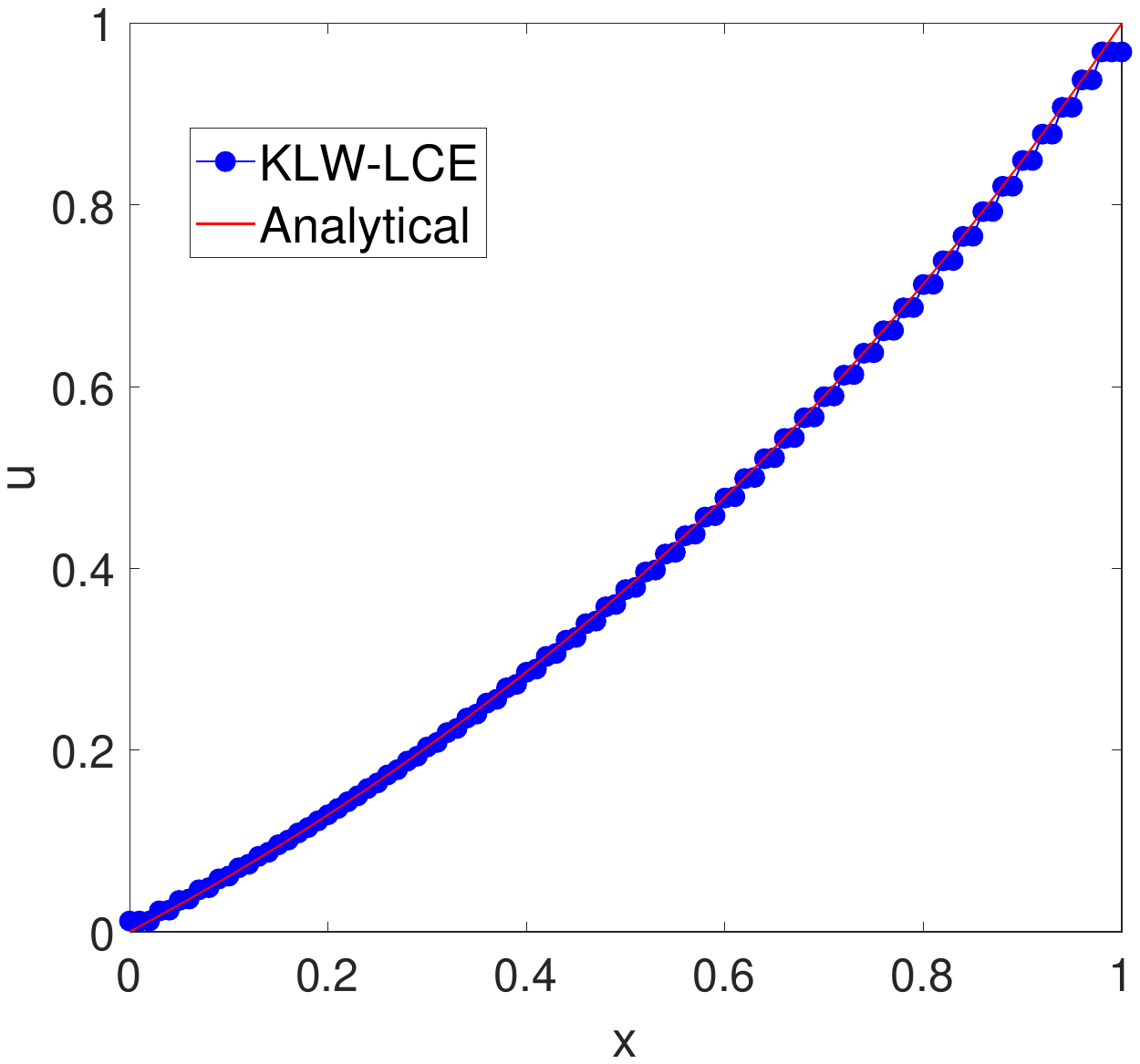} &
\includegraphics[height=3.5cm]{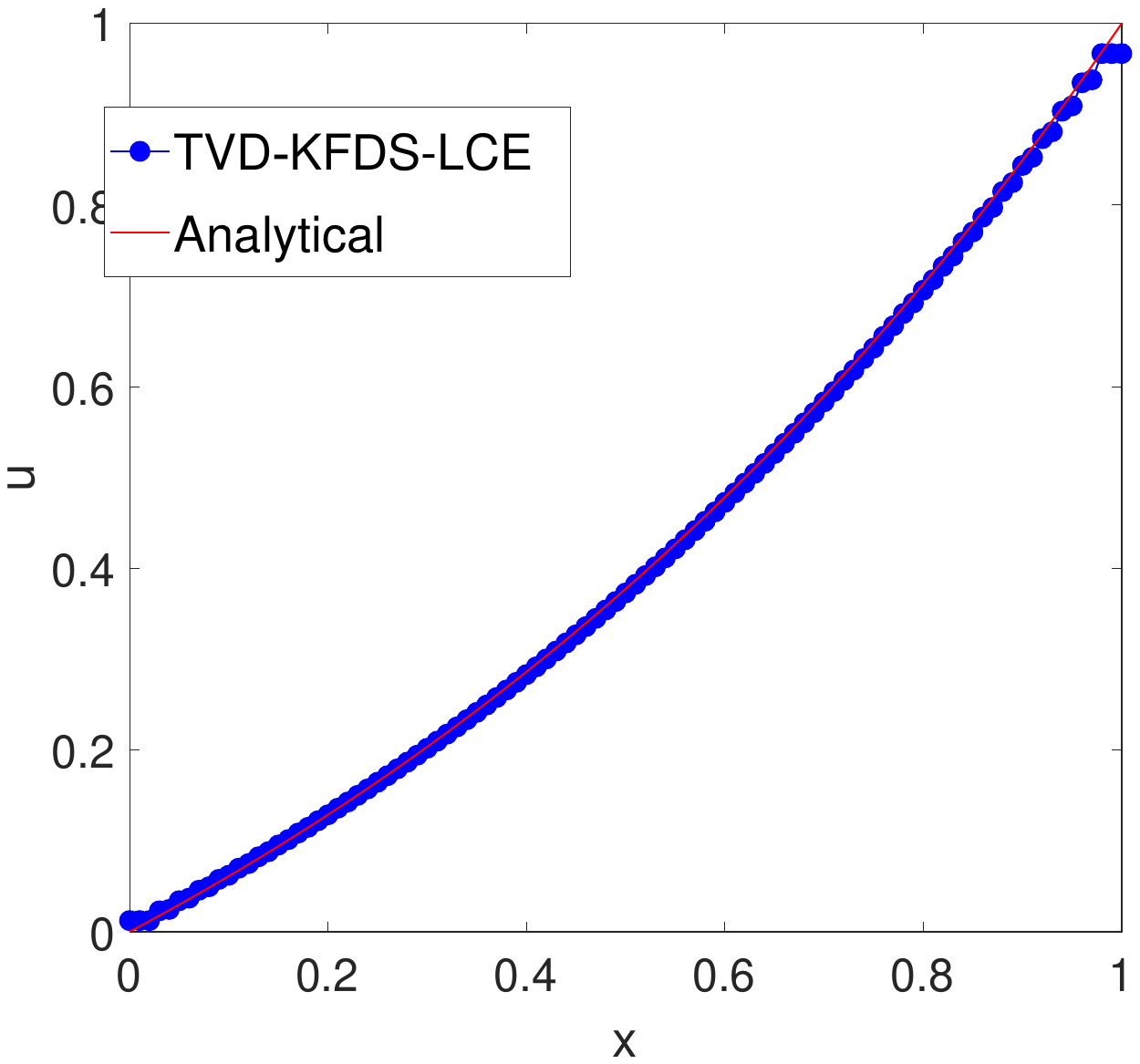} \\
 (a)First order KFDS & (b)Second order KLW & (c)Second order TVD-KFDS
\end{tabular}
\caption{Test Case 2(a) : KFDS, KLW \& TVD schemes in LCE framework with Peclet Number = 1}
\label{TC_BL_LCE_KFDS_P1} 
\end{center} 
\end{figure}

\begin{figure}[!h] 
\begin{center} 
\begin{tabular}{ccc}
\includegraphics[height=3.5cm]{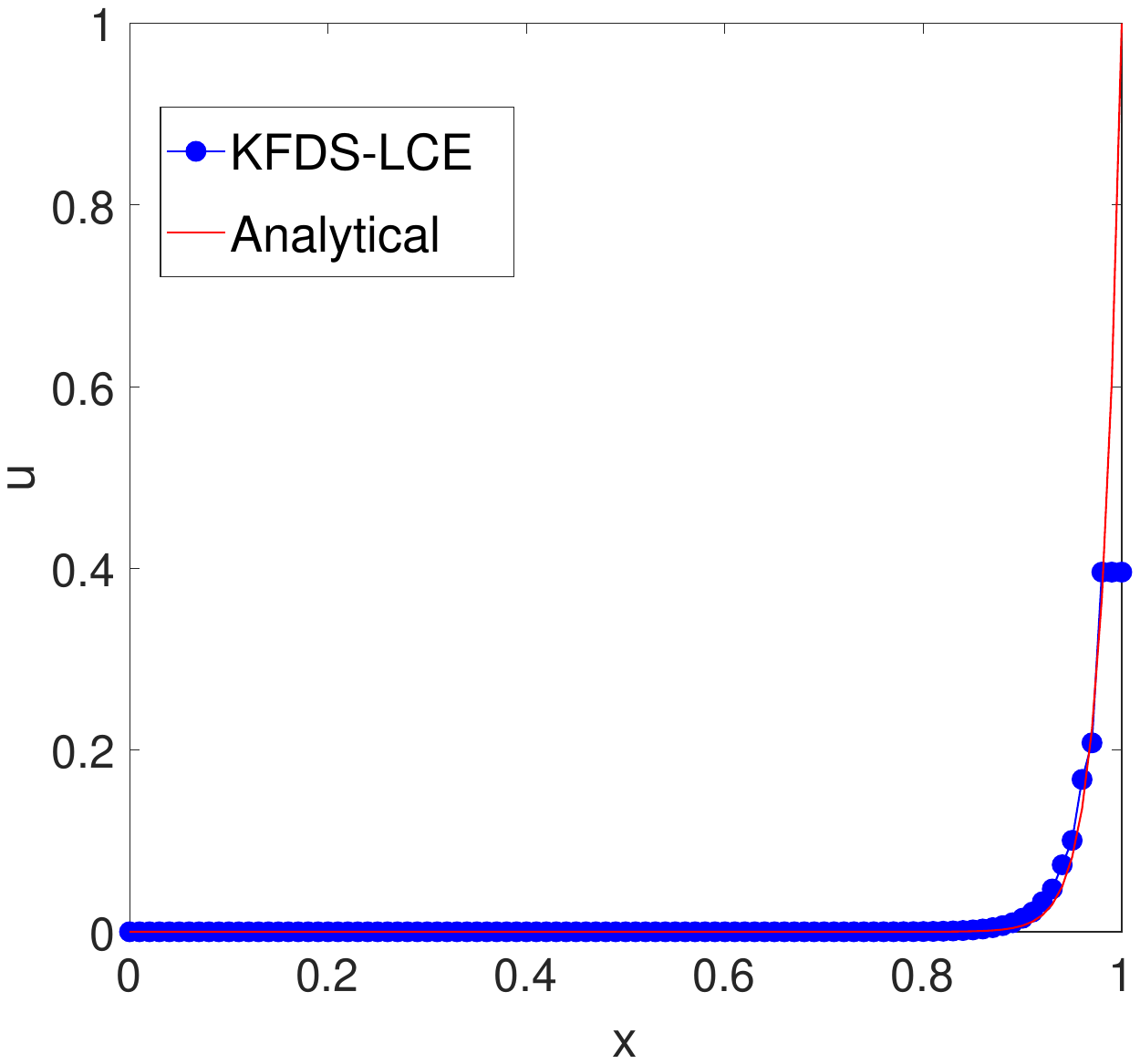} &
\includegraphics[height=3.5cm]{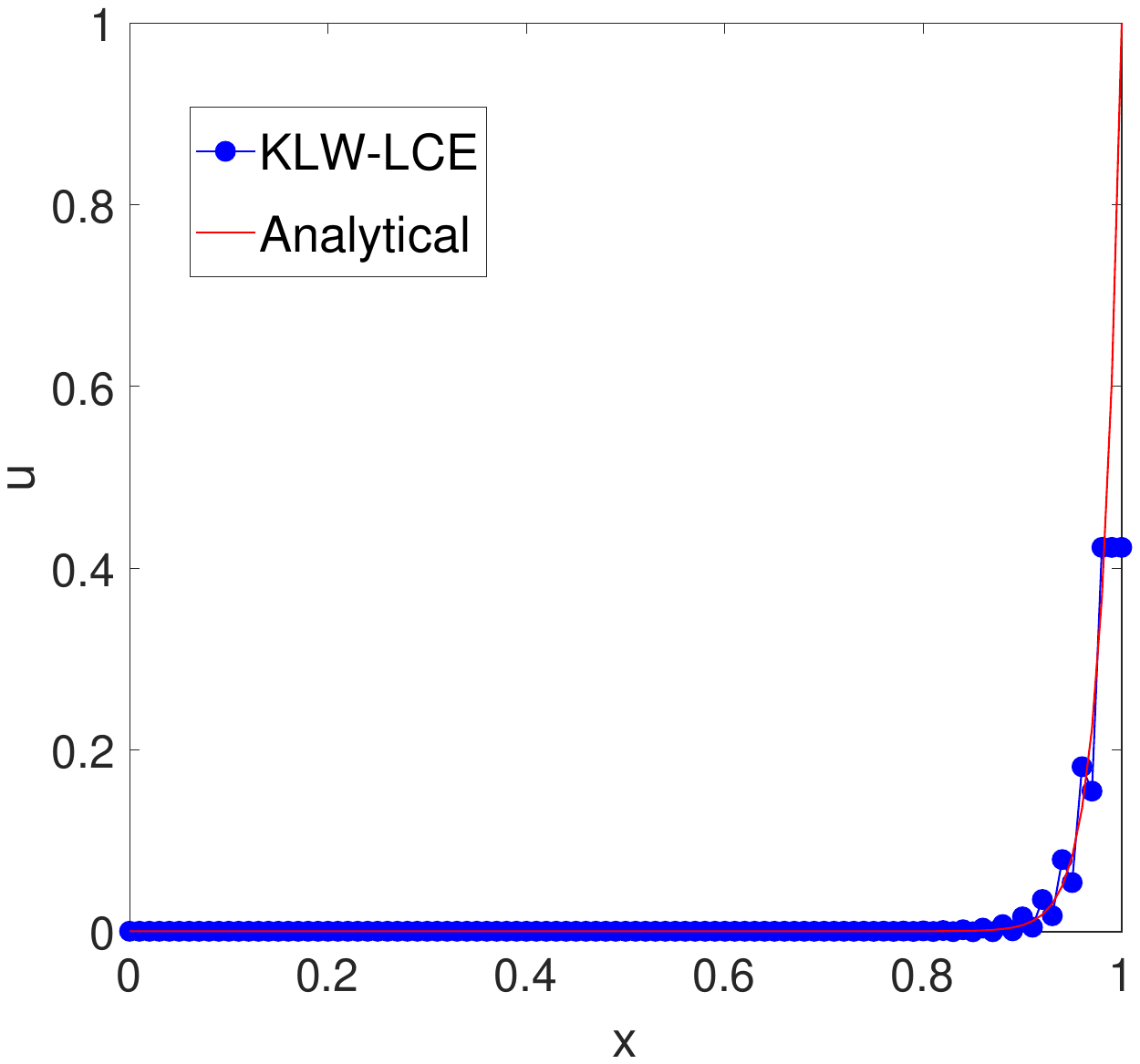} &
\includegraphics[height=3.5cm]{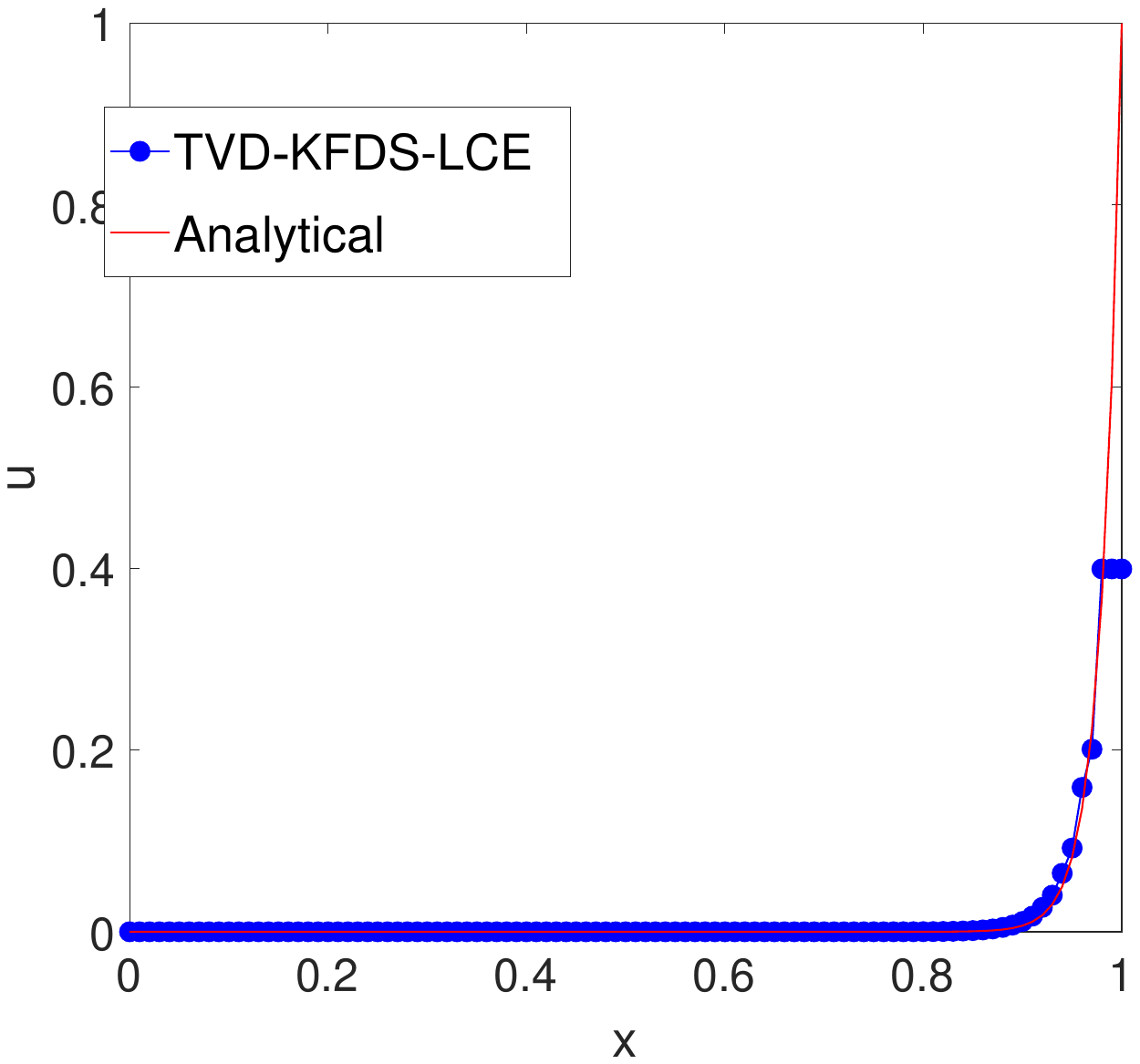} \\
 (a)First order KFDS & (b)Second order KLW & (c)Second order TVD-KFDS
\end{tabular}
\caption{Test Case 2(b) : KFDS, KLW \& TVD schemes in LCE framework with Peclet Number = 50}
\label{TC_BL_LCE_KFDS_P50} 
\end{center} 
\end{figure}

The results for this test case are given in figures (\ref{TC_BL_LCE_KFDS_P1}) and (\ref{TC_BL_LCE_KFDS_P50}). This test case mimics the boundary layer phenomena associated with the viscous flows of the fluids in one dimension.  Peclet number is the ratio of the convection terms to the diffusion terms.  Higher the Peclet number, greater will be domination of convection in the flow in comparison to diffusion.  In this test case, at Peclet Number = 50, we see a large gradient near $x = 1$. Each of the schemes captures the steep gradient with moderate accuracy.  The accuracy increases significantly with an increase in the number of grid points in the gradient region. At Peclet number =1, the convection and diffusion terms are equally dominant and therefore results in a near linear transition of the flow without any steep gradients. Thus, the schemes demonstrate capability to work in both convection dominant and convection-diffusion dominant cases. 

\subsubsection*{Test case 3: Steady shock wave and evolving expansion wave for inviscid Burgers equation}
\bea
u(x,0) = \left\{ \ba{l} 0 \ \textrm{for} \ x < -\fr{1}{3} \\ 
1 \ \textrm{for} \ -\fr{1}{3} \le x \le \fr{1}{3} \\    
-1 \ \textrm{for} \ x > \fr{1}{3} \ea \right. 
\eea 
With these initial conditions, the jump at $-\fr{1}{3}$ will evolve into an expansion fan and the jump at 
$\fr{1}{3}$ will create a steady shock wave (as shock speed is zero from the R-H condition, for the given left and right states here).  The solution for this problem is sought at time $t=0.3$, well before any possible interaction between the shock and expansion waves and well before the waves reach the boundaries.  

\begin{figure}[!h] 
\begin{center} 
\begin{tabular}{ccc}
\includegraphics[height=3.5cm]{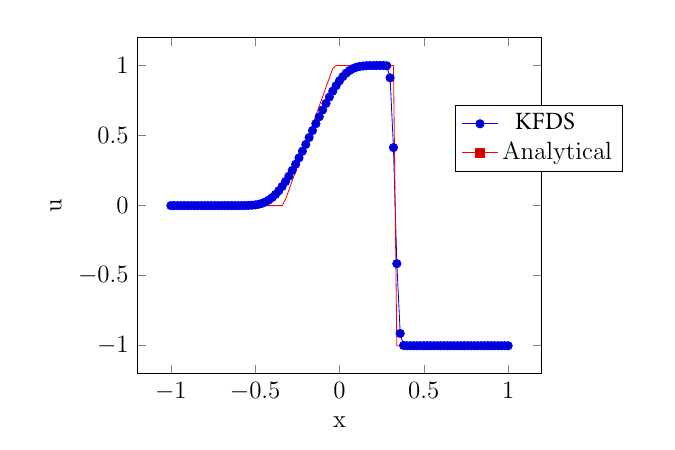} &
\includegraphics[height=3.5cm]{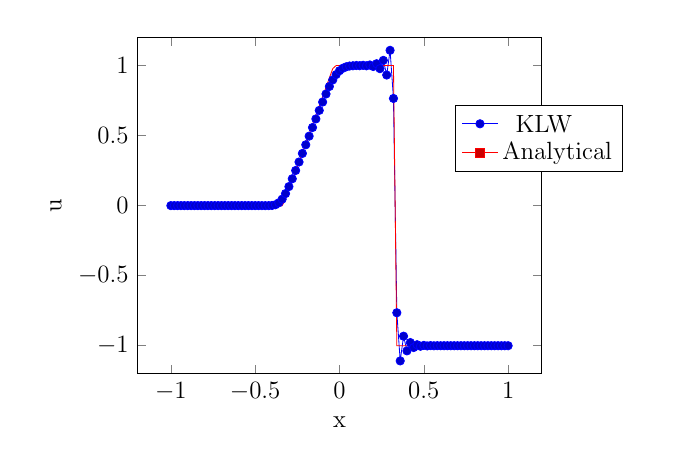} &
\includegraphics[height=3.5cm]{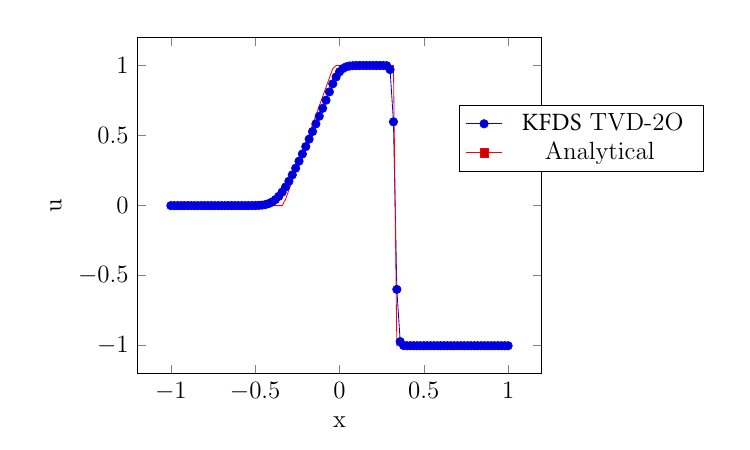} \\
 (a)First order KFDS & (b)Second order KLW & (c)TVD-KFDS\\
\includegraphics[height=3.5cm]{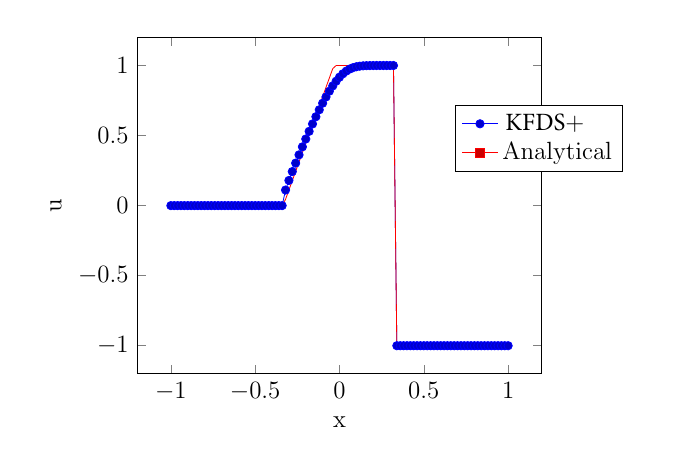} & &
\includegraphics[height=3.5cm]{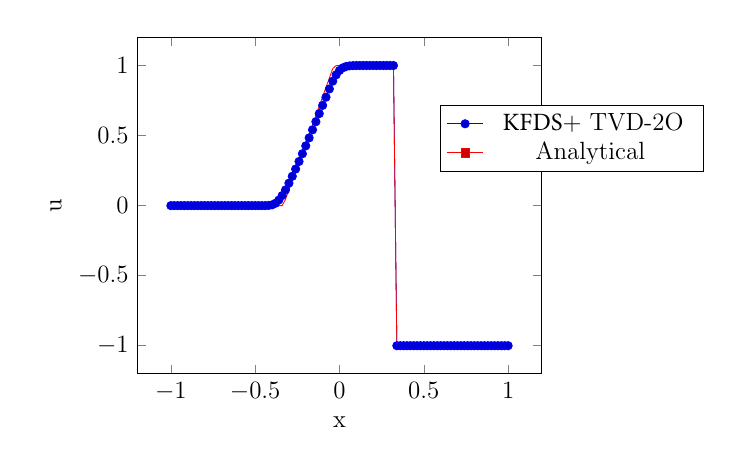} \\
 (d)First order KFDS+ &  & (e)TVD-KFDS+ \\
\end{tabular}
\caption{Test Case 3 : KFDS, KLW, TVD-KFDS, KFDS+ \& TVD-KFDS+  schemes in Burgers framework with 100 points}
\label{TC_1_KFDS} 
\end{center} 
\end{figure} 

The results for this test case are given in figures (\ref{TC_1_KFDS}).  For Test Case 3, as the shock is steady, first order accurate {\em KFDS} scheme captures the shock wave with some numerical diffusion while the first order accurate {\em KFDS+} scheme captures the shock exactly by design, without any numerical diffusion.  Both schemes capture the evolving expansion fan reasonably, with some numerical diffusion.  The second order KLW scheme is more accurate in capturing both waves but generates oscillations near the shock wave.  The high resolution TVD versions of both schemes capture shock waves without any oscillations, with KFDS+ scheme retaining its exact capturing feature. Both schemes capture the expansion fan more accurately the the first order versions.     

\subsubsection*{Test case 4: Unsteady shock wave and evolving expansion wave} 
  In this test case, the initial conditions given below are chosen such that an unsteady shock wave and an evolving expansion wave arise in the solution.   
\bea
u(x,0) = \left\{ \ba{l} 
0 \ \textrm{for} \ x < - \fr{1}{3} \\ 
1 \ \textrm{for} \ \fr{1}{3} \le x \le \fr{1}{3} \\ 
0 \ \textrm{for} \ x > \fr{1}{3} 
\ea \right. 
\eea 
These initial conditions define a square wave, with a jump from $0$ to $1$ at 
$-\fr{1}{3}$ and another jump from $1$ to $0$ at $\fr{1}{3}$.  The jump at 
$-\fr{1}{3}$ will evolve into an expansion fan and the jump at $\fr{1}{3}$ will be a moving shock wave.  The solution for this unsteady problem is sought at time $t=0.3$, well before any possible interaction between the shock and expansion waves and well before the waves reach the boundaries. 

\begin{figure}[!h] 
\begin{center} 
\begin{tabular}{ccc}
\includegraphics[height=3.5cm]{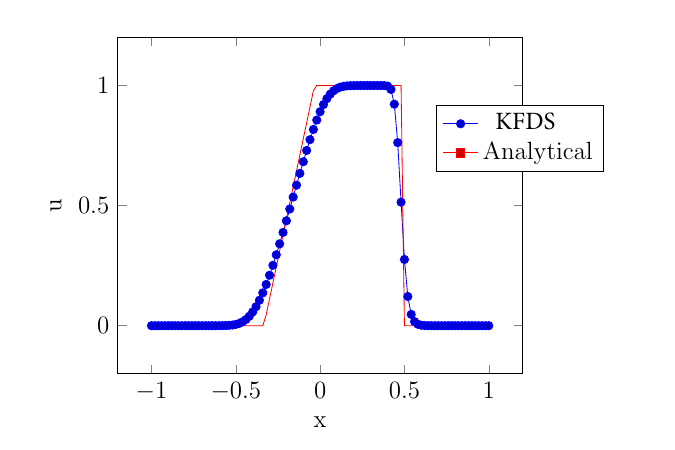} &
\includegraphics[height=3.5cm]{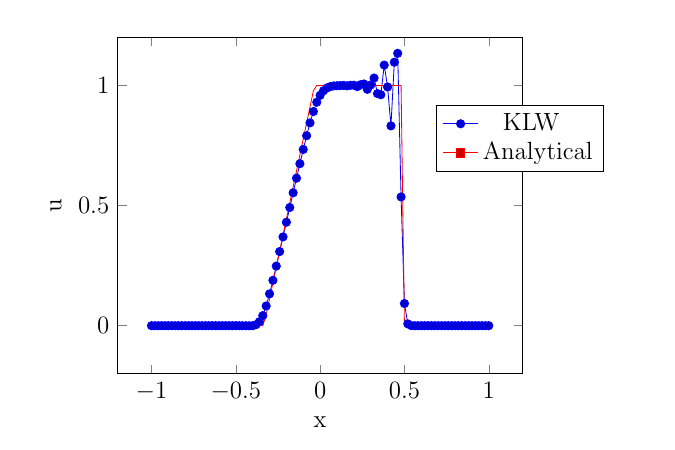} &
\includegraphics[height=3.5cm]{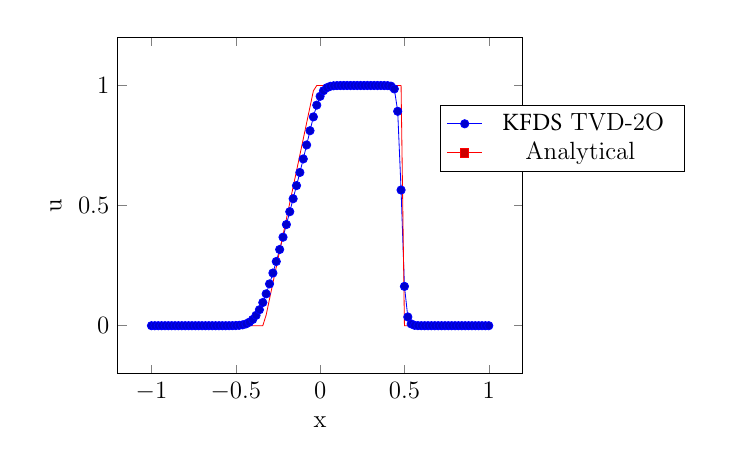} \\
 (a)First order KFDS & (b)Second order KLW & (c)TVD-KFDS\\
\includegraphics[height=3.5cm]{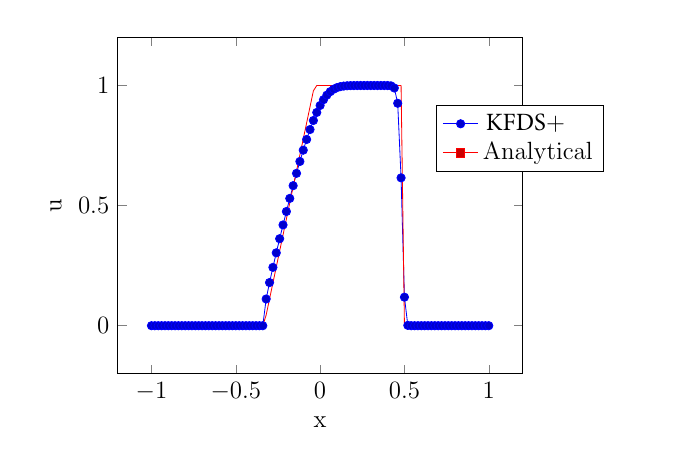} & &
\includegraphics[height=3.5cm]{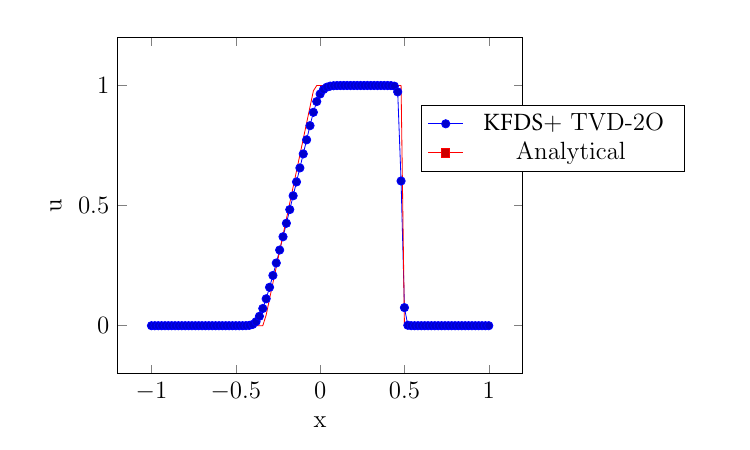} \\
 (d)First order KFDS+ &  & (e)TVD-KFDS+ \\
\end{tabular}
\caption{Test Case 4 : KFDS, KLW, TVD-KFDS, KFDS+ \& TVD-KFDS+  schemes in Burgers framework with 100 points}
\label{TC_2_1O_KFDS} 
\end{center} 
\end{figure}

The results for this test case are given in figures (\ref{TC_2_1O_KFDS}).  The unsteady shock wave and the evolving expansion fan are captured with some numerical diffusion by both first order accurate {\em KFDS} and {\em KFDS+} schemes.  However, {\em KFDS+} scheme captures the shock with less number of points compared to the {\em KFDS} scheme.  The second order accurate {\em KLW} scheme captures both the waves more accurately but introduces oscillations near the shock.  The high resolution TVD versions of both {\em KFDS} and {\em KFDS+} schemes avoid these oscillations, while capturing both waves accurately.  

\subsubsection*{Test case 5: Steady shock wave and evolving expansion wave with a sonic point}
\bea
u(x,0) = \left\{ \ba{l} -1 \ \textrm{for} \ x < -\fr{1}{3} \\ 
1 \ \textrm{for} \ -\fr{1}{3} \le x \le \fr{1}{3} \\    
-1 \ \textrm{for} \ x > \fr{1}{3} \ea \right. 
\eea 
With these initial conditions, the jump at $-\fr{1}{3}$ will evolve into an expansion fan, which contains a sonic point where the wave speed changes sign. The jump at $\fr{1}{3}$ will create a steady shock wave (as shock speed is zero for the given left and right states here).  The solution for this problem is sought at time $t=0.3$, well before any possible interaction between the shock and expansion waves and well before the waves reach the boundaries.

\begin{figure}[!h] 
\begin{center} 
\begin{tabular}{ccc}
\includegraphics[height=3.8cm]{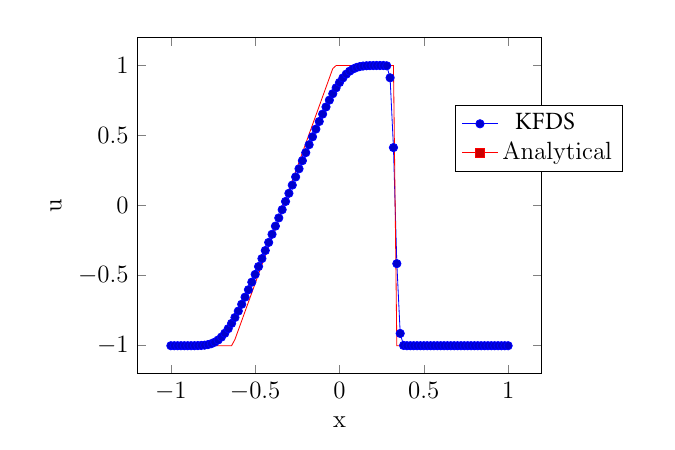} &
\includegraphics[height=3.8cm]{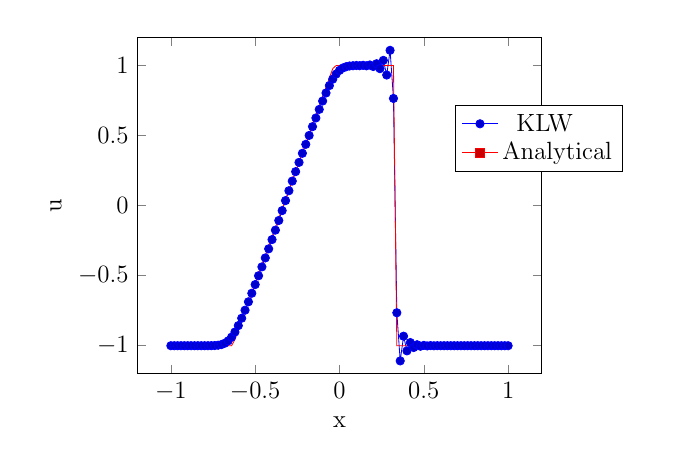} &
\includegraphics[height=3.8cm]{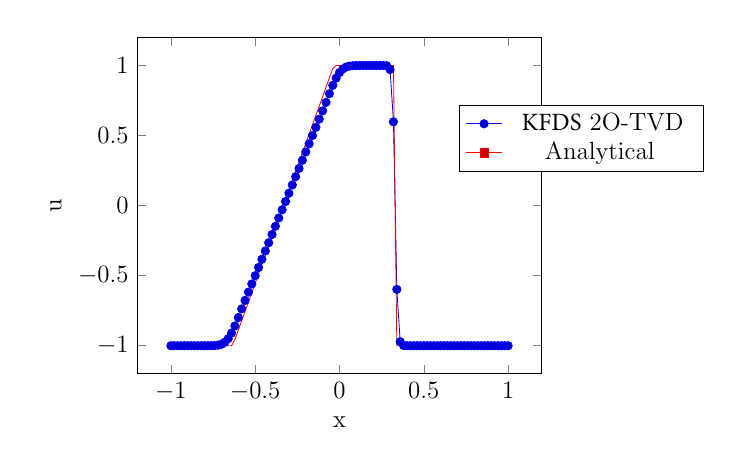} \\
 (a)First order KFDS & (b)Second order KLW & (c)TVD-KFDS\\
\includegraphics[height=3.8cm]{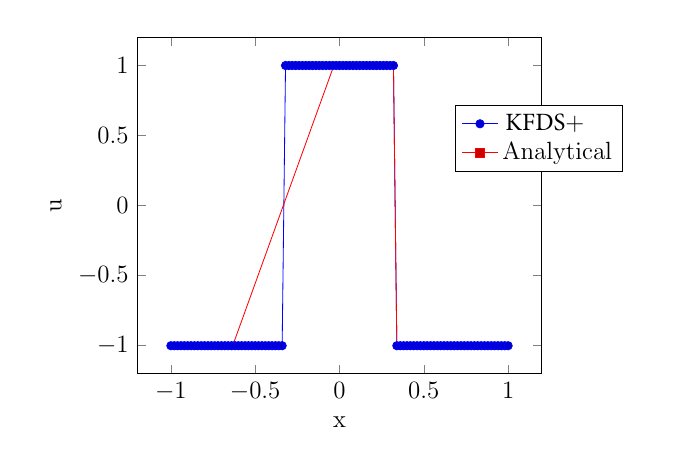} & &
\includegraphics[height=3.8cm]{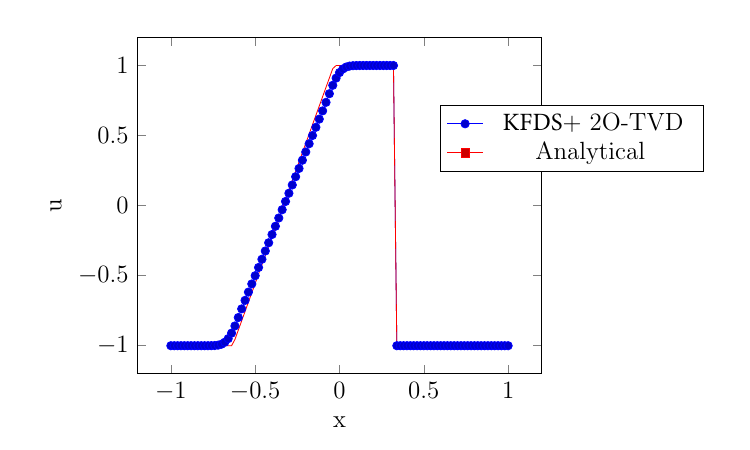} \\
 (d)First order KFDS+ &  & (e)TVD-KFDS+ \\
\end{tabular}
\caption{Test Case 5 : KFDS, KLW, TVD-KFDS, KFDS+ \& TVD-KFDS+  schemes in Burgers framework with 100 points}
\label{TC_3_1O_KFDS}    
\end{center} 
\end{figure} 

Results for this test case are given in figure (\ref{TC_3_1O_KFDS}).  In Test Case 5, as the shock is steady, {\em KFDS+} scheme again captures the shock exactly, without any numerical diffusion and without any intermediate point in the shock.  {\em KFDS} scheme, however, captures the shock with some numerical diffusion.  While the {\em KFDS} scheme captures the expansion fan with a sonic point well, the {\em KFDS+} scheme fails to alter the initial discontinuity and results in an expansion shock.  This is because the numerical diffusion in KFDS+ scheme vanishes at sonic point.  Traditionally, in such cases, an entropy fix is used.  However, we introduce a different technique based on the second order accurate results.  We can see that second order accurate {\em KLW} and {\em KFDS} schemes capture the expansion wave more accurately and have no problem at the sonic point.  Therefore, we introduce the strategy of using Chapman-Enskog condition based $\lambda$ in the expansion fan and R-H condition based $\lambda$ in the shock wave, in the high resolution TVD versions of {\em KFDS+} scheme.  The switching between the two waves is done by recognizing a discontinuity based on information of characteristics.  If $a(u)_{j} > 0$ and $a(u)_{j+1} < 0$ at a cell-interface 
$x_{j+\frac{1}{2}}$, then once can expect a shock wave as the characteristics are converging from left and right sides.  With this hybrid of {\em KFDS} and 
{\em KFDS+} schemes, the steady shock wave is captured exactly and yet there is no trouble at the sonic point in the expansion fan, which is captured accurately.       
			
In general, {\em KFDS+} scheme is less diffusive than {\em KFDS} scheme.  This reduction is numerical diffusion in general, as well as the exact capturing of the shock wave in steady state, is due to the utilization of the Rankine-Hugoniot jump condition in the construction of the scheme.  

\subsubsection*{Test case 6: A high gradient profile}
  This test case \cite{W.Shen} is aimed at testing the robustness, accuracy and time sensitivity of the numerical scheme. The test case is solved on a domain [-2,3] which is divided into 100 points. This test problem is to be simulated for two different cases: $(a)$ a viscosity coefficient $\nu$ = 0.05 and  time $t = 3s$ and $(b)$ a viscosity coefficient $\nu$ = 0.001 and  time $t = 1s$.  The initial conditions in the domain for $case \ (a)$ are given below.  
\bea
u(x,0) = \left\{ \ba{l}-1 \  \textrm{if} \left( \frac { x } { 4 \nu } \right) < -20\\ 
  1 \   \textrm{if} \left( \frac { x } { 4 \nu } \right)  >  0 \\ 
 \frac { 1 } { 2 } \left( 1 - \frac { \exp \left( \frac { x } { 4 \nu } \right) - \exp \left( \frac { - x } { 4 \nu } \right) } { \exp \left( \frac { x } { 4 \nu } \right) + \exp \left( \frac { - x } { 4 \nu } \right) } \right) \ \textrm{ if } \- 20 \leq \frac { x } { 4 \nu } \leq 20 \ea \right. 
\eea 

\begin{figure}[!h] 
\begin{center} 
\begin{tabular}{ccc}
\includegraphics[height=4cm]{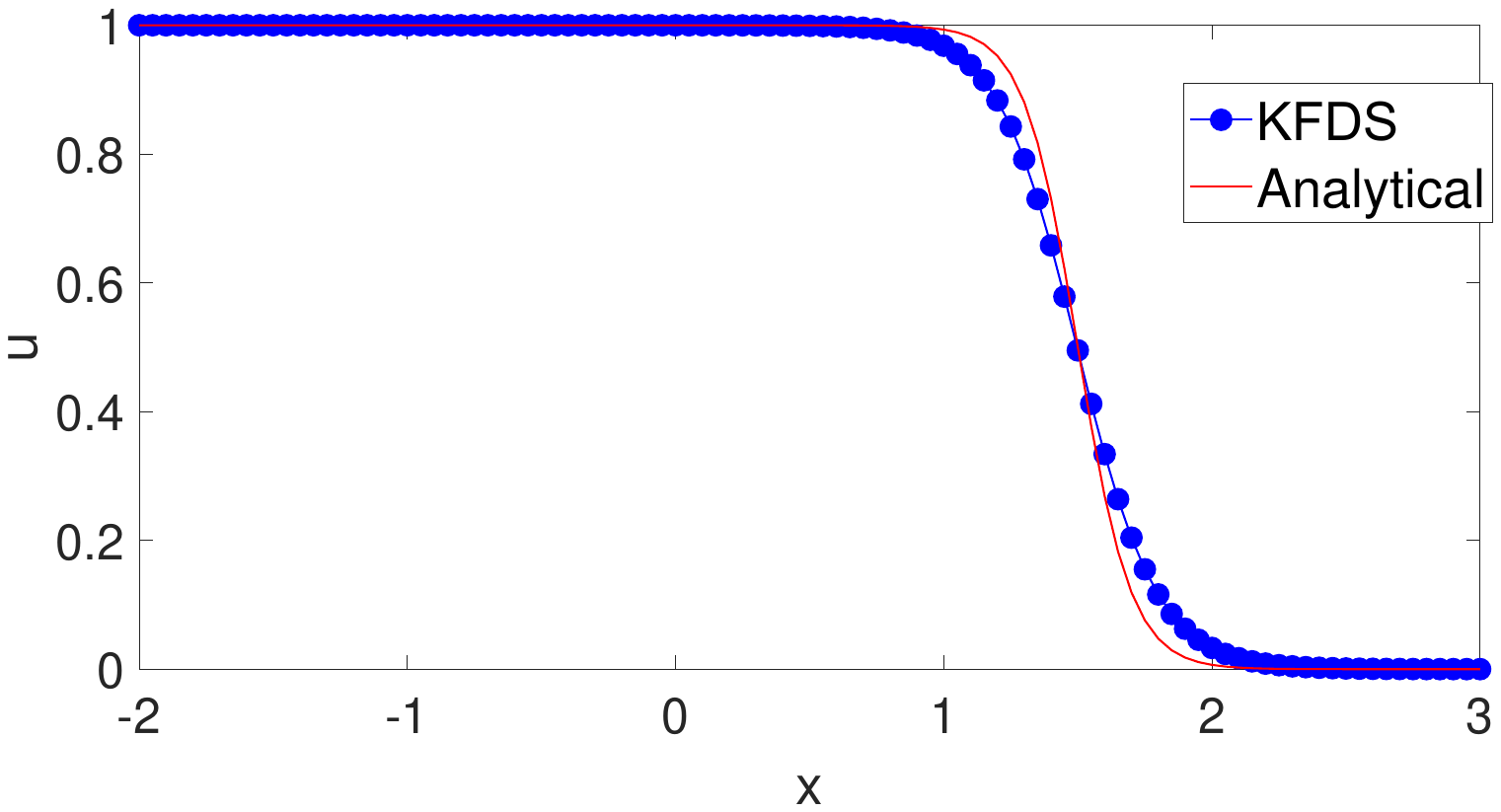} &
\includegraphics[height=4cm]{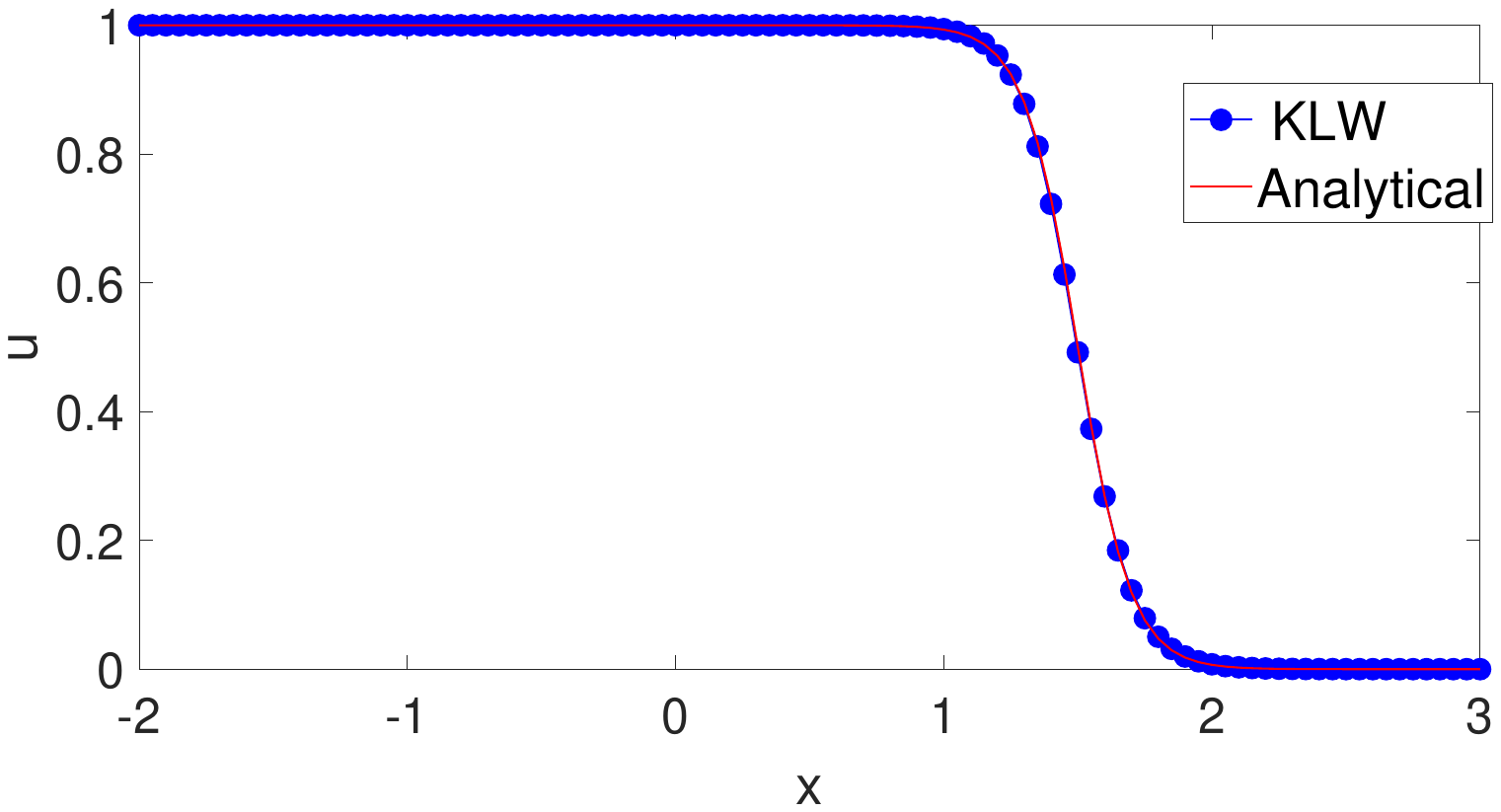} &
\includegraphics[height=4cm]{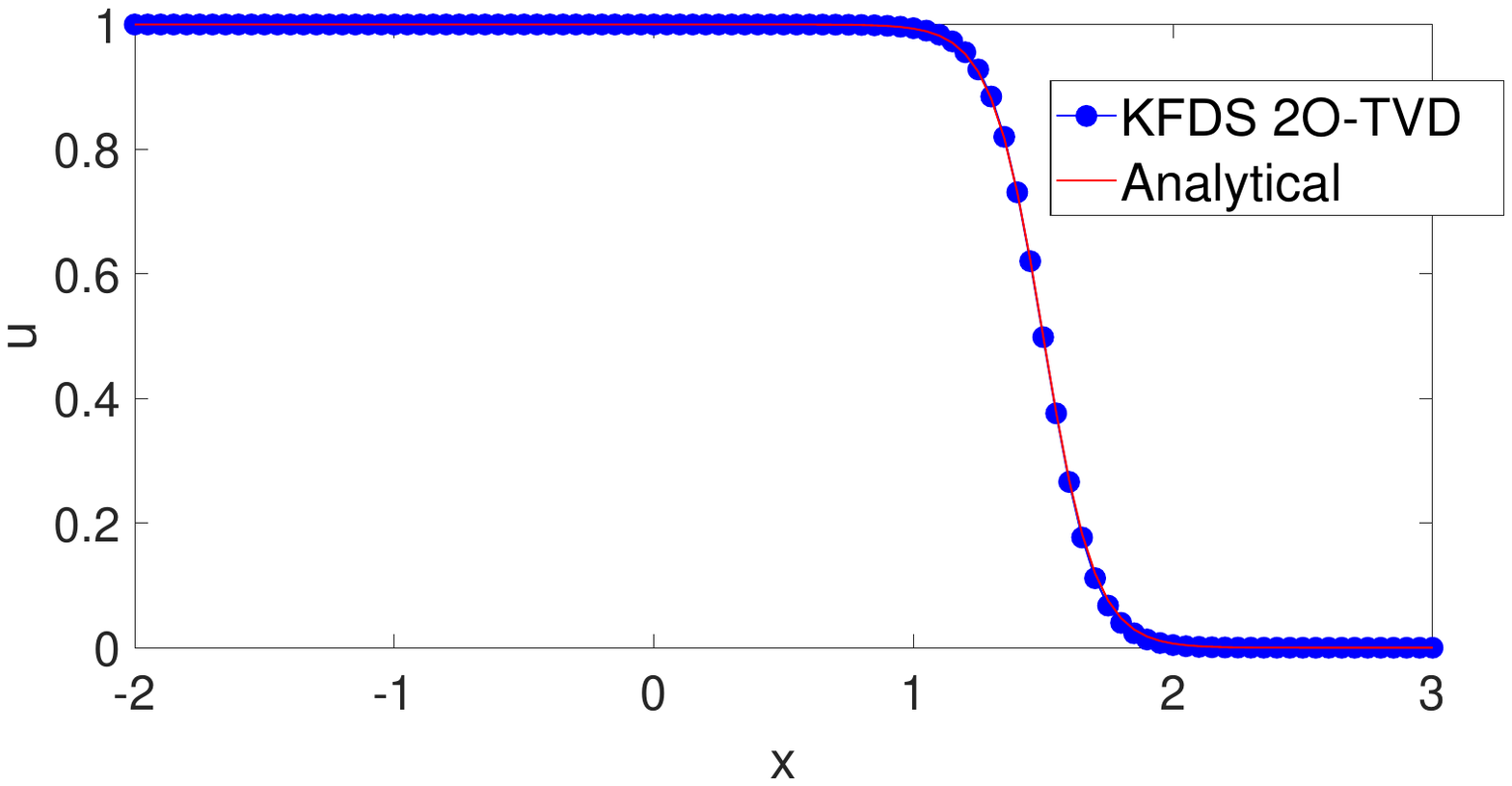} \\
 (a)First order KFDS & (b)Second order KLW & (c)TVD-KFDS\\
\includegraphics[height=4cm]{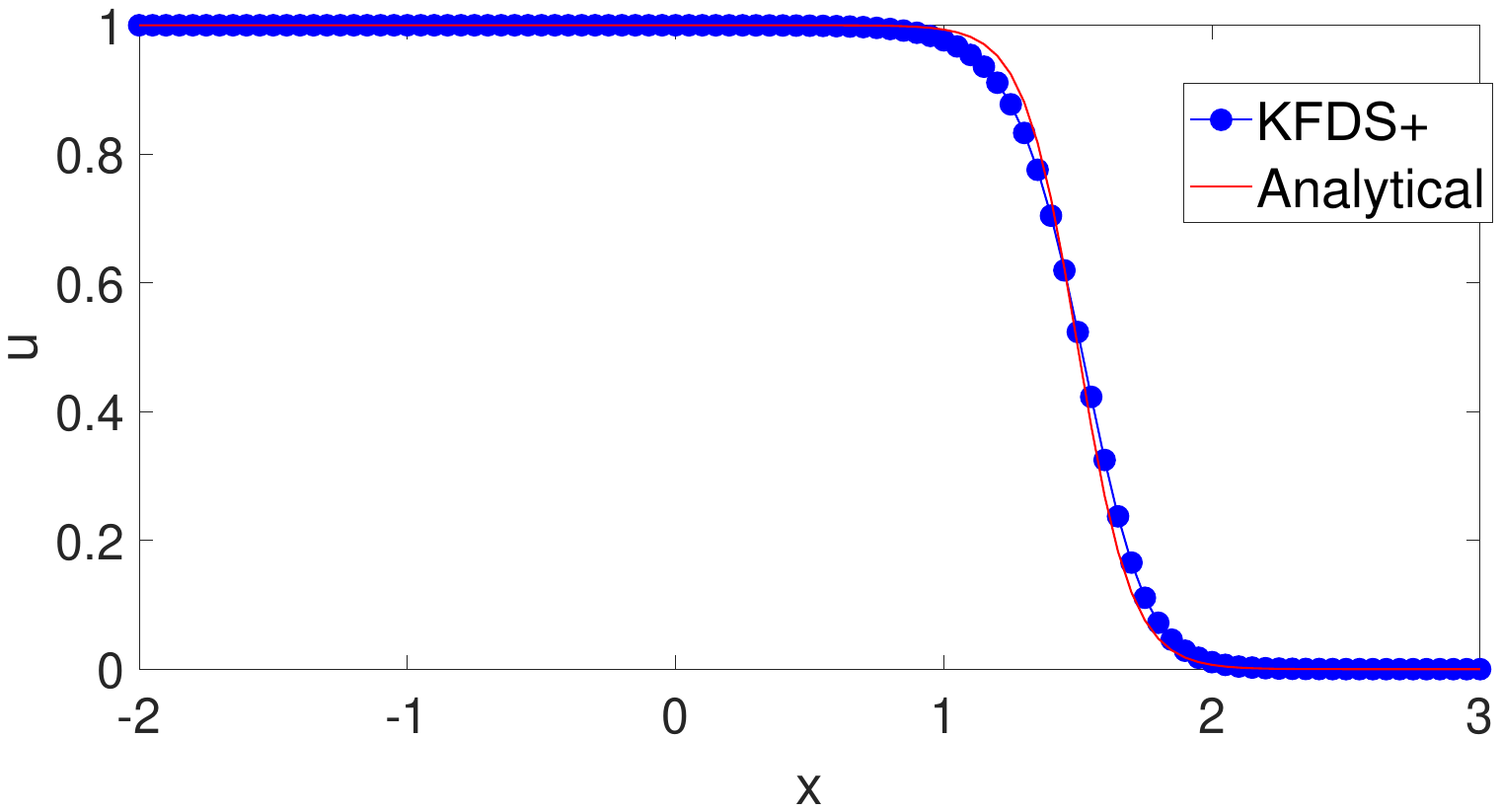} &
\includegraphics[height=4cm]{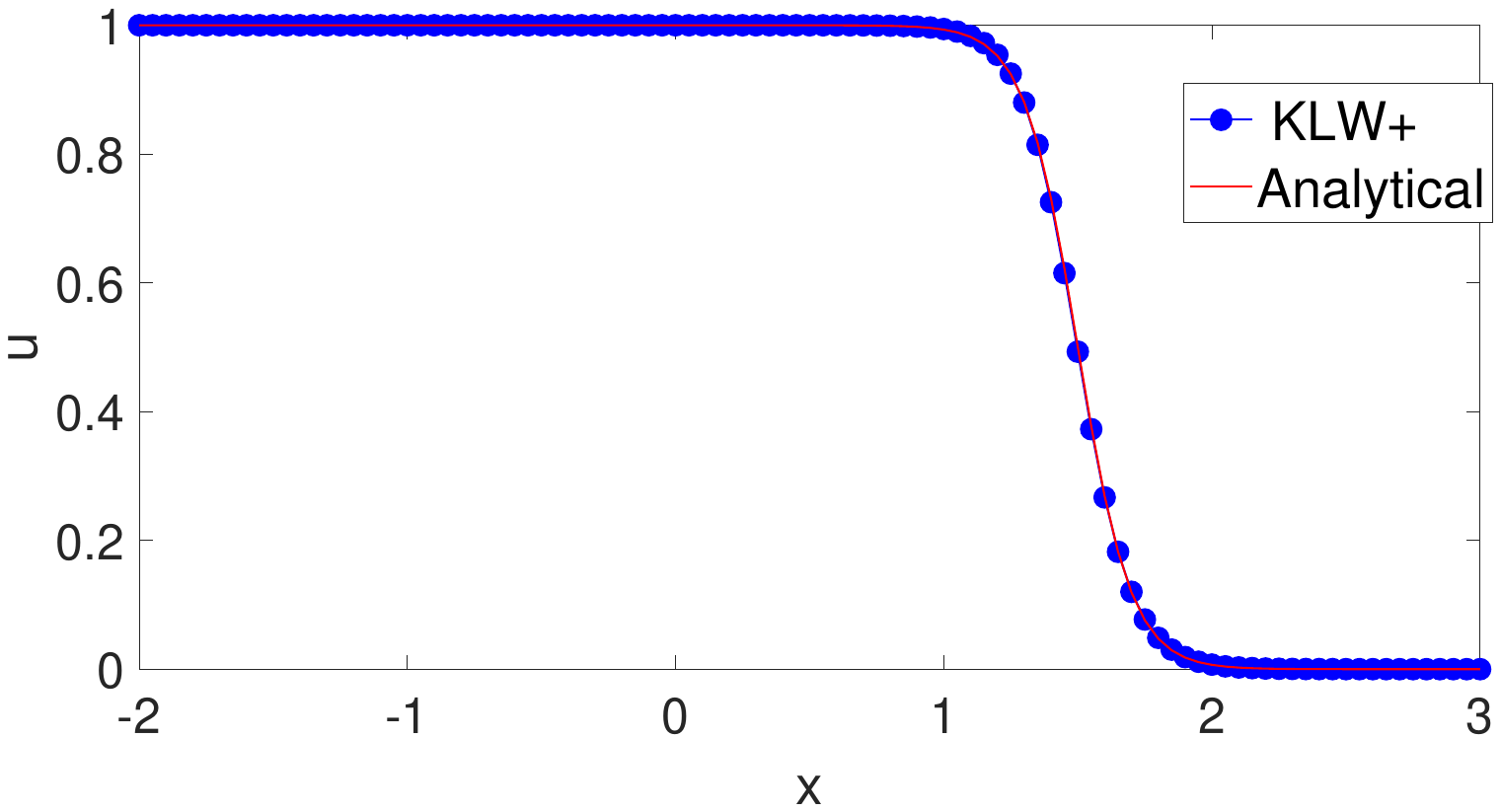} &
\includegraphics[height=4cm]{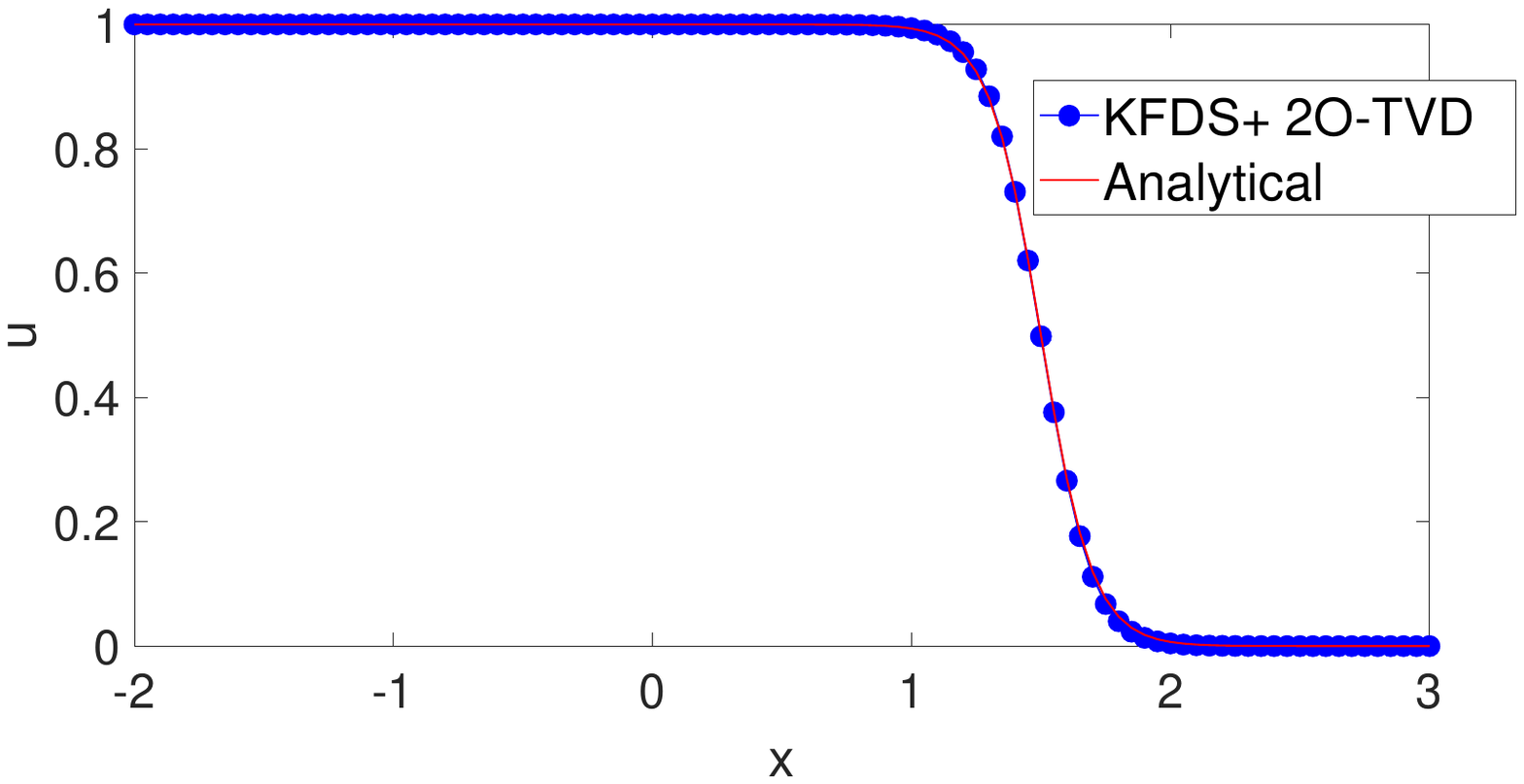} \\
 (d)First order KFDS+ &(e) KLW+  & (f)TVD-KFDS+ \\
\end{tabular}
\caption{Test Case 6(a) : KFDS, KLW \& TVD-KFDS schemes in viscous Burgers equation framework with 100 points}
\label{TC_6a_VISC_KFDS} 
\end{center} 
\end{figure}  
The initial conditions in the domain for $ case \ (b)$ are given below.  
\bea
u(x,0) = \left\{ \ba{l}  1 \ \textrm{ for } \ x \le 0 \\ 
0 \ \textrm{for} \ x  > 0    \ea \right. 
\eea
The exact solution for the test case is given below.  
\bea
 u _ { e x a c t } ( x , t ) = 1 - \frac { 1 } { 2 } \left( 1 - \frac { \exp \left( \frac { 0.5 t - x } { 4 \nu } \right) - \exp \left( \frac { x - 0.5 t } { 4 \nu } \right) } { \exp \left( \frac { 0.5 t - x } { 4 \nu } \right) + \exp \left( \frac { x - 0.5 t } { 4 \nu } \right) } \right)
\eea

\begin{figure}[!h] 
\begin{center} 
\begin{tabular}{ccc}
\includegraphics[height=4cm]{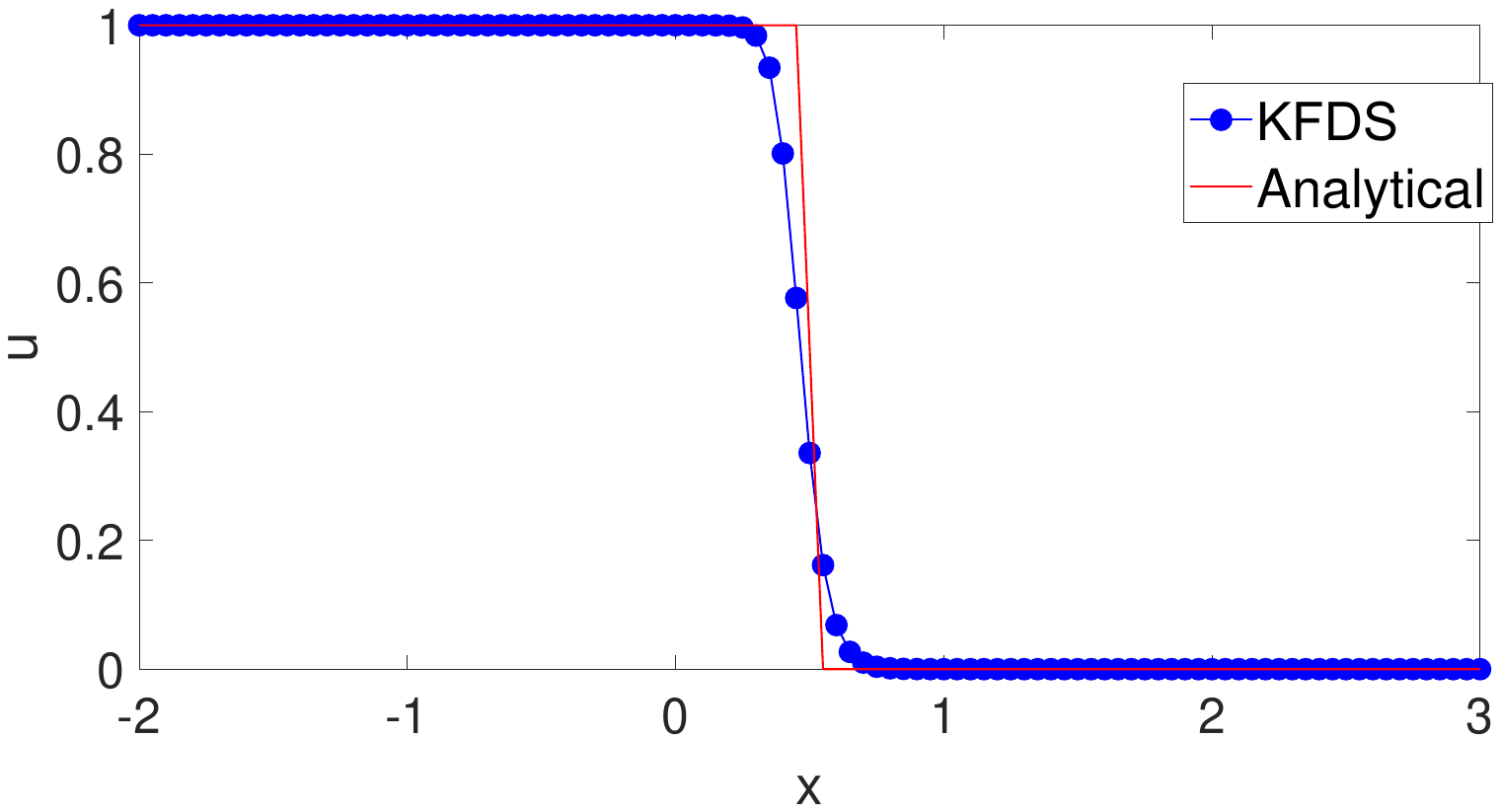} &
\includegraphics[height=4cm]{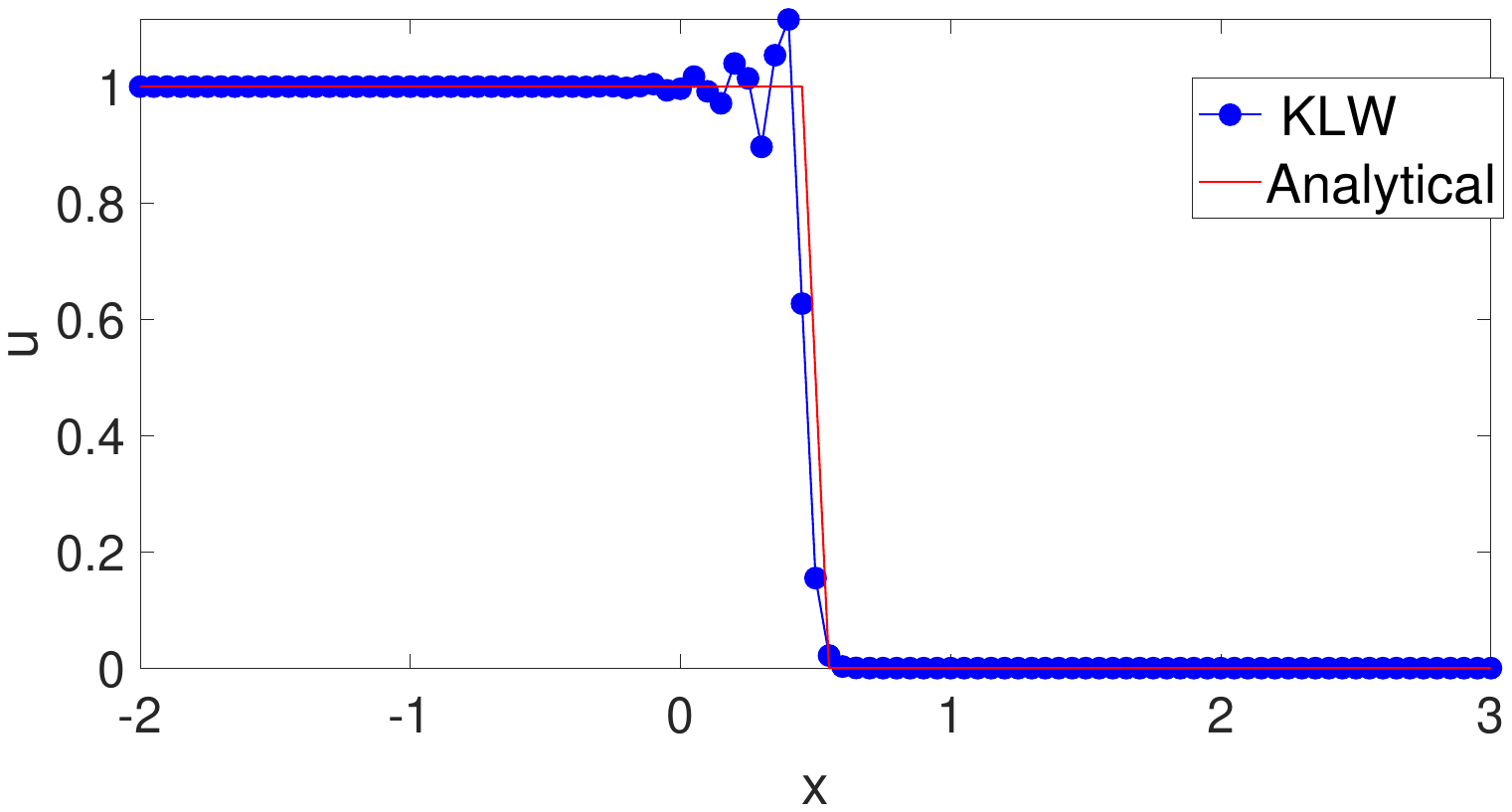} &
\includegraphics[height=4cm]{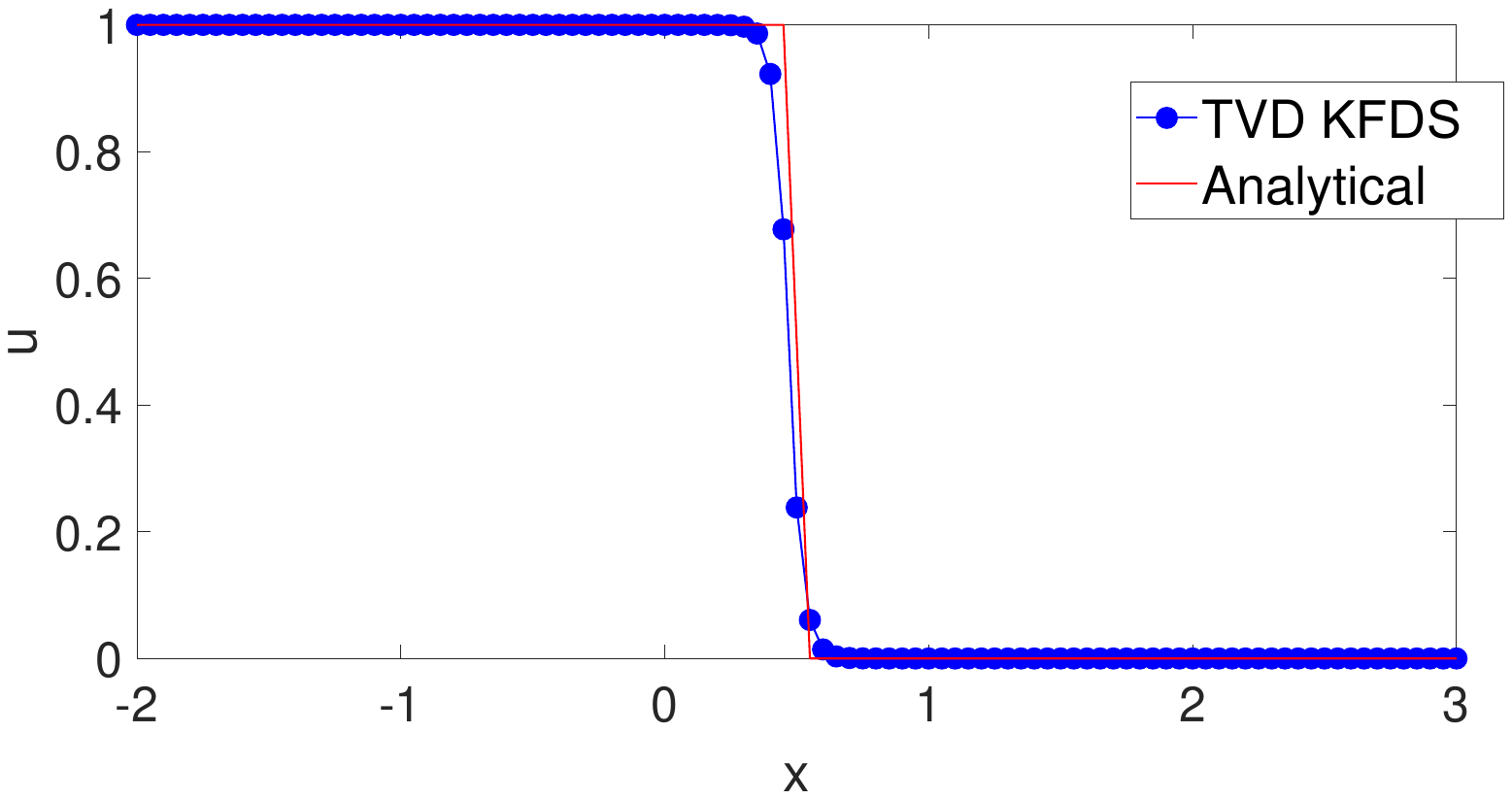} \\
 (a)First order KFDS & (b)Second order KLW & (c)TVD-KFDS\\
\includegraphics[height=4cm]{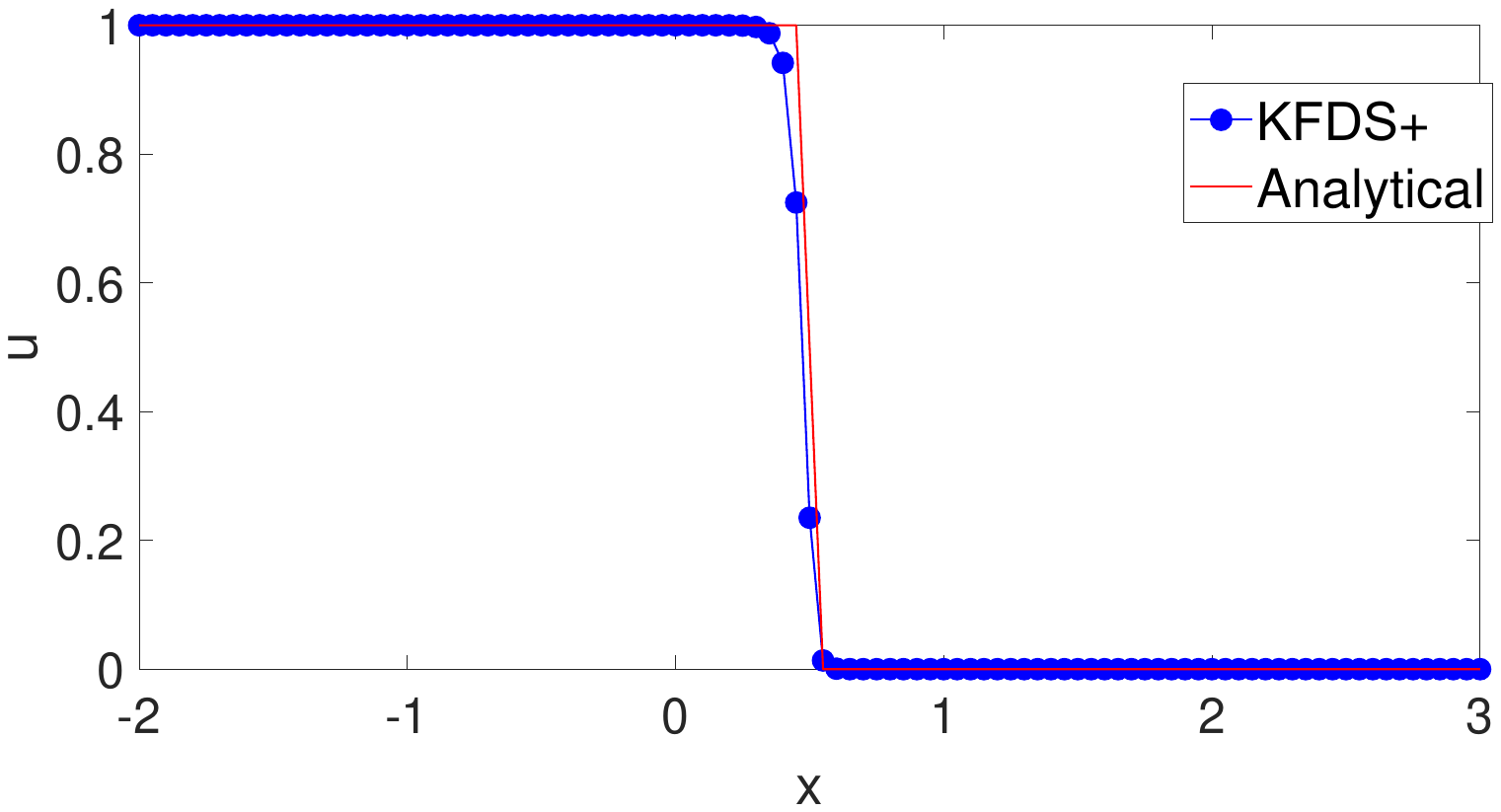} &
\includegraphics[height=4cm]{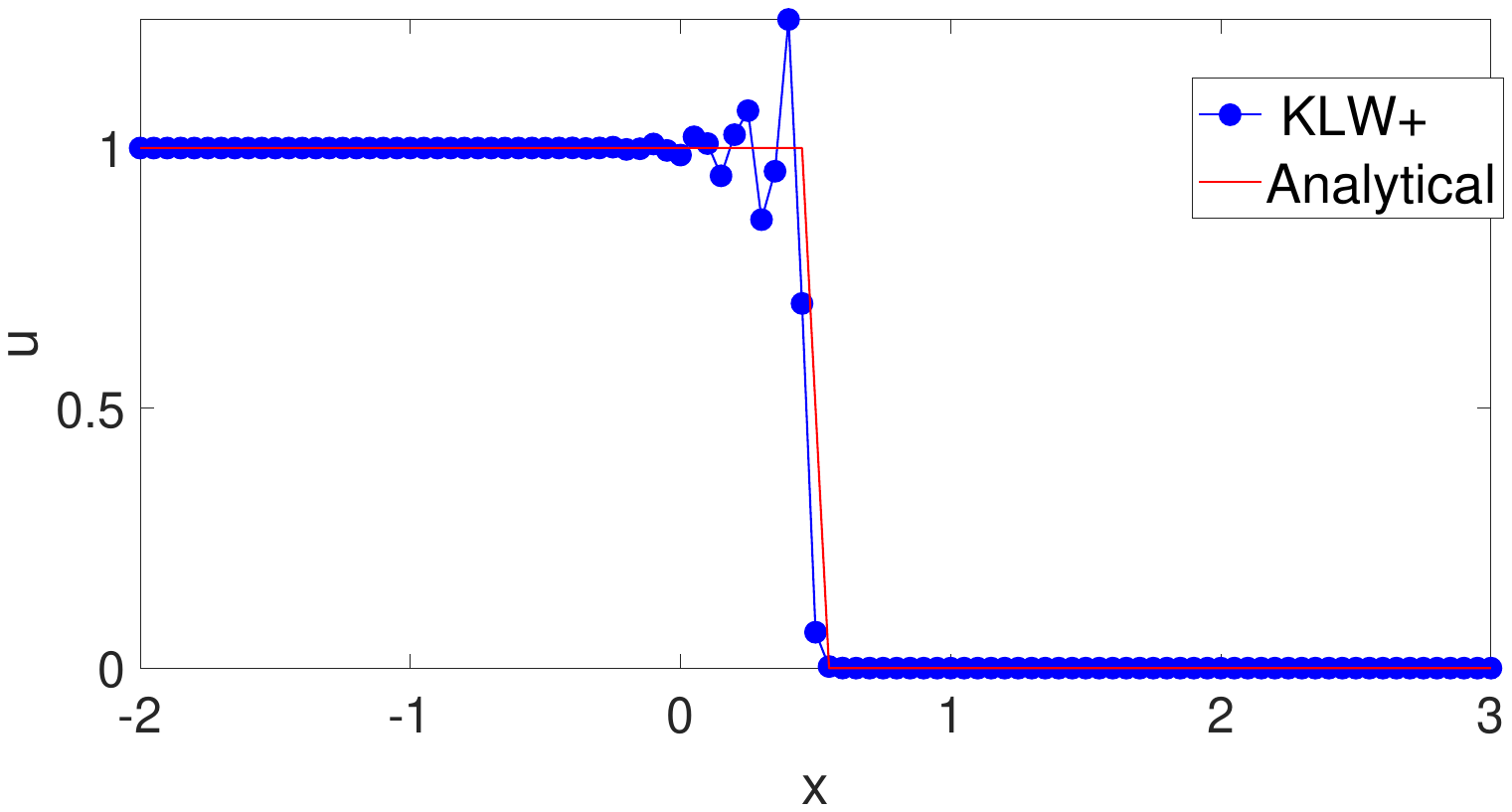} &
\includegraphics[height=4cm]{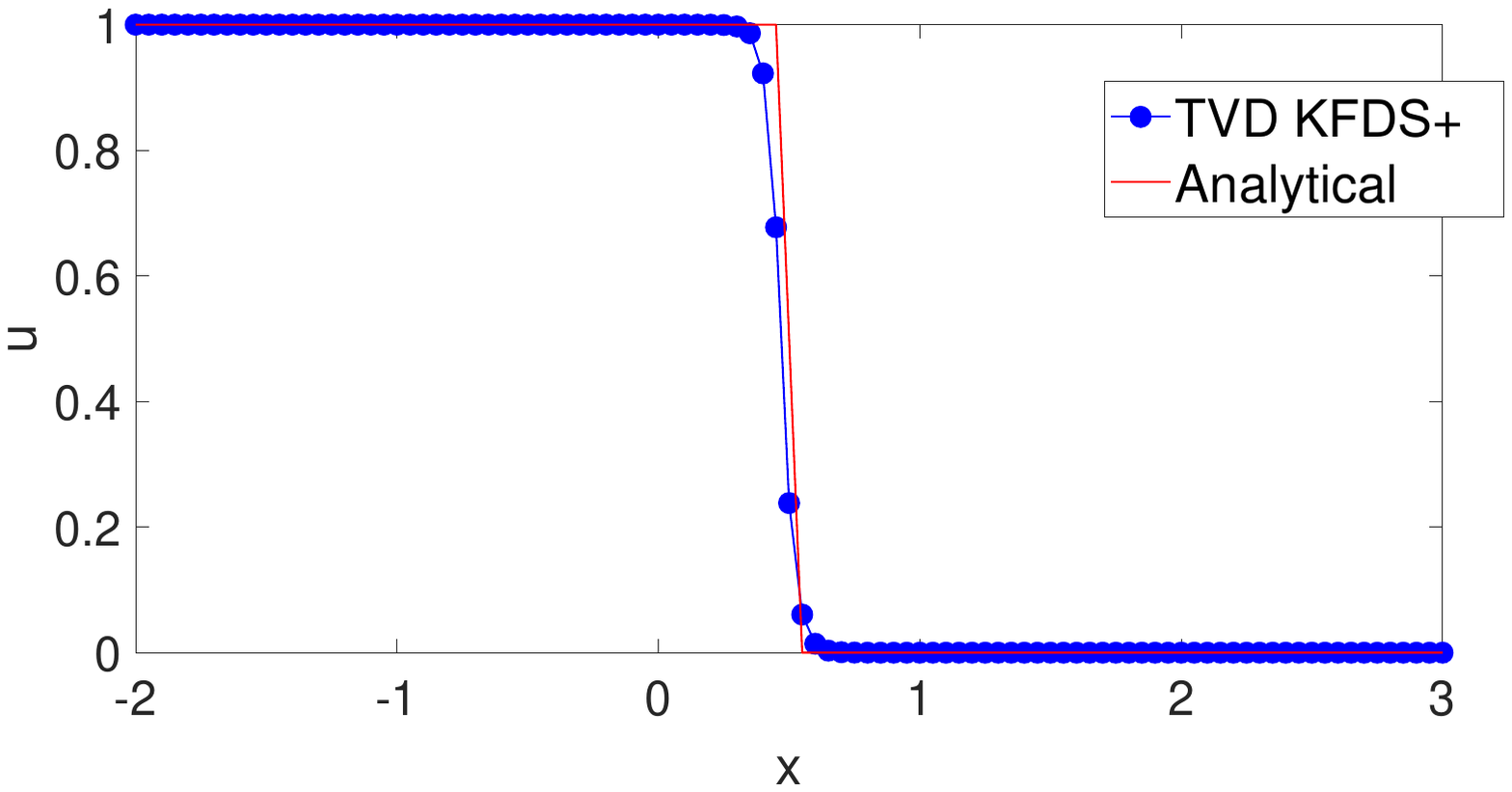} \\
  (d)First order KFDS+ &(e) KLW+  & (f)TVD-KFDS+ \\
\end{tabular}
\caption{Test Case 6(b) : KFDS, KLW \& TVD-KFDS  schemes in Viscous Burgers framework with 100 points}
\label{TC_6b_VISC_KFDS} 
\end{center} 
\end{figure}

The results for the test case $6(a)$  and are given in figure (\ref{TC_6a_VISC_KFDS}). This test case is analogous to a moving shock case for the inviscid Burgers equation. The numerical solutions are computed for $t=3s$ and are compared with the analytical solution. Both the first order KFDS and KFDS+ methods capture the shock position with reasonably well, with KFDS+ being more accurate.  All the higher order versions match the analytical solution with much more closely.  

The results for the test case $6(b)$  and are given in figure (\ref{TC_6b_VISC_KFDS}). This test case is characterised by a low viscosity coefficient of $\nu = 0.001$. The low viscosity coefficient results in the domination of convection fluxes in comparison to the viscous fluxes.  In general, higher order schemes typically exhibit oscillatory behaviour in the regions of large gradients. This is seen clearly in the results of KLW schemes. The TVD version effectively demonstrates the switching of the scheme from KLW to KFDS specifically in regions of large gradients and yields an  oscillation-free higher order solution.

\subsubsection*{Test case 7: Sinusoidal wave problem}
The viscous Burgers equation is solved for a sinusoidal wave \cite{Mund} with initial condition given by $u ( x , 0 ) = - \sin ( \pi x )$. The computational domain $[-1, 1]$ is divided into 100 cells. The simulation is carried out with a viscosity coefficient $\nu$ = 0.1  and time $t = 2.55237s$.  The boundary conditions for this test case are 
$u(-1, t) = 0$  and  $u(1, t) = 0$.  The exact solution for the test case \cite{Benton} is given below.  
\bea
u _ { e x a c t } ( x , t ) = 4 \pi \nu \frac { \sum _ { n = 1 } ^ { \infty } n a _ { n } \exp \left( - \nu n ^ { 2 } \pi ^ { 2 } t \right) \sin ( n \pi x ) } { a _ { 0 } + 2 \sum _ { n = 1 } ^ { \infty } n a _ { n } \exp \left( - \nu n ^ { 2 } \pi ^ { 2 } t \right) \cos ( n \pi x ) }\
\eea
where $a _ { n } = ( - 1 ) ^ { n } I _ { n } \left( \frac { 1 } { 2 \pi \nu } \right)$.    Here $I _ { n } ( \frac { 1 } { 2 \pi \nu } )$ denotes the modified Bessel function of the first kind and order $n$.  

\begin{figure}[!h] 
\begin{center} 
\begin{tabular}{ccc}
\includegraphics[height=4cm]{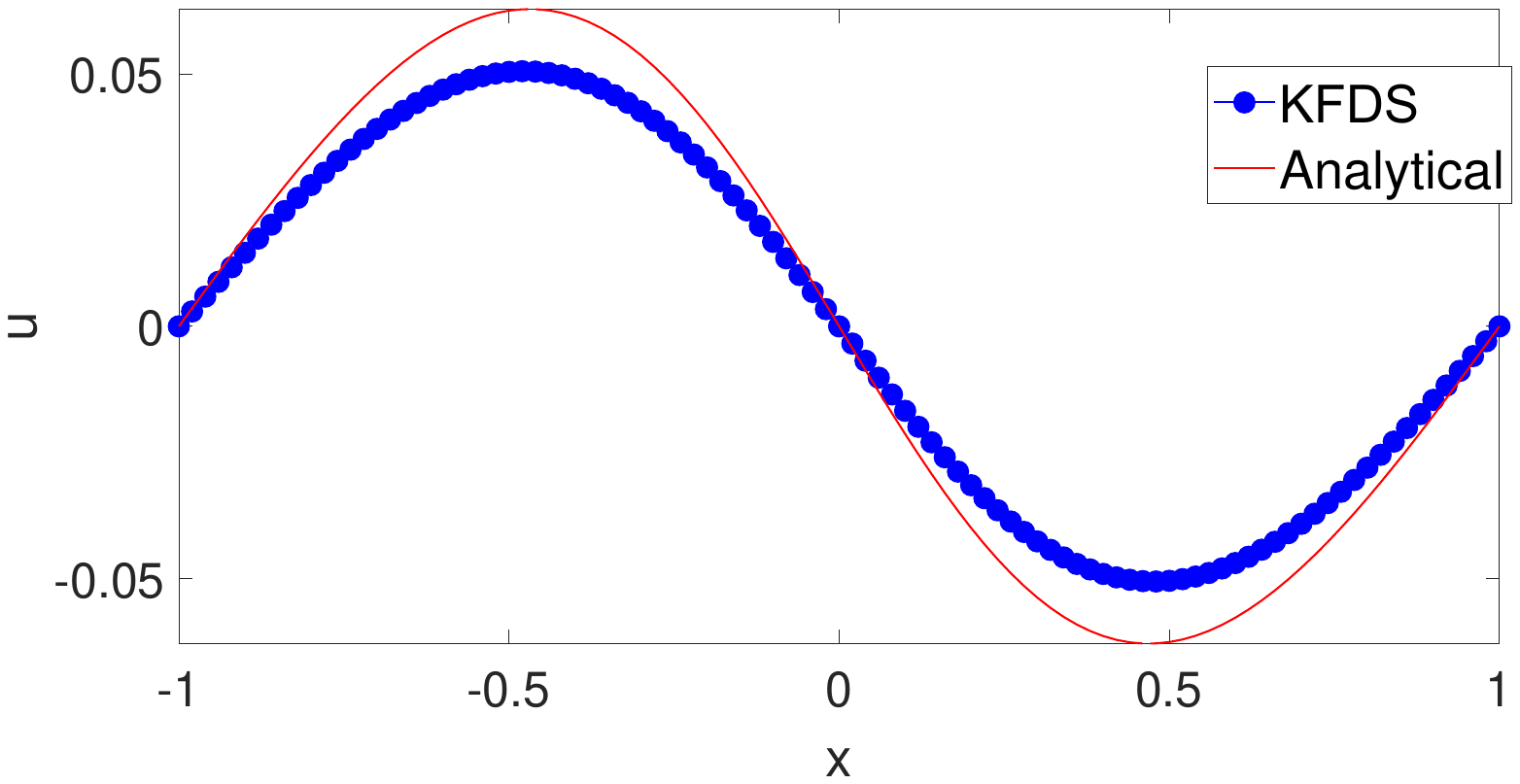} &
\includegraphics[height=4cm]{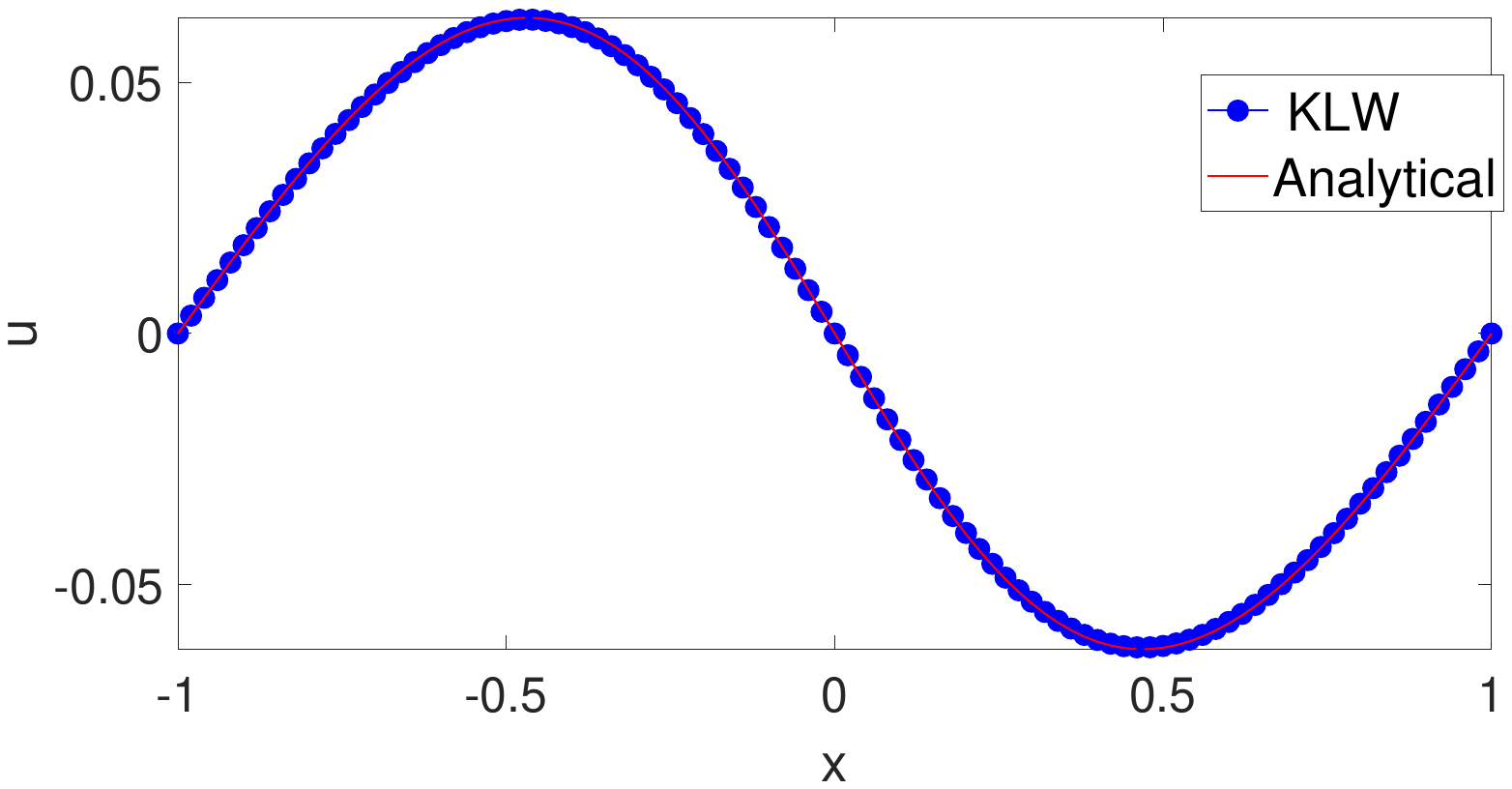} &
\includegraphics[height=4cm]{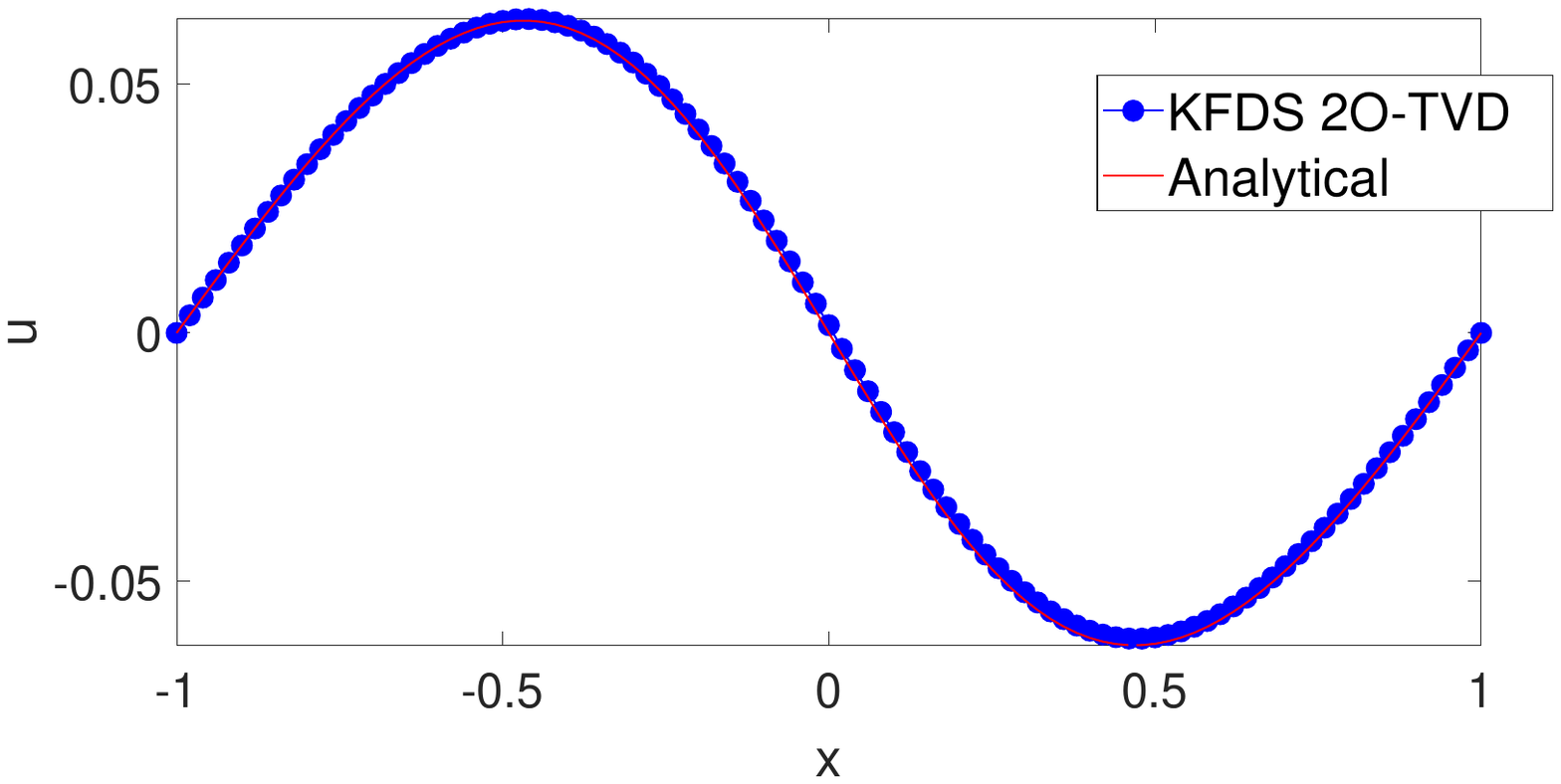} \\
 (a)First order KFDS & (b)Second order KLW & (c)TVD-KFDS\\
\includegraphics[height=4cm]{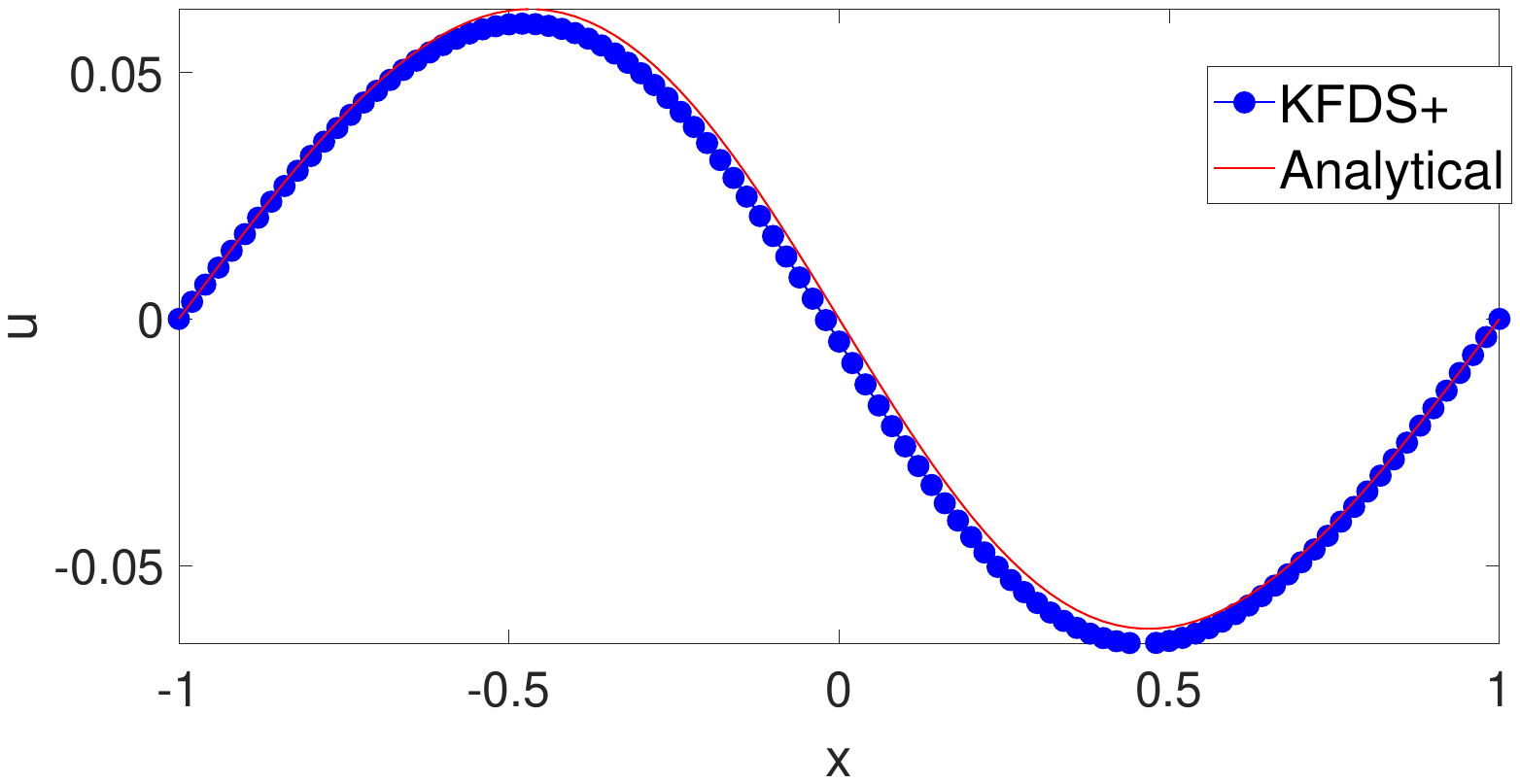} &
\includegraphics[height=4cm]{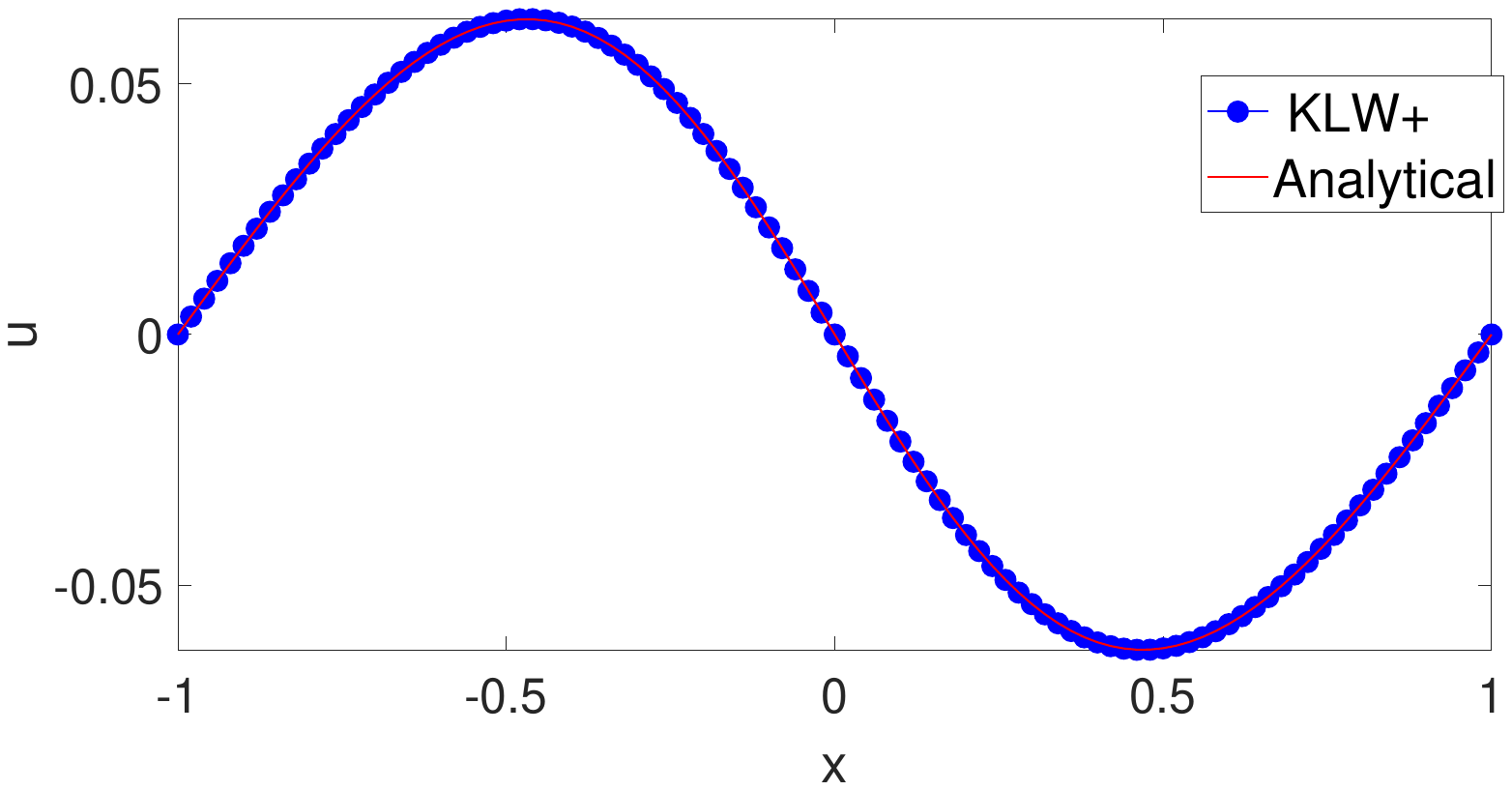} &
\includegraphics[height=4cm]{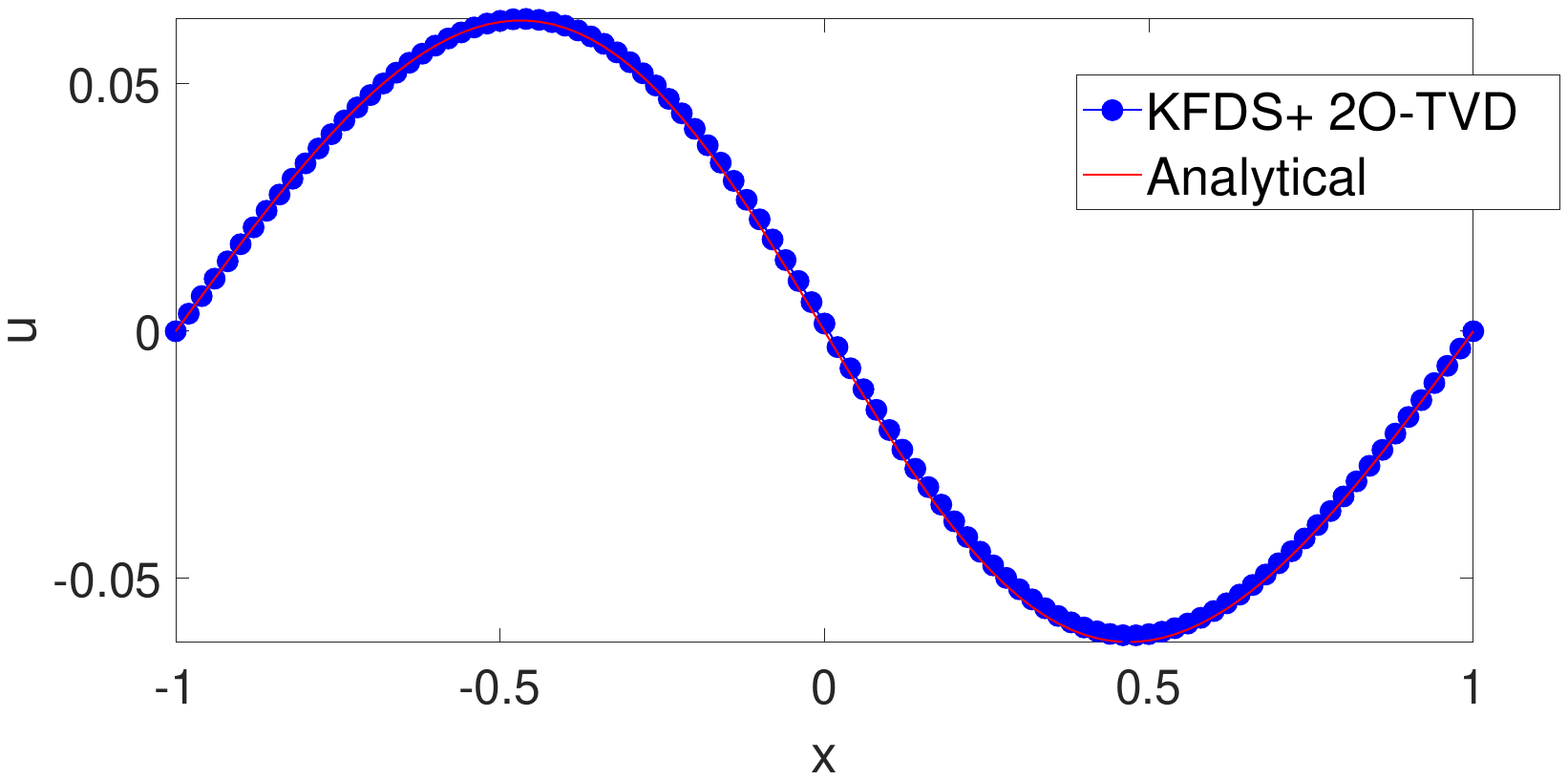} \\
 (d)First order KFDS+ &(e) KLW+  & (f)TVD-KFDS+ \\
\end{tabular}
\caption{Test Case 7 : KFDS, KLW \& TVD-KFDS  schemes in Viscous Burgers framework with 100 points}
\label{TC_7_VISC_KFDS} 
\end{center} 
\end{figure}

The results for this test case are given in figure (\ref{TC_7_VISC_KFDS}). The sinusoidal wave test case effectively helps in evaluating the inherent numerical diffusion in the scheme.  There is a significant difference between the KFDS scheme and the analytical results indicating the scheme generates additional numerical diffusion. However, the KFDS+ scheme has better match with the analytical solution. The second order KLW schemes adhere to the analytical solution closely.  Similarly, the TVD schemes also exhibit good match with the analytical solution.

\subsubsection*{Test case 8: Steady shock test case}
The viscous version of the steady shock test case \cite{Fletcher} is aimed at evaluating the ability of the numerical scheme in capturing the viscous diffusion.  The domain $[-1, 1]$ is divided into 100 cells and the initial conditions are as given below.  
\bea
u(x,0) = \left\{ \ba{l}  1 \ \textrm{ for } \ x \le 0 \\ 
-1 \ \textrm{for} \ x  > 0    \ea \right. 
\eea
The exact solution for this test case \cite{Benton} is given below  
\bea
u _ { e x a c t } ( x , t ) = \frac { 1 } { 2 } \left[ \left( u _ { L } + u _ { R } \right) - \left( u _ { L } - u _ { R } \right) \tanh \left( \frac { x \left( u _ { L } - u _ { R } \right) } { 4 \nu } \right) \right]
\eea
where $u_L$ and $u_R$ correspond to the initial left and right states of the domain. The cases are simulated for a viscosity coefficient of $(a)$  $\nu = 0.1$ and $(b)$ $\nu = 0.001$.  

\begin{figure}[!h] 
\begin{center}
\begin{tabular}{ccc}
\includegraphics[height=4cm]{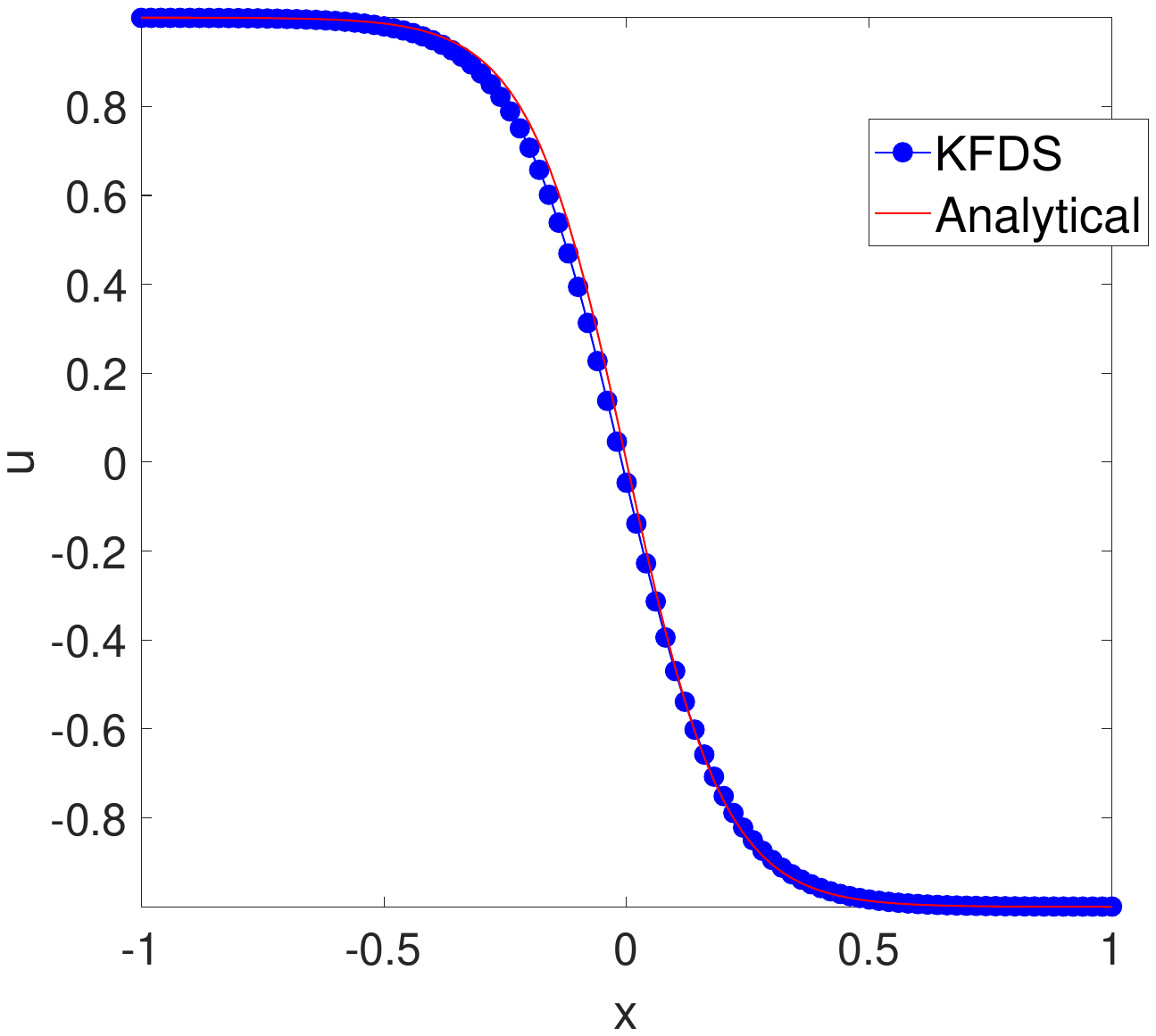} &
\includegraphics[height=4cm]{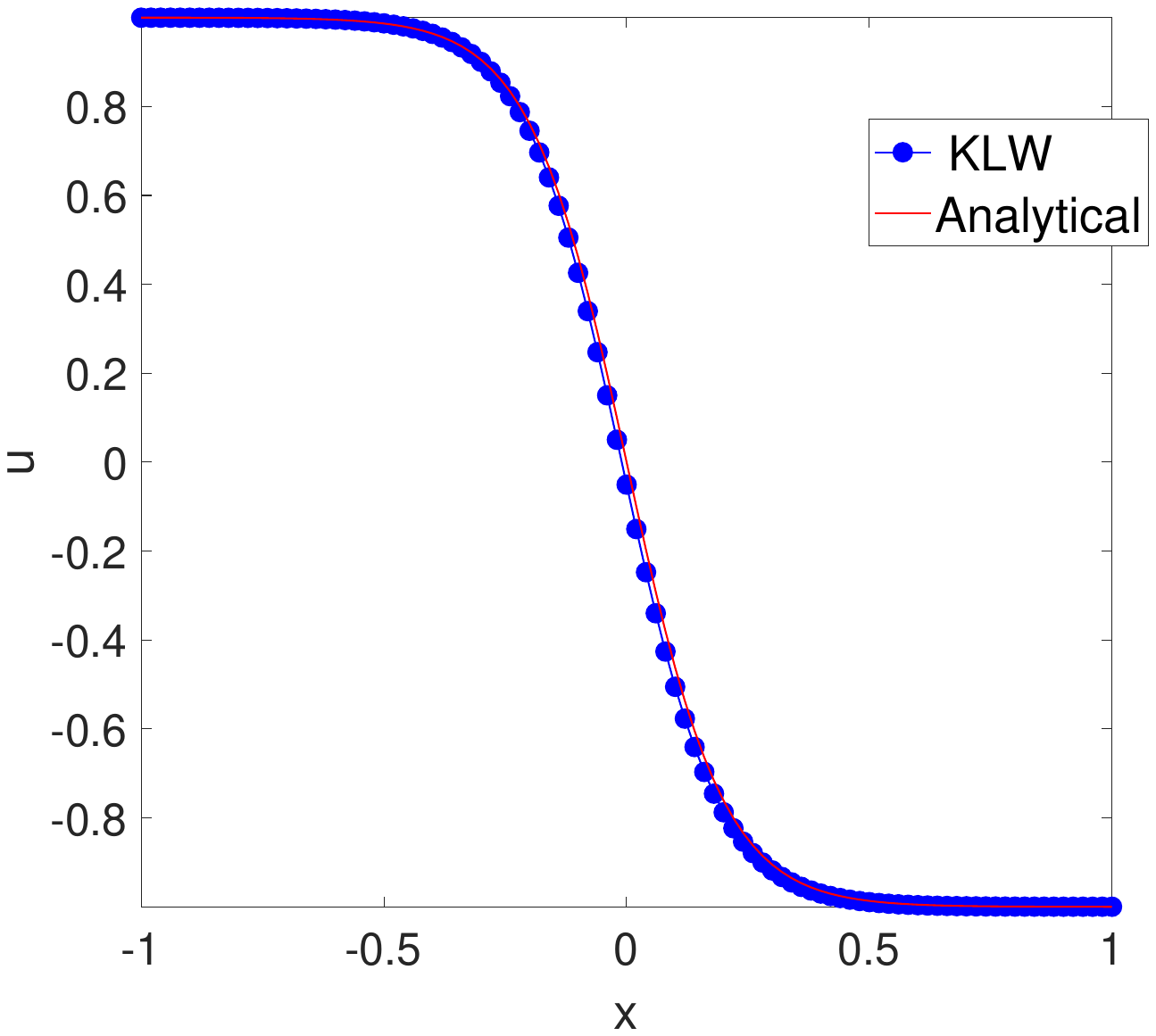} &
\includegraphics[height=4cm]{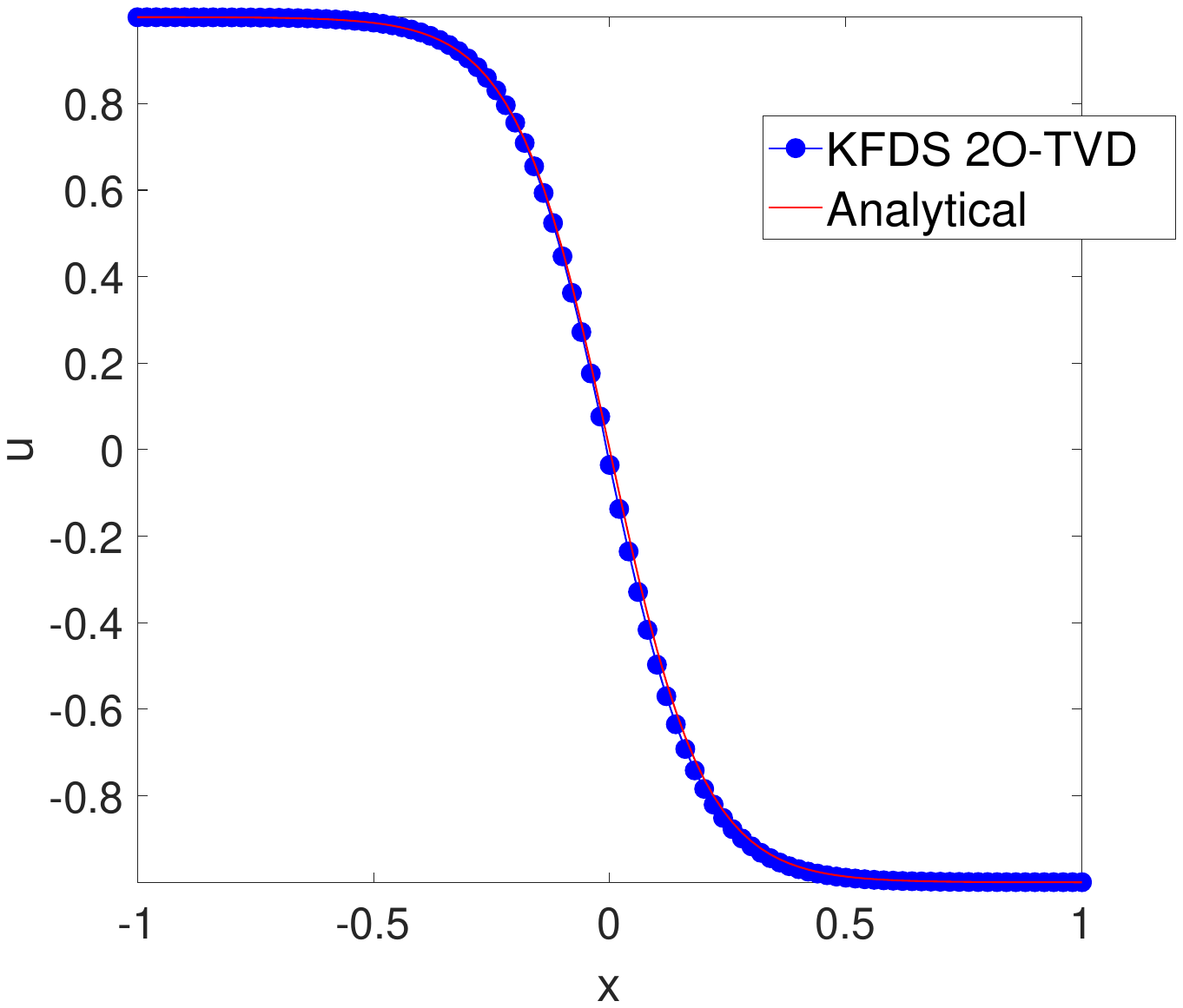} \\
 (a)First order KFDS & (b)Second order KLW & (c)TVD-KFDS\\
\includegraphics[height=4cm]{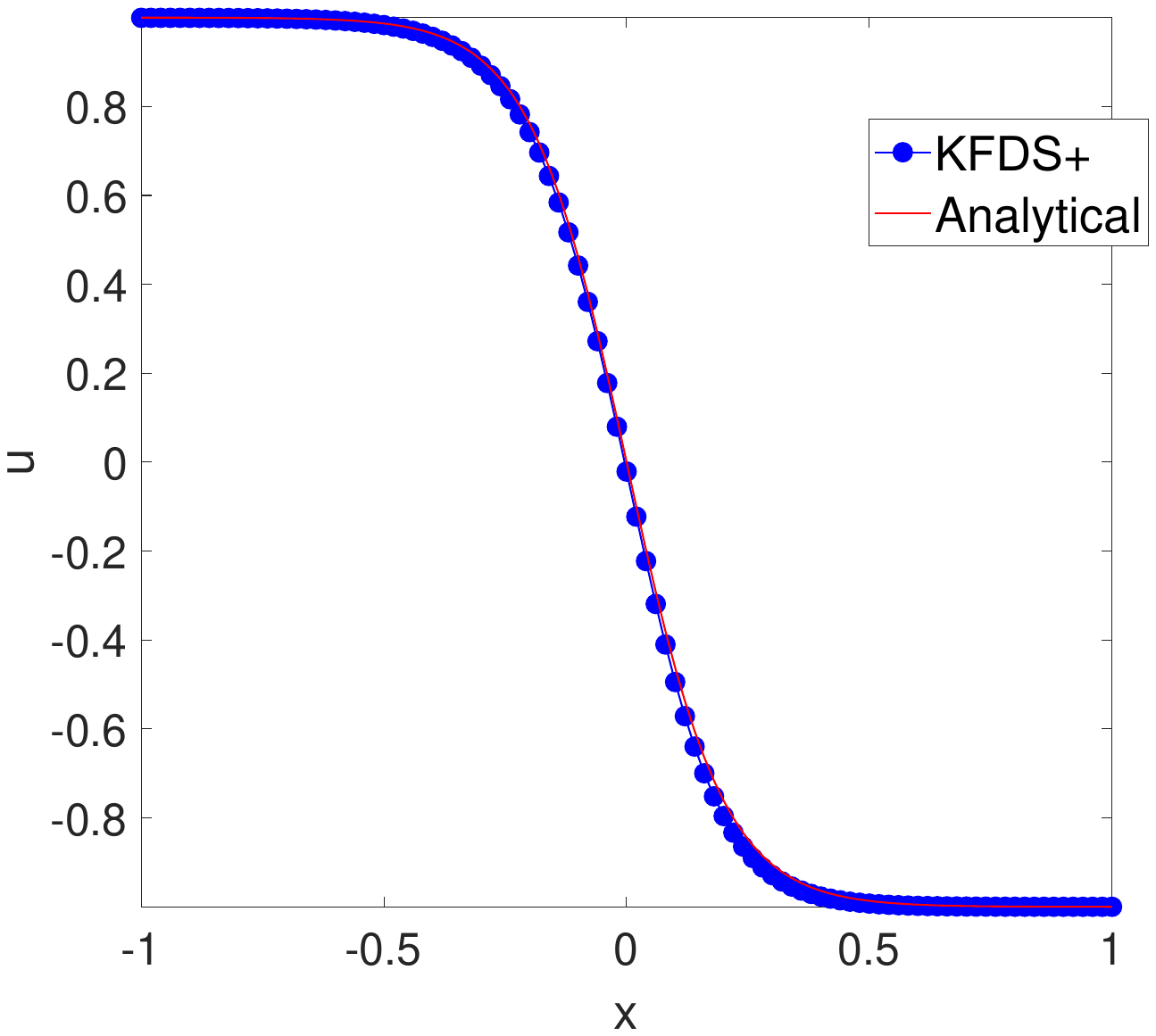} &
\includegraphics[height=4cm]{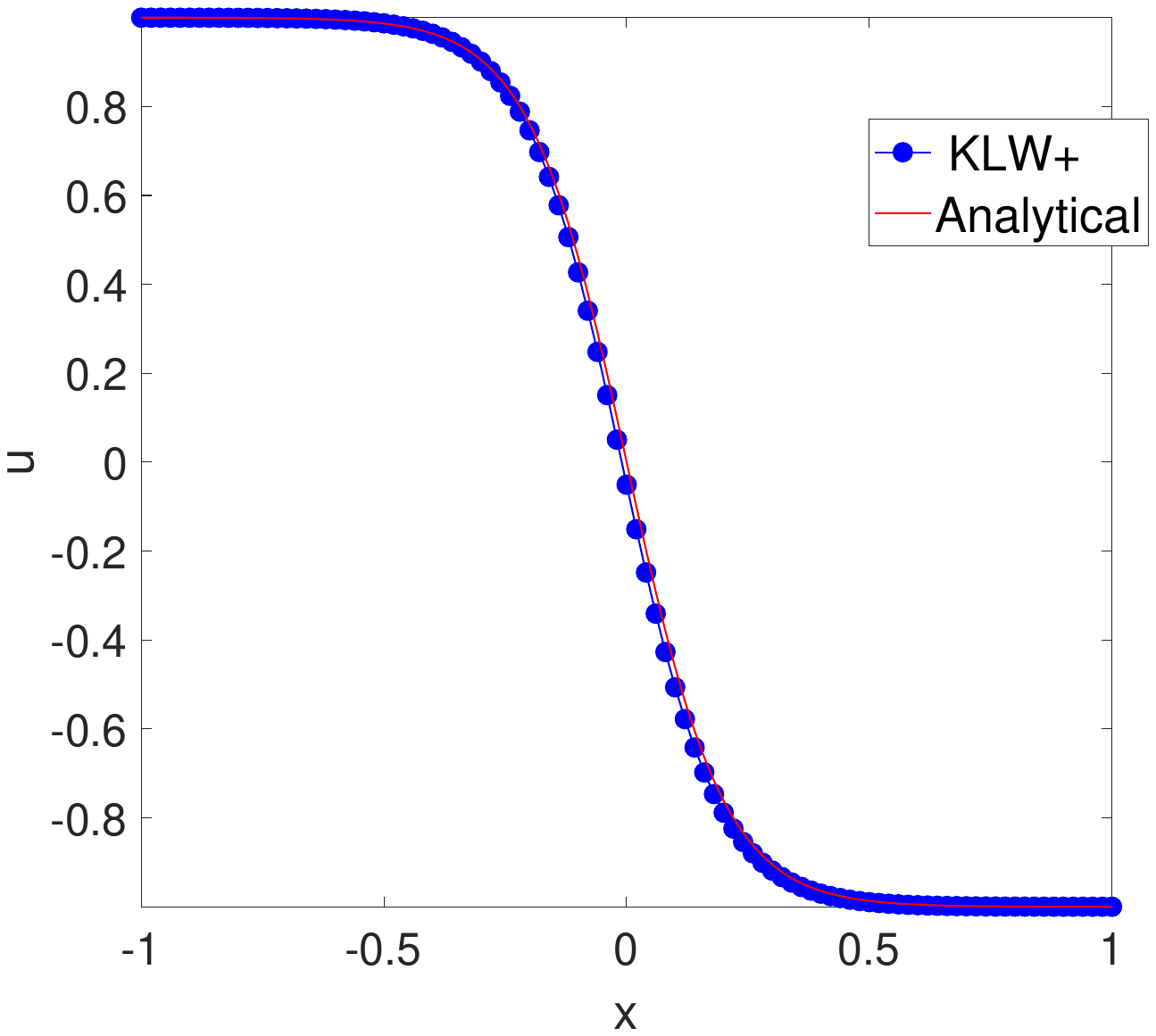} &
\includegraphics[height=4cm]{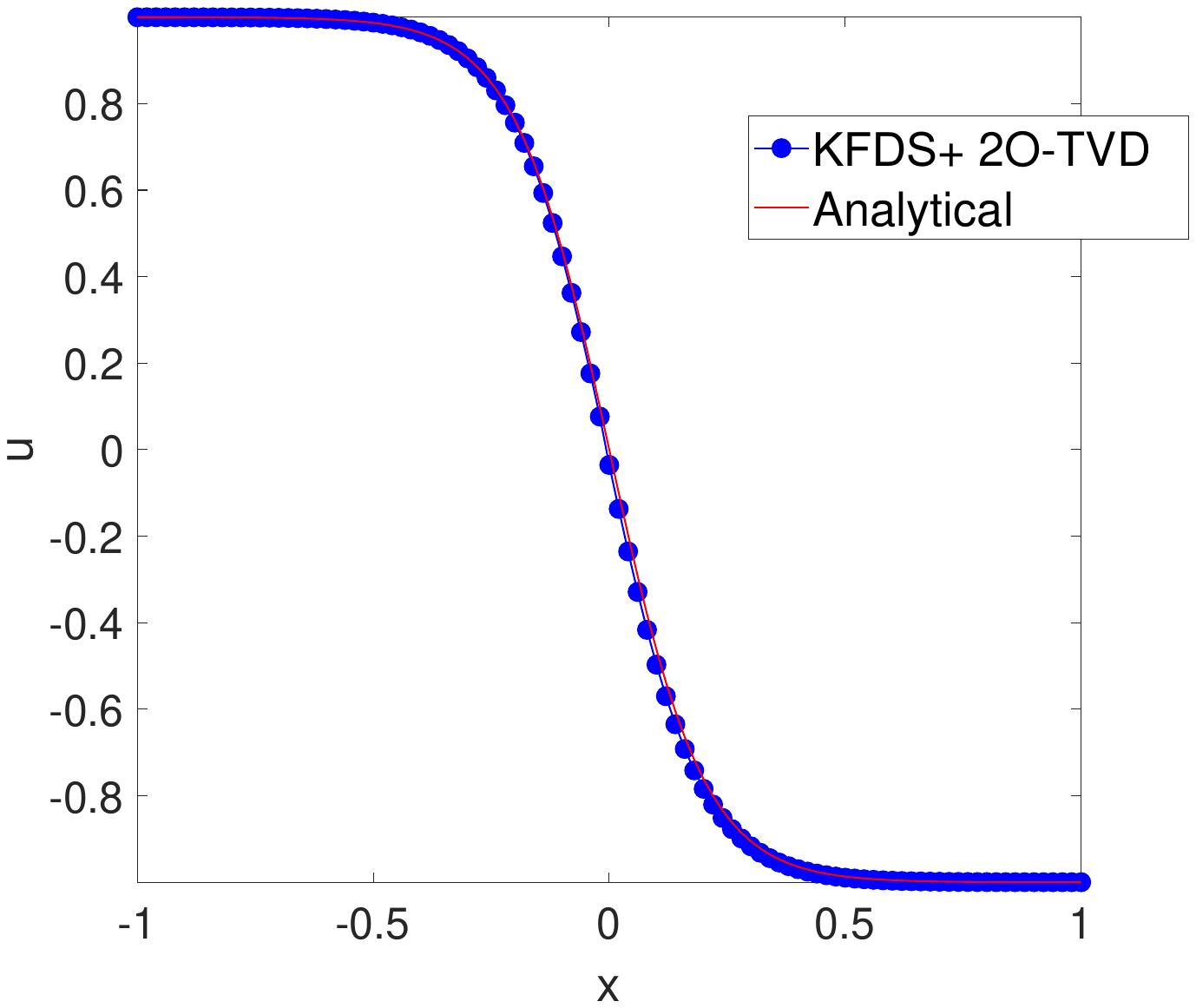} \\
  (d)First order KFDS+ &(e) KLW+  & (f)TVD-KFDS+ \\
\end{tabular}
\caption{Test Case 8(a) : KFDS, KLW \& TVD-KFDS schemes in viscous Burgers equation framework with 100 points}
\label{TC_8a_VISC_KFDS} 
\end{center} 
\end{figure}

\begin{figure}[!h] 
\begin{center} 
\begin{tabular}{ccc}
\includegraphics[height=4cm]{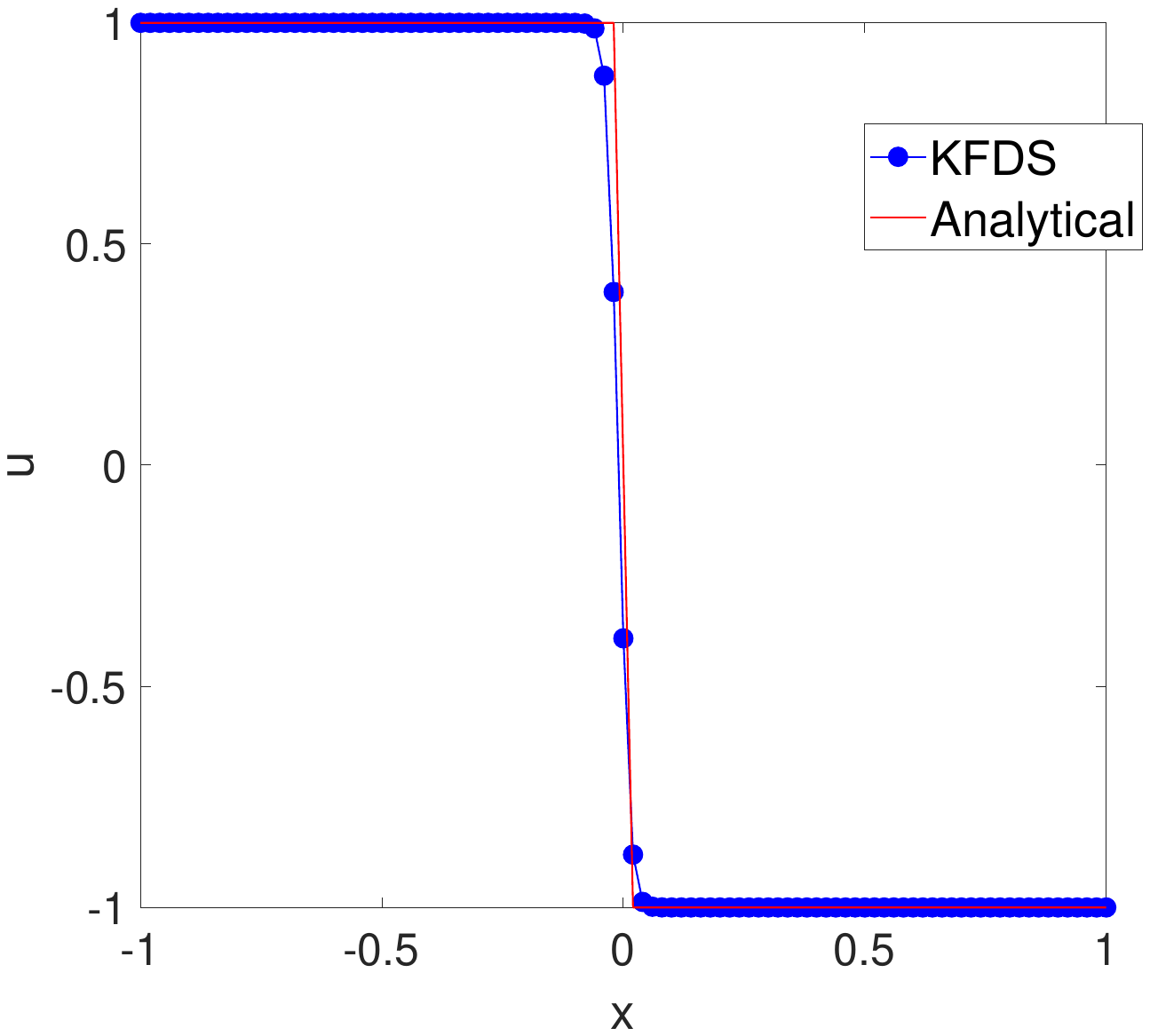} &
\includegraphics[height=4cm]{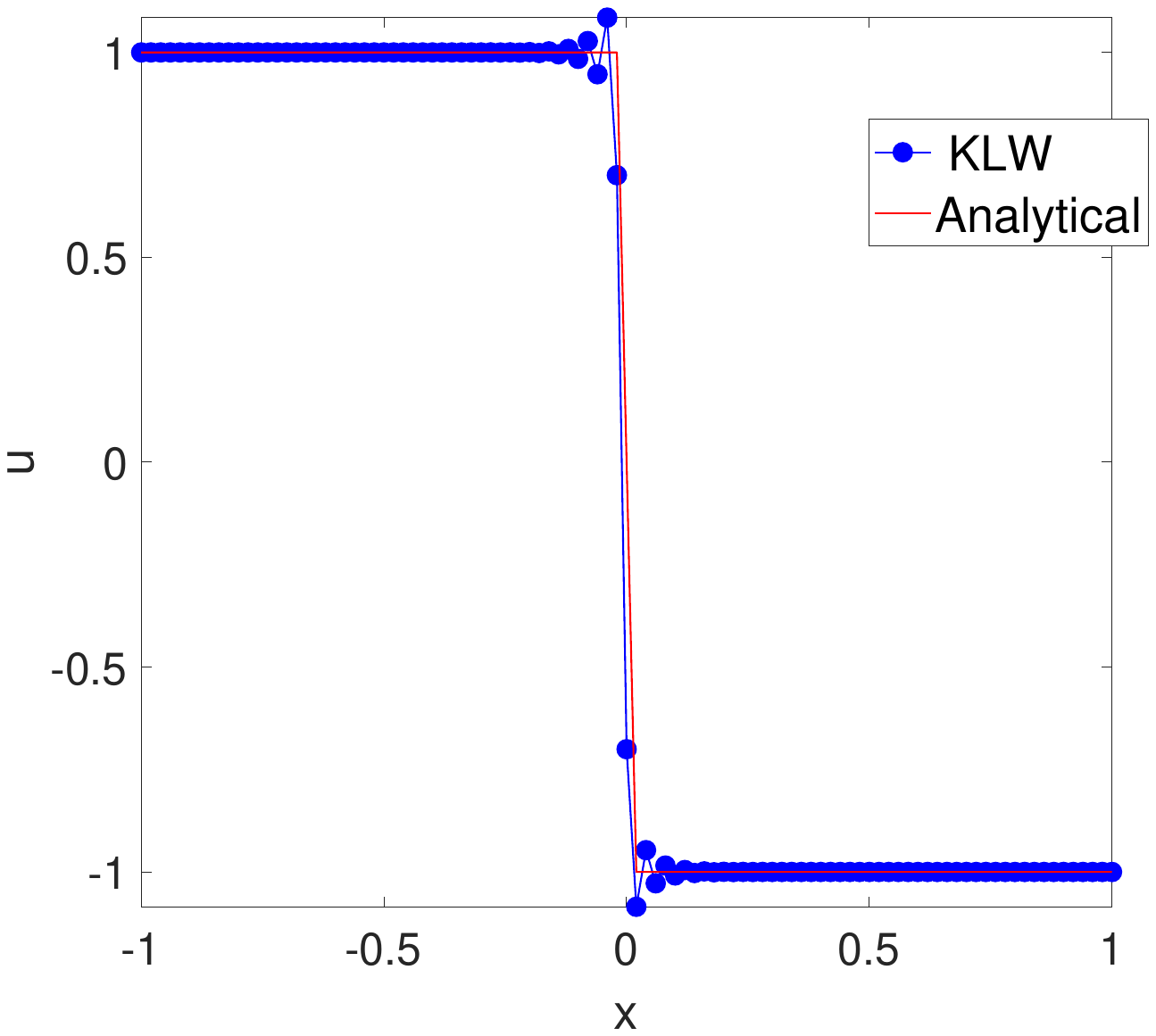} &
\includegraphics[height=4cm]{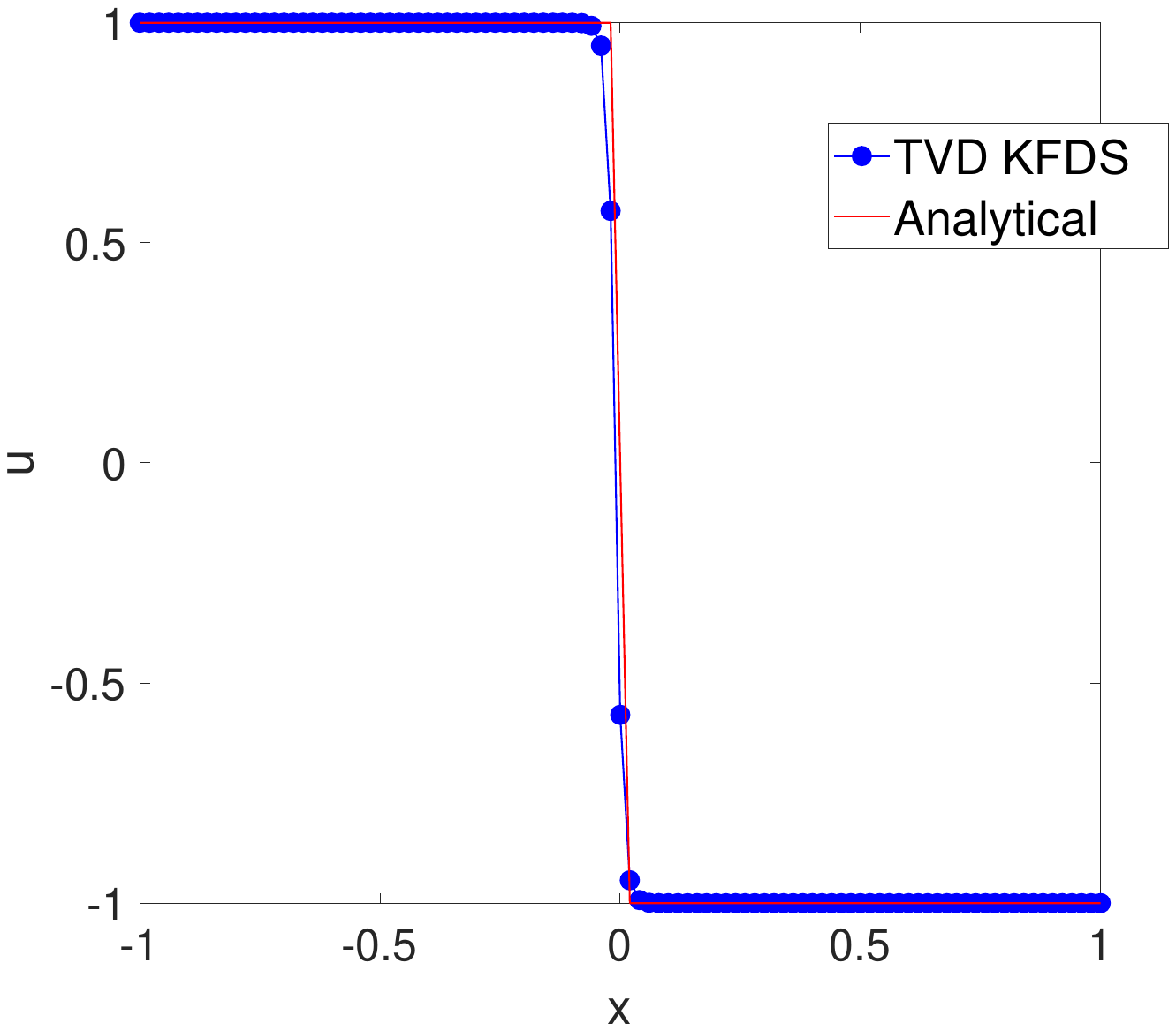} \\
 (a)First order KFDS & (b)Second order KLW & (c)TVD-KFDS\\
\includegraphics[height=4cm]{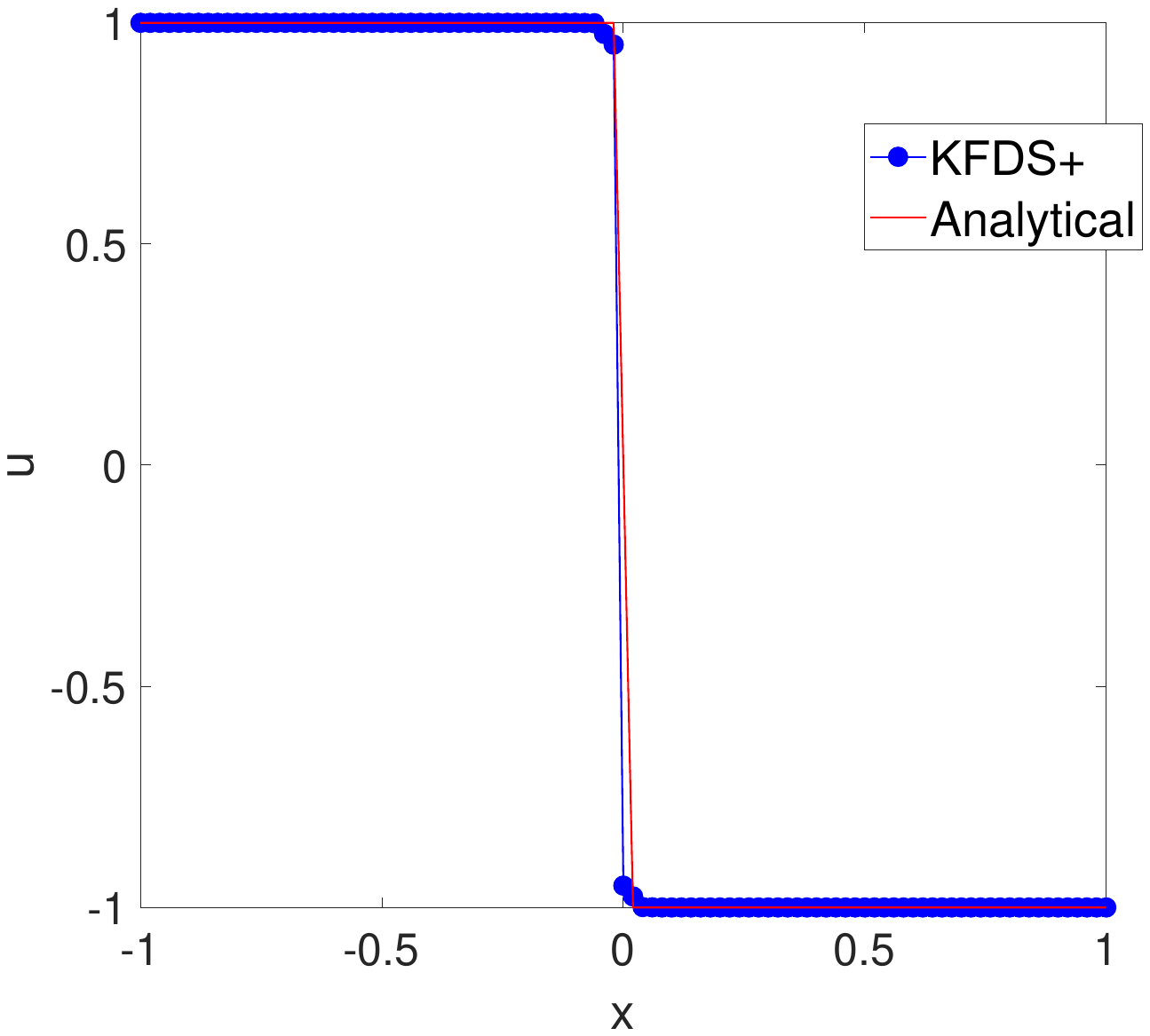} &
\includegraphics[height=4cm]{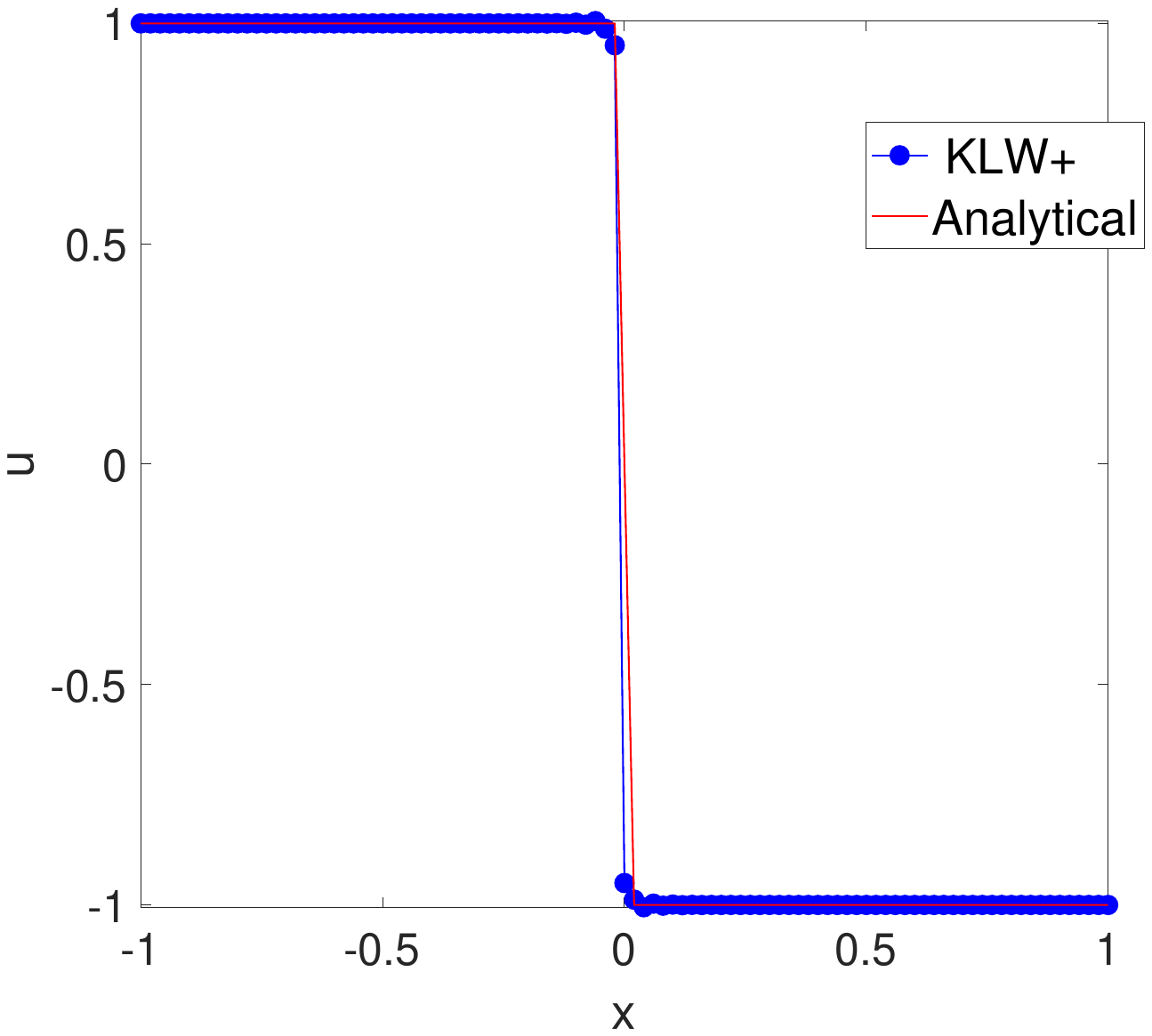} &
\includegraphics[height=4cm]{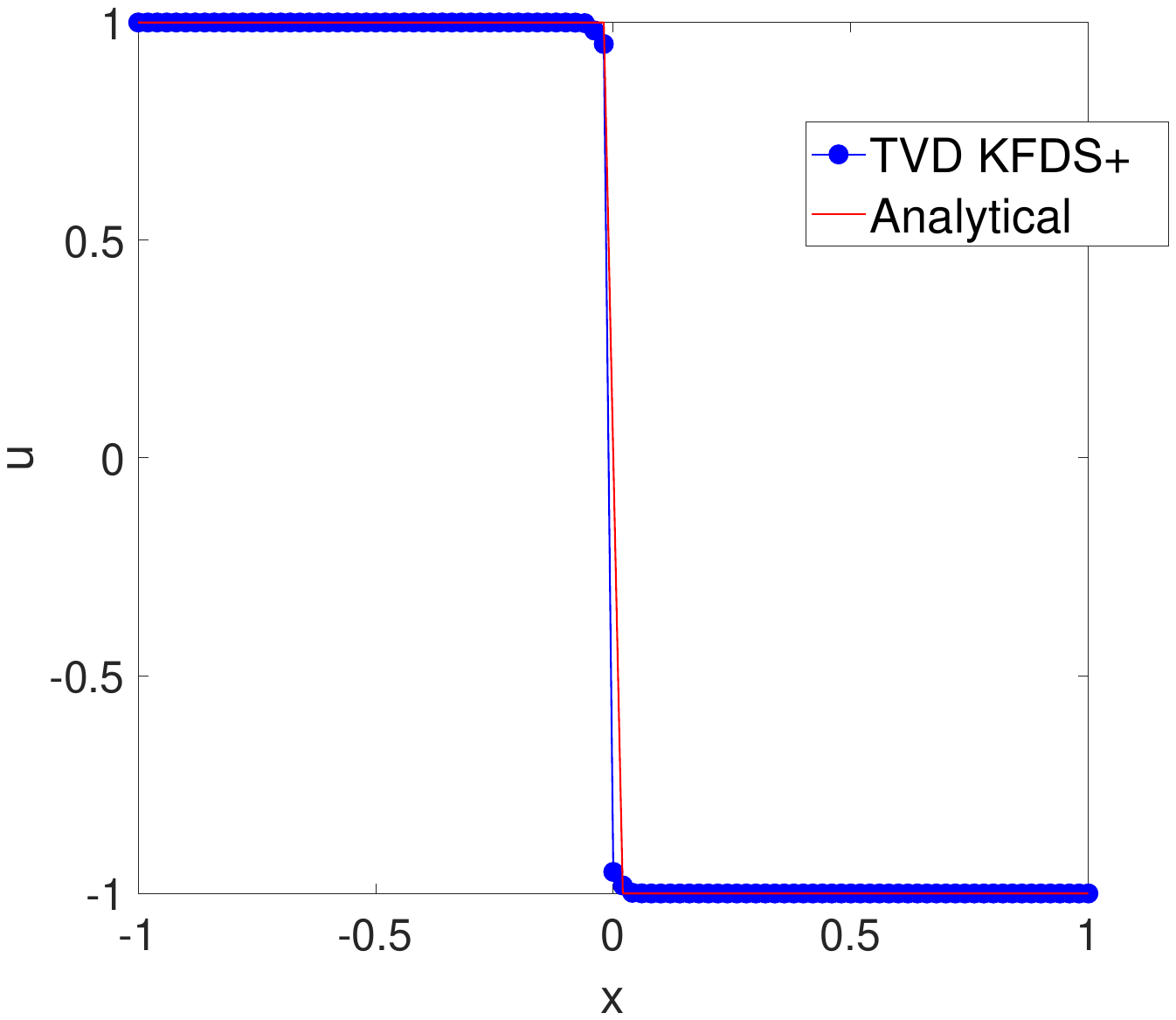} \\
 (d)First order KFDS+ &(e) KLW+  & (f)TVD-KFDS+ \\
\end{tabular}
\caption{Test Case 8(b) : KFDS, KLW \& TVD-KFDS  schemes in Viscous Burgers framework with 100 points}
\label{TC_8b_VISC_KFDS} 
\end{center} 
\end{figure}  
The results for the test case 8(a) are given in figure (\ref{TC_8a_VISC_KFDS}).  In this test case, a steady shock is present at $x=0$ which is gradually getting diffused into the domain due to the viscous flux. In comparison to inviscid test case, the presence of inherent viscous fluxes are handy in eliminating or damping any kind of oscillations that are typical of higher order schemes. The shock is captured very accurately in KFDS+ scheme and all the higher order versions in comparison to the exact solution. There is a very minor difference in the shock thickness in the KFDS scheme which diminishes on using finer grid points.  The results for the test case 8(b) are given in figure (\ref{TC_8b_VISC_KFDS}). Like test case 6(b), this test case also demonstrates the oscillatory nature of KLW scheme at large gradients, with excellent accuracy in smooth regions.  Thus the TVD scheme blended out of KFDS and KLW versions effectively eliminates the oscillations and results in solutions that are very accurate compared the first order case.  

\subsubsection*{Test case 9: 1D Dam break on variable depth bed}
This test problem involving  a variable river bed is considered to test the ability of scheme in resolving flow features with source terms in a one dimensional shallow water system, which is nonlinear and hyperbolic. The computational domain  $[-1, 1]$ is divided into 100 cells. The  initial discontinuity is located at $x= 0.5$.  The final computation time is 
$t = 0.1$.  The river bed profile, the initial conditions and the boundary conditions of the test case are as given below.  \\  
River bed Profile :
\bea
b(x) = \left\{ \ba{l}  \frac{1}{8}\left(\cos \left(10 \pi\left(x-\frac{1}{2}\right)\right)+1\right) \ \textrm{ for } \frac{2}{5} \leq x \leq \frac{3}{5}  \\ 
0 \ \  \   \  \  \ \   \  \  \  \textrm{Otherwise}    \ea \right. 
\eea
The initial condition is given by  $u(x,0) =0$ and the boundary condition is given by 
\bea
h(x,0) = \left\{ \ba{l}  1-b(x) \ \textrm{ if  } \ 0 \leq x \leq \frac{1}{2} \\ 
0.5-b(x) \ \textrm{if} \ \frac{1}{2}<x \leq 1    \ea \right. 
\eea  

\begin{figure}[!h] 
\begin{center} 
\begin{tabular}{cccc}
&\includegraphics[height=4.5cm]{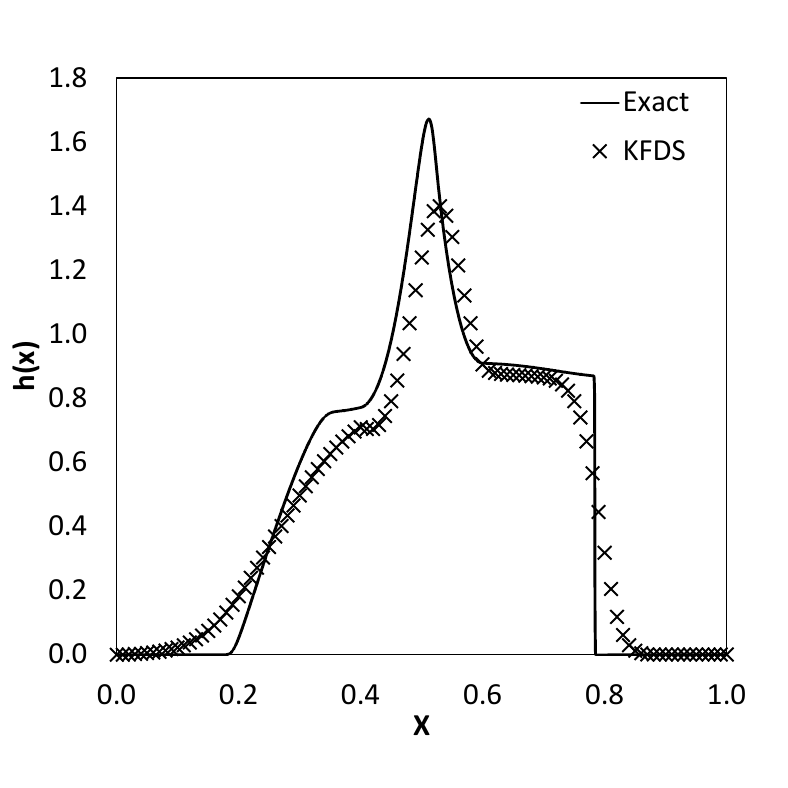} &
\includegraphics[height=4.5cm]{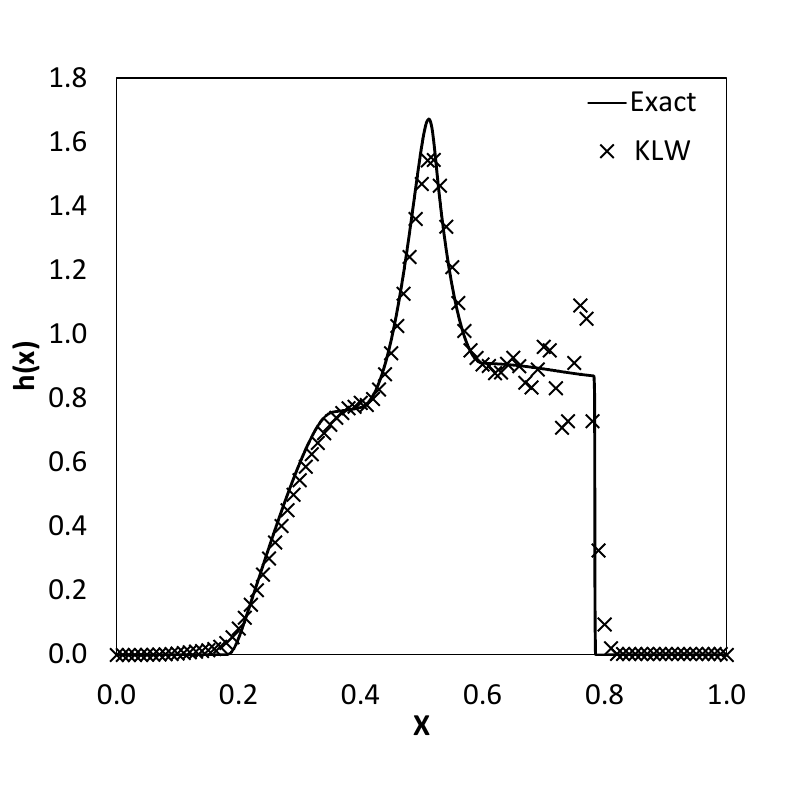} &
\includegraphics[height=4.5cm]{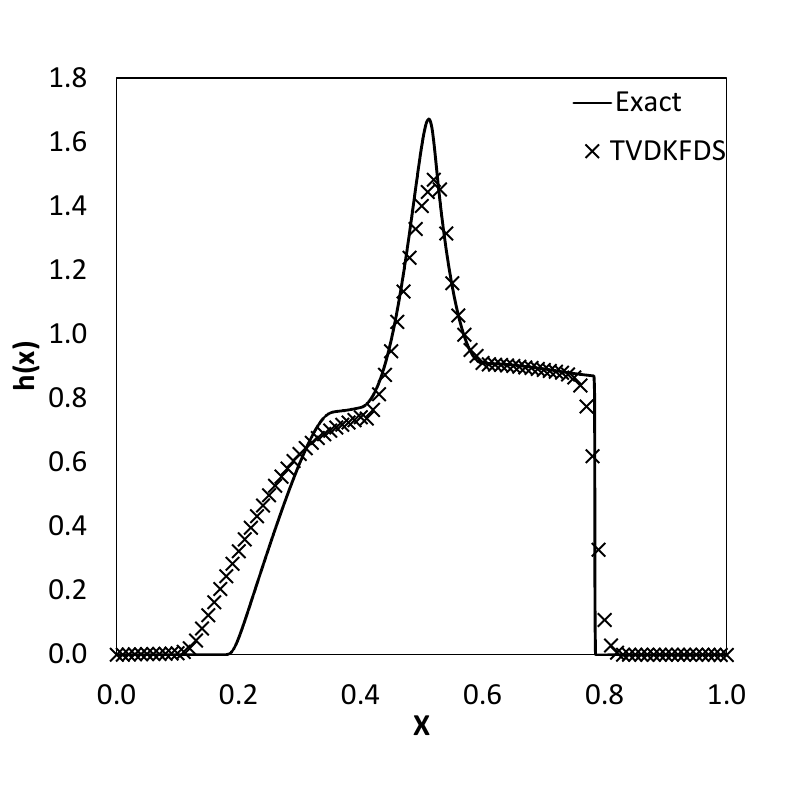} \\
 &(a)First order KFDS & (b)Second order KLW & (c)TVD-KFDS\\
&\includegraphics[height=4.5cm]{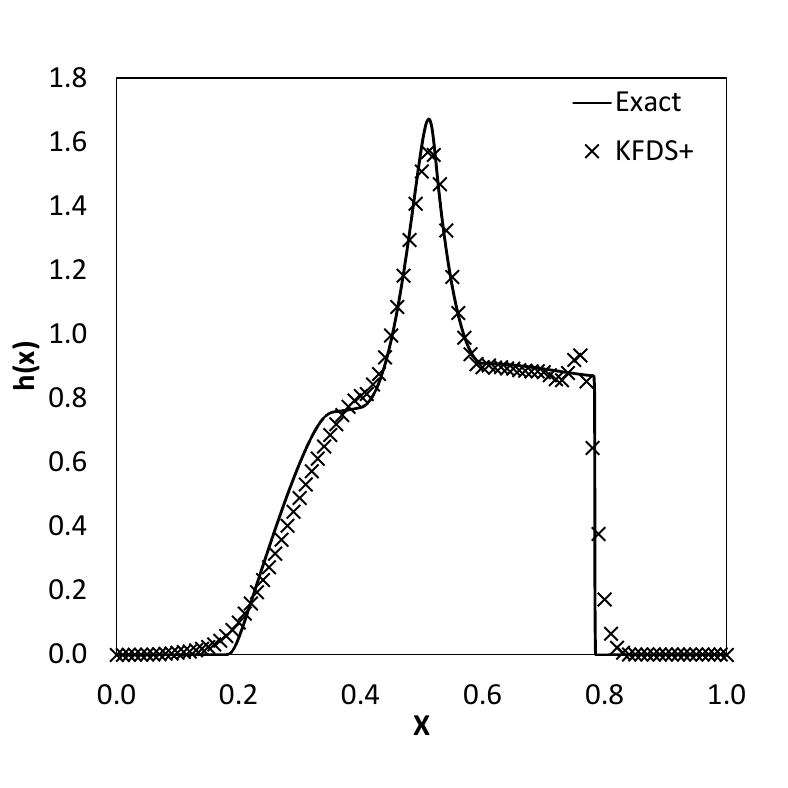} & &
\includegraphics[height=4.5cm]{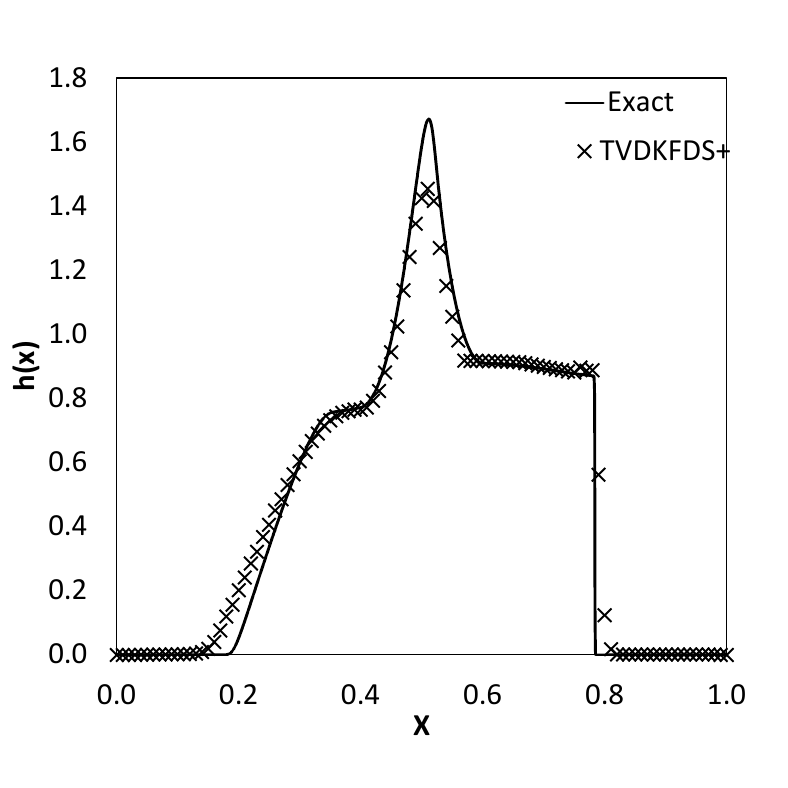} \\
& (d)First order KFDS+  &  & (e)TVD-KFDS+ \\
\end{tabular}
\caption{Test Case 9 : KFDS, KLW \& TVD-KFDS  schemes in Shallow Water framework  : Water Height}
\label{SWE_1D_KFDS_H} 
\end{center} 
\end{figure}

\begin{figure}[!h] 
\begin{center} 
\begin{tabular}{ccc}
\includegraphics[height=4.5cm]{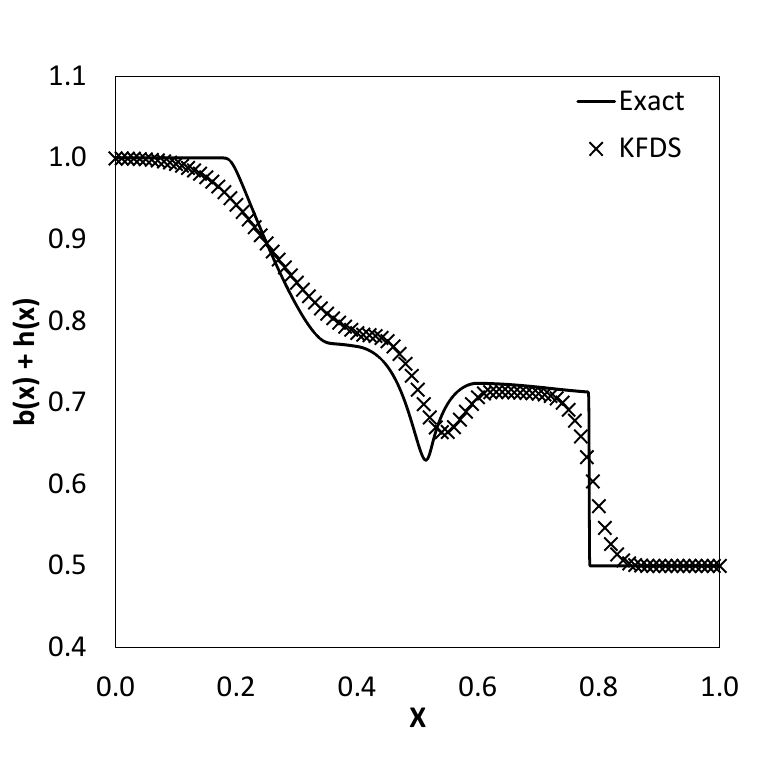} &
\includegraphics[height=4.5cm]{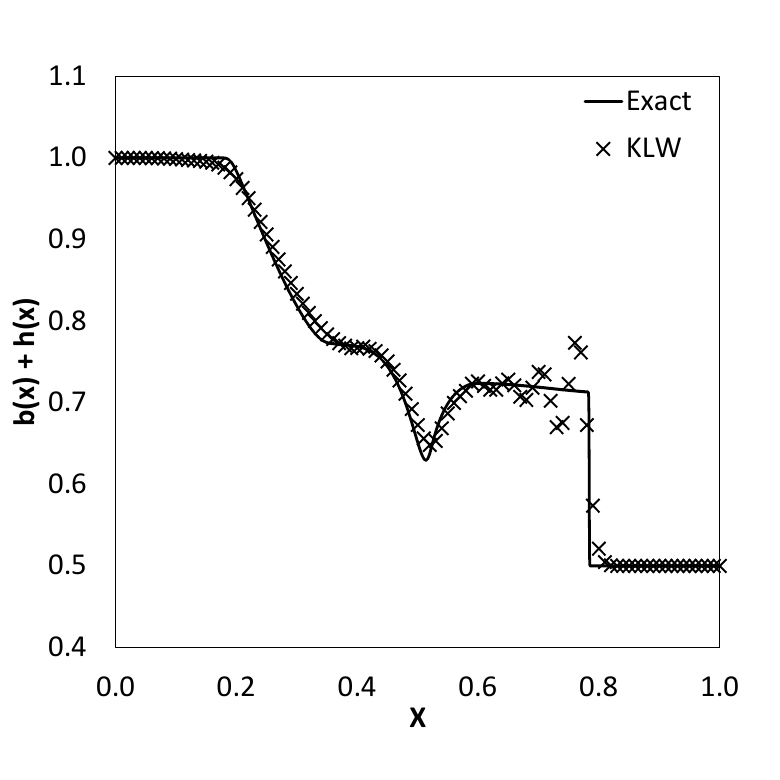} &
\includegraphics[height=4.5cm]{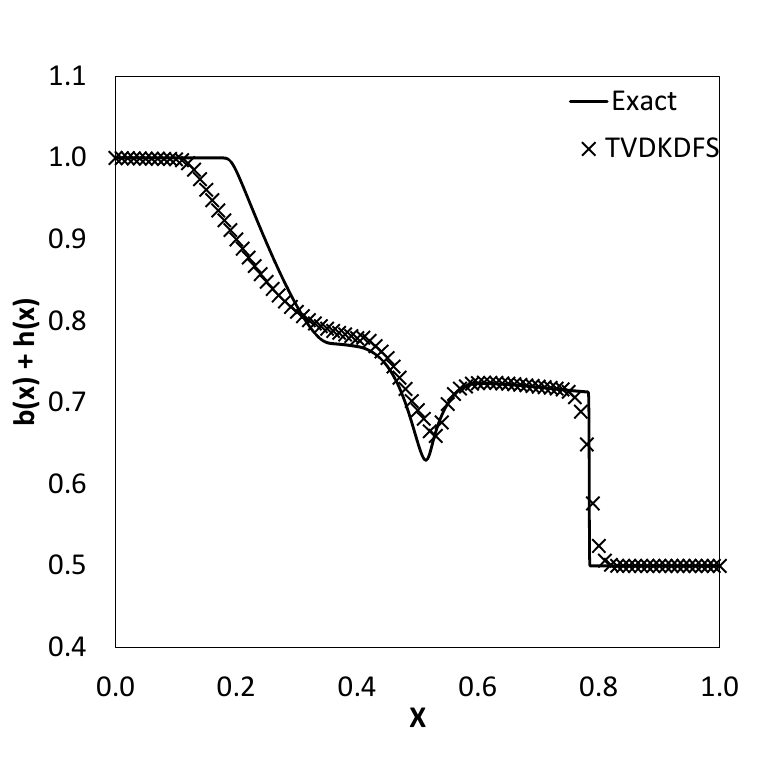} \\
 (a)First order KFDS & (b)Second order KLW & (c)TVD-KFDS\\
\includegraphics[height=4.5cm]{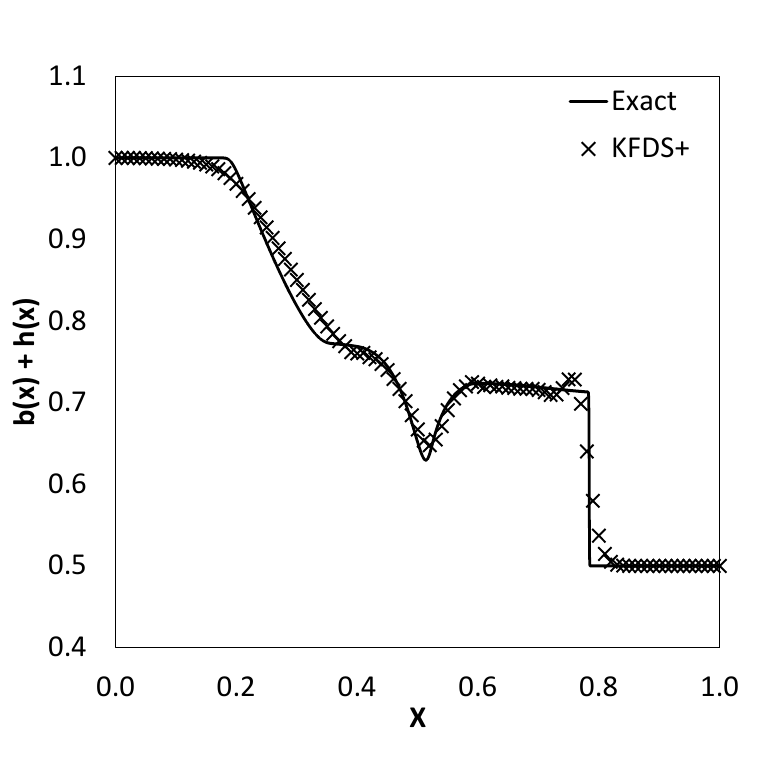} & &
\includegraphics[height=4.5cm]{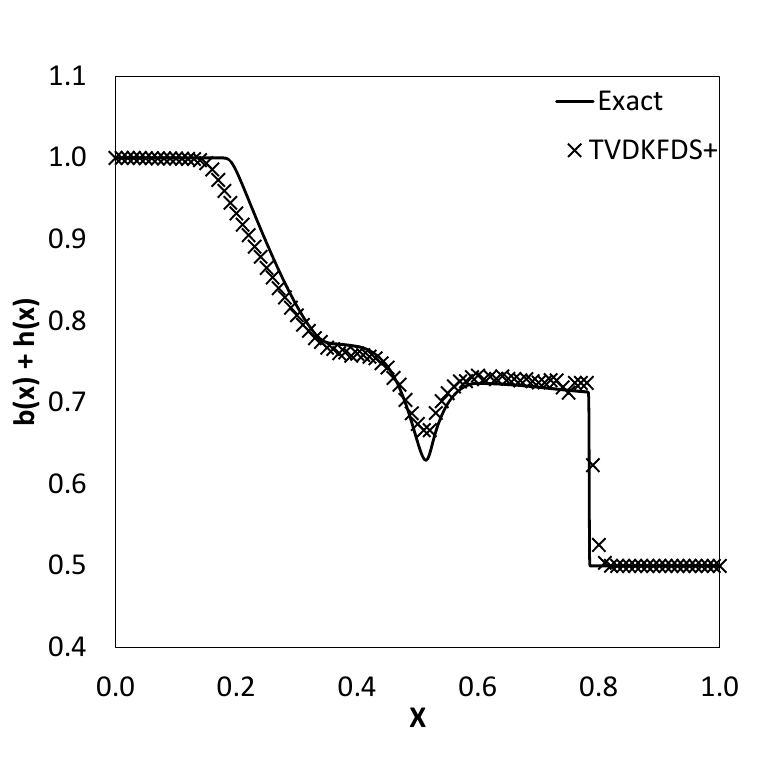} \\
 (d)First order KFDS+  &  & (e)TVD-KFDS+ \\
\end{tabular}
\caption{Test Case 9 : KFDS, KLW \& TVD-KFDS  schemes in Shallow Water framework : Water + Base Height}
\label{SWE_1D_KFDS_hb} 
\end{center} 
\end{figure}

This test case is for simulating the shallow water system of equations. The results are given in Figures (\ref{SWE_1D_KFDS_H}) \& (\ref{SWE_1D_KFDS_hb}). The presence of the variation in the river bed adds complexity in simulation of this test case in the form of an additional source term. The first order KFDS scheme, although diffusive,  gives a reasonable capture of the flow gradients including the shock region. The solution improves with the KFDS+ version. However a marginal overshoot is observed in the pre-shock region.  Although the second order KLW scheme captures the gradients in the flow better than the KFDS scheme, the solution evolves with a pre-shock osciilation.  Both the TVD schemes show significant improvement in the solution as compared to KFDS variants and are devoid of undue osillations. 

The subsequent test cases (Test cases 10 to 13) are two dimensional in nature and are to be solved on a 64x64 grid for a standard comparison. 

\subsubsection*{Test case 10: 2D  Linear Convection Equation - Diagonal Discontinuity} 
This test case is designed for 2D linear convection equation \cite{Spekreijse} given by 
\be
\frac{\partial u}{\partial t}+a \frac{\partial u}{\partial x}+b \frac{\partial u}{\partial y}=0;  \     \   a = cos\phi , \     \  b = sin\phi
\ee 
The test case is to be solved in a square domain  $x:[0,1]$, $y:[0,1]$ with $\phi = 45^{o}$. The boundary conditions for this test case are given as  
\bea
\left\{ \begin{array} { l l } { u (0,y) = 1 , } & {0 < y < 1 } \\ { u(x,0) = 0 , } & {0 < x <  1 }  \end{array} \right.
\eea 
The exact solution for this test case is given by 
\bea
\left\{ \begin{array} { l l } { u_{exact} (x,y) = 1 , } & {if  \ bx - ay < 0 } \\ { u_{exact} (x,y) = 0 , } & {if \ bx - ay > 0 }  \end{array} \right.
\eea

\begin{figure}[!h] 
\begin{center} 
\begin{tabular}{cccc}
\includegraphics[height=4.2cm]{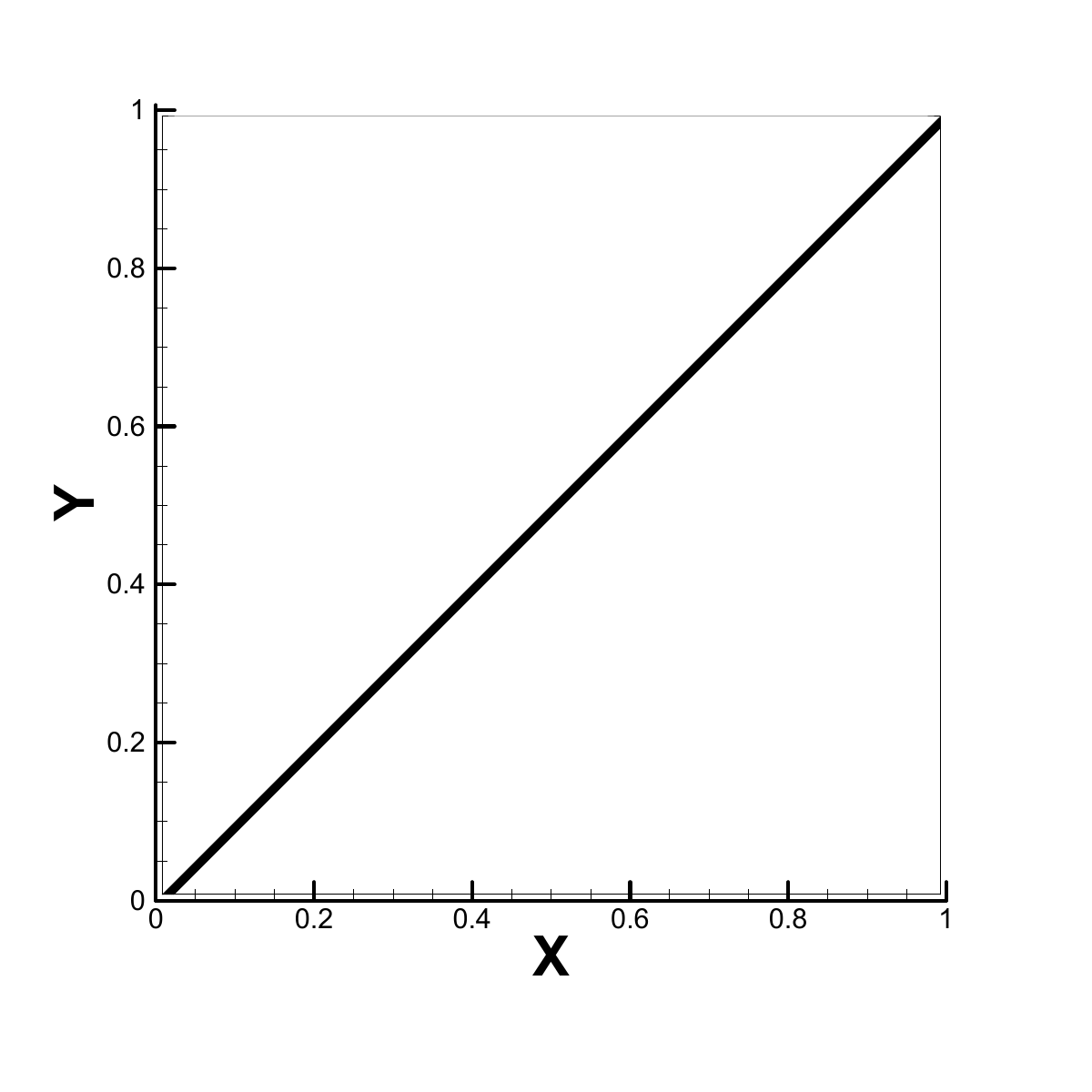} &
\includegraphics[height=4.2cm]{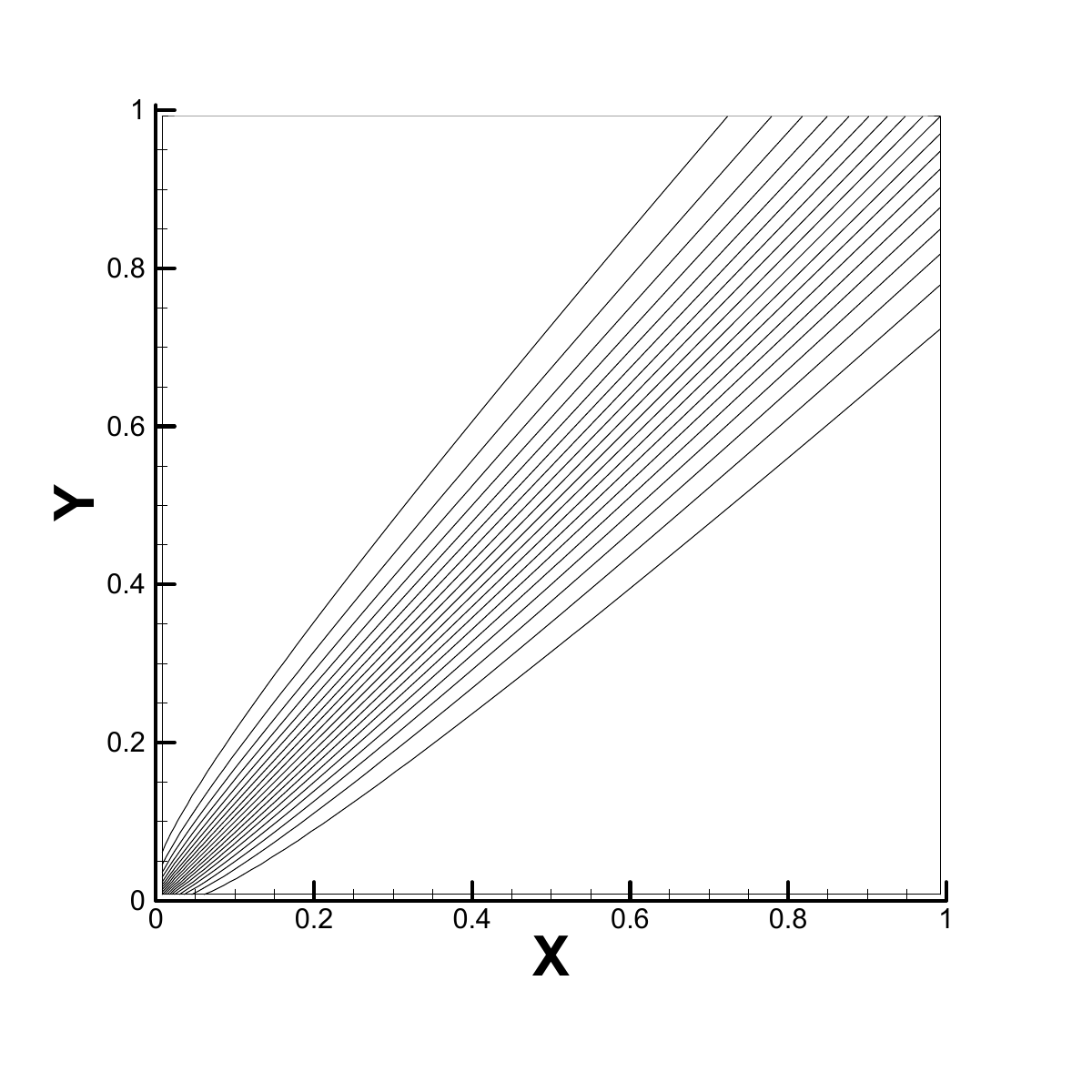} &
\includegraphics[height=4.2cm]{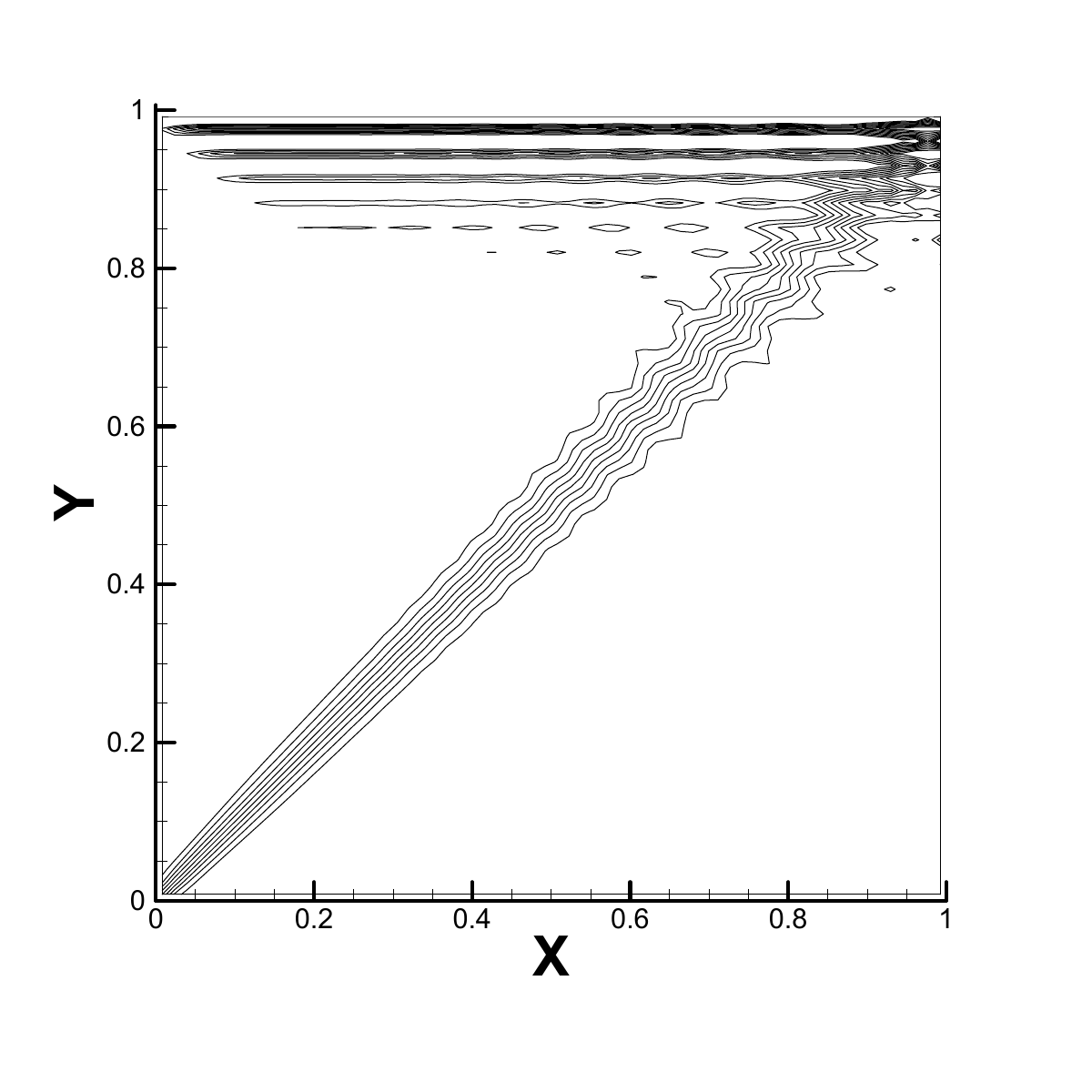} &
\includegraphics[height=4.2cm]{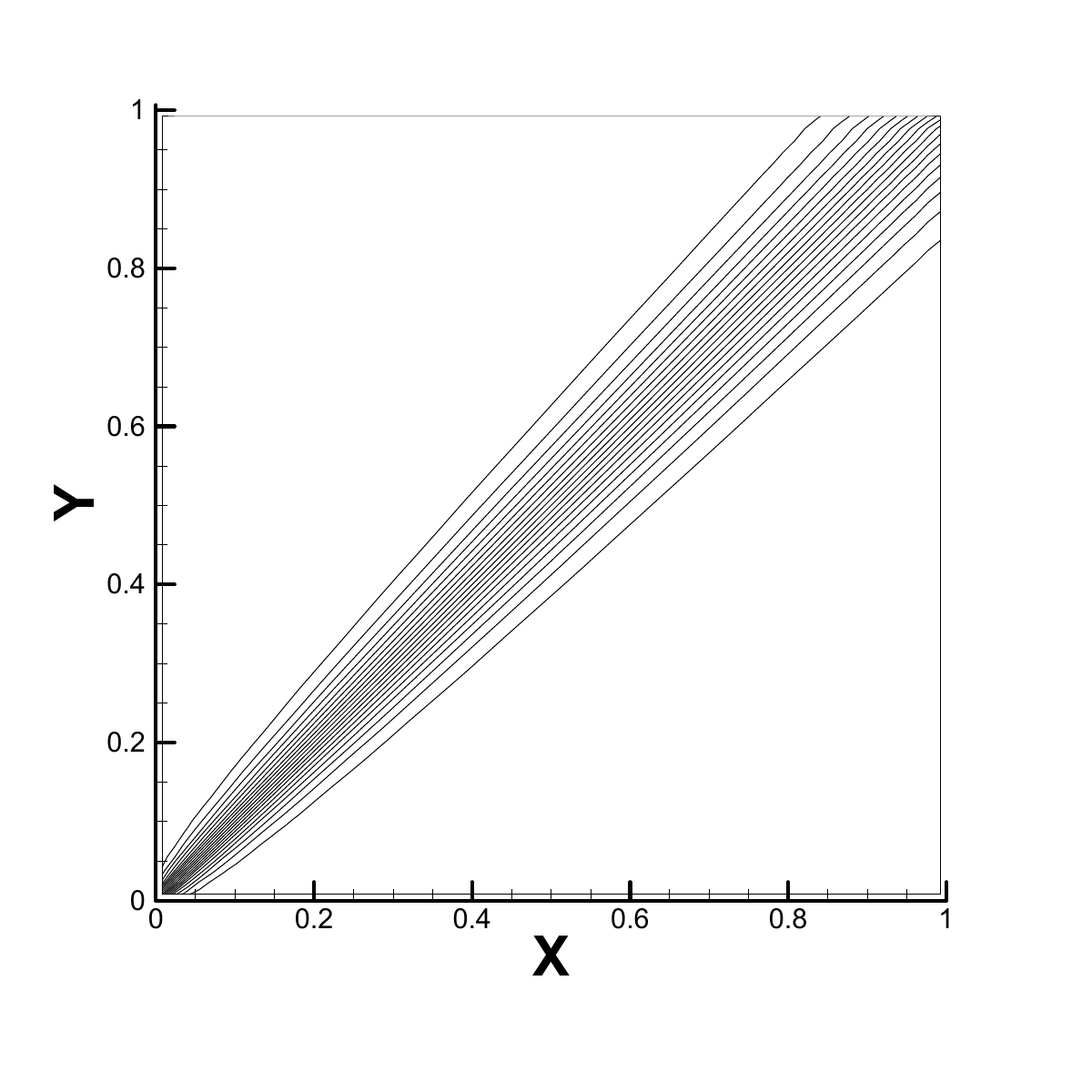} \\
(a)Exact solution &    (b)First order KFDS & (c)KLW  &   (d)TVD-KFDS \\
\end{tabular}
\caption{Test Case 10 : KFDS, KLW \& TVD-KFDS  schemes in LCE framework with 64x64 points}
\label{TC_2D_LCE_KFDS} 
\end{center} 
\end{figure}

The results for this test case are given in figure (\ref{TC_2D_LCE_KFDS}). The given initial conditions for this 2-D linear convection equation results in a steady discontinuity along the diagonal of the domain.  The test cases have been simulated on a $64 \times 64$ grid. The KFDS scheme captures the position of the contact discontinuity with some numerical diffusion.  The KLW scheme demonstrates lower diffusion in comparison to KFDS scheme but exhibits oscillations. 
The TVD scheme blends the ability of both KFDS and KLW schemes and provides an oscillation-free solution better than KFDS scheme.  

\subsubsection*{Test case 11: 2D Burgers equation - steady normal shock and a smooth region}

This test case is designed for 2D inviscid Burgers equation \cite{Spekreijse}. The initial conditions as well as the boundary conditions for this test case are as given below. 
\bea
\left\{ \begin{array} { l l } { u ( 0 , y ) = 1 , } & { 0 < y < 1 } \\ { u ( 1 , y ) = - 1 , } & { 0 < y < 1 } \\ { u ( x , 0 ) = 1 - 2 x , } & { 0 < x < 1 } \end{array} \right.
\eea
The test case is solved in a square domain $x:[0,1]$, $y:[0,1]$.  The exact solution for this test case \cite{Spekreijse} consists of a normal shock in the middle of the domain, from the top boundary to point (0.5,0.5) and a smooth variation of the solution below that point till the bottom boundary.   

\begin{figure}[!h] 
\begin{tabular}{cccc}
\includegraphics[height=4.2cm]{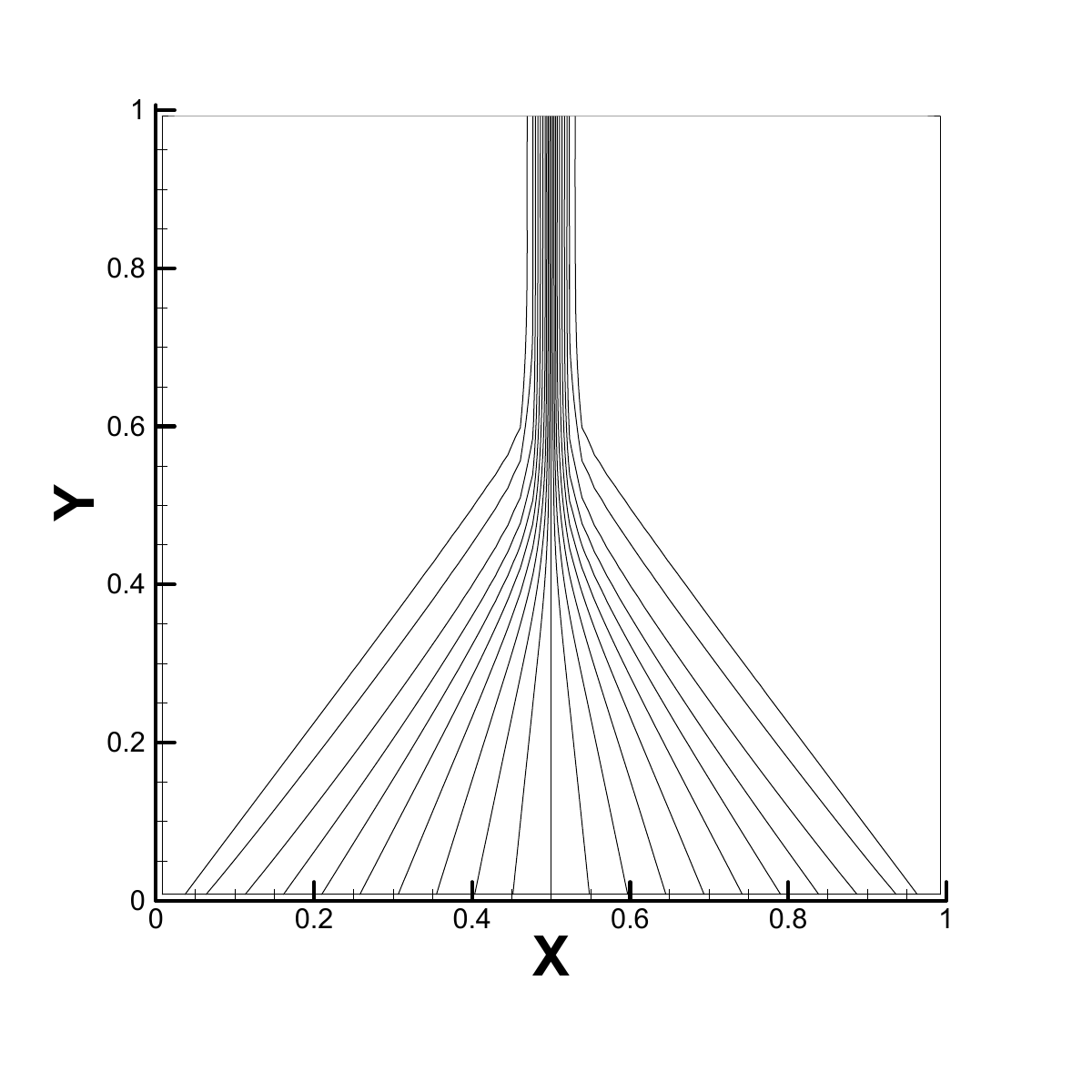} &
\includegraphics[height=4.2cm]{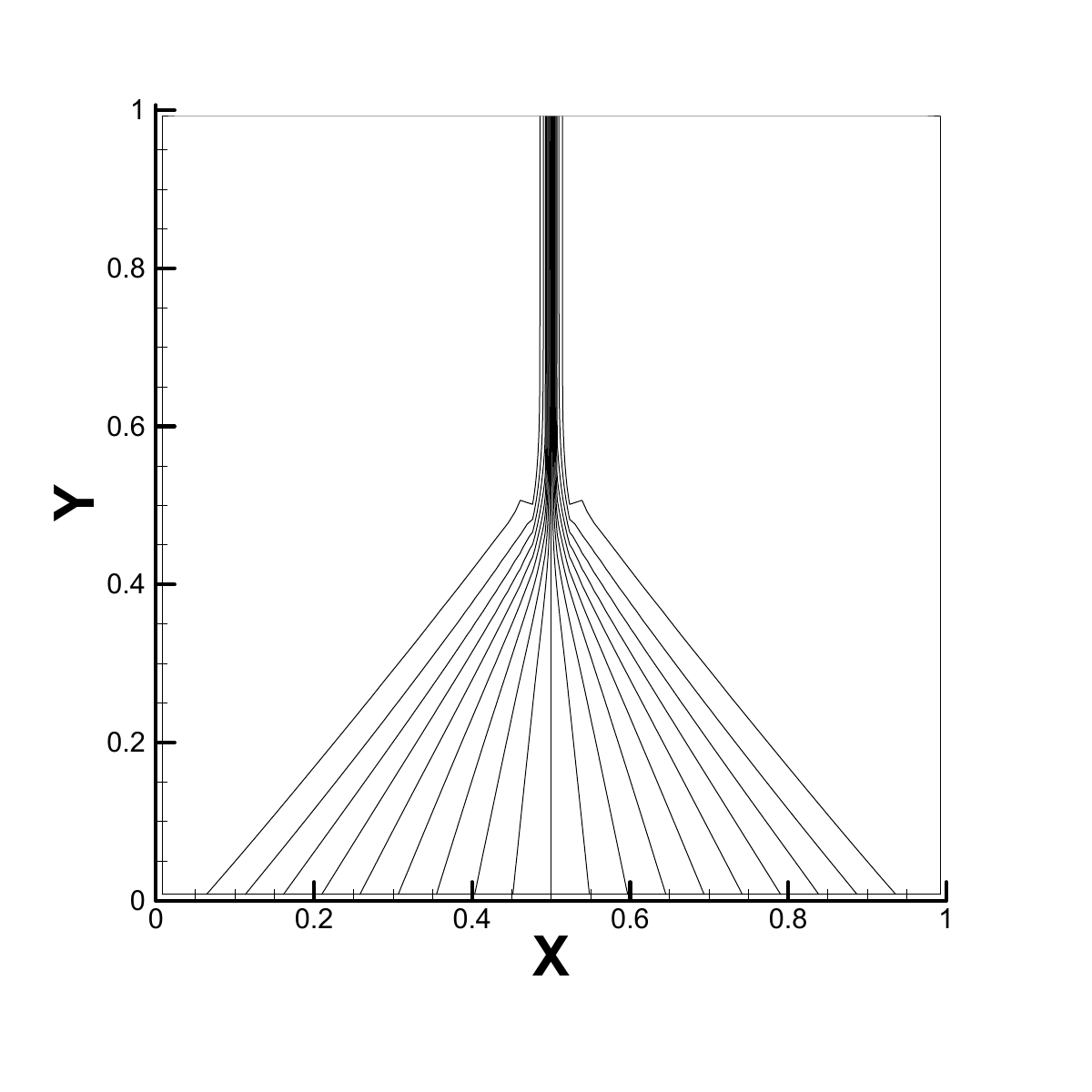} &
\includegraphics[height=4.2cm]{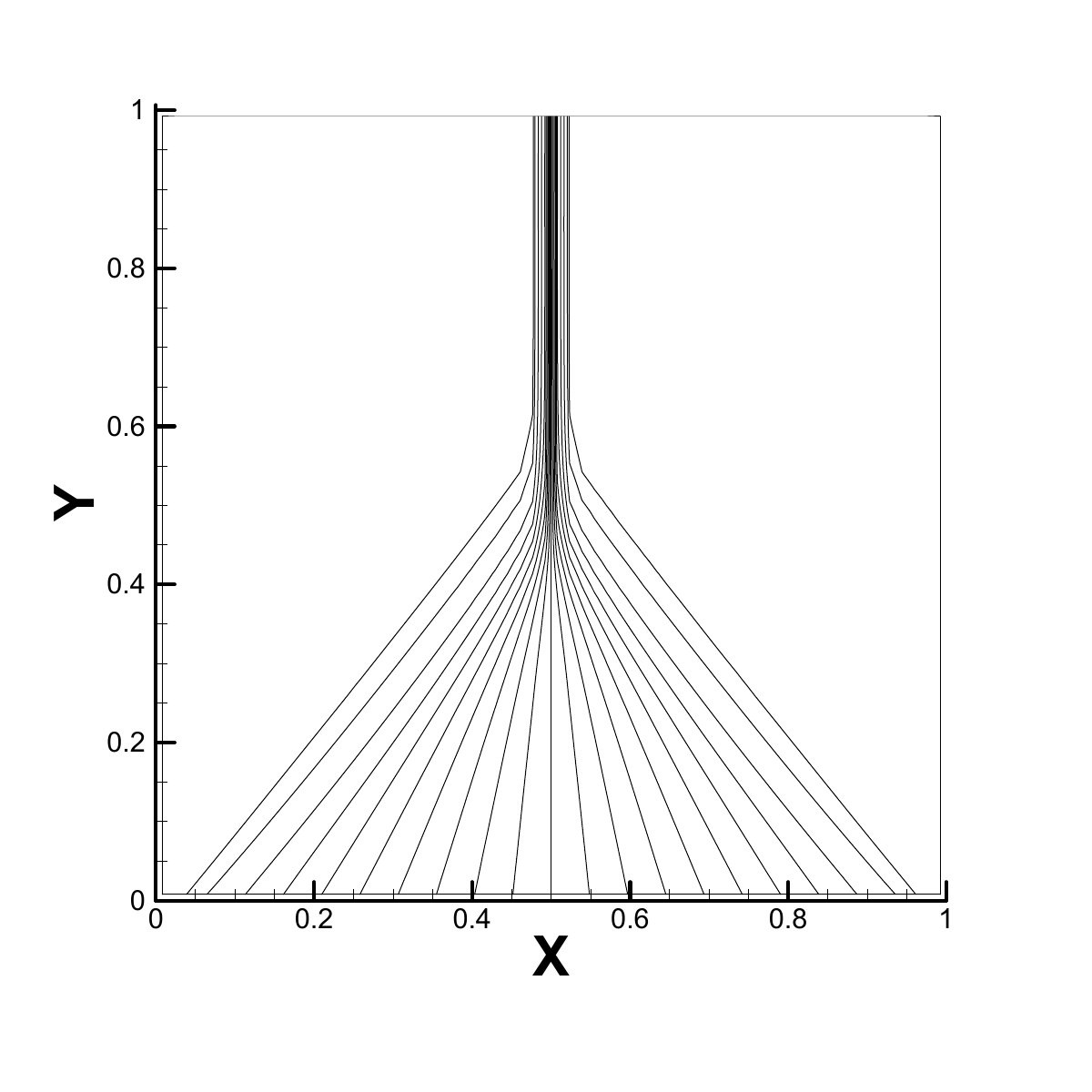}&
\includegraphics[height=4.2cm]{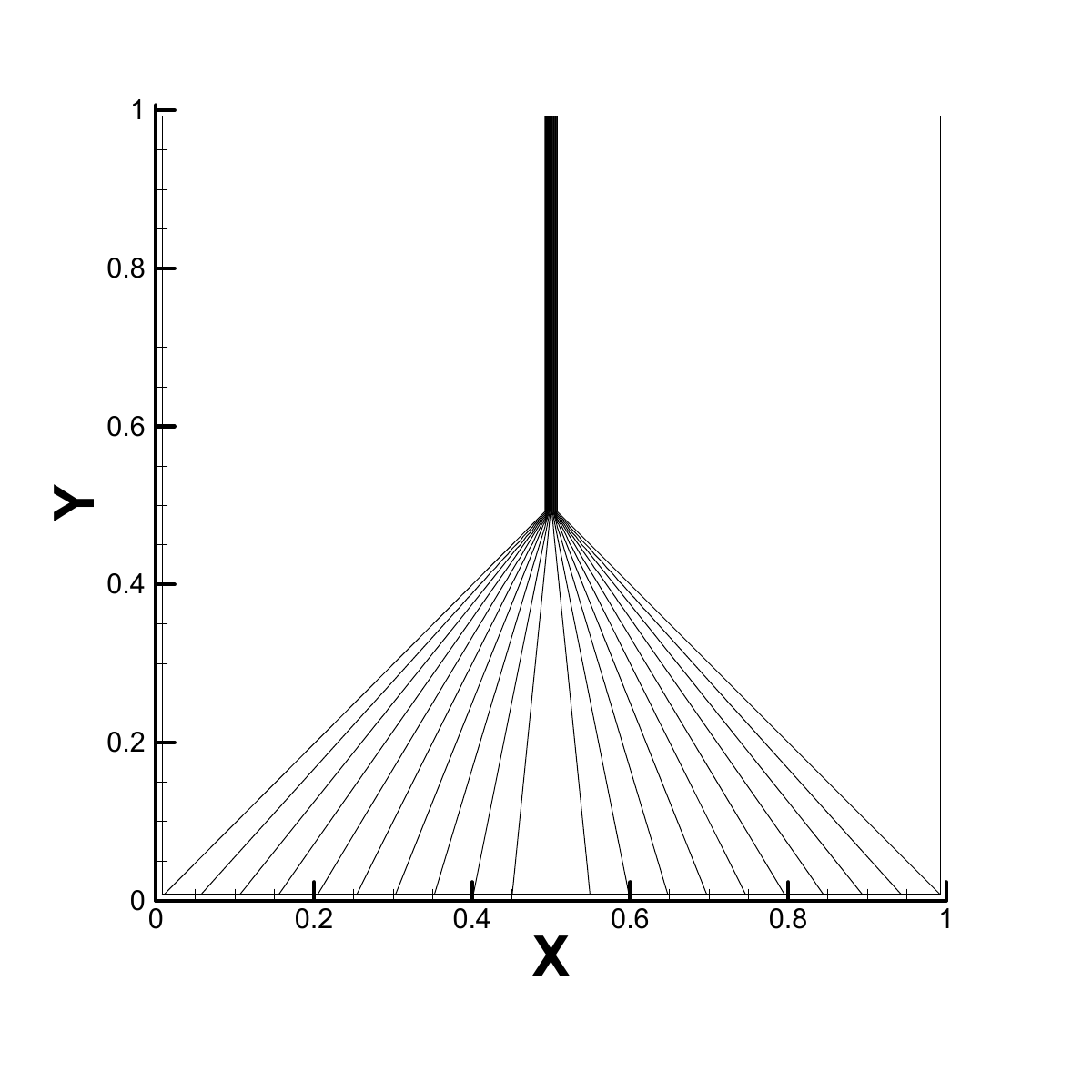}  \\
 (a)First order KFDS & (b)Second order KLW & (c)TVD-KFDS & (d) Exact solution \\
\includegraphics[height=4.2cm]{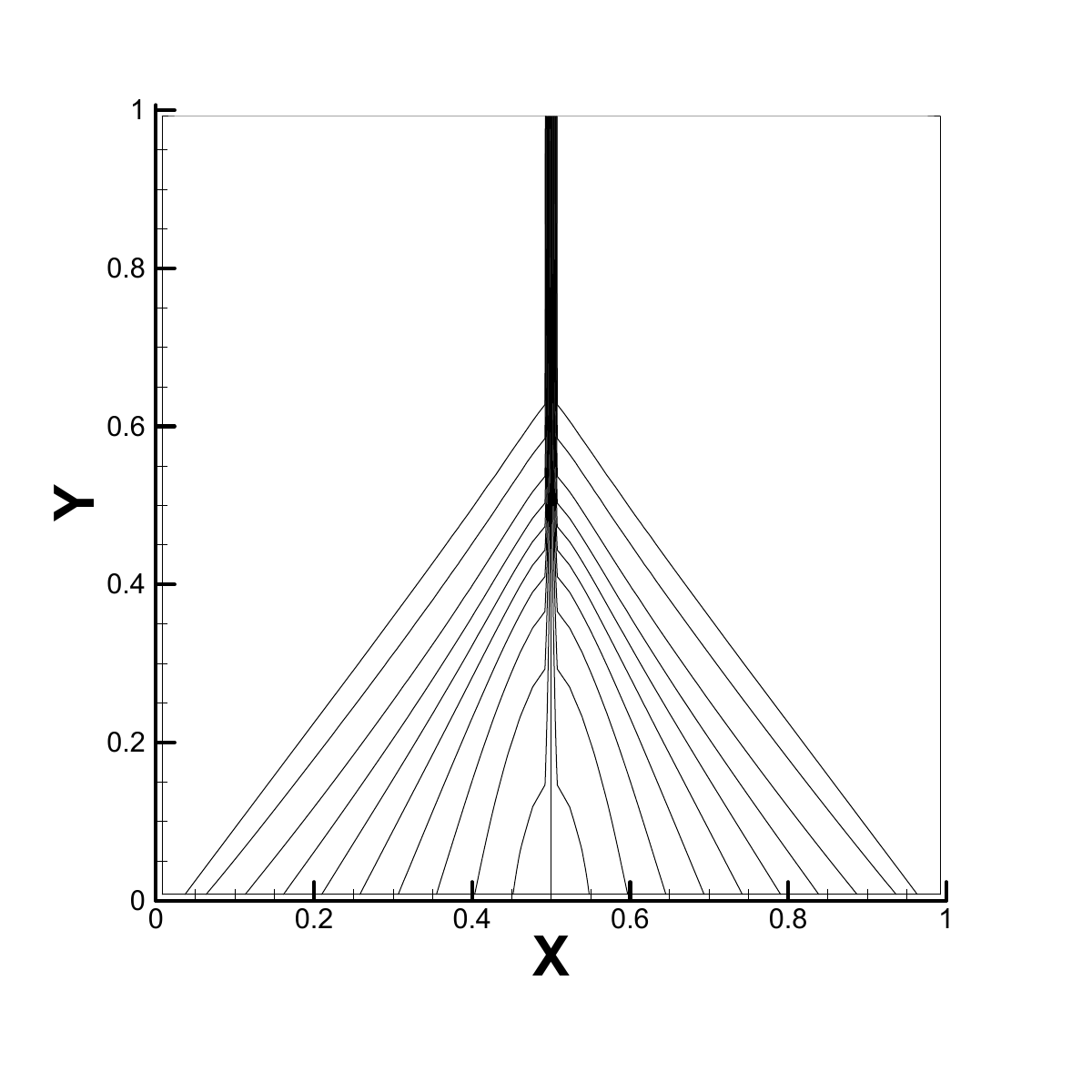} &
\includegraphics[height=4.2cm]{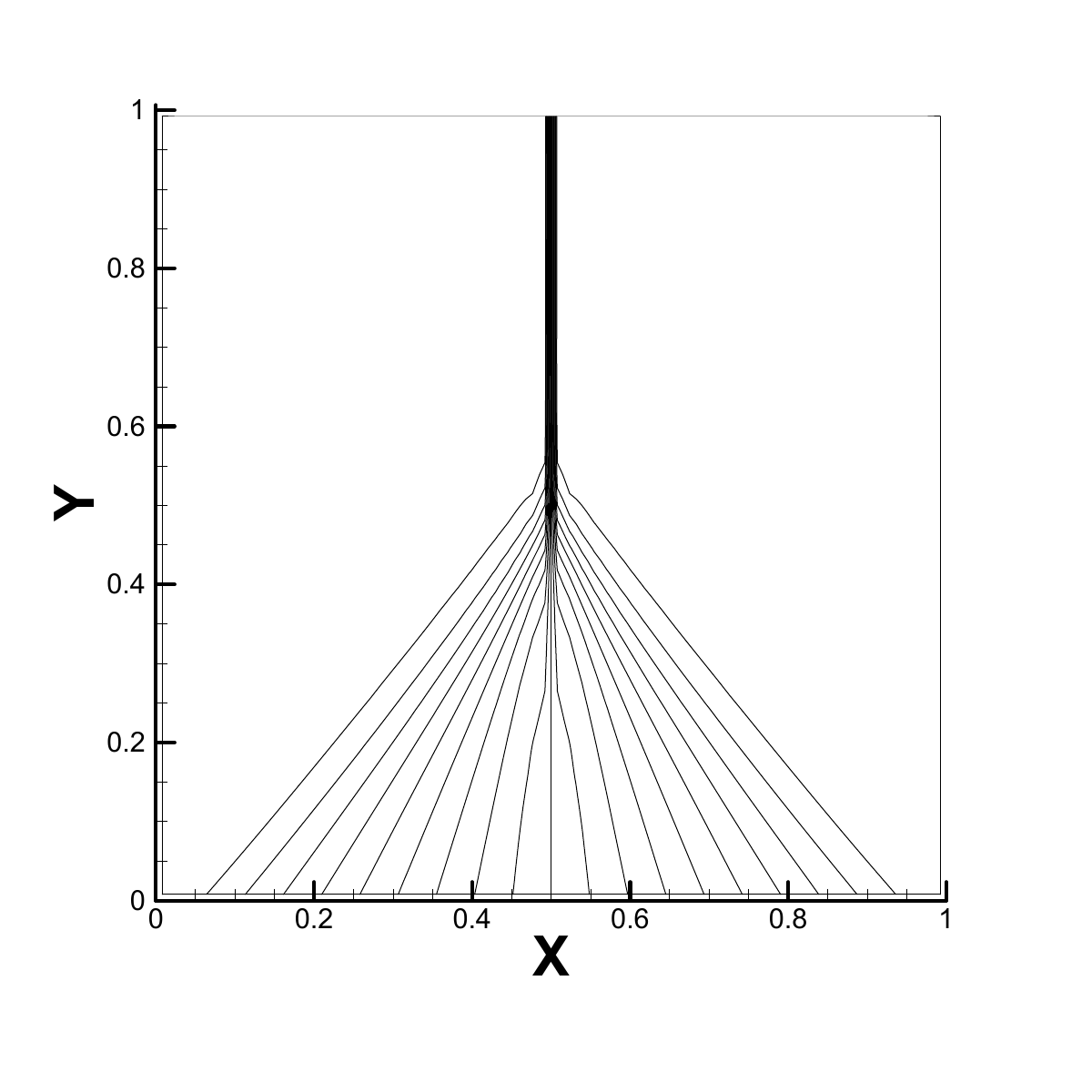} &
\includegraphics[height=4.2cm]{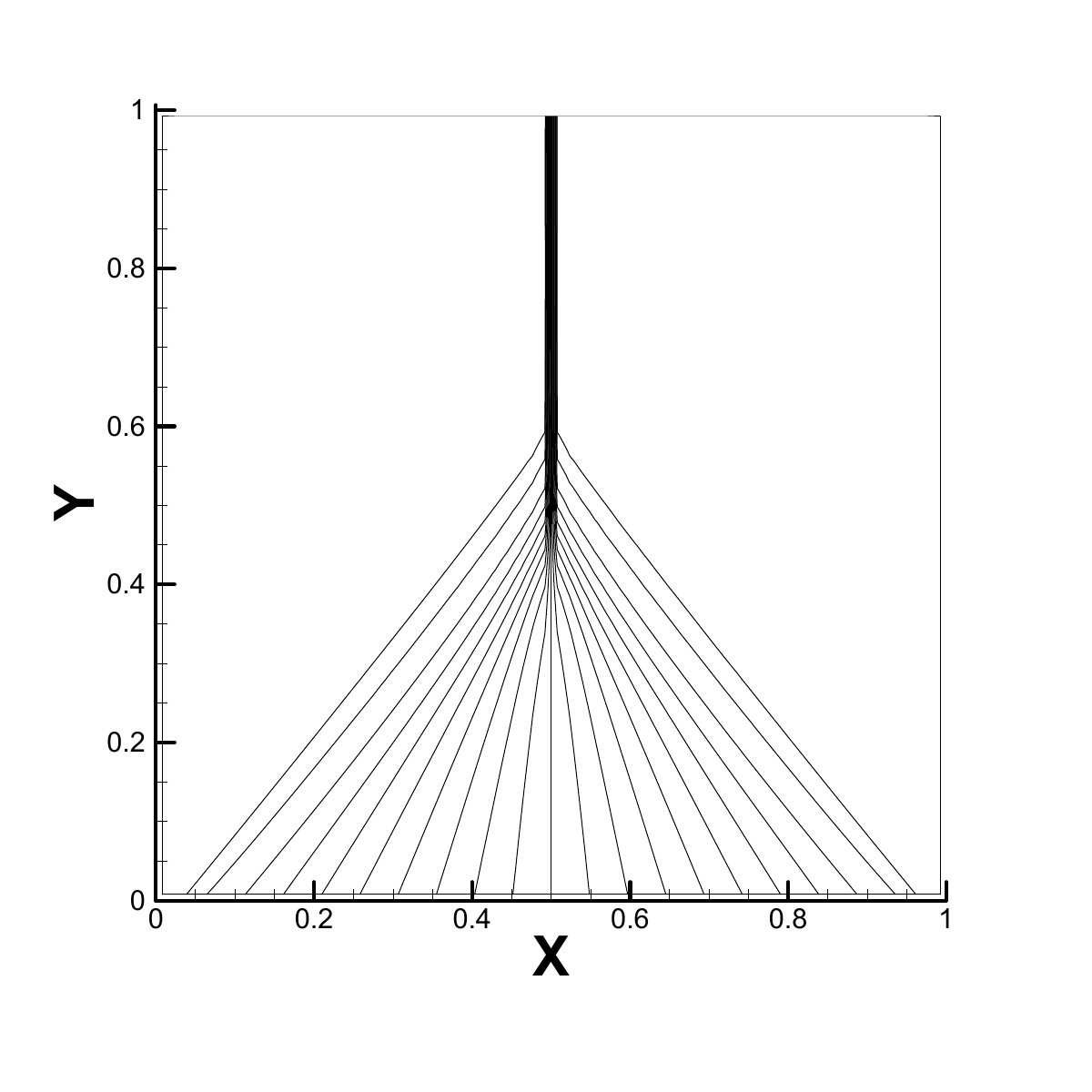} &\\
 (d)First order KFDS+ &(e) KLW+  & (f)TVD-KFDS+& \\
\end{tabular}
\caption{Test Case 11 : KFDS, KLW \& TVD-KFDS  schemes in 2D Burgers framework with 64 x64  points}
\label{TC_11_2D_KFDS} 
\end{figure}  

The results for the test case 11 are given in figure (\ref{TC_11_2D_KFDS}). The initial conditions and the boundary conditions for this test case lead to the formation of a normal shock located between (0.5,0.5) to (0.5,1.0) and a smooth variation located right beneath it. The exact solution for this test case consists of a normal shock and a smooth variation. The first order schemes, KFDS and KFDS+, are able to capture the stationary normal shock along with the expansion fan.  However, the shock appears diffused in the KFDS scheme in comparison to the exact solution. The KFDS+ scheme is able to capture the shock with good accuracy but generates an expansion shock in the smooth region. The second order schemes, KLW and KLW+, are able to capture both the shock and expansion waves with reasonable accuracy. However, they do exhibit oscillations near the shock. The TVD schemes demonstrate the ability to maintain the good features of both KFDS and KLW schemes and exhibit good capturing of both the shock and smooth regions sans oscillations.

\subsubsection*{Test case 12: 2D Burgers equation - Steady oblique shock and a smooth region}
This test case is designed for 2D inviscid Burgers equation \cite{Spekreijse}. The boundary conditions for this test case are as given below. 
\bea
\left\{ \begin{array} { l l } { u ( 0 , y ) = 1.5 , } & { 0 < y < 1 } \\ { u ( 1 , y ) = - 0.5 , } & { 0 < y < 1 } \\ { u ( x , 0 ) = 1.5 - 2 x , } & { 0 < x < 1 } \end{array} \right.
\eea
The exact solution for this test case \cite{Spekreijse} is similar to the previous test case, except that the shock is oblique to the grid, starting from the top-right boundary to the point (0.75,0.5).  

\begin{figure}[!h] 
\begin{center} 
\begin{tabular}{cccc}
\includegraphics[height=4.2cm]{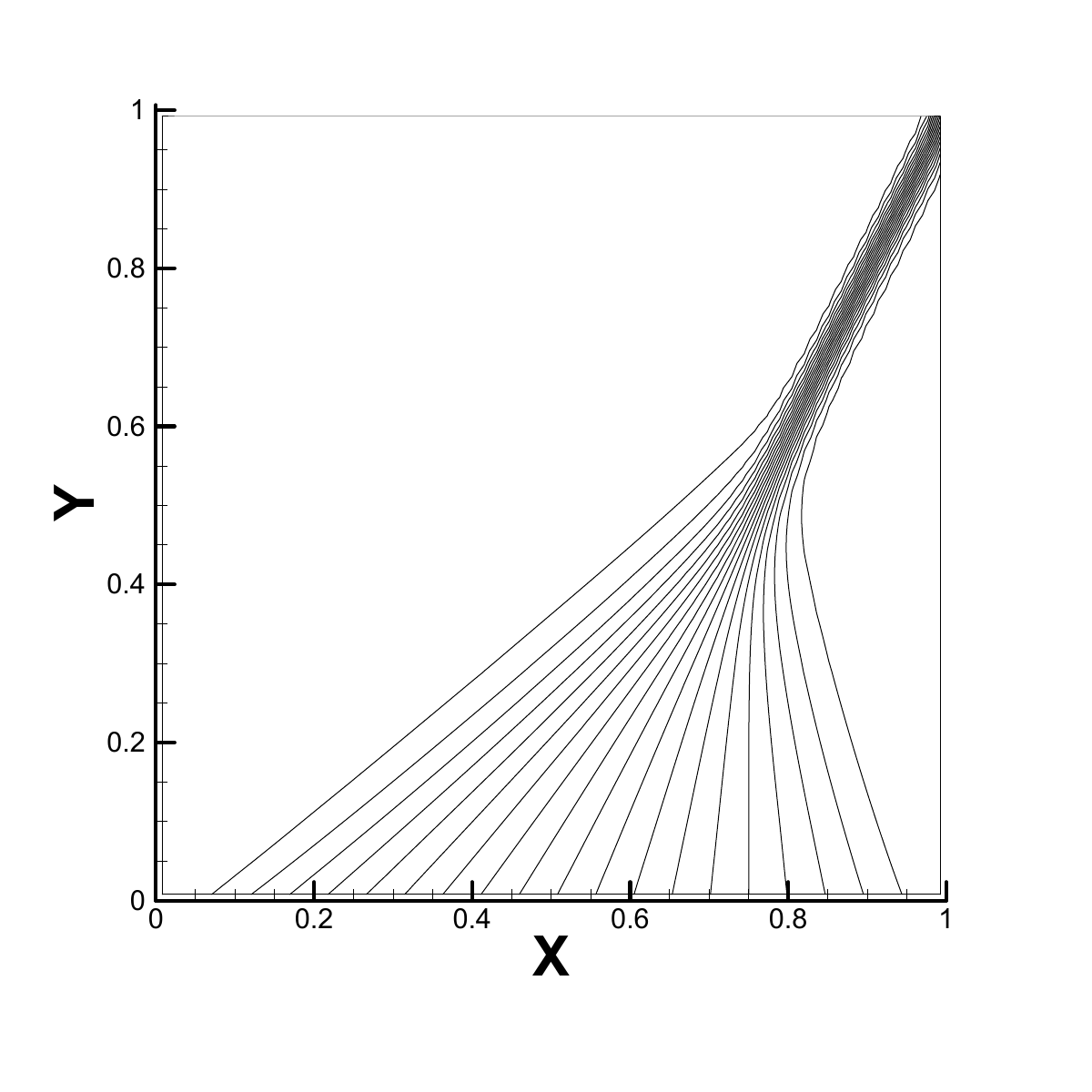} &
\includegraphics[height=4.2cm]{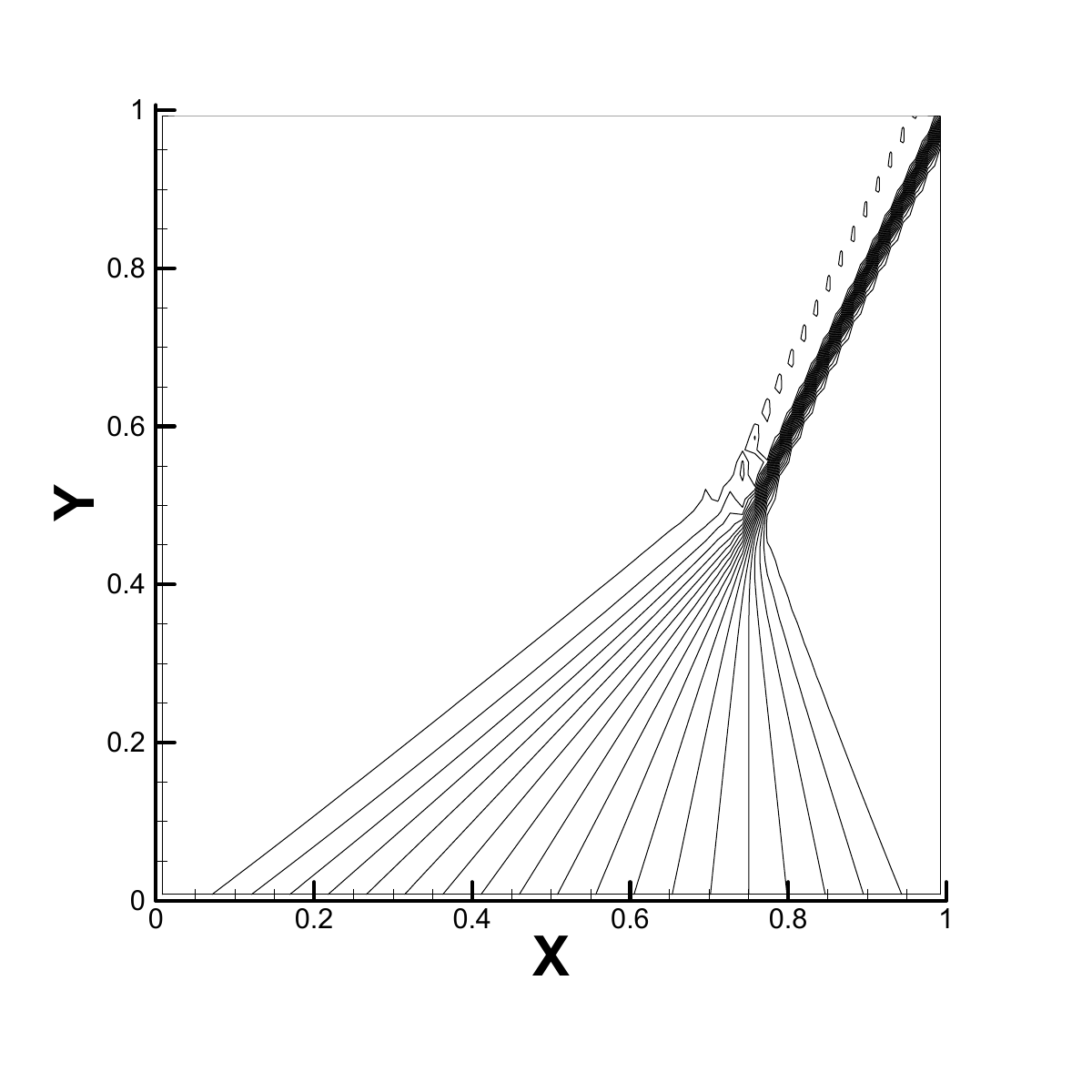} &
\includegraphics[height=4.2cm]{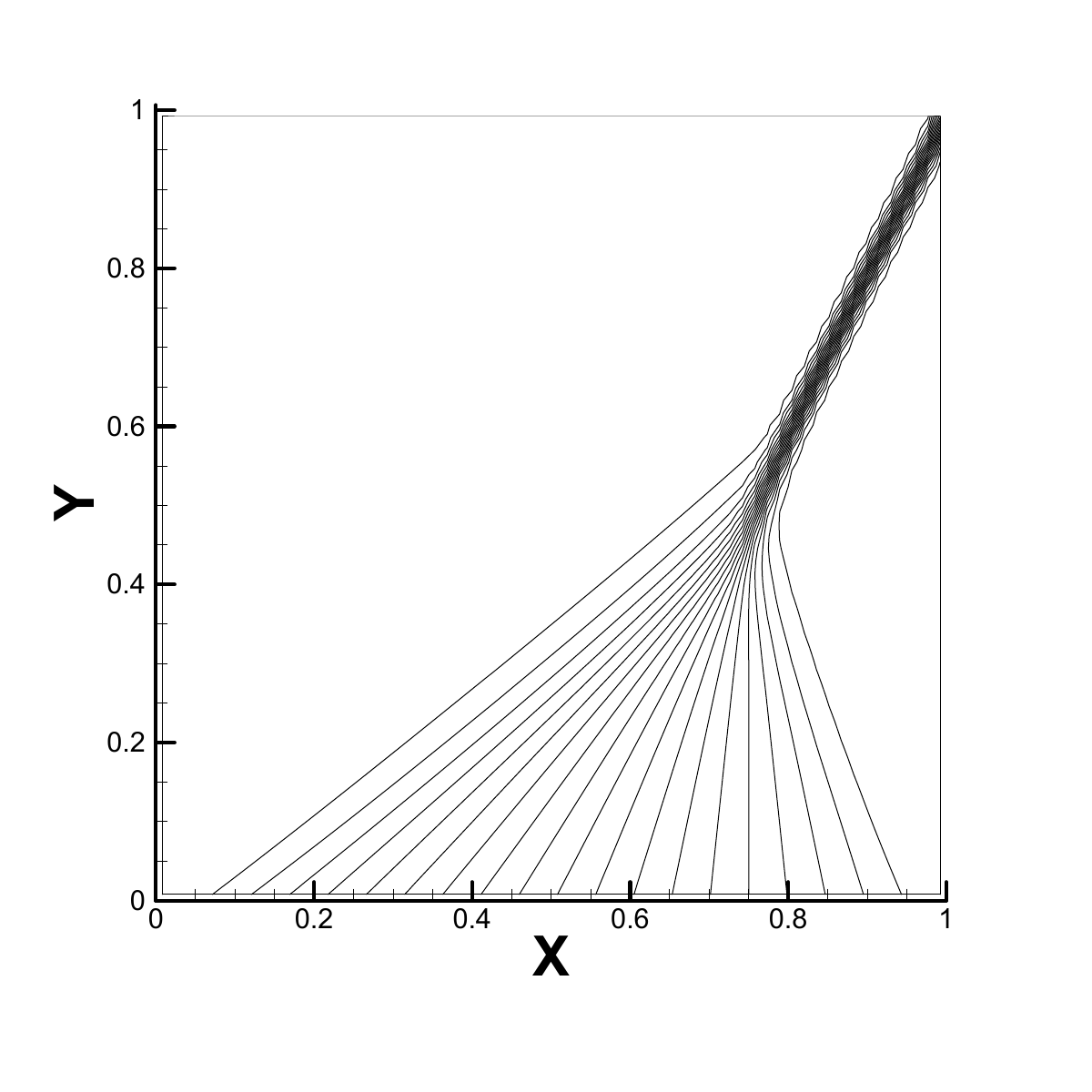}&
 \includegraphics[height=4.2cm]{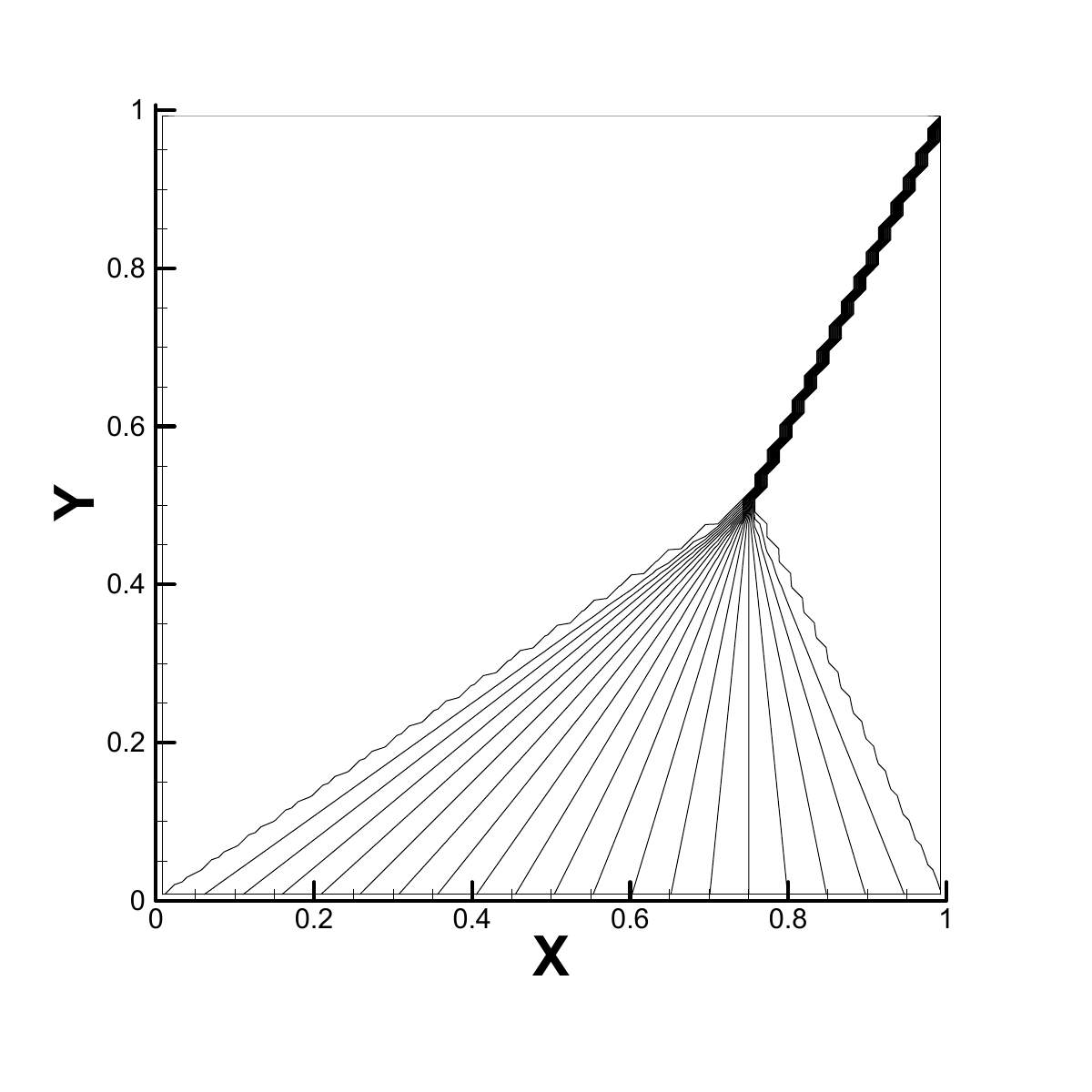}\\
 (a)First order KFDS & (b)Second order KLW & (c)TVD-KFDS & (d) Exact solution \\
\includegraphics[height=4.2cm]{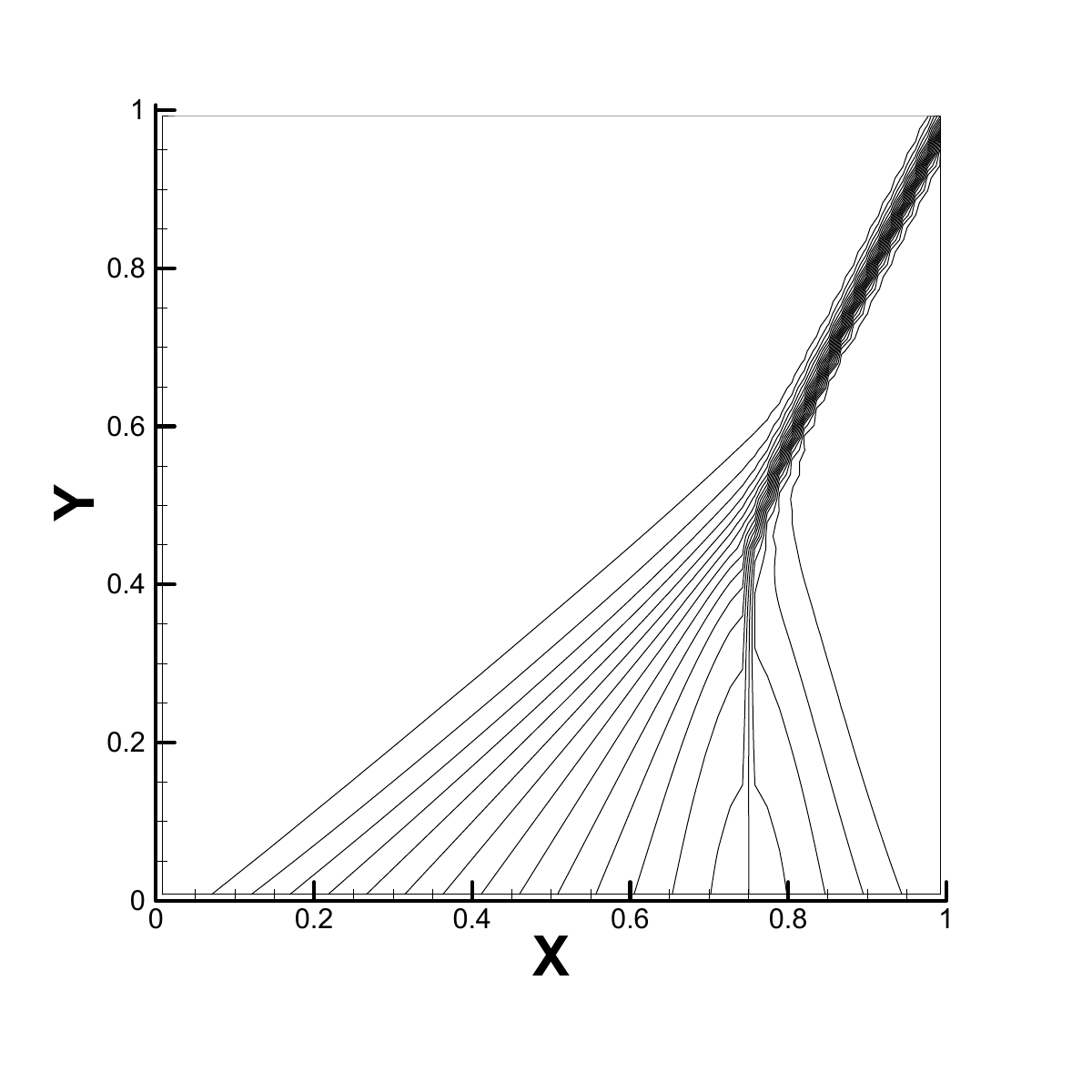} &
\includegraphics[height=4.2cm]{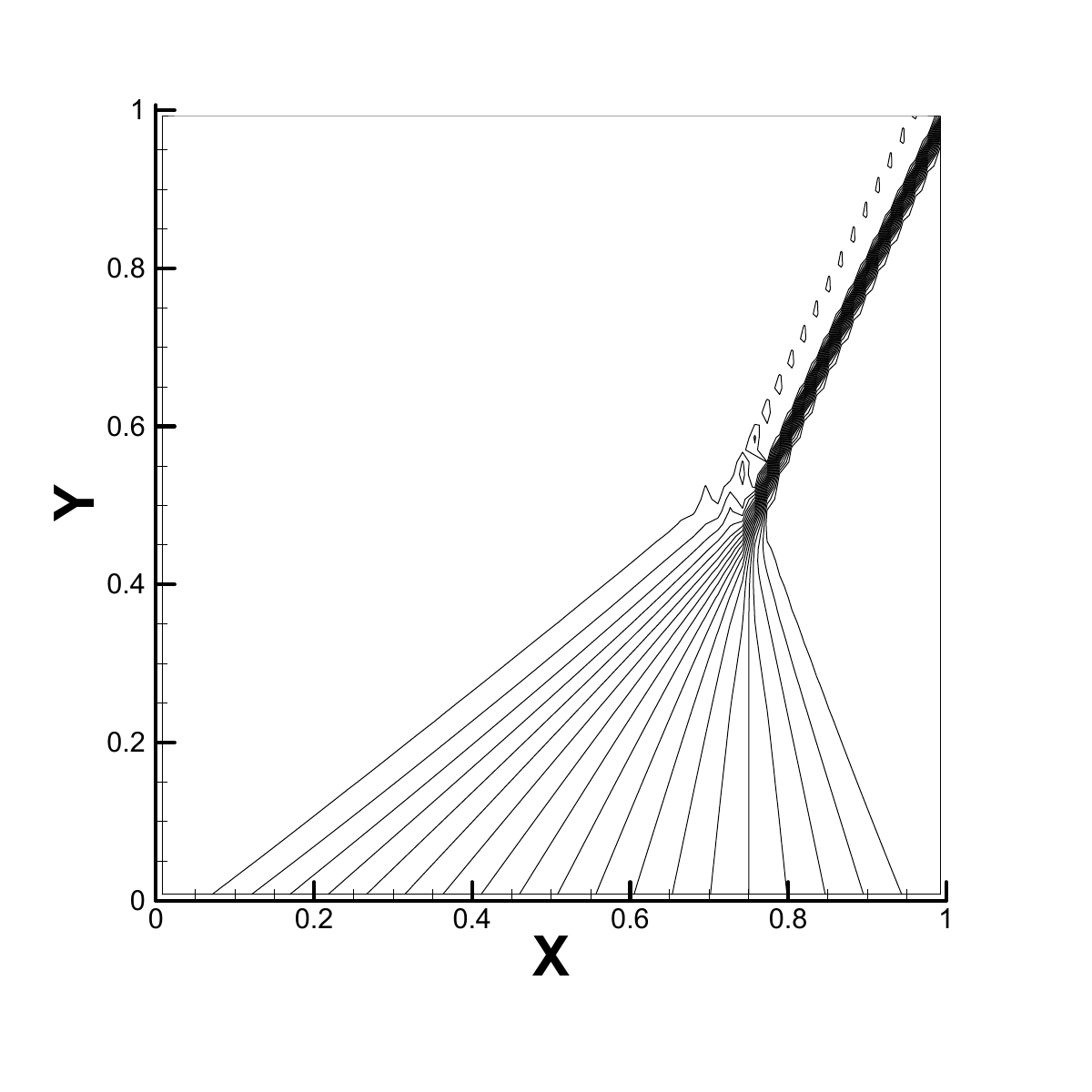} &
\includegraphics[height=4.2cm]{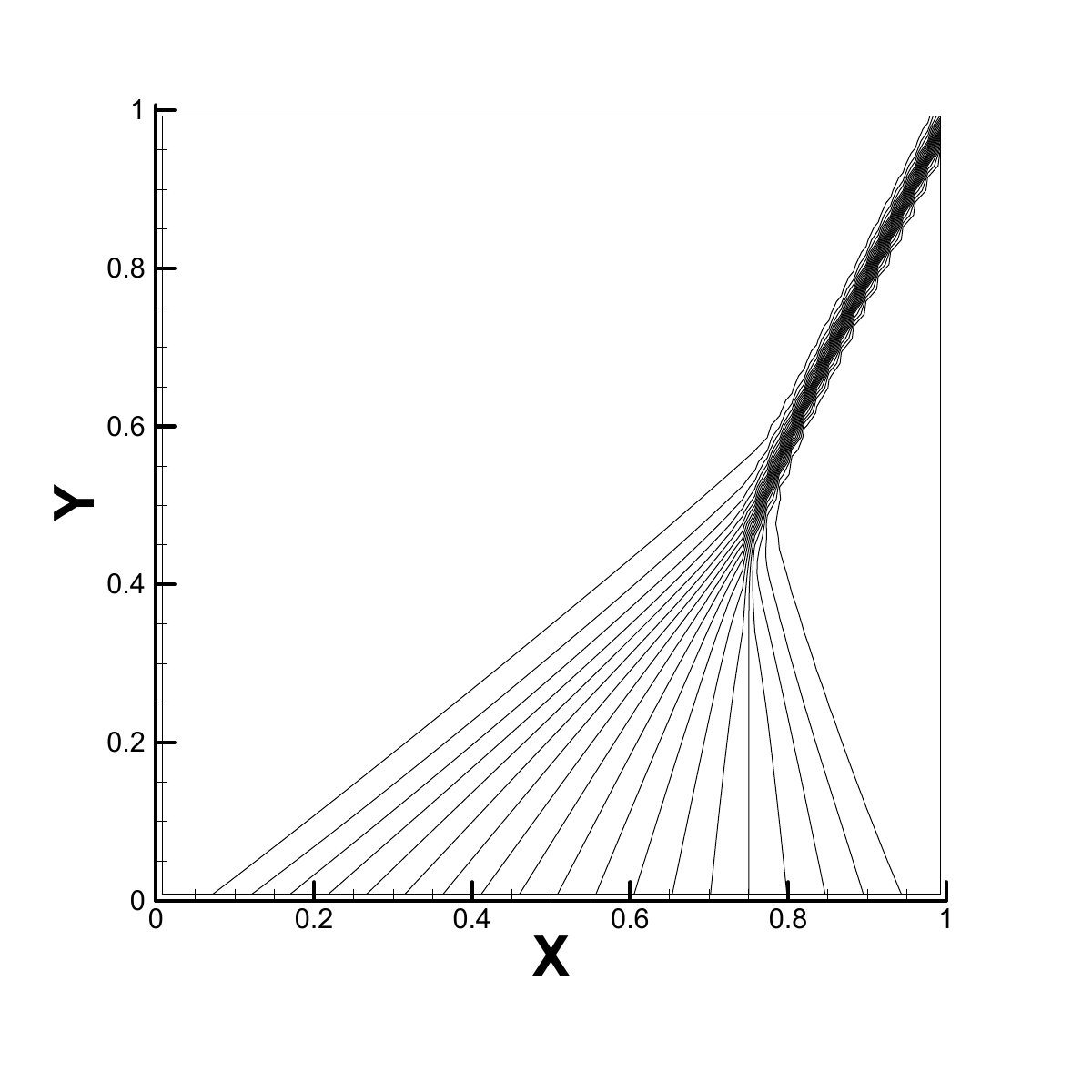} \\
 (d)First order KFDS+ &(e) KLW+  & (f)TVD-KFDS+ \\
\end{tabular}
\caption{Test Case 12 : KFDS, KLW \& TVD-KFDS  schemes in 2D Burgers framework with 64 x64  points}
\label{TC_12_2D_KFDS} 
\end{center} 
\end{figure}

The results for the test case 12 are given in figure(\ref{TC_12_2D_KFDS}). Here the  initial conditions and the boundary conditions for this test case lead to the formation of an oblique shock located between (0.75,0.5) to (1.0,1.0) and a  smooth region located beneath it.  The observations are similar to the test case 11 where TVD schemes demonstrate ability to capture the essential features of the test case with good accuracy and without oscillations.  

\subsubsection*{Test case 13:}
In this test case \cite{Zhang}, the 2-D viscous Burgers equation is solved on a square domain $x: [-0.5, 0.5]$, 
$y: [-0.5, 0.5]$ for the time period of $0.1s$ and with $\nu = 0.01$.  The initial conditions as well as the initial conditions are deduced from the exact solution, given below.   
\be
u ( x , y , t ) = 0.5 - \tanh \left( \frac { x + y - t } { 2 \nu } \right)
\ee

\begin{figure}[!h] 
\begin{center} 
\begin{tabular}{cccc}
\includegraphics[height=4cm]{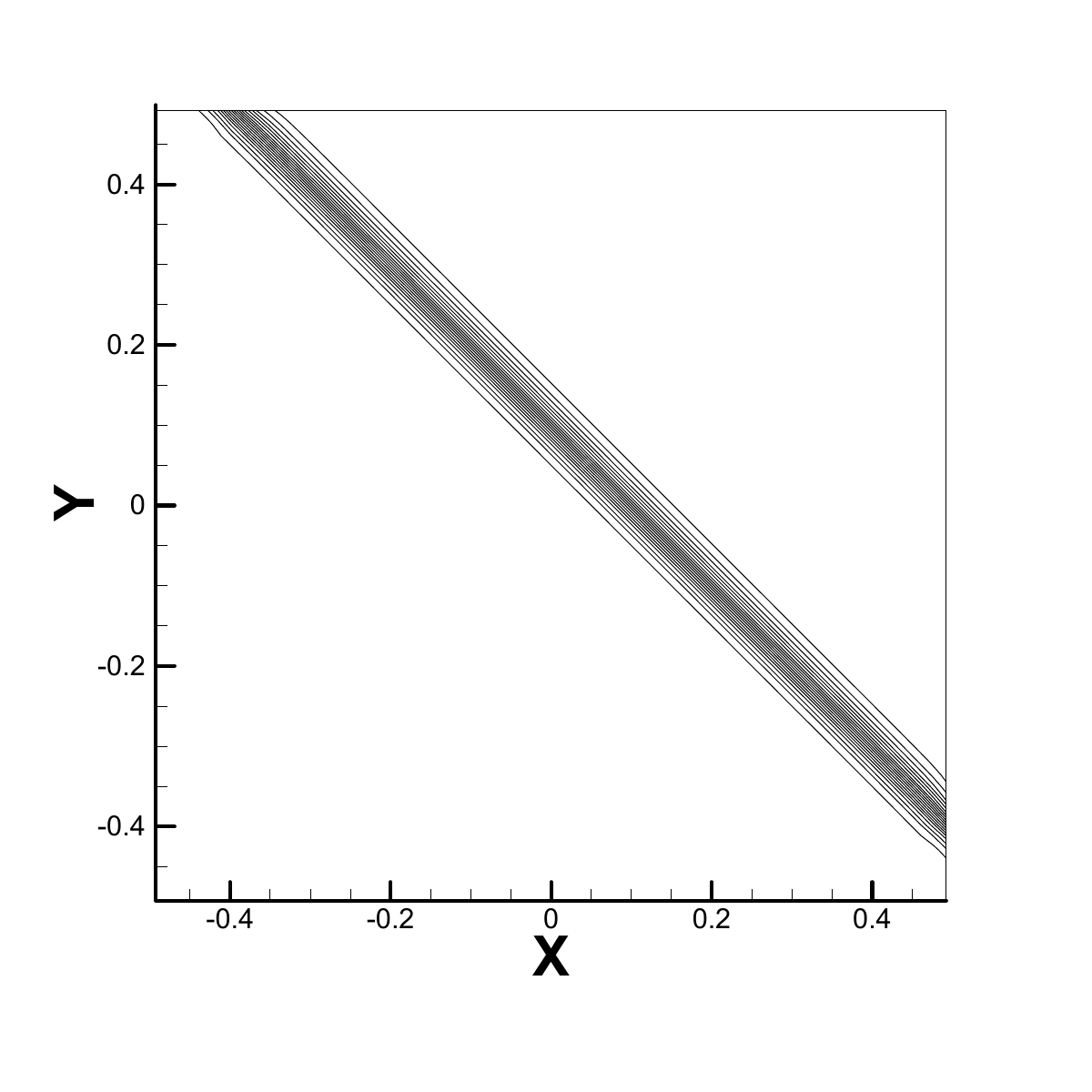} &
\includegraphics[height=4cm]{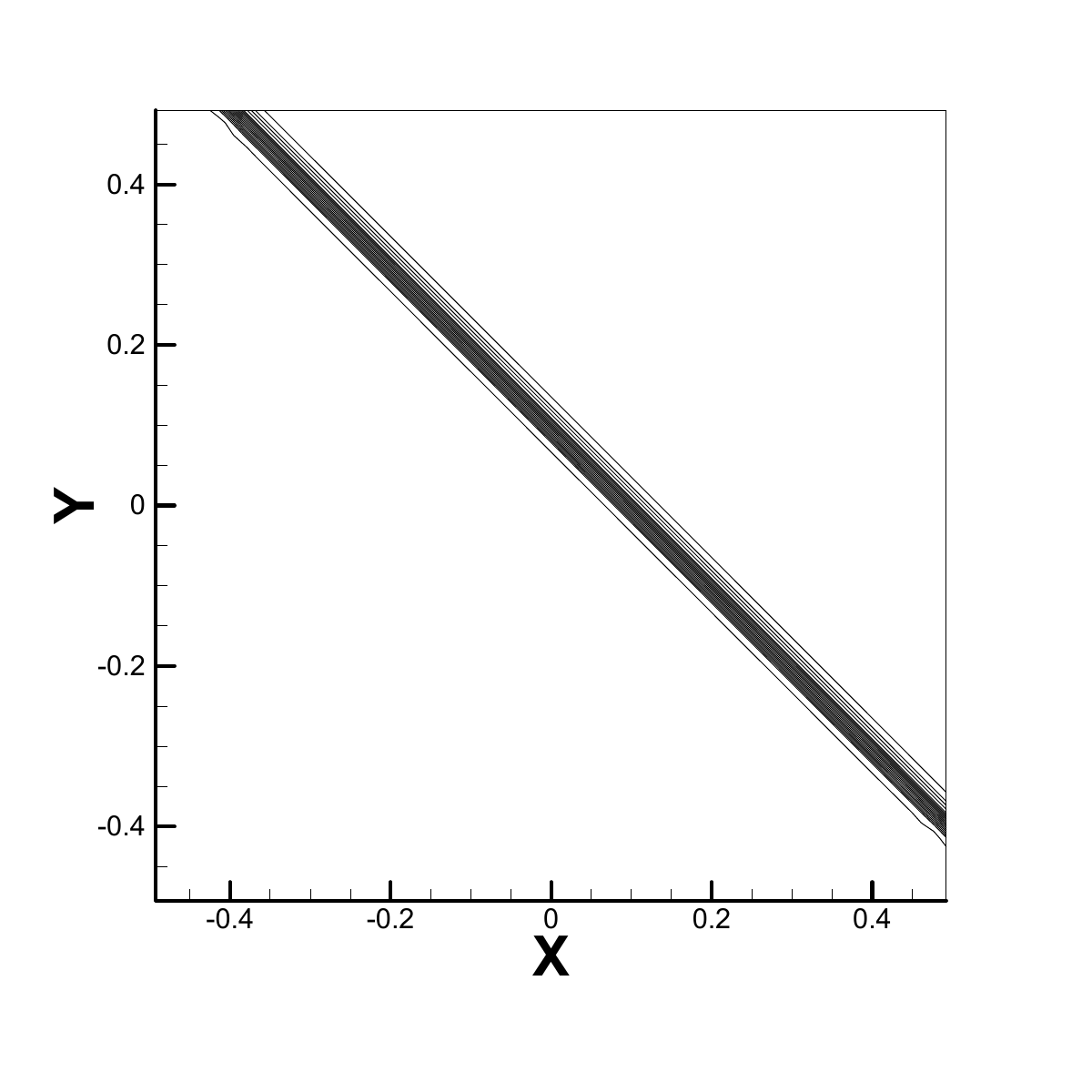} &
\includegraphics[height=4cm]{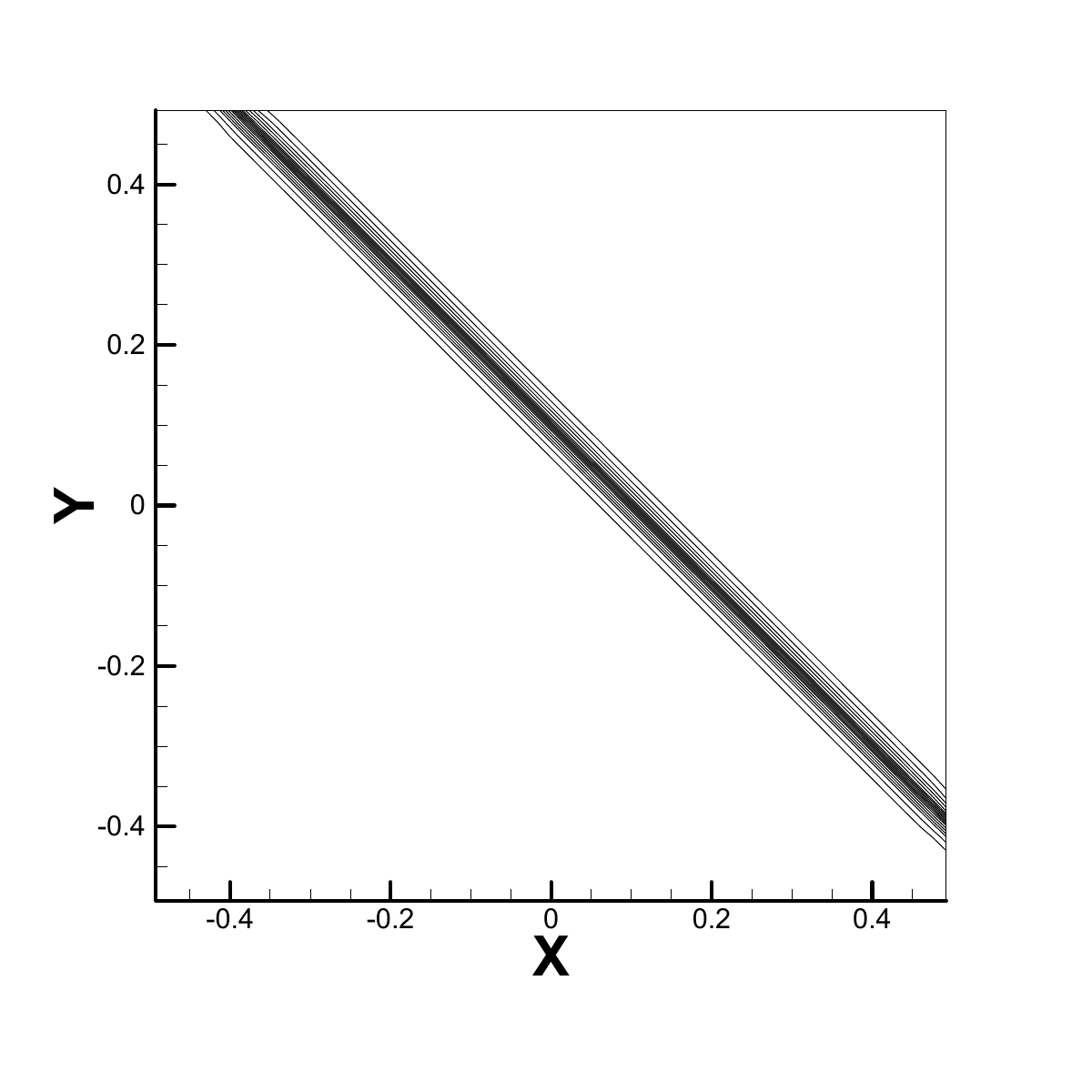} &
\includegraphics[height=4cm]{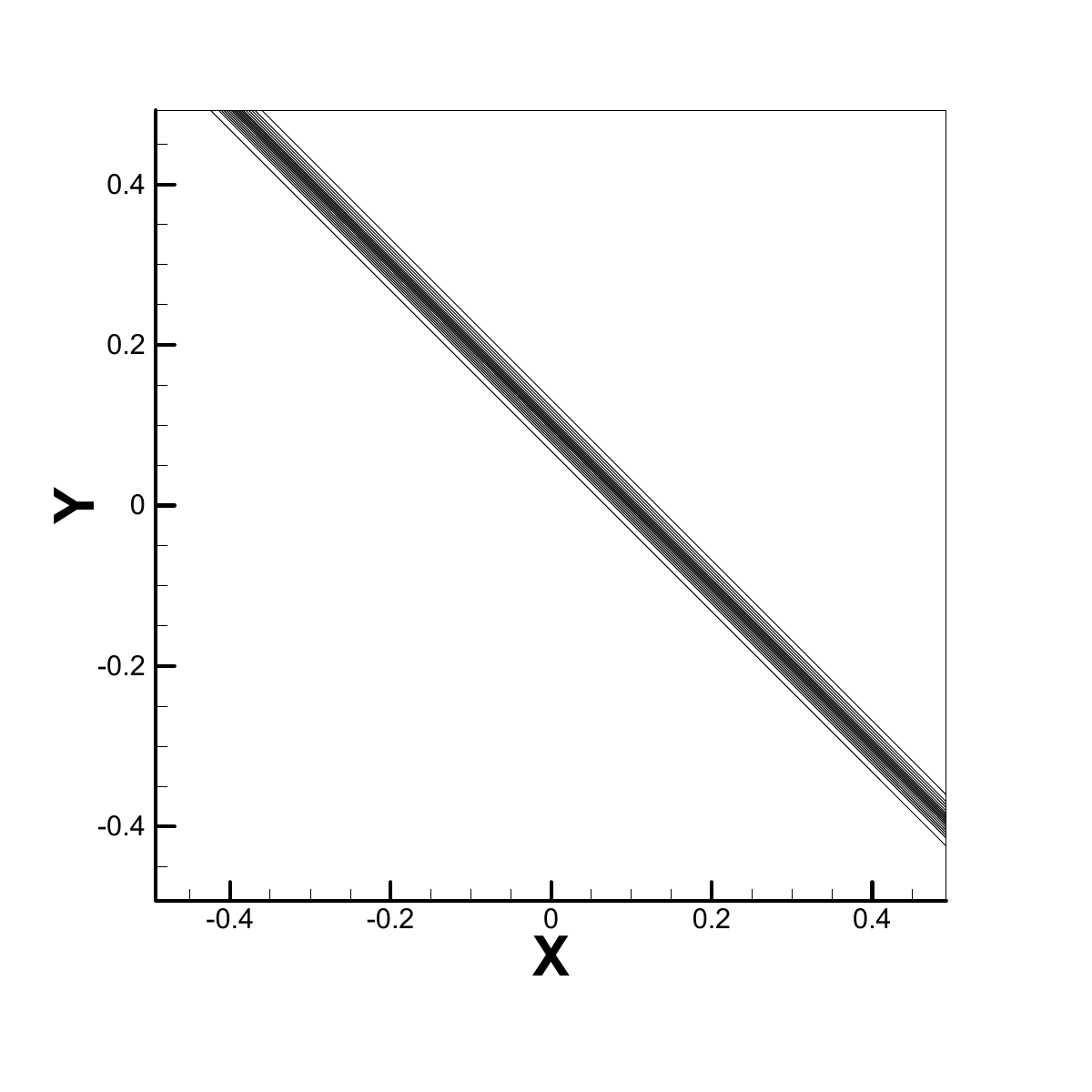} \\
(a) First order KFDS & (b)Second order KLW & (c)TVD-KFDS & (d)Exact Solution \\
\includegraphics[height=4cm]{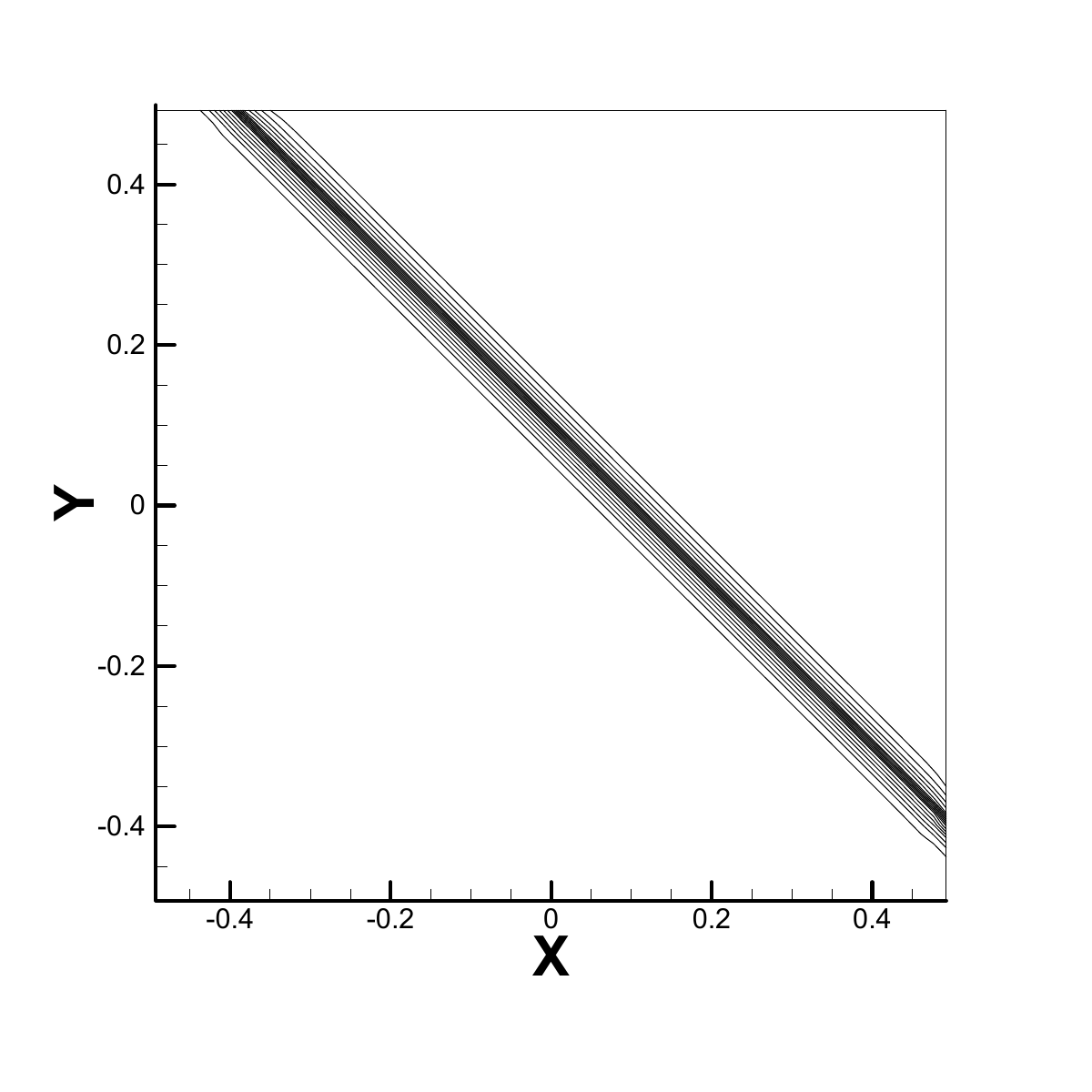} &
\includegraphics[height=4cm]{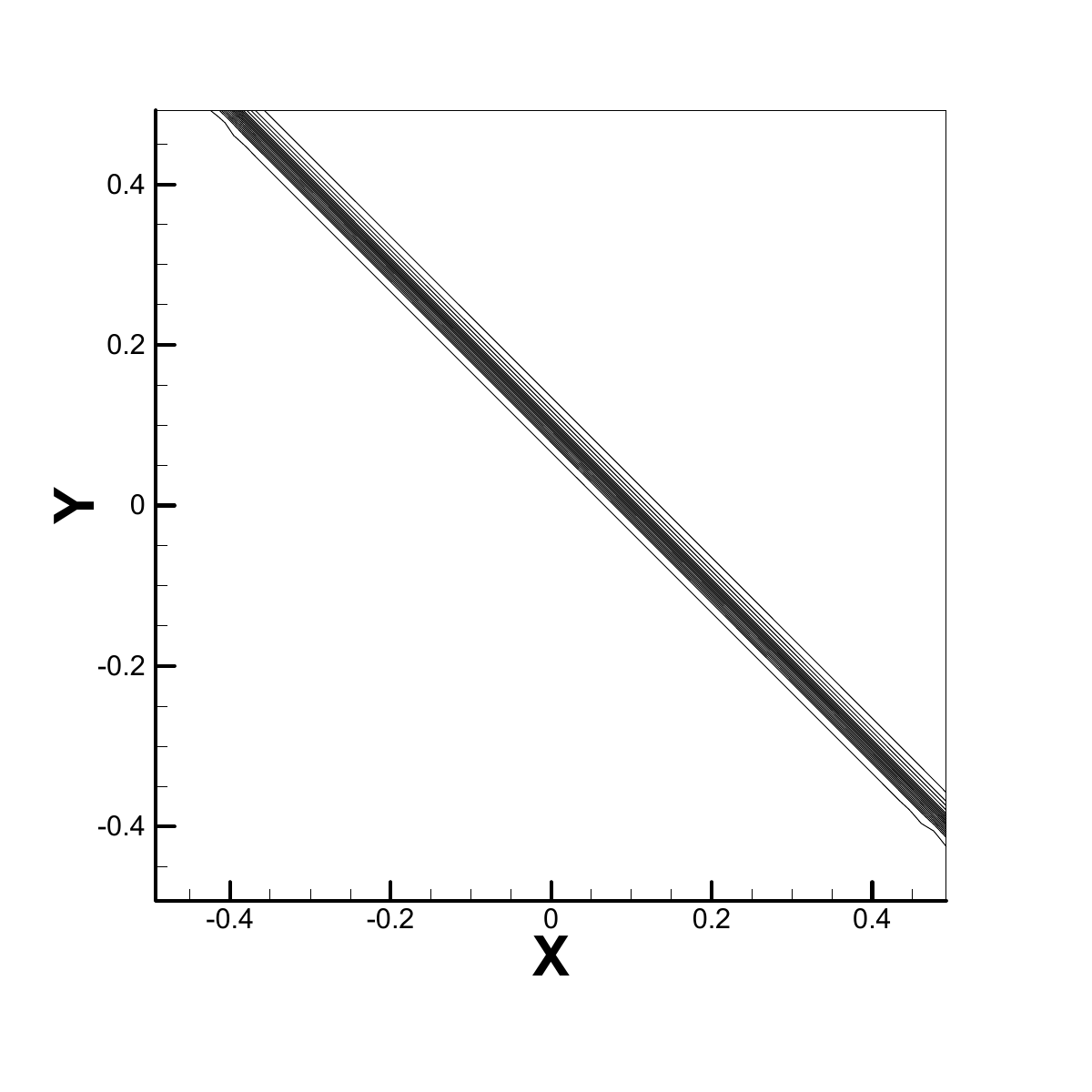} &
\includegraphics[height=4cm]{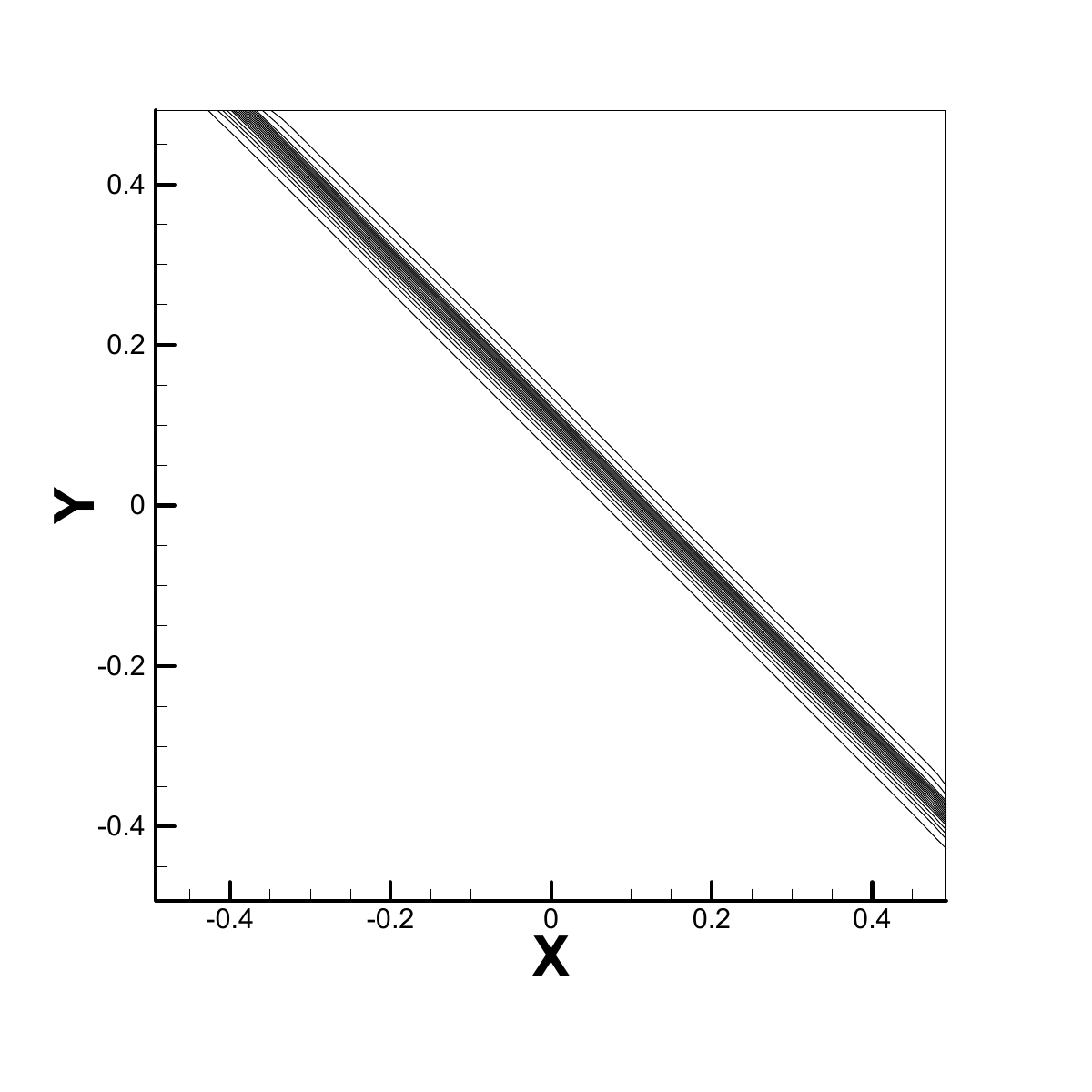} & \\
 (e)First order KFDS+ &(f) KLW+  & (g)TVD-KFDS+ & \\
\end{tabular}
\caption{Test Case 13 : KFDS, KLW \& TVD-KFDS  schemes in 2D Burgers framework with 64 x64  points}
\label{TC_13_2D_KFDS} 
\end{center} 
\end{figure}

The results for the test case 13 are given in figure (\ref{TC_13_2D_KFDS}).  In this test case there is an initial steep gradient located diagonally  from $(0,1)$ to $(1,0)$ which is diffused due to the viscous term. The exact solution for this test case at time $T = 0.1s$ is given in figure (\ref{TC_13_2D_KFDS}). All the schemes have been able to demonstrate capturing of the discontinuity at the right location in comparison to the exact solution. As expected, the second order schemes show accurate capturing of the discontinuity in comparison to the first order schemes.  

\subsubsection*{Test case 14:}
In this test case \cite{Zhang}, 2-D viscous Burgers equation is solved with the initial condition given by 
\bea
u ( x , y , 0 ) = \sin ( 2 \pi x ) \cos ( 2 \pi y )
\eea
The test case is solved in a computational domain of $x:[0, 1]$, $y:[0, 1]$, with $\nu= 0.01$ and time $t = 0.1s$.   Periodic boundary conditions are applied to all the boundaries.  

\begin{figure}[!h] 
\begin{center} 
\begin{tabular}{ccc}
\includegraphics[height=4cm]{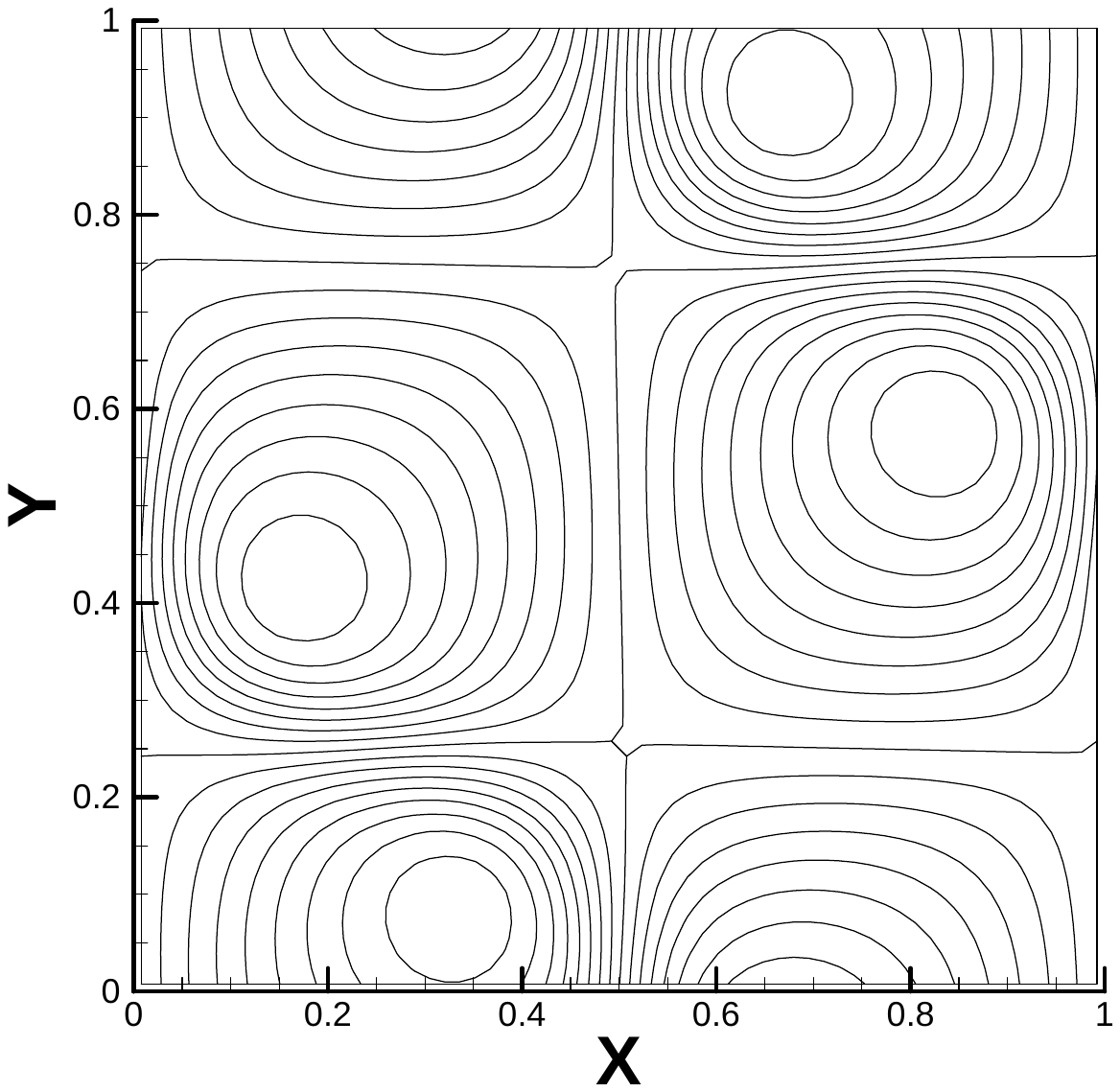} &
\includegraphics[height=4cm]{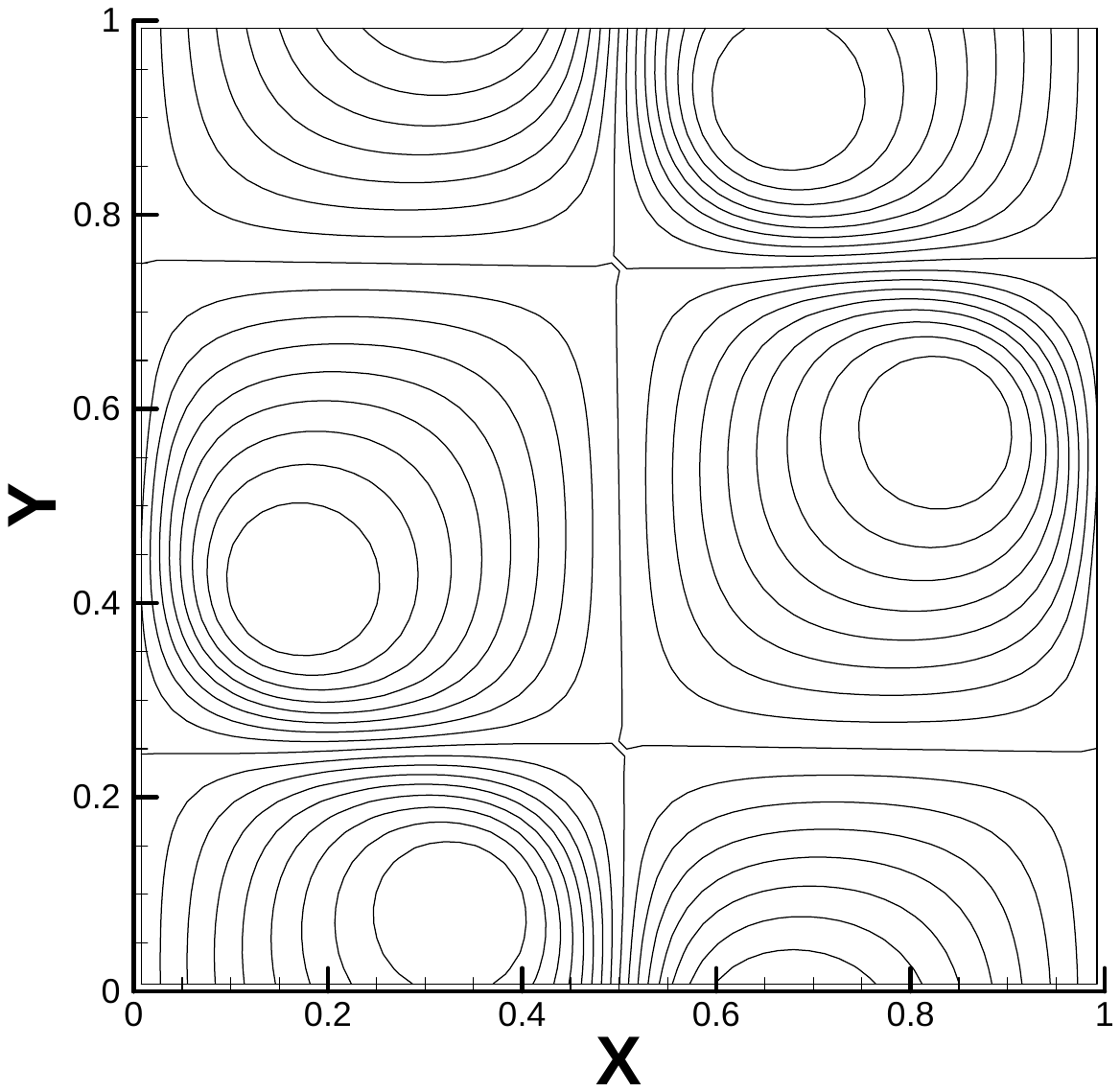} &
\includegraphics[height=4cm]{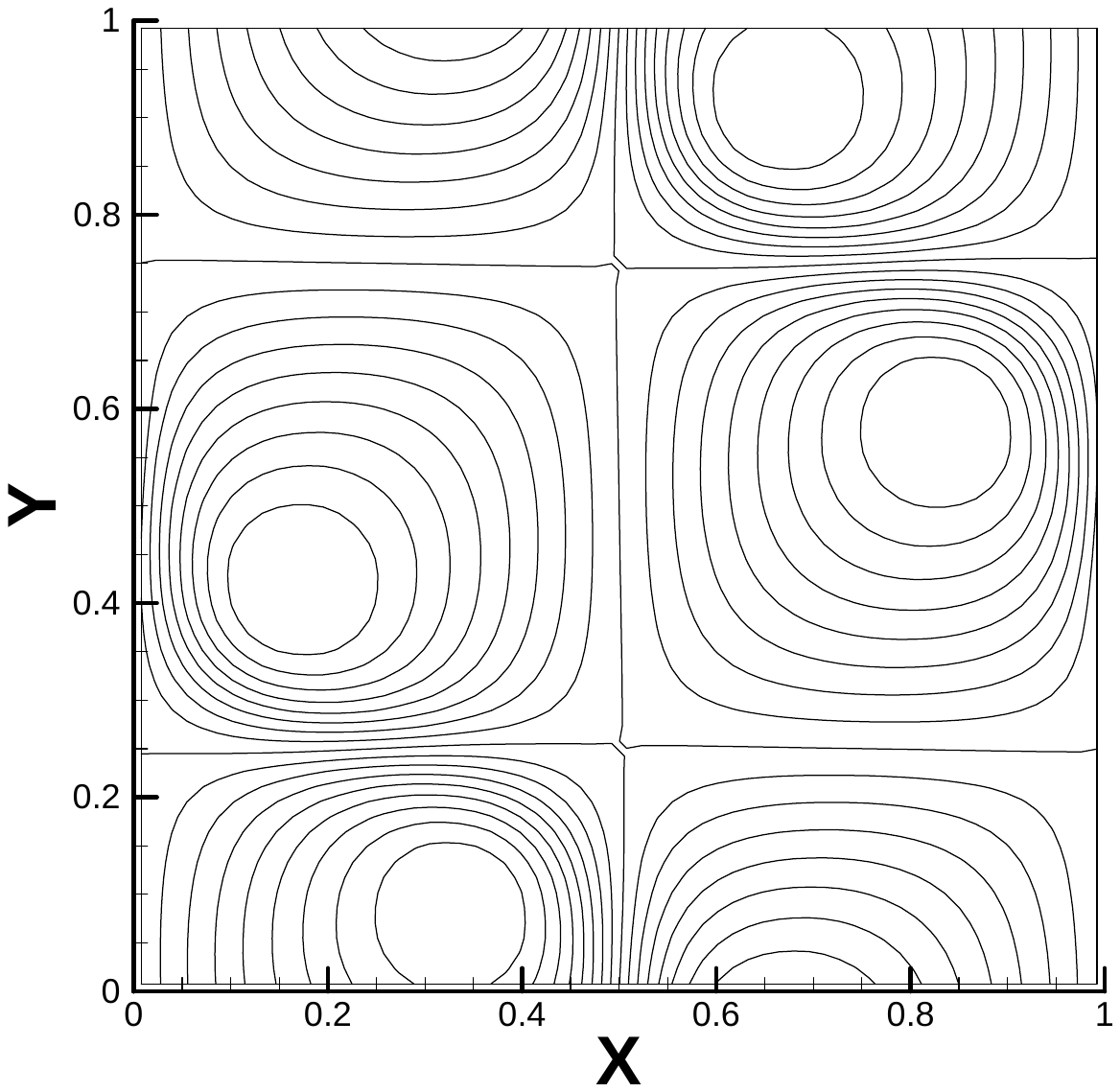} \\
 (a)First order KFDS & (b)Second order KLW & (c)TVD-KFDS\\
\includegraphics[height=4cm]{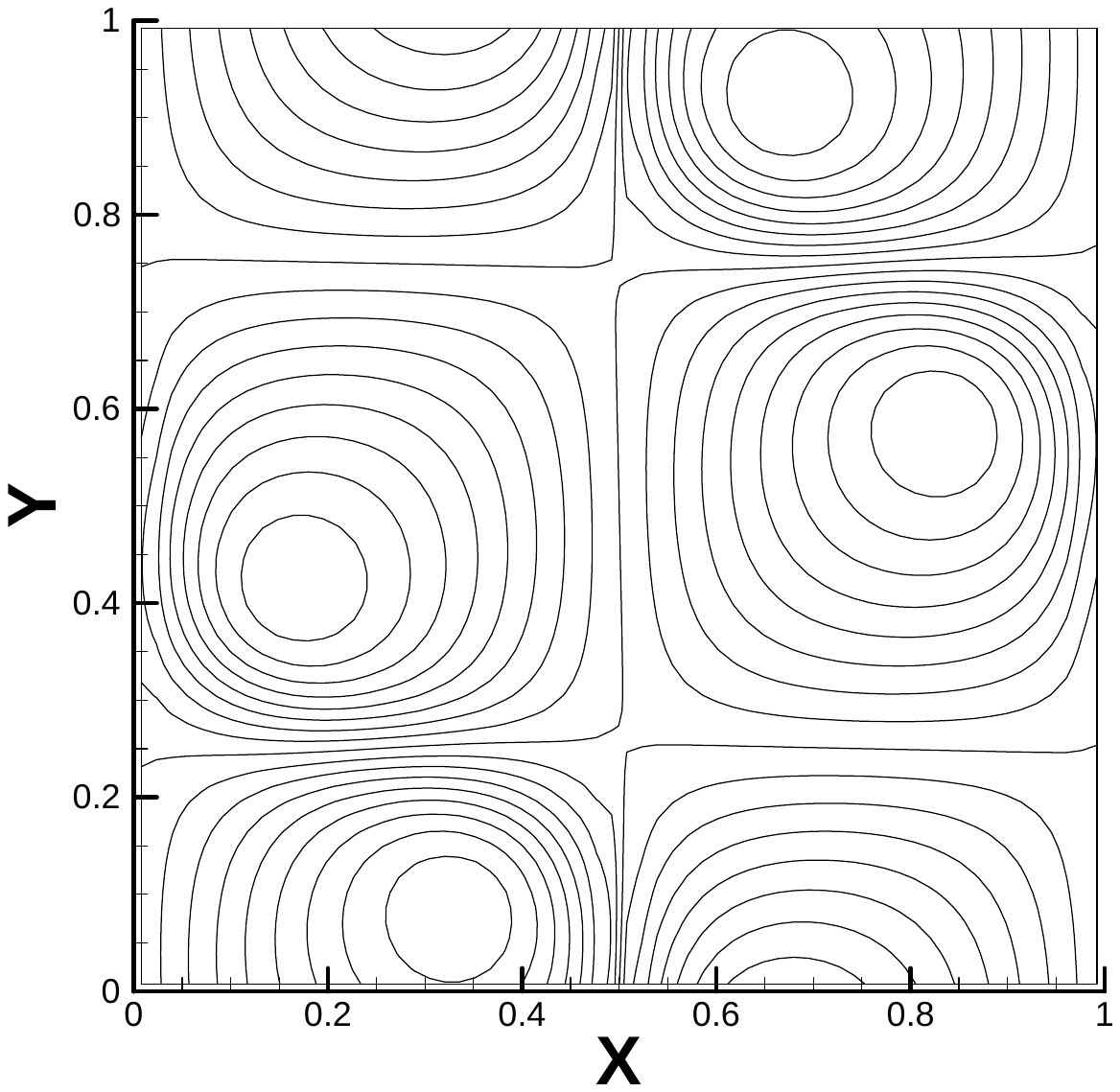} &
\includegraphics[height=4cm]{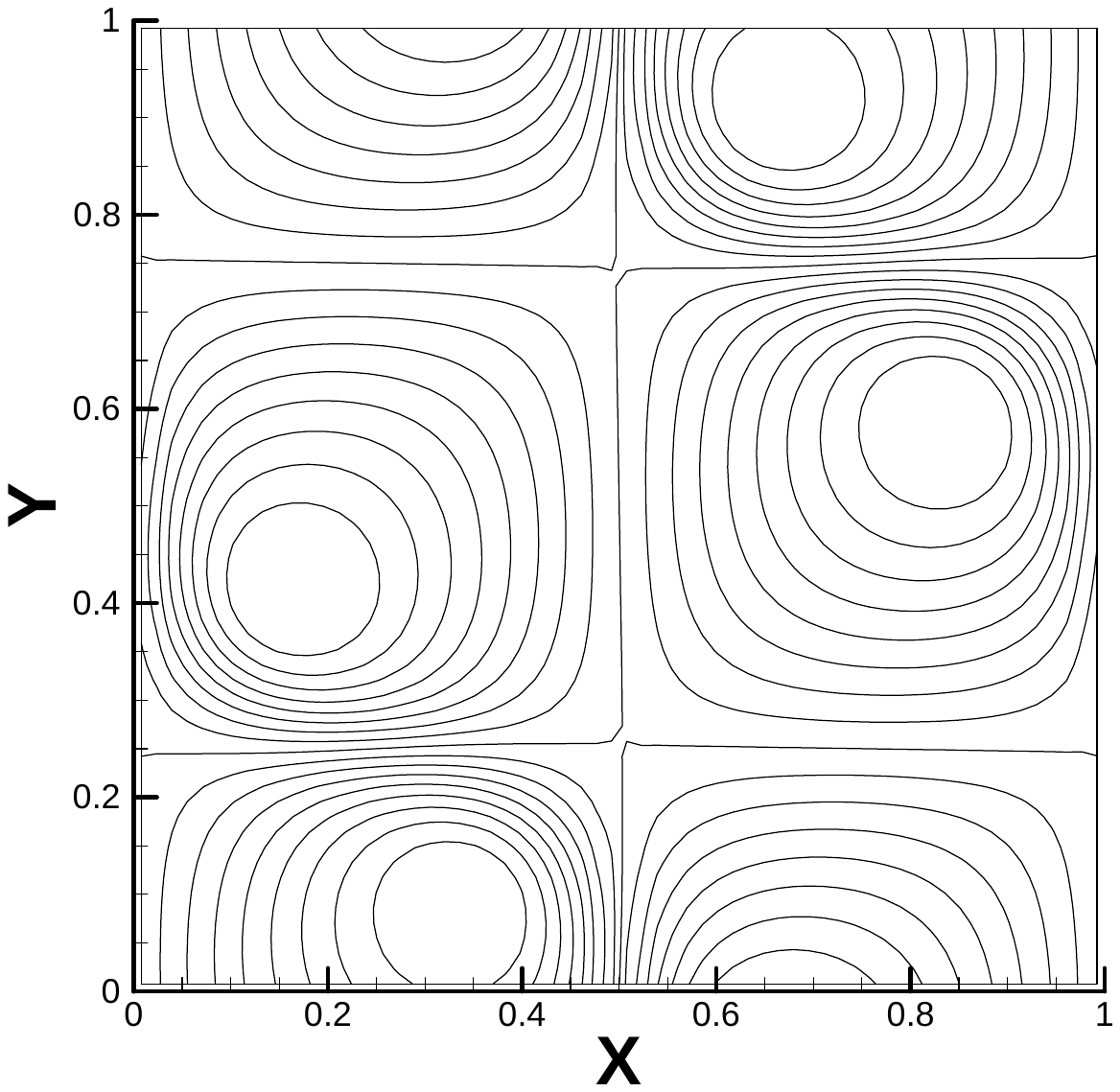} &
\includegraphics[height=4cm]{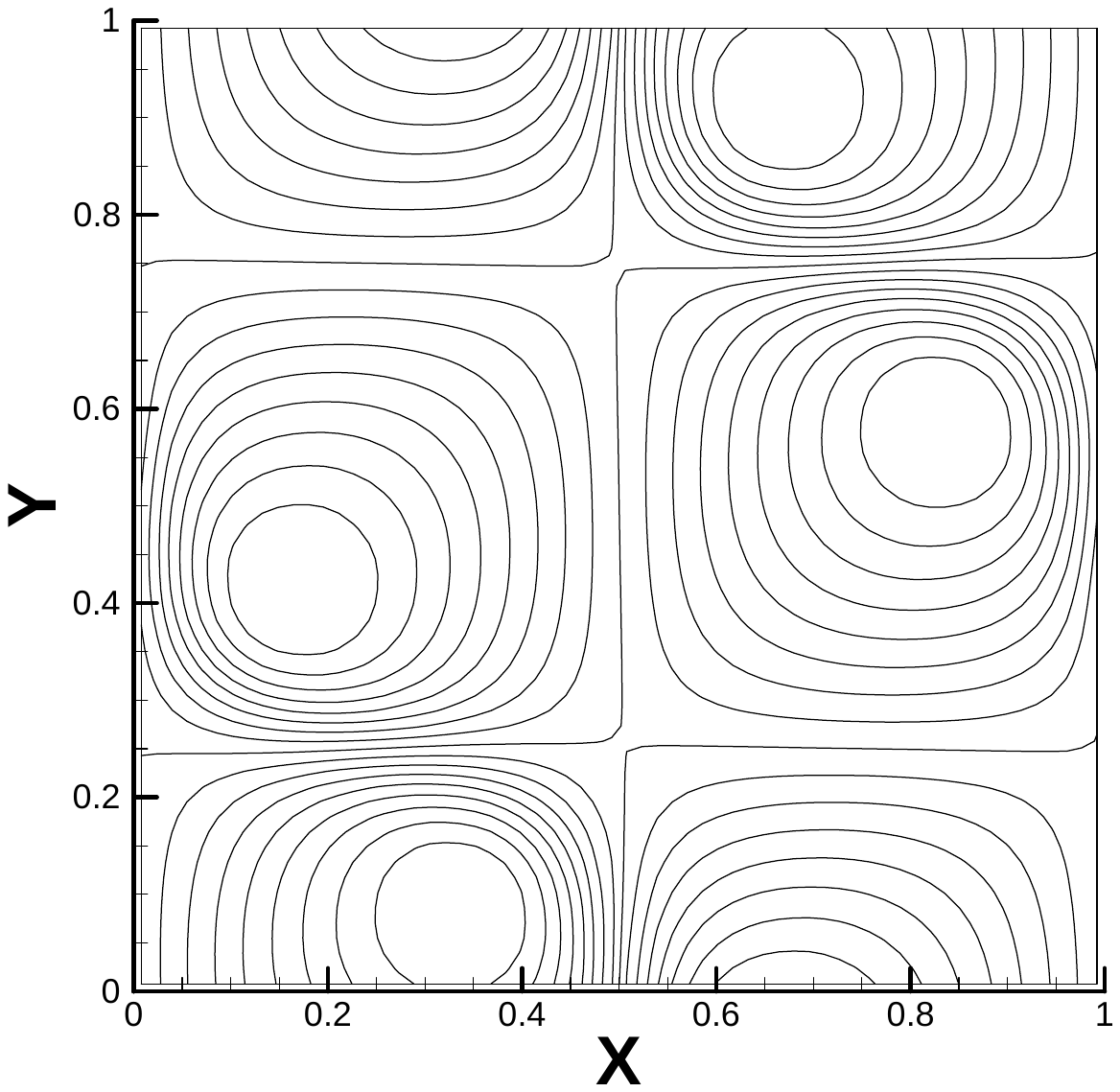} \\
 (d)First order KFDS+ &(e) KLW+  & (f)TVD-KFDS+ \\
\end{tabular}
\caption{Test Case 14 : KFDS, KLW \& TVD-KFDS  schemes in 2D Burgers framework with 64 x64  points}
\label{TC_14_2D_KFDS} 
\end{center} 
\end{figure}

The results for the test case 14 are given in figures (\ref{TC_14_2D_KFDS}). This test case is a scalar mimic of a pattern of vortices which propagate through the domain.  In the absence of an exact solution, the results are compared with  those provided in \cite{Zhang}.  The transportation of the vortices and their positions are reasonably captured in all the schemes in comparison to the solutions provided in \cite{Zhang}.  As there are no discontinuities in this test case, there is not much difference between different versions of the schemes.  

\subsubsection*{Test case 15:}
This test case involves an instantaneous breakage of an  idealised circular dam, modelled by 2-D shallow water equations. The case involves a circular dam of $11m$ radius centred inside rectangular computational domain of $x=50m$ and $y=50m$. The domain is discretised with $40 \times 40$ rectangular grid. Computations are obtained at a final time of $t = 0.69s$. 
		
\begin{figure}[!h] 
\begin{center} 
\begin{tabular}{cc}
\includegraphics[height=5.5cm]{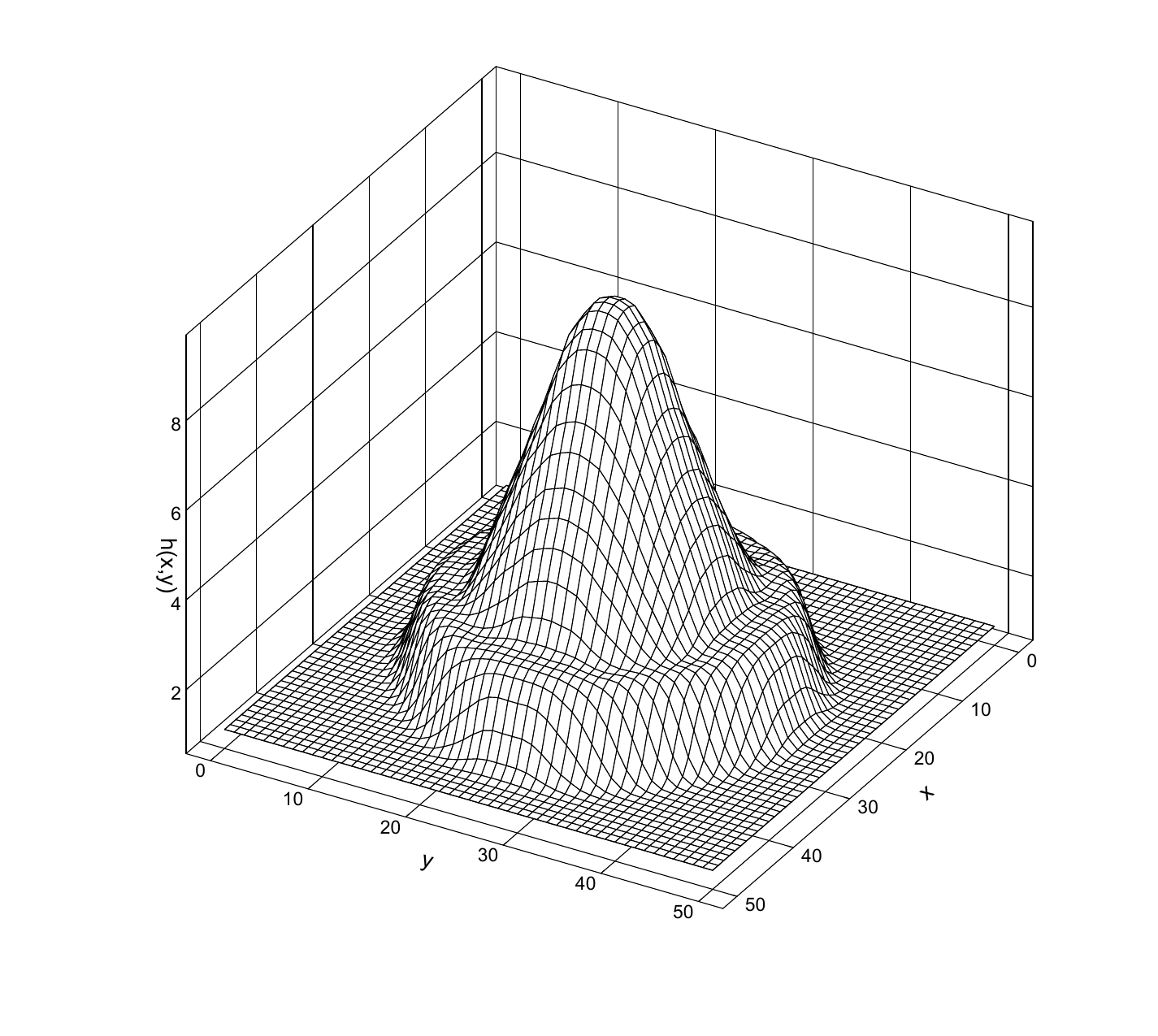} &
\includegraphics[height=5.5cm]{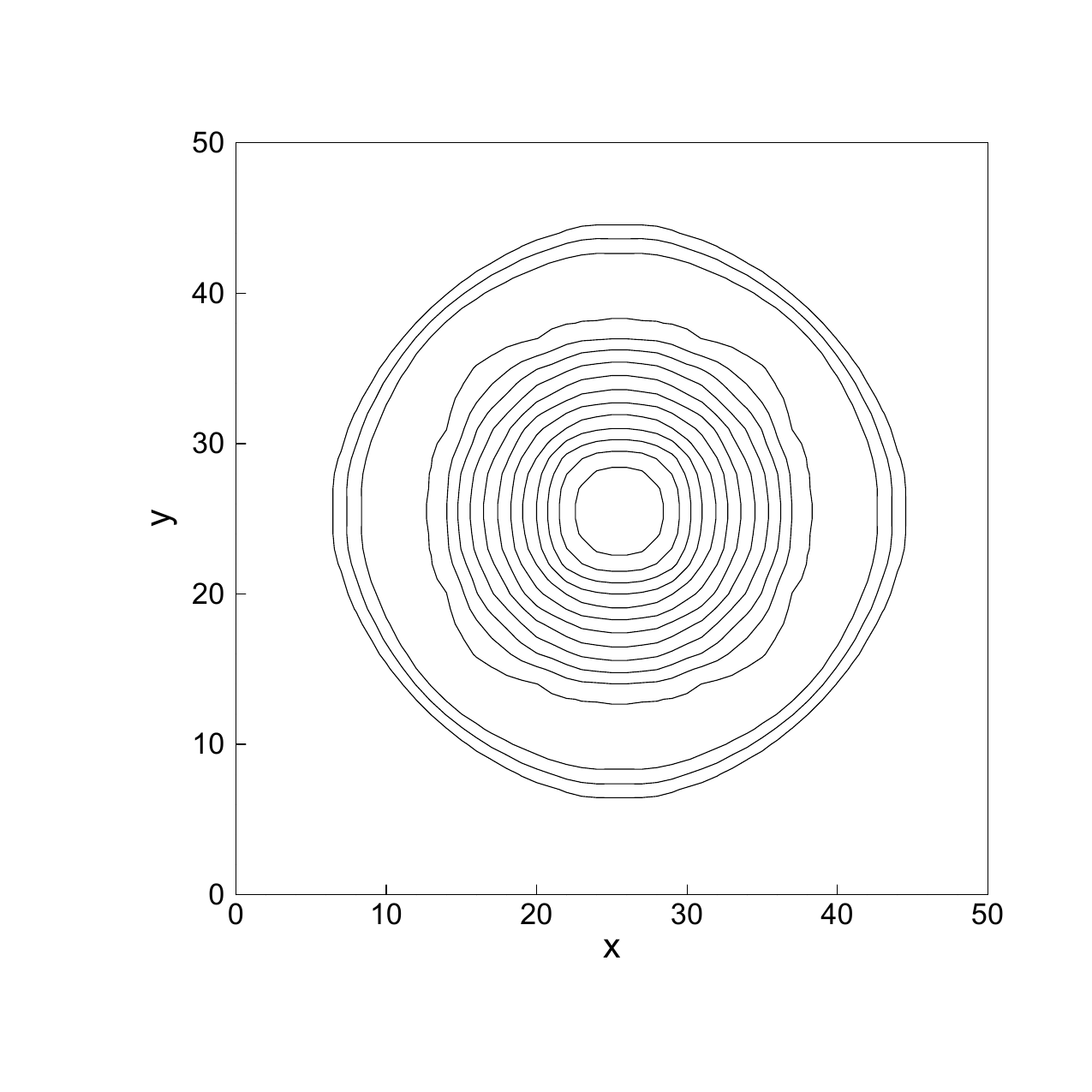}  \\
  & (a)First order KFDS \\
\includegraphics[height=5.5cm]{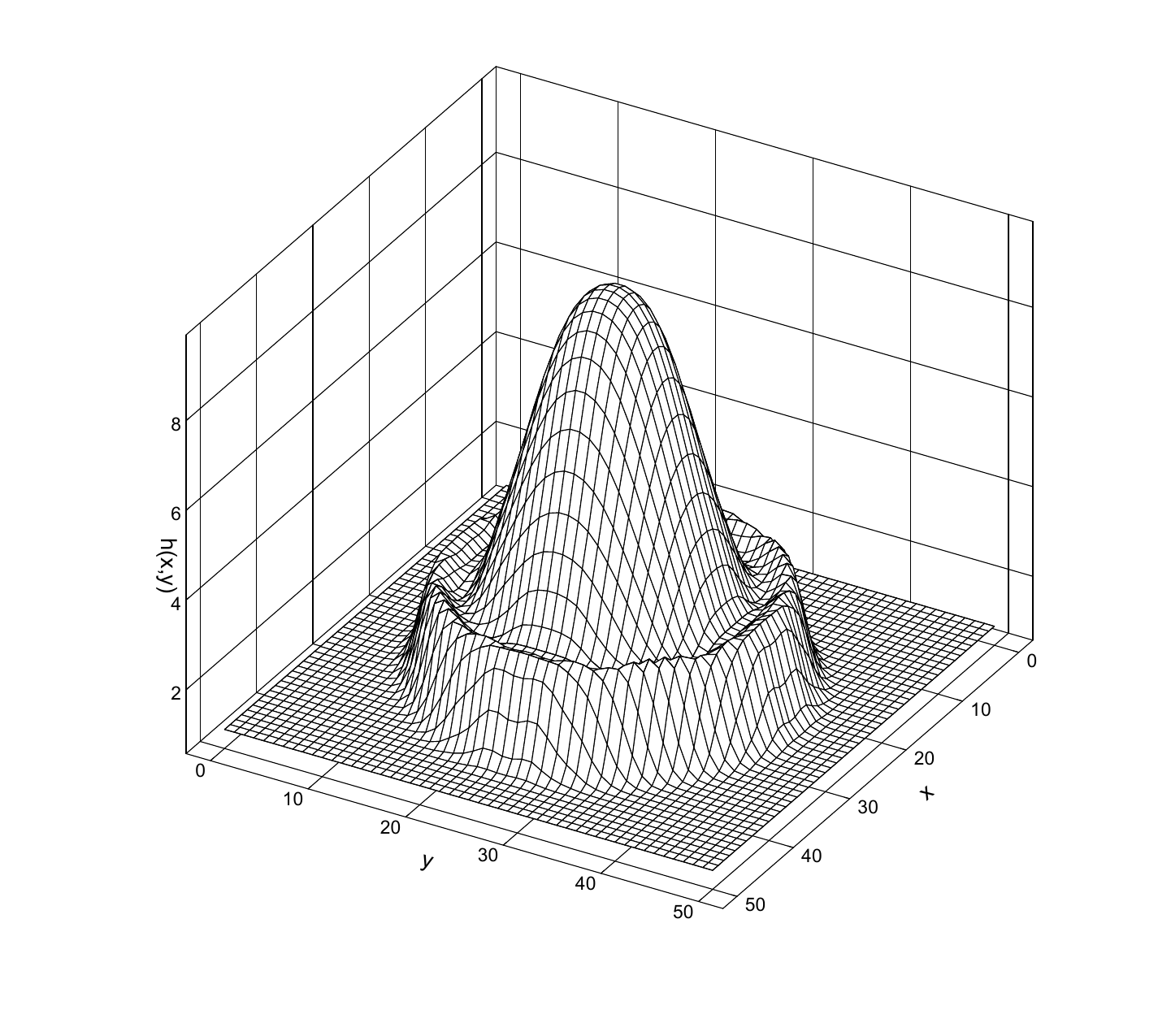} &
\includegraphics[height=5.5cm]{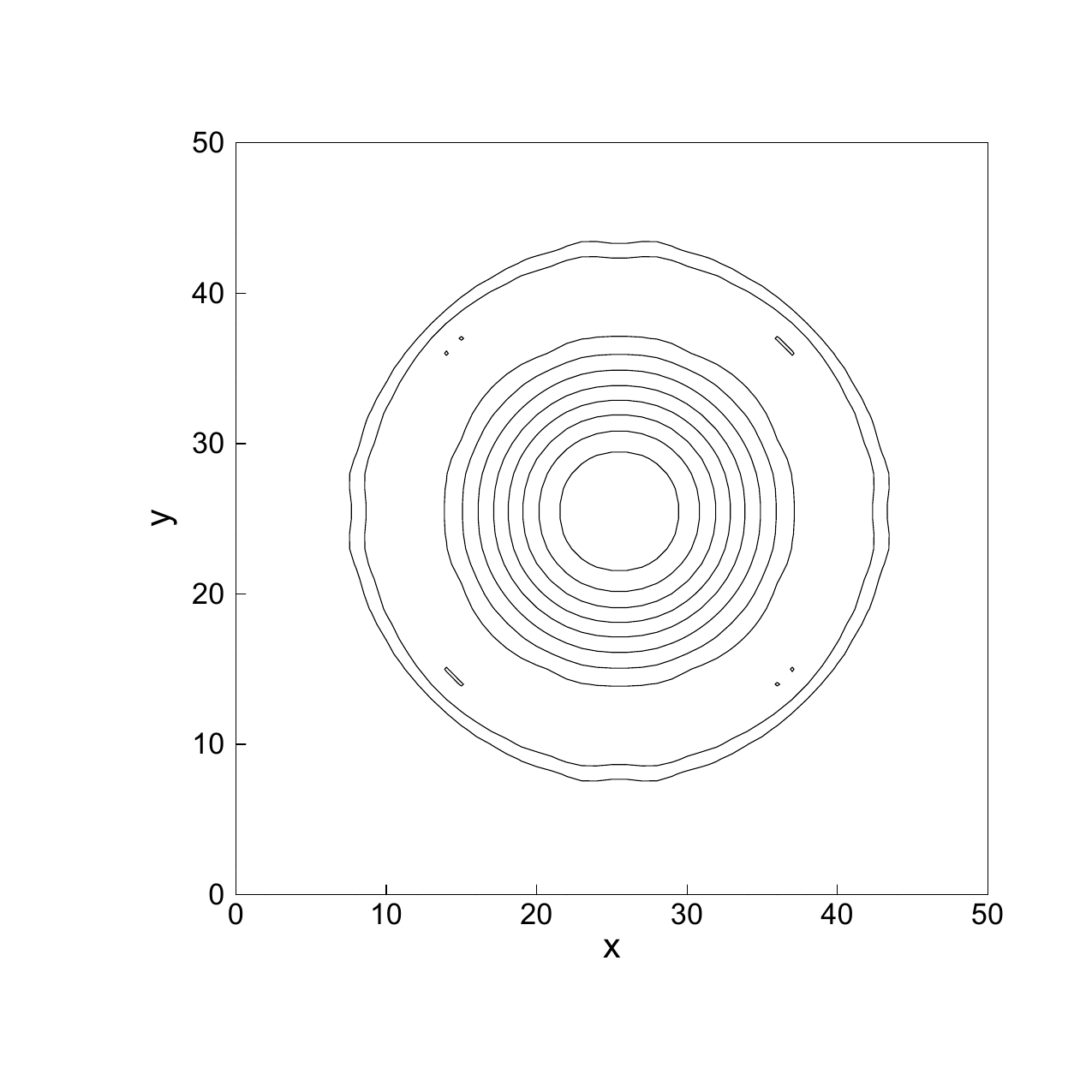}\\
 & (b)Second order KLW \\
\includegraphics[height=5.5cm]{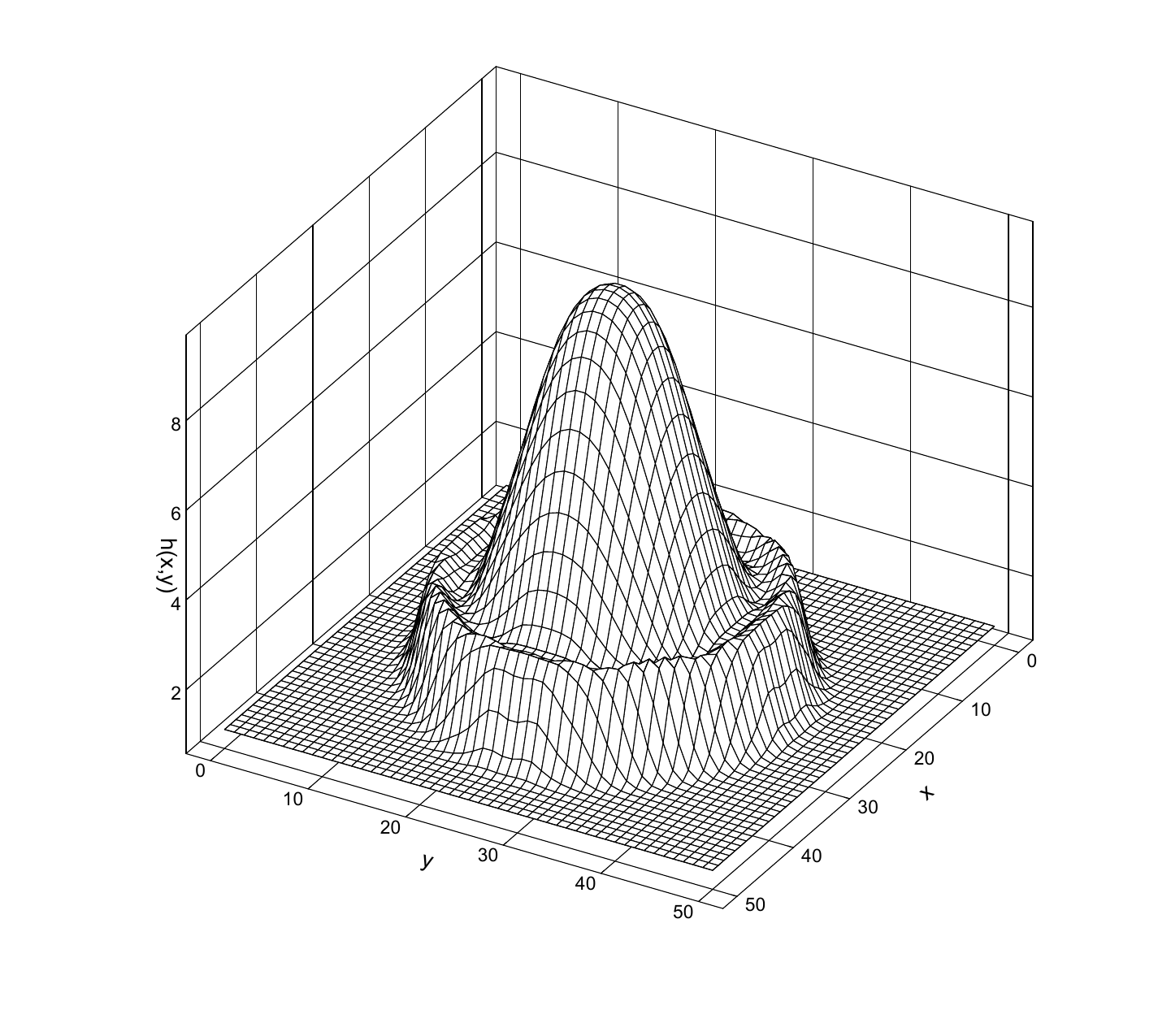} &
\includegraphics[height=5.5cm]{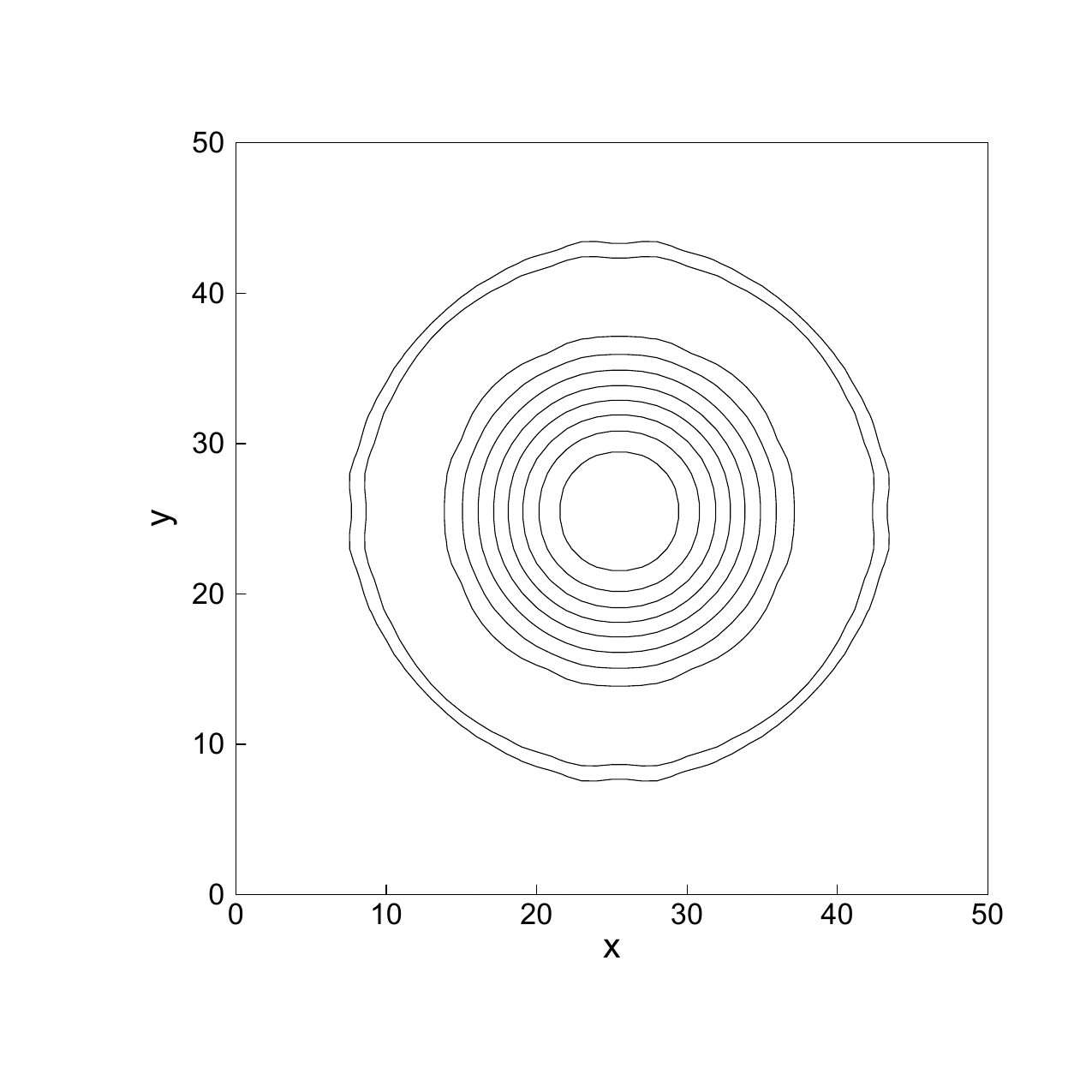} \\
  & (c)TVD-KFDS\\
\end{tabular}
\caption{Test Case 15 : KFDS, KLW \& TVD-KFDS  schemes in Shallow water framework with 40 x40  points}
\label{SWE_2D_CDAM_KFDS} 
\end{center} 
\end{figure} 

The results for the test case 15 are given in figures (\ref{SWE_2D_CDAM_KFDS}). The hydraulic jump arising out of the instantaneous dam break radiates outwards in all directions while the contraction wave travels radially inwards towards the centre. Three dimensional and two dimensional plots of the hydraulic jump is provided in the figures. Each of the numerical scheme captures these variations reasonably. The numerical results have good correlation with those available in the literature.  

\section{Conclusions} 
A firm foundation for the flexible velocity Boltzmann equation is established by deriving it from the classical Boltzmann equation, using Dirac delta functions for the equilibrium distribution.  Based on the so derived flexible velocity Boltzmann system, upwind schemes are introduced for solving the nonlinear convection equation.  The flexible velocity also happens to be the coefficient of numerical diffusion in the finite volume framework and hence the freedom in fixing it is utilized in deriving an accurate shock capturing scheme in steady state, which is otherwise low in numerical diffusion.  A generalized kinetic Lax-Wendroff method is introduced by applying the classical strategy to the flexible velocity Boltzmann framework.  To get rid of the oscillations near shocks for the higher order accurate versions, total variation diminishing technique is introduced within the discrete kinetic framework, thus obtaining new and efficient kinetic schemes for solving the nonlinear convection-diffusion equations.  Discrete Chapman-Enskog distributions are used in the relaxation step for obtaining schemes for the viscous cases.  These new schemes are free of any Riemann solvers and enjoy a potential for further applications in convection dominated flows.   

\section*{Acknowledgments} 
The first author thanks the University of South Africa (UNISA) for providing funding and Indian Institute of Science for granting sabbatical leave, during which this work was conceived.  The third author thanks Hindustan Aeronautics Limited (HAL), Bangalore for funding the author, which facilitated this research activity.

\clearpage 
\begin{appendices} 
{\bf \large Appendix A: Useful Integrals and Functions} 
\be 
J_{n} = \int_{-\infty}^{\infty} w^{n} e^{-w^{2}} dw, \ n=0,1,\cdots 
\ee
\be 
J^{+}_{n} = \int_{0}^{\infty} w^{n} e^{-w^{2}} dw, \ n=0,1,\cdots   
\ee
\be 
J^{-}_{n} = \int_{-\infty}^{0} w^{n} e^{-w^{2}} dw, \ n=0,1,\cdots    
\ee
Note that $J^{+}_{n} + J^{-}_{n} = J_{n}$.   
\begin{center} 
{\bf Table of Integrals} 
\end{center} 
\bea 
\ba{cc} 
J_{0} = \sqrt{\pi} & J^{+}_{0} = \fr{1}{2} \sqrt{\pi} \\[2mm]  
J_{1} = 0 & J^{+}_{1} = \fr{1}{2} \\[2mm]  
J_{2} = \fr{1}{2} \sqrt{\pi} & J^{+}_{2} = \fr{1}{4} \sqrt{\pi} \\[2mm]  
J_{1} = 0 & J^{+}_{1} = \fr{1}{2} \\[2mm]  
J_{2} = \fr{3}{4} \sqrt{\pi} & J^{+}_{2} = \fr{3}{8} \sqrt{\pi} \
\ea 
\eea  

\vspace*{0.5cm} 

\noindent 
{\bf Definition of the Dirac Delta Function:}  
\bea  
\delta \left( x - \alpha \right) = \left\{ \ba{l} 0 \ \textrm{if} \ x \ne \alpha \\ 
\infty \ \textrm{if} \ x = \alpha \ea \right. 
\eea 
{\bf Sifting or Selector property:}   
\be
\int_{-\infty}^{\infty} f(x) \delta \left( x - \alpha \right) dx = f(\alpha)  
\ee
{\bf Other properties:} 
\be
\int_{-\infty}^{\infty} \delta \left( x - \alpha \right) dx = 1 \ \ \ \ \ (\textrm{Normalization}) 
\ee
\be 
\delta \left( a x - \alpha \right) = \fr{1}{|a|} \delta \left( x - \fr{\alpha}{a} \right)  
\ee

\newpage 
\noindent 
{\bf \large Appendix B: Split fluxes in 1-D} 
\be 
g^{+}(u) = {\bf P} \Lambda^{+} {\bf f}^{eq}
\ee
Therefore 
$$ \ba{rcl} g^{+}(u) & = & \left[ 1 \ 1 \right] \left[ \ba{cc} \lambda_{+} & 0 \\ 0 & 0 \ea \right] \left[ \ba{c} {\bf f_{+}}^{eq} \\ {\bf f_{-}}^{eq} \ea \right] \\[2mm]   
& = & \lambda_{+} {\bf f_{+}}^{eq} \\[2mm] 
& = & \lambda \left[ \fr{1}{2} u + \fr{1}{2 \lambda} g(u) \right] \\[2mm] 
\ea 
$$ 
Therefore 
\be 
g^{+}(u) = \fr{1}{2} g(u) + \fr{\lambda}{2} u 
\ee 
Similarly, starting with $g^{-}(u) = {\bf P} \Lambda^{-} {\bf f}^{eq}$, we get

$$ \ba{rcl} g^{-}(u) & = & \left[ 1 \ 1 \right] \left[ \ba{cc} 0 & 0 \\ 
0 & \lambda_{-} \ea \right] \left[ \ba{c} {\bf f_{+}}^{eq} \\ {\bf f_{-}}^{eq} \ea \right] \\[2mm]   
& = & - \fr{\lambda}{2} u + \fr{1}{2} g(u)    
\ea 
$$ 

Thus 
\be 
g^{\pm}(u) = \fr{1}{2} g(u) \pm \fr{\lambda}{2} u 
\ee 
\newpage 

\noindent 
{\bf \large Appendix C} 
\section*{Chapman-Enskog type expansion for the flexible velocity Boltzmann equation} 
The flexible velocity Boltzmann equation (FVBE) is 
\be 
\fr{\del {\bf f}}{\del t} + \Lambda \fr{\del {\bf f}}{\del x} 
= - \fr{1}{\epsilon} \left[ {\bf f} - {\bf f}^{eq} \right]
\ee
Let us start with an asymptotic expansion of the distribution function as 
\be 
{\bf f} = {\bf f}^{(0)} + \epsilon {\bf f}^{(1)} 
+ \mathcal{O}\left( \epsilon^{2} \right)
\ee
and substitute it in FVBE to obtain  
$$ \fr{\del }{\del t} 
\left[ {\bf f}^{(0)} + \epsilon {\bf f}^{(1)} 
+ \mathcal{O}\left( \epsilon^{2} \right) \right] 
+ \Lambda \fr{\del }{\del x} 
\left[ {\bf f}^{(0)} + \epsilon {\bf f}^{(1)} 
+ \mathcal{O}\left( \epsilon^{2} \right) \right] 
= - \fr{1}{\epsilon} \left[ 
{\bf f}^{(0)} + \epsilon {\bf f}^{(1)} 
+ \mathcal{O}\left( \epsilon^{2} \right) - {\bf f}^{eq} 
\right] $$ 
or 
$$ \epsilon \left[ \fr{\del {\bf f}^{(0)}}{\del t} 
+ \Lambda \fr{\del {\bf f}^{(0)}}{\del x} \right] 
+ \epsilon^{2} \left[ \fr{\del {\bf f}^{(1)}}{\del t} 
+ \Lambda \fr{\del {\bf f}^{(1)}}{\del x} \right] 
+ \mathcal{O}\left( \epsilon^{3} \right)
= - \left[ {\bf f}^{(0)} + \epsilon {\bf f}^{(1)} 
+ \mathcal{O}\left( \epsilon^{2} \right) - {\bf f}^{eq} \right] $$   
Now, let us collect terms of different orders of $\epsilon$.  

\subsubsection*{\bf $\epsilon^{0}$ terms:} 

$$ 0 = - {\bf f}^{(0)} + {\bf f}^{eq}$$ 
or  
\be 
{\bf f}^{(0)} = {\bf f}^{eq} 
\ee  

\subsubsection*{Zeroth approximation:} 
If we use the zeroth approximation ${\bf f} = {\bf f}^{(0)} = {\bf f}^{eq}$ in FVBE and take moments, we obtain 
$$ {\bf P} \left[ \fr{\del {\bf f}^{eq}}{\del t} + \Lambda \fr{\del {\bf f}^{eq}}{\del x} 
= - \fr{1}{\epsilon} \left[{\bf f}^{eq} - {\bf f}^{eq} \right] \right] $$ 
or 
$$ \fr{\del \left({\bf P}{\bf f}^{eq}\right)}{\del t} 
+ \fr{\del \left({\bf P} \Lambda {\bf f}^{eq}\right)}{\del x} = 0 $$    
or 
\be  
\fr{\del u}{\del t} + \fr{\del g(u)}{\del x} = 0 
\ee 
which is the original convection equation that is modeled.  

\subsubsection*{\bf $\epsilon^{1}$ terms:} 

$$ \fr{\del {\bf f}^{(0)}}{\del t} + \Lambda \fr{\del {\bf f}^{(0)}}{\del x} 
= - {\bf f}^{(1)} $$ 
which, after using ${\bf f}^{(0)} = {\bf f}^{eq}$, leads to  
$$ {\bf f}^{(1)} = - \left[ \fr{\del {\bf f}^{eq}}{\del t} 
+ \Lambda \fr{\del {\bf f}^{eq}}{\del x} \right] $$ 
Using it in the expression ${\bf f} = {\bf f}^{(0)} + \epsilon {\bf f}^{(1)}$, we obtain 
\be 
{\bf f} = {\bf f}^{eq} - \epsilon \left[ \fr{\del {\bf f}^{eq}}{\del t} 
+ \Lambda \fr{\del {\bf f}^{eq}}{\del x} \right] 
\ee

\subsubsection*{Chapman-Enskog approximation} 
Let us denote the first approximation as the Chapman-Enskog approximation and write 
\be 
{\bf f}^{CE} = {\bf f}^{(0)} + \epsilon {\bf f}^{(1)} 
= {\bf f}^{eq} - \epsilon \left[ \fr{\del {\bf f}^{eq}}{\del t} 
+ \Lambda \fr{\del {\bf f}^{eq}}{\del x} \right] 
\ee 
To evaluate the derivatives of the equilibrium distribution, let us start with 
$$ {\bf f}^{eq} = \left[ \ba{c} \fr{1}{2} u + \fr{1}{2 \lambda} g(u) \\[4mm]  
\fr{1}{2} u - \fr{1}{2 \lambda} g(u) \ea \right] 
= \fr{1}{2} \left[ \ba{c} u \\ u \ea \right] + \fr{1}{2 \lambda} \left[ \ba{c} g(u) \\ - g(u) \ea \right] $$
Therefore 
$$ \fr{\del {\bf f}^{eq}}{\del t} 
= \fr{1}{2} \left[ \ba{c} \fr{\del u}{\del t} \\[4mm] \fr{\del u}{\del t} \ea \right] 
+ \fr{1}{2 \lambda} \left[ \ba{c} \fr{\del g(u)}{\del t} \\[4mm] - \fr{\del g(u)} {\del t} \ea \right] $$ 
Now, using the original convection equation, which is already derived as the zeroth approximation, we can write 
$$ \fr{\del u}{\del t} = - \fr{\del g(u)}{\del x} = - \fr{\del g(u)}{\del u} \fr{\del u}{\del x} 
= - a(u) \fr{\del u}{\del x} $$  
and 
$$ \fr{\del g(u)}{\del t} = \fr{\del g(u)}{\del u} \fr{\del u}{\del t} = a(u) \fr{\del u}{\del t} 
= a(u) \left[ - a(u) \fr{\del u}{\del x} \right] = - \left(a(u)\right)^{2} \fr{\del u}{\del x} $$ 
Using the above conversion of time derivatives to space derivatives, which is typical of the Chapman-Enskog expansion, we obtain 
$$ \fr{\del {\bf f}^{eq}}{\del t} 
= \fr{1}{2} \left[ \ba{c} - a(u) \fr{\del u}{\del x} \\[4mm] - a(u) \fr{\del u}{\del x} \ea \right] 
+ \fr{1}{2 \lambda} \left[ \ba{c} - \left(a(u)\right)^{2} \fr{\del u}{\del x} \\[4mm] 
\left(a(u)\right)^{2} \fr{\del u}{\del x}  \ea \right] 
= \fr{1}{2} \fr{\del u}{\del x} \left[ \ba{c} - a(u) \\ - a(u) \ea \right] 
+ \fr{1}{2 \lambda} a(u) \fr{\del u}{\del x} \left[ \ba{c} -a(u) \\ a(u) \ea \right] $$ 
Now 
$$ \ba{rcl} 
\Lambda \fr{\del {\bf f}^{eq}}{\del x} 
& = & \Lambda \left[
\fr{1}{2} \left[ \ba{c} \fr{\del u}{\del x} \\[4mm] \fr{\del u}{\del x} \ea \right] 
+ \fr{1}{2 \lambda} \left[ \ba{c} \fr{\del g(u)}{\del x} \\[4mm] - \fr{\del g(u)}{\del x} \ea \right]
\right] \\[10mm]  
& = & \fr{1}{2} \fr{\del u}{\del x} 
\left[ \ba{c} \lambda \\ - \lambda \ea \right] 
+ \fr{1}{2 \lambda} a(u) \fr{\del u}{\del x} 
\left[ \ba{c} \lambda \\ \lambda \ea \right]   
\ea $$ 
Therefore 
$$ \ba{rcl} 
\fr{\del {\bf f}^{eq}}{\del t} + \Lambda \fr{\del {\bf f}^{eq}}{\del x} 
& = & \fr{1}{2} \fr{\del u}{\del x} \left[ \ba{c} - a(u) \\ - a(u) \ea \right] 
+ \fr{1}{2 \lambda} a(u) \fr{\del u}{\del x} \left[ \ba{c} -a(u) \\ a(u) \ea \right] 
+ \fr{1}{2} \fr{\del u}{\del x} 
\left[ \ba{c} \lambda \\ - \lambda \ea \right] 
+ \fr{1}{2 \lambda} a(u) \fr{\del u}{\del x} 
\left[ \ba{c} \lambda \\ \lambda \ea \right] \\[5mm] 
& = & \fr{1}{2} \fr{\del u}{\del x} 
\left[ \ba{c} \lambda - \fr{\left(a(u)\right)^{2}}{\lambda} \\ 
- \lambda + \fr{\left(a(u)\right)^{2}}{\lambda} \ea \right] \\ 
& = & \fr{1}{2} \fr{\del u}{\del x} 
\left[ \lambda^{2} - \left(a(u)\right)^{2} \right] 
\left[ \ba{c} \fr{1}{\lambda} \\[4mm] - \fr{1}{\lambda} \ea \right] 
\ea \\ $$ 
Therefore, from the earlier expression,  
$ {\bf f}^{CE} = {\bf f}^{eq} - \epsilon \left[ \fr{\del {\bf f}^{eq}}{\del t} 
+ \Lambda \fr{\del {\bf f}^{eq}}{\del x} \right]$,   
we obtain 
\be 
{\bf f}^{CE} = {\bf f}^{eq} - \epsilon \fr{1}{2} \fr{\del u}{\del x} 
\left[ \lambda^{2} - \left(a(u)\right)^{2} \right] 
\left[ \ba{c} \fr{1}{\lambda} \\[4mm] - \fr{1}{\lambda} \ea \right]  
\ee   
We now define 
\be 
\nu = \epsilon \left[ \lambda^{2} - \left(a(u)\right)^{2} \right] 
\ee
Therefore 
\be 
{\bf f}^{CE} = {\bf f}^{eq} -  \fr{1}{2} \nu \fr{\del u}{\del x} 
\left[ \ba{c} \fr{1}{\lambda} \\[4mm] - \fr{1}{\lambda} \ea \right]  
\ee 
which we write as 
\be 
{\bf f}^{CE} = {\bf f}^{eq} - {\bf f}^{eq}_{v} 
\ee 
where 
\be 
{\bf f}^{eq}_{v} = \fr{1}{2} \nu \fr{\del u}{\del x} 
\left[ \ba{c} \fr{1}{\lambda} \\[4mm] - \fr{1}{\lambda} \ea \right]  
\ee
\subsubsection*{Some useful moments for the discrete Chapman-Enskog distribution} 
\be 
{\bf f}^{CE} = {\bf f}^{eq} - {\bf f}_{v}^{eq}
\ee 
where 
\be
{\bf f}_{v}^{eq} = \fr{1}{2} \nu \fr{\del u}{\del x} 
\left[ \ba{c} \fr{1}{\lambda} \\[4mm] - \fr{1}{\lambda} \ea \right], \ 
\nu = \epsilon \left[\lambda^{2} - \left(a(u)\right)^{2} \right] 
\ee 
Let us evaluate some useful moments of ${\bf f}_{v}^{eq}$.  
$$ \ba{rcl} 
{\bf P}{\bf f}_{v}^{eq} & = & \fr{1}{2} \nu \fr{\del u}{\del x} {\bf P} 
\left[ \ba{c} \fr{1}{\lambda} \\[4mm] - \fr{1}{\lambda} \ea \right] \\[10mm] 
& = & \fr{1}{2} \nu \fr{\del u}{\del x} 
\left[ 1 \ 1 \right] 
\left[ \ba{c} \fr{1}{\lambda} \\[4mm] - \fr{1}{\lambda} \ea \right] \\[5mm]  
& = & 0 
\ea $$ 
or 
\be 
{\bf P}{\bf f}_{v}^{eq} = 0 
\ee 
$$ \ba{rcl} 
{\bf P} \Lambda {\bf f}_{v}^{eq} & = & \fr{1}{2} \nu \fr{\del u}{\del x} 
{\bf P} \Lambda 
\left[ \ba{c} \fr{1}{\lambda} \\[4mm]  - \fr{1}{\lambda} \ea \right] \\[10mm] 
& = & \fr{1}{2} \nu \fr{\del u}{\del x} 
\left[ 1 \ 1 \right] 
\left[\ba{cc} \lambda & 0 \\ 0 & - \lambda \ea \right] 
\left[ \ba{c} \fr{1}{\lambda} \\[4mm] - \fr{1}{\lambda} \ea \right] \\[10mm] 
& = & \nu \fr{\del u}{\del x} 
\ea $$ 
Therefore 
\be 
{\bf P} \Lambda {\bf f}_{v}^{eq} =  \nu \fr{\del u}{\del x}  = g_{v}(u) 
\ee 
$$ \ba{rcl} 
{\bf P}\Lambda ^{2}{\bf f}_{v}^{eq} & = & \fr{1}{2} \nu \fr{\del u}{\del x} 
\left[ 1 \ 1 \right] \left[\ba{cc} \lambda & 0 \\ 0 & - \lambda \ea \right]  
\left[\ba{cc} \lambda & 0 \\ 0 & - \lambda \ea \right] 
\left[ \ba{c} \fr{1}{\lambda} \\[4mm] - \fr{1}{\lambda} \ea \right] \\[5mm] 
& = & 0 
\ea $$ 
Thus 
\be 
{\bf P}\Lambda ^{2}{\bf f}_{v}^{eq} = 0 
\ee 
$$ \ba{rcl} 
{\bf P} \Lambda^{+} {\bf f}_{v}^{eq} & = & \fr{1}{2} \nu \fr{\del u}{\del x} 
\left[ 1 \ 1 \right] \left[ \ba{cc} \lambda & 0 \\ 0 & 0 \ea \right] 
\left[ \ba{c} \fr{1}{\lambda} \\[4mm] - \fr{1}{\lambda} \ea \right] \\[10mm] 
& = & \fr{1}{2} \nu \fr{\del u}{\del x} 
\left[ 1 \ 1 \right] \left[ \ba{c} 1 \\ 0 \ea \right] \\[2mm]  
& = &   \fr{1}{2}g_{v}(u)
\ea $$ 
$$ \ba{rcl} 
{\bf P} \Lambda^{-} {\bf f}_{v}^{eq} & = & \fr{1}{2} \nu \fr{\del u}{\del x} 
\left[ 1 \ 1 \right] \left[ \ba{cc} 0 & 0 \\ 0 & - \lambda \ea \right] 
\left[ \ba{c} \fr{1}{\lambda} \\[4mm] - \fr{1}{\lambda} \ea \right] \\[10mm] 
& = & \fr{1}{2} \nu \fr{\del u}{\del x} 
\left[ 1 \ 1 \right] \left[ \ba{c} 0 \\ 1 \ea \right] \\[2mm]  
& = &   \fr{1}{2}g_{v}(u)
\ea $$ 
Therefore 
\be 
{\bf P} \Lambda^{\pm} {\bf f}_{v}^{eq} =  \fr{1}{2}g_{v}(u) 
\ee 
Thus 
\bea
{\bf P}{\bf f}^{CE} & = & {\bf P} {\bf f}^{eq}  - {\bf P} {\bf f_{v}}^{eq}  = u\\[2mm]
{\bf P} \Lambda {\bf f}^{CE} & = & {\bf P} \Lambda {\bf f}^{eq} -{\bf P} \Lambda {\bf f_{v}}^{eq} = g(u) -  g_{v}(u) 
\eea

\subsubsection*{FVBE as vanishing diffusion approximation for the inviscid case} 
Let us now use the derived Chapman-Enskog distribution function in FVBE and obtain the moment.  
$$ {\bf P} \left[ \fr{\del {\bf f}^{CE}}{\del t} 
+ \Lambda \fr{\del {\bf f}^{CE}}{\del x} = - \fr{1}{\epsilon} 
\left[ {\bf f}^{CE} - {\bf f}^{eq} \right] \right] $$  
or 
$$ \fr{\del \left({\bf P} {\bf f}^{CE}\right)}{\del t} + 
\fr{\del \left({\bf P} \Lambda {\bf f}^{CE} \right)}{\del x} 
= - \fr{1}{\epsilon} \left[ {\bf P} {\bf f}^{CE} - {\bf P} {\bf f}^{eq} \right] $$   
or 
$$ \fr{\del u}{\del t} + \fr{\del }{\del x} 
\left[ g(u) - g_{v}(u) \right]
= - \fr{1}{\epsilon} \left[ u - u \right] = 0 $$ 
or 
\be 
\fr{\del u}{\del t} + \fr{\del g(u)}{\del x} = \fr{\del g_{v}(u)}{\del x} 
\ee 
Thus, the use of Chapman-Enskog distribution leads to the viscous Burgers equation.  However, if our concern is only with the inviscid convection equation, the Chapman-Enskog type expansion leads to the following expression, which can be interpreted as the FVBE leading to a vanishing diffusion approximation to the nonlinear convection equation. 
\be 
\fr{\del u}{\del t} + \fr{\del g(u)}{\del x} = \epsilon \fr{\del}{\del x} 
\left[ \left\{ \lambda^{2} - \left( a(u) \right)^{2} \right\} \fr{\del u}{\del x} \right] 
\ee

Therefore, the Chapman-Enskog type expansion of the flexible velocity Boltzmann equation has yielded a convection equation augmented with vanishing diffusion, as the right hand side contains a second derivative.  For the model to be stable, the diffusion must be non-negative.  Therefore 
\be 
\lambda^{2} - \left(a(u)\right)^{2} \ge 0 
\ee 
which can be used to fix the value of $\lambda$.    
\newpage 

\noindent 
{\bf \large Appendix D} 
\section*{Integration of Boltzmann equation to obtain shallow water equations}

Consider 1D Boltzmann equations including a force term as given by
\be
\fr{\del  f}{\del t} + v \fr{\del f}{\del x} +\mathbb{F}\fr{\del f}{\del v} = - \fr{1}{\epsilon} \left[  f - f^{eq} \right] 
\ee
where $\mathbb{F} = \mathbb{F}(x)$ is the Force acting on the fluid. The equilibrium distribution function, $f^{eq}$ is defined as 
\be 
f^{eq} = \phi \sqrt{\fr{\beta}{\pi}} e^{- {\beta}\left( v - {u} \right)^{2} } , \ \ \  \textrm{where} \ \beta = \fr{1}{\phi} \ \textrm{which \ leads \ to} \ f^{eq} = \sqrt{\fr{\phi}{\pi}} e^{- \fr{\left( v - u \right)^{2}}{\phi}} 
\ee 
Here $\phi$ is the geopotential and is defined by $\phi=gh$.  To obtain the 1D shallow water equations given in [\ref{1D_SWE}] and [\ref{1D_SWE_COMP}], we multiply [\ref{1D_BOLTZ}] by a vector $\psi$ such as 
\be  \label{1D_BOLTZ3}
\psi \left \{ \fr{\del f}{\del t} + v \fr{\del f}{\del x} +\mathbb{F}\fr{\del f}{\del v} = - \fr{1}{\epsilon} \left[  f -  f^{eq} \right] \right \} , \ \ where \ \psi = \left[ \begin{array}{c} 1 \\ v \end{array}\right].  
\ee 
and take moments. The moment thus taken for each of the term is explained below. \\

{( i ) $ \int\limits_{- \infty}^{\infty}{ {\psi} \fr{\del f}{\del t}}dv$}\\
\be
\int\limits_{- \infty}^{\infty}{ {\psi} \fr{\del f}{\del t}}dv =   \fr{\del }{\del t}\int\limits_{- \infty}^{\infty}{\psi {f}_{eq}}dv -\int\limits_{- \infty}^{\infty}\bcancelto{{f}_{eq}\fr{\del \psi }{\del t}}{0}dv = \ \ \fr{\del }{\del t}\int\limits_{- \infty}^{\infty}{\psi {f}_{eq}}dv
\ee

\bea
\psi = 1 :  \fr{\del }{\del t}\int\limits_{- \infty}^{\infty}{{f}_{eq}}dv \ = \fr{\del \phi}{\del t}\\
\psi = v :  \fr{\del }{\del t}\int\limits_{- \infty}^{\infty}{v {f}_{eq}}dv \ = \fr{\del (\phi u)}{\del t} 
\eea

{( ii ) $\int\limits_{- \infty}^{\infty} {\psi} v \fr{\del f}{\del x}dv$}\\
\be
\int\limits_{- \infty}^{\infty} {\psi} v \fr{\del f}{\del x}dv = \fr{\del }{\del x}\int\limits_{- \infty}^{\infty}{\psi v {f}_{eq}}dv -\int\limits_{- \infty}^{\infty}\bcancelto{\psi {f}_{eq}\fr{\del v }{\del x}}{0}dv - \int\limits_{- \infty}^{\infty}{v {f}_{eq}\fr{\del \psi }{\del x}}dv
\ee
\bea
\psi = 1 :  \int\limits_{- \infty}^{\infty}  v \fr{\del f}{\del x}dv \ = \fr{\del }{\del x}\int\limits_{- \infty}^{\infty}{v {f}_{eq}}dv - \int\limits_{- \infty}^{\infty}{v {f}_{eq}\fr{\del (1) }{\del x}}dv = \fr{\del (\phi u)}{\del x}\\
\psi = v :  \int\limits_{- \infty}^{\infty} {v}{v}\fr{\del f}{\del x}dv \ =  \fr{\del }{\del x}\int\limits_{- \infty}^{\infty}{{v}^{2} {f}_{eq}}dv - \int\limits_{- \infty}^{\infty}{v {f}_{eq}\fr{\del (v) }{\del x}}dv = \fr{\del }{\del x} (\fr{{\phi}^{2}}{2}+\phi {u}^{2})
\eea

{( iii ) $\int\limits_{- \infty}^{\infty} {\psi} \mathbb{F}\fr{\del f}{\del v}dv$}\\
Applying integration by parts the integrand can be expanded as 
\be
\int\limits_{- \infty}^{\infty} {\psi} \mathbb{F}\fr{\del f}{\del v}dv =  \fr{\del }{\del v}\int\limits_{- \infty}^{\infty}{\psi \mathbb{F} {f}_{eq}}dv- \int\limits_{- \infty}^{\infty}{\psi {f}_{eq}\fr{\del \mathbb{F} }{\del v}}dv - \int\limits_{- \infty}^{\infty}{\mathbb{F} {f}_{eq}\fr{\del \psi }{\del v}}dv.
\ee
To evaluate the first term of the expansion we use divergence theorem and convert the volume integral into a surface integral in the velocity space given by $\int{\psi \mathbb{F} {f}_{eq}}dv = \oint {\psi \mathbb{F} {f}_{eq}}d{S}_{v}$. The surface boundary is positioned at infinity where the distribution function goes to zero. Therefore the first term in the integral vanishes. The second term in the integral vanishes too as, $ \mathbb{F}$ being a conservative force and is independent of the velocity space. Therefore the integral can be evaluated as 

\bea
\int\limits_{-\infty}^{\infty} {\psi} \mathbb{F}\fr{\del f}{\del v}dv =  \fr{\del }{\del v}\int\limits_{- \infty}^{\infty}\bcancelto{\psi \mathbb{F} {f}_{eq}}{0}dv- \int\limits_{- \infty}^{\infty}\bcancelto{\psi {f}_{eq}\fr{\del \mathbb{F} }{\del v}}{0}dv - \int\limits_{- \infty}^{\infty}{\mathbb{F} {f}_{eq}\fr{\del \psi }{\del v}}dv. \\
\psi = 1 : \int\limits_{- \infty}^{\infty} {\psi} \mathbb{F}\fr{\del f}{\del v}dv = - \int\limits_{- \infty}^{\infty}{\mathbb{F} {f}_{eq}\fr{\del \psi }{\del v}}dv = - \int\limits_{- \infty}^{\infty}{\mathbb{F} {f}_{eq}\fr{\del (1) }{\del v}}dv = 0 \\
\psi = v :  \int\limits_{- \infty}^{\infty} {\psi} \mathbb{F}\fr{\del f}{\del v}dv =  - \int\limits_{- \infty}^{\infty}{\mathbb{F} {f}_{eq}\fr{\del (v) }{\del v}}dv = -\mathbb{F} \phi
\eea  
Therefore, the obtained conservation equations can be written as 
\be
\fr{\partial U}{\partial t}+\fr{\partial G(U)}{\partial x}+ S(U) = 0 .
\ee
where $U$ is the conserved variable vector, $G(U)$ the Flux vector and $S(U)$ is the source term vector. These are defined by
\be \label{1D_SWE_COMP2}
U=\left[\begin{array}{c} \phi \\ \phi u \end{array}\right], \quad G(U)=\left[\begin{array}{c} \phi u \\ \phi u^{2}+\frac{1}{2} {\phi}^{2} \end{array}\right] \quad and \quad S(U)=\left[\begin{array}{c} 0 \\ -\mathbb{F} \phi \end{array}\right] 
\ee 
by substituting the definition of $\phi$  and $\mathbb{F} = g{b}_{x}$ we get
\be
U=\left[\begin{array}{c} h \\ h u \end{array}\right], \quad G(U)=\left[\begin{array}{c} h u \\ h u^{2}+\frac{1}{2} g h^{2} \end{array}\right] \quad and \quad S(U)=\left[\begin{array}{c} 0 \\ -g h b_{x} \end{array}\right] 
\ee 
which form the definition of the conserved variable vector and flux vector for the shallow water equations, together with the source term.

%

\end{appendices}

\end{document}